\documentclass{amsart}
\usepackage{amssymb}
\usepackage{mathrsfs}
\usepackage{stmaryrd}
\usepackage{bbm}
\usepackage{oldgerm}
\usepackage[english]{babel}
\usepackage[T1]{fontenc}
\usepackage[latin1]{inputenc}
\usepackage[all]{xy}
\usepackage{hyperref}
\usepackage{upref}
\usepackage{xcolor}
\usepackage{tikz-cd}

\hypersetup{
 colorlinks,
 linkcolor={red!50!black},
 citecolor={blue!50!black},
 urlcolor={blue!80!black}
}

\newtheorem{teo}[subsection]{Theorem}
\newtheorem{prop}[subsection]{Proposition}
\newtheorem{cor}[subsection]{Corollary}
\newtheorem{lem}[subsection]{Lemma}

\theoremstyle{definition}

\newtheorem{defi}[subsection]{Definition}
\newtheorem{rema}[subsection]{Remark}

\newtheorem{exemple}[subsection]{Example}

\newtheorem{por}[subsection]{Porism}
\newtheorem{notat}[subsection]{Notation}

\numberwithin{equation}{subsection}

\newcommand{\adm}{\mathrm{adm}}
\newcommand{\abx}{\mathrm{Ab}(S_\mathrm{fppf})}
\newcommand{\Aut}{\mathrm{Aut}}
\newcommand{\bk}{\overline{k}}
\newcommand{\colim}{\mathrm{colim}}
\newcommand{\cyc}{\mathrm{cyc}}
\newcommand{\End}{\mathrm{End}}
\newcommand{\ev}{\mathrm{ev}}
\newcommand{\Ext}{\mathrm{Ext}}
\newcommand{\fppf}{\mathrm{fppf}}
\newcommand{\F}{\mathcal{F}}
\newcommand{\fin}{\mathrm{fin}}
\newcommand{\Frob}{\mathrm{Frob}}
\newcommand{\Fet}{\mathrm{F}\acute{\mathrm{e}}\mathrm{t}}
\newcommand{\G}{\mathcal{G}}
\newcommand{\Ga}{\mathbb{G}}
\newcommand{\Hl}{\mathcal{H}}
\newcommand{\GL}{\mathrm{GL}}
\newcommand{\Gal}{\mathrm{Gal}}
\newcommand{\Hom}{\mathrm{Hom}}
\newcommand{\Ind}{\mathrm{Ind}}
\newcommand{\id}{\mathrm{id}}

\newcommand{\Lc}{\mathcal{L}}
\newcommand{\Mac}{\mathcal{M}}
\newcommand{\Loc}{\mathrm{Loc}}
\newcommand{\m}{\mathfrak{m}}
\newcommand{\Ow}{\mathcal{O}}

\newcommand{\Pic}{\mathrm{Pic}}

\newcommand{\pr}{\mathrm{pr}}

\newcommand{\Rep}{\mathrm{Rep}}
\newcommand{\Res}{\mathrm{Res}}
\newcommand{\rk}{\mathrm{rk}}
\newcommand{\Sh}{\mathrm{Sh}}
\newcommand{\sw}{\mathrm{sw}}
\newcommand{\Spec}{\mathrm{Spec}}
\newcommand{\Sym}{\mathrm{Sym}}

\newcommand{\Tr}{\mathrm{Tr}}
\newcommand{\Trip}{\mathrm{Trip}}
\newcommand{\V}{\mathcal{V}}
\newcommand{\Ver}{\mathrm{Ver}}

\newcommand{\Z}{\mathbb{Z}}

\title{Geometric local $\varepsilon$-factors}
\author{Quentin Guignard}
\address{Institut des
Hautes \'Etudes Scientifiques, 35 route de Chartres, 91440 Bures-sur-Yvette, France}
\address{
\'Ecole Normale Sup\'erieure, 45 rue d'Ulm, 75005 Paris, France}
\email{quentin.guignard@ens.fr}

\begin{document}

\begin{abstract}
Inspired by the work of Laumon on local $\varepsilon$-factors and by Deligne's $1974$ letter to Serre, we give an explicit cohomological definition of $\varepsilon$-factors for $\ell$-adic Galois representations over henselian discrete valuation fields of positive equicharacteristic $p \neq \ell$, with (not necessarily finite) perfect residue fields. These geometric local $\varepsilon$-factors are completely characterized by an explicit list of purely local properties, such as an induction formula and the compatibility with geometric class field theory in rank $1$, and satisfy a product formula for $\ell$-adic sheaves on a curve over a perfect field of characteristic $p$.
\end{abstract}

\maketitle
\tableofcontents

\section{Introduction}

\subsection{\label{chap310.0}} The theory of local $\varepsilon$-factors over local fields with finite residue fields originated from Tate's thesis in rank $1$, and was brought to its current form by works of Dwork \cite{Dw}, Langlands \cite{L}, Deligne \cite{De73} and Laumon \cite{La87}. Central motivations for these developments were the problem of decomposing the constants of functional equations of Artin's $L$-functions, or of Weil's $L$-functions, as a product of local contributions, and the applications of such a decomposition to Langlands program through Deligne's recurrence principle, cf. (\cite{La87}, 3.2.2). Inspired by the work of Laumon \cite{La87} and by Deligne's $1974$ letter to Serre (\cite{bloch}, Appendix), we provide in this text an explicit cohomological construction of $\varepsilon$-factors for $\ell$-adic Galois representations over equicharacteristic henselian discrete valuation fields, with (not necessarily finite) perfect residue fields of positive characteristic, such as the field of Laurent series $k((t))$ for any perfect field $k$ of positive characteristic $p$. As it turns out, these geometric local $\varepsilon$-factors fit into a product formula for the determinant of the cohomology of an $\ell$-adic sheaf on a curve over a perfect field of characteristic $p$.

After the writing of this text, we were informed by Takeshi Saito that results similar to ours were stated in a work of Seidai Yasuda \cite{Ya3}, itself relying on previous works \cite{Ya1} and \cite{Ya2} by the same author. It appears that the approach used by Yasuda is different from ours: in order to construct local $\varepsilon$-factors and to establish their key properties, he first considers the case of finite coefficients, which enables him to use a spreading out argument in order to be able to rely on previous results regarding the finite residue field case, while we adopted instead a self-contained apprach by first giving an explicit cohomological construction of local $\varepsilon$-factors, cf. \ref{chap310.3}, and then by establishing directly its main properties, most notably the induction formula.


\subsection{\label{chap310.5}} Let us first recall the classical theory for local fields with finite residue fields. We restrict to the equicharacteristic case, and we give a slightly non standard presentation as a preparation for our extension to the case of a general perfect residue field (cf. \ref{chap310.2}). Let us fix an algebraic closure $\overline{\mathbb{F}}_p$ of $\mathbb{F}_p$. Let $\ell$ be a prime number distinct from $p$ and let $\psi : \mathbb{F}_p \rightarrow \overline{\mathbb{Q}}_{\ell}^{\times}$ be a non trivial homomorphism. Let us consider quadruples $(T,\F,\omega, \overline{s})$ where $T$ is a henselian trait of equicharacteristic $p$, whose closed point $s$ is finite over $\mathbb{F}_p$, where $\F$ is a constructible \'etale $ \overline{\mathbb{Q}}_{\ell}$-sheaf on $T$, where $\omega$ is a non zero meromorphic $1$-form on $T$ (cf. \ref{chap32.12}), and where $\overline{s} : \Spec(\overline{\mathbb{F}}_p) \rightarrow T$ is a morphism of schemes.

A \textit{theory of $\ell$-adic local $\varepsilon$-factors over $\mathbb{F}_p$, with respect to $\psi$,} is a rule $\varepsilon$ which assigns to any such quadruple $(T,\F,\omega,\overline{s})$ a homomorphism $\varepsilon_{\overline{s}}(T,\F,\omega)$ from the Galois group $\Gal(\overline{s}/s)$ to $ \overline{\mathbb{Q}}_{\ell}^{\times}$, and which satisfies the following axioms:
\begin{enumerate}
\item the homomorphism $\varepsilon_{\overline{s}}(T,\F,\omega)$ depends only on the isomorphism class of the quadruple $(T,\F,\omega,\overline{s})$;
\item there exists a finite extension $E$ of $\mathbb{Q}_{\ell}$ contained in $ \overline{\mathbb{Q}}_{\ell}$, depending on $\F$, such that $\varepsilon_{\overline{s}}(T,\F,\omega)$ is a continuous homomorphism from $\Gal(\overline{s}/s)$ to $E^{\times}$;
\item for any exact sequence
$$
0 \rightarrow \F' \rightarrow \F \rightarrow \F'' \rightarrow 0,
$$
of constructible \'etale $ \overline{\mathbb{Q}}_{\ell}$-sheaves on $T$, we have 
$$
\varepsilon_{\overline{s}}(T,\F,\omega) = \varepsilon_{\overline{s}}(T,\F',\omega) \varepsilon_{\overline{s}}(T,\F'',\omega);
$$
\item if $\F$ supported on the closed point of $T$, then $\varepsilon_{\overline{s}}(T,\F,\omega)$ is the $\ell$-adic character of $\Gal(\overline{s}/s)$ corresponding to the $1$-dimensional representation $\det \left( \F_{\overline{s}} \right)^{-1}$;
\item \label{chap3fr5}for each finite generically \'etale extension $f : T' \rightarrow T$ of henselian traits, there exists a homomorphism $\lambda_f(\omega)$ from the Galois group $\Gal(\overline{s}/s)$ to $ \overline{\mathbb{Q}}_{\ell}^{\times}$ such that
$$
\varepsilon_{\overline{s}}(T,f_*\F,\omega) = \lambda_f(\omega)^{\rk(\F)} \delta_{s'/s}^{a(T',\F,f^* \omega)} \Ver_{s'/s} \left( \varepsilon_{\overline{s}'}(T',\F,f^*\omega) \right), 
$$
for any constructible \'etale $ \overline{\mathbb{Q}}_{\ell}$-sheaf $\F$ on $T'$, of generic rank $\rk(\F)$, where the verlagerung $\Ver_{s'/s}$ and the signature $\delta_{s'/s}$ are defined in \ref{chap3verlagerung}, and the conductor $a(T,\F,f^* \omega)$ is defined in \ref{chap32.11};
\item \label{chap3fr6}if $j : \eta \rightarrow T$ is the inclusion of the generic point of $T$ and if $\F$ is a lisse \'etale $ \overline{\mathbb{Q}}_{\ell}$-sheaf of rank $1$ on $\eta$, then we have
$$
(-1)^{a(T,j_* \F)} \varepsilon_{ \overline{s}}(T,j_* \F,\omega) (\Frob_s) = \varepsilon(\chi_{\F}, \Psi_{\omega}),
$$
where $\Frob_s$ is the geometric Frobenius element of $\Gal(\overline{s}/s)$, where the conductor $a(T,j_*\F)$ is defined in \ref{chap32.11}, where $\Psi_{\omega} : k(\eta) \rightarrow \Lambda^{\times}$ is the additive character given by $z \mapsto \psi(\Tr_{k/\mathbb{F}_p}(z \omega))$, where $\chi_{\F}$ is the character of $k(\eta)^{\times}$ associated to $\F$ by local class field theory, and where $ \varepsilon(\chi_{\F}, \Psi_{\omega})$ is the automorphic $\varepsilon$-factor of the pair $(\chi_{\F}, \Psi_{\omega})$, cf. (\cite{La87}, 3.1.3.2).

\end{enumerate}

\begin{teo}\label{chap3teo1} For any prime number $p$, any prime number $\ell$ distinct from $p$ and any non trivial homomorphism $\psi : \mathbb{F}_p \rightarrow \overline{\mathbb{Q}}_{\ell}^{\times}$, there exists a unique theory of $\ell$-adic local $\varepsilon$-factors over $\mathbb{F}_p$, with respect to $\psi$.
\end{teo}

Since the Galois group $\Gal(\overline{s}/s)$ is procyclic, the $\ell$-adic character $\varepsilon_{ \overline{s}}(T,j_* \F,\omega)$ is completely determined by its value at the geometric Frobenius element of $\Gal(\overline{s}/s)$. Actually, the rule which associates the quantity 
$$
(-1)^{a(T,\F)} \varepsilon_{\overline{s}}(T,\F,\omega)(\Frob_s),
$$
to a quadruple $(T,\F,\omega, \overline{s})$, where $a(T,\F)$ is the conductor of the pair $(T,\F)$ (cf. \ref{chap32.11}), satisfies the properties listed in (\cite{La87}, 3.1.5.4). Thus Theorem \ref{chap3teo1} is a reformulation of the theorem of Langlands \cite{L} and Deligne \cite{De73} regarding the existence and uniqueness of local $\varepsilon$-factors. 

The proof of existence by Deligne and Langlands in the finite field case is somewhat indirect: starting with the prescribed values (\ref{chap3fr6}) of the local $\varepsilon$-factors in rank $1$, local $\varepsilon$-factors are defined in arbitrary rank by Brauer's theory and by the induction property (\ref{chap3fr5}), and the main problem is then to prove that the resulting factors are independent of the choices made. Our approach is different: we first give a simple cohomological definition of local $\varepsilon$-factors in arbitrary rank (cf. \ref{chap3localepsfact}) using the theory of Gabber-Katz extensions (cf. \ref{chap310.2} below), and we use Brauer's theory only to establish the main properties of these local $\varepsilon$-factors. 

Laumon gave a cohomological formula for local $\varepsilon$-factors over local fields with finite residue fields (\cite{La87}, 3.5.1.1). If $\F$ is supported on the generic point $\eta$ of $T$, then Laumon's formula takes the following form:
\begin{align}\label{chap3laumonform}
\varepsilon_{\overline{s}}(T,\F,d \pi) = \det(F_{\pi}^{(0,\infty')}(\F)) \circ \sigma_{\pi},
\end{align}
where $\pi$ is a uniformizer of $k(\eta)$, where Laumon's \emph{local Fourier transform} $F_{\pi}^{(0,\infty')}(\F)$ is an $\ell$-adic representation of $\Gal(\overline{\eta}/\eta)$, cf. (\cite{La87}, 2.4.1), and where $\sigma_{\pi} : \Gal(\overline{s}/s) \rightarrow \Gal(\overline{\eta}/\eta)^{\mathrm{ab}}$ is the section of the natural homomorphism $\Gal(\overline{\eta}/\eta)^{\mathrm{ab}} \rightarrow \Gal(\overline{s}/s)$ corresponding by local class field theory to the unique section of the valuation homomorphism $k(\eta)^{\times} \rightarrow \mathbb{Z}$ sending the element $1$ of $\mathbb{Z}$ to $\pi$.

It is straightforward to extend (\ref{chap3laumonform}) to a rule $\varepsilon_{\mathrm{Lau}}$ satisfying the properties $(1),(2),(3)$ and $(4)$ of a theory of $\ell$-adic local $\varepsilon$-factors over $\mathbb{F}_p$. Moreover, the normalization in rank $1$, namely property $(6)$, can be proved directly for $\varepsilon_{\mathrm{Lau}}$ by using Laumon's $\ell$-adic stationary phase method from \cite{La87}. Unfortunately, there seems to be no direct proof that $\varepsilon_{\mathrm{Lau}}$ satisfies the property $(5)$, namely the induction formula, and it is therefore not possible to take Laumon's formula as a definition of local $\varepsilon$-factors. However, the $\ell$-adic stationary phase method yields that the rule $\varepsilon_{\mathrm{Lau}}$ produced from Laumon's formula (\ref{chap3laumonform}) coincides with our own definition (cf. \ref{chap310.3}) in the finite field case (cf. \ref{chap3xteolaumon}). Thus our main Theorem \ref{chap3teo2} below, in conjunction with the $\ell$-adic stationary phase method, proves the induction formula for $\varepsilon_{\mathrm{Lau}}$.
 
By using (\ref{chap3laumonform}) and the $\ell$-adic stationary phase method, Laumon proved the following product formula: 

\begin{teo}[\cite{La87}, Th. 3.2.1.1]\label{chap3fr} Let $X$ be a connected smooth projective curve of genus $g$ over a finite field $k$, let $\bk$ be an algebraic closure of $k$, let $\omega$ be a non zero global meromorphic differential $1$-form on $X$ and let $\F$ be a constructible $\overline{\mathbb{Q}}_{\ell}$-sheaf on $X$ of generic rank $\rk(\F)$. The $\ell$-adic character $\varepsilon_{\bk}(X, \F)$ of $\Gal(\bk/k)$ associated to the $1$-dimensional representation $\det(R\Gamma(X_{\bk},\F))^{-1}$ (cf. \ref{chap3globaleps}) admits the following decomposition:
$$
\varepsilon_{\bk}(X, \F) = \chi_{\cyc}^{N(g-1) \rk(\F)} \prod_{x \in |X|} \delta_{x/k}^{a(X_{(x)},\F_{|X_{(x)}})} \Ver_{x/k} \left( \varepsilon_{\overline{x}}(X_{(x)},\F_{|X_{(x)}},\omega_{|X_{(x)}}) \right),
$$
where $N$ is the number of connected components of $X_{\bk}$, where $|X|$ is the set of closed points of $X$, where $X_{(x)}$ is the henselization of $X$ at a closed point $x$, and where $\chi_{\cyc}$ is the $\ell$-adic cyclotomic character of $k$. All but finitely many terms in this product are identically equal to $1$.
\end{teo}

The formulation of the product formula in Theorem \ref{chap3fr} differs from Laumon's (\cite{La87}, Th. 3.2.1.1), but yields an equivalent formula. Indeed, if $k$ is of cardinality $q$, then evaluating the product formula in Theorem \ref{chap3fr} at the geometric Frobenius $\Frob_k$ yields that the determinant $\det \left( \Frob_k \ | \ R\Gamma(X_{\bk},\F) \right)^{-1}$ is equal to
$$
q^{N(1-g(X))\rk(\F)} \prod_{x \in |X|} (-1)^{([k(x):k]-1)a(X_{(x)},\F_{|X_{(x)}})} \varepsilon_{\overline{x}}(X_{(x)},\F_{|X_{(x)}},\omega_{|X_{(x)}})(\Frob_x),
$$
and $\sum_{x \in X} [k(x):k] a(X_{(x)},\F_{|X_{(x)}})$ has the same parity as the Euler characteristic $-\chi(X_{\bk},\F)$ by the Grothendieck-Ogg-Shafarevich formula, hence the product formula asserts that the quantity
$$
 q^{N(1-g(X))\rk(\F)} \prod_{x \in |X|} (-1)^{a(X_{(x)},\F_{|X_{(x)}})} \varepsilon_{\overline{x}}(X_{(x)},\F_{|X_{(x)}},\omega_{|X_{(x)}})(\Frob_x),
$$
coincides with the determinant $\det \left( - \Frob_k \ | \ R\Gamma(X_{\bk},\F) \right)^{-1}$, as in Laumon's formulation (\cite{La87}, Th. 3.2.1.1).

For $\ell$-adic sheaves with finite geometric monodromy, the product formula in Theorem \ref{chap3fr} reduces by Brauer's induction theorem to the rank $1$ case, and the latter follows from Tate's thesis, cf. (\cite{La87}, 3.2.1.7). A geometric proof of the product formula in rank $1$ was given by Deligne in his $1974$ letter to Serre (\cite{bloch}, Appendix), using geometric class field theory. Deligne's proof in the rank $1$ case, which we review in Section \ref{chap3productsection1}, extends to the case of an arbitrary perfect base field $k$, and constitutes an important ingredient in the proof of the main theorem \ref{chap3teo2} below.
%
%

\subsection{\label{chap310.1}} Let us consider a quadruple $(T,\F,\omega,\overline{s})$, where $T$ is an equicharacteristic henselian trait, with perfect residue field of positive characteristic $p$, equipped with an algebraic closure $\overline{s}$ of its closed point $s$, where $\omega$ is a non zero meromorphic $1$-form on $T$ (cf. \ref{chap32.12}) and where $\F$ is a constructible \'etale $ \overline{\mathbb{Q}}_{\ell}$-sheaf on $T$.

Let us assume for simplicity that $\F$ is irreducible, with vanishing fiber at $s$, and that $k(s)$ is the perfection of a finitely generated extension of $\mathbb{F}_p$, so that $\F$ has finite geometric monodromy by Grothendieck's local monodromy theorem. If one wishes to construct an $\varepsilon$-factor $\varepsilon_{\overline{s}}(T,\F,\omega)$ by using Brauer's theorem from finite group theory, in order to reduce through additivity and induction to the rank $1$ case, we need $\F$ to have finite monodromy, rather than merely having finite geometric monodromy. When $k$ is finite, the Galois group $\Gal(\overline{s}/s)$ is procyclic, hence some twist of $\F$ by a geometrically constant $ \overline{\mathbb{Q}}_{\ell}$-sheaf of rank $1$ has finite monodromy, cf. (\cite{De73}, 4.10), and this allows Deligne and Langlands to reduce to the finite monodromy case. 

In the general case, the Galois group $\Gal(\overline{s}/s)$ is not procyclic, nor abelian, and twisting by geometrically constant $ \overline{\mathbb{Q}}_{\ell}$-sheaves of rank $1$ is not enough to reduce to the finite monodromy case from the finite geometric monodromy case. However, it is possible to allow for such a reduction by considering more general twists. More precisely, we can reduce to the finite monodromy case at the following costs (cf. \ref{chap3dec}, \ref{chap3irrfinite}):
\begin{itemize}
\item[(a)] considering $ \overline{\mathbb{Q}}_{\ell}$-sheaves on $T$ twisted by a $ \overline{\mathbb{Q}}_{\ell}^{\times}$-valued $2$-cocycle on $\Gal(\overline{s}/s)$, rather than merely $ \overline{\mathbb{Q}}_{\ell}$-sheaves,
\item[(b)] allowing twists by higher rank (twisted) geometrically constant sheaves, rather than rank $1$ such sheaves.

\end{itemize}

The notion of twisted sheaf is recalled in \ref{chap31.3}. Let us simply describe here the corresponding notion of twisted $ \overline{\mathbb{Q}}_{\ell}$-representation. If $\eta$ is the generic point of $T$ and if $\overline{\eta}$ is a separable closure of $\eta_{\overline{s}}$, with Galois group $\Gal(\overline{\eta}/\eta)$ endowed with the natural homomorphism $r : \Gal(\overline{\eta}/\eta) \rightarrow \Gal(\overline{s}/s)$, then a $ \overline{\mathbb{Q}}_{\ell}$-representation of $\Gal(\overline{\eta}/\eta)$ twisted by a $ \overline{\mathbb{Q}}_{\ell}^{\times}$-valued $2$-cocycle $\mu$ on $\Gal(\overline{s}/s)$, is a continuous map
$$
\rho : \Gal(\overline{\eta}/\eta) \rightarrow \GL(V),
$$
where $V$ is a finite dimensional vector space over some finite extension of $\mathbb{Q}_{\ell}$ contained in $ \overline{\mathbb{Q}}_{\ell}$, which satisfies
$$
\rho(g) \rho(h) = \mu(r(g),r(h)) \rho(gh),
$$
for all $g,h$ in $\Gal(\overline{\eta}/\eta)$. When $\mu=1$ is the trivial cocycle, a twisted $ \overline{\mathbb{Q}}_{\ell}$-representation of $\Gal(\overline{\eta}/\eta) $ is simply a $ \overline{\mathbb{Q}}_{\ell}$-Galois representation over $\eta$. The preliminary section \ref{chap3preliminaries} is devoted to a more thorough discussion of twisted representations. 

\subsection{\label{chap310.2}} Let $k$ be a perfect field of positive characteristic $p$, with algebraic closure $\bk$, and let $\ell, \psi$ be as in \ref{chap310.0}. Let $\Lambda$ be either $ \overline{\mathbb{F}}_{\ell}$, $ \overline{\mathbb{Q}}_{\ell}$ or the integral closure $ \overline{\mathbb{Z}}_{\ell}$ of $\mathbb{Z}_{\ell}$ in $ \overline{\mathbb{Q}}_{\ell}$.

 Let us consider quadruples $(T,\F,\omega,\overline{s})$ where $T$ is a henselian trait over $k$, whose closed point $s$ is finite over $k$, equipped with a $k$-morphism $\overline{s} : \Spec(\bk) \rightarrow T$, where $\omega$ is a non zero meromorphic $1$-form on $T$ (cf. \ref{chap32.12}) and where $\F$ is a constructible \'etale $\Lambda$-sheaf on $T$ twisted (cf. \ref{chap31.3}) by some \textit{unitary} $2$-cocycle on $\Gal(\overline{s}/s)$, i.e. a $2$-cocycle which is continuous with values in a finite subgroup of $\Lambda^{\times}$ (cf. \ref{chap3unitary}). 


A \textit{theory of twisted $\ell$-adic local $\varepsilon$-factors over $k$, with respect to $\psi$,} is a rule $\varepsilon$ which assigns to any such quadruple $(T,\F,\omega,\overline{s})$ a map $\varepsilon_{\overline{s}}(T,\F,\omega)$ from the Galois group $\Gal(\overline{s}/s)$ to $ \Lambda^{\times}$, and which satisfies the following axioms:
\begin{itemize}

\item[$(i)$] the map $\varepsilon_{\overline{s}}(T,\F,\omega)$ depends only on the isomorphism class of the quadruple $(T,\F,\omega,\overline{s})$;
\item[$(ii)$] there exists a sub-$\mathbb{Z}_{\ell}$-algebra $\Lambda_0$ of $\Lambda$ of finite type, depending on $\F$, such that $\varepsilon_{\overline{s}}(T,\F,\omega)$ is a continuous map from $\Gal(\overline{s}/s)$ to $\Lambda_0^{\times}$;
\item[$(iii)$](cf. \ref{chap3multiplicativity}) for any exact sequence
$$
0 \rightarrow \F' \rightarrow \F \rightarrow \F'' \rightarrow 0,
$$
of constructible \'etale $\Lambda$-sheaves on $T$ twisted by the same unitary $2$-cocycle, we have 
$$
\varepsilon_{\overline{s}}(T,\F,\omega) = \varepsilon_{\overline{s}}(T,\F',\omega) \varepsilon_{\overline{s}}(T,\F'',\omega);
$$
\item[$(iv)$](cf. \ref{chap3onapoint}) if $\F$ supported on the closed point of $T$, then the value of $\varepsilon_{\overline{s}}(T,\F,\omega)$ at an element $g$ of $\Gal(\overline{s}/s)$ is given by
$$
\varepsilon_{\overline{s}}(T,\F,\omega)(g) = \det \left( g \ | \ \F_{\overline{s}} \right)^{-1};
$$
\item[$(v)$](cf. \ref{chap3inductionformula}) for each finite generically \'etale extension $f : T' \rightarrow T$ of henselian traits $T'$ and $T$ over $k$, there exists a homomorphism $\lambda_f(\omega)$ from the Galois group $\Gal(\overline{s}/s)$ to $ \overline{\mathbb{Z}}_{\ell}^{\times}$ such that
$$
\varepsilon_{\overline{s}}(T,f_*\F,\omega) = \lambda_f(\omega)^{\rk(\F)} \delta_{s'/s}^{a(T',\F,f^* \omega)} \Ver_{s'/s} \left( \varepsilon_{\overline{s}'}(T',\F,f^*\omega) \right), 
$$
for any constructible \'etale $\Lambda$-sheaf $\F$ on $T'$, of generic rank $\rk(\F)$, twisted by some unitary $2$-cocycle on $\Gal(\overline{s}/s)$, where the signature $\delta_{s'/s}$ and the \textit{verlagerung} or \textit{transfer} $ \Ver_{s'/s}$ are defined in \ref{chap31.4.6};
\item[$(vi)$] if the fiber of $\F$ at the closed point $s$ of $T$ vanishes and if $\F$ is generically of rank $1$, with Swan conductor $\nu - 1$, then the value of $\varepsilon_{\overline{s}}(T,\F,\omega)$ at an element $g$ of $\Gal(\overline{s}/s)$ is prescribed as follows:
$$
\varepsilon_{\overline{s}}(T,\F,\omega)(g) = \det \left( g \ | \ H_c^{\nu} \left(\Pic^{\nu + v(\omega)}(T,\nu s)_{\overline{s}},\chi_{\F} \otimes \Lc_{\psi} \lbrace \Res_{\omega} \rbrace(-v(\omega)) \right) \right),
$$
where $v(\omega)$ is the valuation of $\omega$ (cf. \ref{chap32.12}), where $ \Pic^{\nu + v(\omega)}(T,\nu s)$ is the component of degree $\nu + v(\omega)$ of the local Picard group (cf. \ref{chap3localpic}), where $\chi_{\F}$ is the multiplicative local system on the group $\Pic(T,\nu s)$ naturally associated to $\F$ by twisted local geometric class field theory (cf. \ref{chap3lgcft2twisted}), and where $\Lc_{\psi} \lbrace \Res_{\omega} \rbrace$ is the Artin-Schreier local system associated to the residue morphism $\Res_{\omega}$, cf \ref{chap32.13} for details;
\item[$(vii)$](cf. \ref{chap3unramtwist}) if $\G$ is a geometrically constant $ \Lambda$-sheaf on $T$, twisted by some unitary $2$-cocycle on $\Gal(\overline{s}/s)$ (possibly different from the $2$-cocycle by which $\F$ is twisted), then we have 
$$
\varepsilon_{\overline{s}}(T,\F \otimes \G,\omega) = \det(\G_{\overline{s}})^{a(T,\F,\omega)} \varepsilon_{\overline{s}}(T,\F,\omega)^{\rk(\G)},
$$
where the conductor $a(T,\F,\omega)$ is defined in \ref{chap32.11}.
\item[$(viii)$] if $\F$ is a twisted $ \overline{\mathbb{Z}}_{\ell}$-sheaf, then the diagram
\begin{center}
 \begin{tikzpicture}[scale=1]

\node (A) at (0,2) {$\Gal(\overline{s}/s)$};
\node (B) at (5,4) {$\overline{\mathbb{Q}}_{\ell}^{\times}$};
\node (C) at (5,2) {$\overline{\mathbb{Z}}_{\ell}^{\times}$};
\node (D) at (5,0) {$\overline{\mathbb{F}}_{\ell}^{\times}$};

\path[->,font=\scriptsize]
(A) edge[bend left] node[below right]{$\varepsilon_{\overline{s}}(T,\F \otimes_{\overline{\mathbb{Z}}_{\ell}} \overline{\mathbb{Q}}_{\ell},\omega)$} (B)
(C) edge  (B)
(A) edge[bend right] node[above right]{$\varepsilon_{\overline{s}}(T,\F \otimes_{\overline{\mathbb{Z}}_{\ell}} \overline{\mathbb{F}}_{\ell},\omega)$} (D)
(A) edge node[above]{$\varepsilon_{\overline{s}}(T,\F,\omega)$} (C)
(C) edge  (D);
\end{tikzpicture} 
\end{center}
is commutative.
\end{itemize}

Our main result can then be stated as follows:

\begin{teo}[cf. \ref{chap35.8}]\label{chap3teo2} Let $k$ be a perfect field of positive characteristic $p > 0$. Then for any prime number $\ell$ distinct from $p$ and any non trivial homomorphism $\psi : \mathbb{F}_p \rightarrow \overline{\mathbb{Q}}_{\ell}^{\times}$, there exists a unique theory of twisted $\ell$-adic local $\varepsilon$-factors over $k$, with respect to $\psi$. Moreover, we have the following properties:
\begin{itemize}
\item[$(viii)$] the map $\varepsilon_{\overline{s}}(T,\F,\omega)$ does not depend on the subfield $k$ of $k(s)$.
\item[$(ix)$](cf. \ref{chap3sanity5}) for any quadruple $(T,\F,\omega, \overline{s})$ over $k$, where $\F$ is twisted by a unitary $2$-cocyle $\mu$ on $\Gal(\overline{s}/s)$, we have
$$
\varepsilon_{\overline{s}}(T,\F,\omega)(g) \varepsilon_{\overline{s}}(T,\F,\omega)(h) = \mu(g,h)^{a(T,\F,\omega)}\varepsilon_{\overline{s}}(T,\F,\omega)(gh)
$$
for any elements $g,h$ of $\Gal(\overline{s}/s)$, where $s$ is the closed point of $T$ and where $a(T,\F,\omega)$ is the conductor defined in \ref{chap32.11}. In particular the $2$-cocycle $\mu^{a(T,\F,\omega)}$ is a coboundary (cf. \ref{chap30.11}). 
\end{itemize}
\end{teo}

In the untwisted case, the property $(vii)$ is a consequence of $(i)-(vi)$, since one can assume that the twist $\G$ is of rank $1$, and one can then use Brauer's induction theorem together with $(i)-(v)$ to reduce to the case where $\F$ is also of rank $1$, so that $(vii)$ then follows from $(vi)$, cf. (\cite{La87}, 3.1.5.6).

When $k$ is finite, then for any finite extension $k(s)$ of $k$ contained in $\bk$, any unitary $2$-cocycle on $\Gal(\bk/k(s))$ is a coboundary, and thus the theory of twisted $\ell$-adic local $\varepsilon$-factors over $k$ is not more general than the classical theory of Deligne and Langlands. Actually, we have:

\begin{teo}[cf. \ref{chap35.9}]\label{chap3teo4} Let $(T,\F,\omega,\overline{s})$ be a quadruple as in \ref{chap310.2} over a finite field $k$, where $\F$ is untwisted, i.e. twisted by the trivial cocycle. Then the quantity
$$
(-1)^{a(T,\F)} \varepsilon_{\overline{s}}(T,\F,\omega)(\Frob_s),
$$
$\Frob_s$ is the geometric Frobenius in $\Gal(\overline{s}/s)$, where $s$ is the closed point of $T$ and where $a(T,\F)$ is the conductor of $(T,\F)$ (cf. \ref{chap32.11}), coincides with the classical local $\varepsilon$-factor, normalized as in \cite[Th. 3.1.5.4]{La87}.
\end{teo}

This result will be deduced from the normalization $(vi)$ of geometric local $\varepsilon$-factors and from the Grothendieck-Lefschetz trace formula (cf. \ref{chap3compafini}).

As in the case of a finite base field (cf. \ref{chap310.1}), we have a product formula:

\begin{teo}[cf. \ref{chap3productform5}]\label{chap3teo3}  Let $\Lambda$ be either $ \overline{\mathbb{F}}_{\ell}$ or $ \overline{\mathbb{Q}}_{\ell}$. Let $X$ be a connected smooth projective curve of genus $g(X)$ over a perfect field $k$, let $\omega$ be a non zero global meromorphic differential $1$-form on $X$ and let $\F$ be a constructible $\Lambda$-sheaf on $X$ of generic rank $\rk(\F)$, twisted by some unitary $2$-cocycle on $\Gal(\bk/k)$ (cf. \ref{chap31.3}).
Then the trace function $\varepsilon_{\bk}(X, \F)$ on $\Gal(\bk/k)$ associated to the twisted $1$-dimensional representation $\det(R\Gamma(X_{\bk},\F))^{-1}$ (cf. \ref{chap3globaleps}) admits the following decomposition:
$$
\varepsilon_{\bk}(X, \F) = \chi_{\cyc}^{N(g(X)-1) \rk(\F)} \prod_{x \in |X|} \delta_{x/k}^{a(X_{(x)},\F_{|X_{(x)}})} \Ver_{x/k} \left( \varepsilon_{\overline{x}}(X_{(x)},\F_{|X_{(x)}},\omega_{|X_{(x)}}) \right),
$$
where $N$ is the number of connected components of $X_{\bk}$, where $|X|$ is the set of closed points of $X$ and $\chi_{\cyc}$ is the $\ell$-adic cyclotomic character of $k$. All but finitely many terms in this product are identically equal to $1$.
\end{teo}

We first prove Theorem \ref{chap3teo3} in the case of (twisted) $\Lambda$-sheaves with finite geometric monodromy, cf. \ref{chap3productform2}, and we then prove the general case in Section \ref{chap3productsection3} by using Laumon's $\ell$-adic stationary phase method. An important ingredient of the proof is the following extension of Laumon's formula (\ref{chap3laumonform}):

\begin{teo}[cf. \ref{chap3xteolaumon}]\label{chap3teo5} Let $T$ be an henselian trait with closed point $s$, such that $k(s)$ is a perfect field of positive characteristic $p$, and let $\F$ be a constructible $\Lambda$-sheaf on $T$ with vanishing fiber at $s$, twisted by some unitary $2$-cocycle on $\Gal(\bk/k)$ (cf. \ref{chap31.3}). Let $\pi$ be a uniformizer on $T$, and let $\chi_{\det(\mathrm{F}_{\pi}^{(0,\infty')}(\F))}$ be the multiplicative $\Lambda$-local system associated to $\det(\mathrm{F}_{\pi}^{(0,\infty')}(\F))$ by geometric class field theory (cf. \ref{chap3lgcft2twisted}), where $\mathrm{F}_{\pi}^{(0,\infty')}(\F)$ is Laumon's local Fourier transform, cf. (\cite{La87}, 2.4.1). Then the trace map of the stalk of $\chi_{\det(\mathrm{F}_{\pi}^{(0,\infty')}(\F))}$ at $\pi^{-1}$ coincides with $\varepsilon_{\overline{k}}(T, \F, d \pi)$.
\end{teo}

Laumon's proof of this result when $k$ is finite starts with a reduction to the tamely ramified case (\cite{La87}, 3.5.3.1), and then resort to a computation in the latter case (\cite{La87}, 2.5.3.1). Instead of adapting Laumon's proof to the general case, we choose to avoid these steps in our treatment of Theorem \ref{chap3teo5} : we give a direct proof by using the $\ell$-adic stationary phase method (cf. \ref{chap3xteostat}) and by specializing Theorem \ref{chap3teo2} to the case where the base field $k$ is (the perfection of) a henselian discretely valued field of equicharacteristic $p$.

%

\subsection{\label{chap310.3}} Let us briefly describe our definition of twisted $\ell$-adic local $\varepsilon$-factors over a perfect field $k$ of positive characteristic $p$ (cf. \ref{chap310.2}). Let $(T,\F,\omega,\overline{s})$ be a quadruple over $k$ as in \ref{chap310.2}, and let $s$ be the closed point of $T$. We fix a uniformizer $\pi$ of $\Ow_T$, and denote by $\pi$ as well the morphism
$$
\pi : T \rightarrow \mathbb{A}^1_s,
$$
corresponding to the unique morphism $k(s)[t] \rightarrow \Ow_T$ of $k(s)$-algebras which sends $t$ to $\pi$. The theory of Gabber-Katz extensions, originating from \cite{katz} and reviewed in Section \ref{chap3GKsection}, ensures the existence of a (twisted) $\Lambda$-sheaf $\pi_{\diamondsuit} \F$ on $\mathbb{A}^1_s$, unique up to isomorphism, such that:
\begin{enumerate}
\item the pullback $\pi^{-1} \pi_{\diamondsuit} \F$ is isomorphic to $\F$;
\item the $\Lambda$-sheaf $\pi_{\diamondsuit} \F$ is tamely ramified at infinity;
\item the restriction of $\pi_{\diamondsuit} \F$ to $\mathbb{G}_{m,s}$ is a local system whose geometric monodromy group has a unique $p$-Sylow.
\end{enumerate}
We then simply define
$$
\varepsilon_{\overline{s}}(T,\F,d \pi) =\det \left( R \Gamma_c(\mathbb{A}^1_{\overline{s}}, \pi_{\diamondsuit }\F \otimes \Lc_{\psi}^{-1}) \right)^{-1},
$$
where $\Lc_{\psi}$ is the Artin-Schreier sheaf on the affine line associated to $\psi$. Using geometric class field theory, we then define $\varepsilon_{\overline{s}}(T,\F,\omega)$ for arbitrary meromorphic $1$-forms $\omega$ on $T$, cf. \ref{chap3localepsfact}.

We then show in Section \ref{chap3gfar}, using a variant of Brauer's induction theorem (cf. \ref{chap3brauer2}), that the resulting local $\varepsilon$-factor is independent of the choice of $\pi$ (cf. \ref{chap3indepunif}) and that it satisfies the properties $(i)-(ix)$ listed in \ref{chap310.2} and in Theorem \ref{chap3teo2}. The most notable of these properties is the induction formula $(v)$, which is proved using generalized Gabber-Katz extensions (cf. \ref{chap3GK3}) and the product formula (cf. \ref{chap3teo2}) in generic rank $1$, proved by Deligne in his $1974$ letter to Serre, the latter being published as an appendix in \cite{bloch} and reviewed in Section \ref{chap3productsection1}.

\subsection{\label{chap310.4}} We now describe the organization of this paper. Section \ref{chap3preliminaries} contains preliminary definitions and results on representations of groups twisted by a $2$-cocycle. It notably includes an extension of Brauer's induction theorem to this context, namely Theorem \ref{chap3brauer2}, and a useful decomposition of a twisted representation according to its restriction to a finite normal subgroup in Proposition \ref{chap3dec}.

Section \ref{chap3twistedsection} is devoted to basic definitions and results regarding $\ell$-adic sheaves and their twisted counterparts. 

We review in Section \ref{chap3GKsection} the theory of Gabber-Katz extensions, following the exposition by Katz in \cite{katz}. We provide mild generalizations of the results found in the latter article, namely an extension to twisted $\ell$-adic sheaves on arbitrary Gabber-Katz curves.

Section \ref{chap3gcftsection} is devoted to geometric class field theory, in both of its global and local incarnations. Since this topic is of independent interest, we choose to present more material than what is strictly necessary in order to prove the main results of this text. We discuss in particular the relations between different formulations of geometric local class field theory, namely those of Serre \cite{JPS}, of Contou-Carr\`ere \cite{CC} and Suzuki \cite{TS}, or of Gaitsgory. We also prove local-global compatibility in geometric class field theory (cf. \ref{chap32.6}), as well as functoriality with respect to the norm homomorphism (cf. \ref{chap32.19}). 

In Section \ref{chap3extensions}, we perform a series of computations aiming at describing multiplicative local systems, namely the geometric analog of characters of abelian groups, on certain groups schemes, such as the additive group $\Ga_a$ or the group of Witt vectors of length $2$ over $\mathbb{F}_2$. All of these computation can be considered as being part of the proof of the main proposition \ref{chap3rankcohom} in Section \ref{chap3lgf1}, which describes the cohomology groups appearing in our definition \ref{chap3localfact1} of geometric local $\varepsilon$-factors in generic rank $1$.

In Section \ref{chap3productsection1}, we review Deligne's $1974$ letter to Serre on $\varepsilon$-factors, where the product formula \ref{chap3teo3} is proved in generic rank $1$ by using geometric class field theory. We also provide a mild generalization to the context of twisted $\ell$-adic sheaves.

Section \ref{chap3gfar} is devoted to the proofs of the main results of this text, namely that the twisted $\ell$-adic $\varepsilon$-factors defined with Gabber-Katz extensions as in \ref{chap310.3} are independent of the choice of uniformizer and satisfy the properties $(i)-(ix)$ listed in \ref{chap310.2} and in Theorem \ref{chap3teo2}. Our main tools are the reduction to the rank $1$ case allowed by the results of Section \ref{chap3preliminaries}, and the product formula in generic rank $1$ from Section \ref{chap3productsection1}. We also prove Theorem \ref{chap3teo4} in this section.

Finally, we prove the product formula for (twisted) $\ell$-adic sheaves of arbitrary rank, first under various finiteness hypotheses in section \ref{chap3productsection2}, by using the results of Section \ref{chap3preliminaries} to reduce to the rank $1$ case handled in Section \ref{chap3productsection1}, and then in the general case in Section \ref{chap3productsection3}, by following closely Laumon's proof in the finite field case.

\subsection*{Acknowledgements} This work is part of the author's PhD dissertation and it was prepared at the Institut des Hautes \'Etudes Scientifiques and the \'Ecole Normale Sup\'erieure while the author benefited from their hospitality and support. The author is indebted to Ahmed Abbes for his numerous suggestions and remarks. Further thanks go to Fabrice Orgogozo and Fu Lei for their corrections and remarks, to Dennis Gaitsgory for a fruitful discussion on geometric local class field theory, and to Ofer Gabber for suggesting that the ``potentially unipotent'' condition appearing in an earlier version of this text could be removed. 

\subsection{Conventions and notation}\label{chap3conv} We fix a perfect field $k$ of positive characteristic $p$, and we denote by $\bk$ a fixed algebraic closure of $k$. We denote by $G_k = \Gal(\bk/k)$ the Galois group of the extension $\bk/k$. For any $k$-scheme $X$, and for any $k$-algebra $k'$, we denote by $X_{k'}$ the fiber product of $X$ and $\Spec(k')$ over $\Spec(k)$. The group $G_k$ acts on the left on $\bk$, and thus acts on the right on $X_{\bk}$.

We fix as well a prime number $\ell$ different from $p$, and we denote by $C$ an algebraic closure of the field $\mathbb{Q}_{\ell}$ of $\ell$-adic numbers, endowed with the topology induced by the $\ell$-adic valuation. We denote by $\Z_{\ell}(1)$ the invertible $\Z_{\ell}$-module consisting of sequences $(\zeta_n)_{n \geq 0}$ of elements of $\bk$ such that $\zeta_0 = 1$ and $\zeta_{n+1}^{\ell} = \zeta_n$ for each $n$, endowed with the natural action of $G_k$. For each integer $\nu$, we denote by $\Z_{\ell}(\nu)$ the invertible $\Z_{\ell}$-module $\Z_{\ell}(1)^{\otimes \nu}$. More generally, for any $\ell$-adic sheaf $\F$ on a $k$-scheme (cf. \ref{chap31.0.4}) we denote by $\F(\nu)$ the tensor product of $\F$ and $\Z_{\ell}(\nu)$ over $\Z_{\ell}$. We also denote by
\begin{align*}
\chi_{\cyc} : G_k &\rightarrow \Z_{\ell}^{\times} \\
g &\rightarrow \Tr(g \ | \ \Z_{\ell}(1) ),
\end{align*}
the character associated to the $\ell$-adic representation $\Z_{\ell}(1)$ of $G_k$.

\section{Preliminaries on representations of twisted groups \label{chap3preliminaries}}


\subsection{\label{chap30.0.0.0.1}} Let $G$ be a profinite topological group, and let $f : G \rightarrow X$ be a continuous map onto a finite set $X$ endowed with the discrete topology. The open normal subgroups of $G$ form a basis of open neighbourhoods at the unit element of $G$. Hence, for each element $g$ of $G$, there exists an open normal subgroup $I_g$ such that the coset $g I_g$ is contained in the open subset $f^{-1}(f(g))$. Since $G$ is compact, there exists a finite family $(g_j)_{j \in J}$ of elements of $G$ such that the open subsets $(g_j I_{g_j})_{j \in J}$ form a cover of $G$. Thus, if $I$ is the intersection of the open normal subgroups $(I_{g_j})_{j \in J}$, then $I$ is itself an open normal subgroup of $G$ and $f$ is both left and right $I$-invariant.

\subsection{\label{chap30.0.0.1}} An \textit{admissible $\ell$-adic ring} is a commutative topological ring which is isomorphic to one of the following:
\begin{enumerate}
\item a finite local $\Z / \ell^n$-algebra for some integer $n$, endowed with the discrete topology,
\item the ring of integers in a finite extension of $\mathbb{Q}_{\ell}$, endowed with the topology defined by the $\ell$-adic valuation,
\item a finite extension of $\mathbb{Q}_{\ell}$, endowed with the topology defined by the $\ell$-adic valuation.
\end{enumerate}
An \textit{$\ell$-adic coefficient ring} is a commutative topological ring $\Lambda$ such that any finite subset of $\Lambda$ is contained in a subring of $\Lambda$, which is an admissible $\ell$-adic ring for the subspace topology. In particular, any admissible $\ell$-adic ring is an $\ell$-adic coefficient ring as well.

\begin{exemple} The ring $C$ (cf. \ref{chap3conv}) endowed with the topology induced by the $\ell$-adic valuation, is an $\ell$-adic coefficient ring.
\end{exemple}

\begin{rema}\label{chap3filtered} Any admissible $\ell$-adic ring is a finitely presented $\Z_{\ell}$-agebra. In particular, the set of admissible $\ell$-adic subrings of an $\ell$-adic coefficient ring is filtered when ordered by inclusion.
\end{rema}

\subsection{\label{chap30.0.0.0}} Let $T$ be a topological space and let $\Lambda$ be an $\ell$-adic coefficient ring (cf. \ref{chap30.0.0.1}). A map $f : T \rightarrow \Lambda$ is said to be \textit{$\Lambda$-admissible} if it is continuous and if its image is contained in an admissible $\ell$-adic subring of $\Lambda$. 

Similarly, if $V$ is a free $\Lambda$-module of finite rank, a map $f : T \rightarrow \Aut_{\Lambda}(V)$ (resp. $f : T \rightarrow \Aut_{\Lambda}(V)/\Lambda^{\times}$) is said to be \textit{$\Lambda$-admissible} if it is continuous and if there is an admissible $\ell$-adic subring $\Lambda_0$ of $\Lambda$ and a $\Lambda_0$-form $V_0$ of $V$ such that $f$ factors through $\Aut_{\Lambda_0}(V_0)$ (resp. $\Aut_{\Lambda_0}(V_0)/\Lambda_0^{\times}$).
%

\subsection{\label{chap30.11}} Let $G$ be a topological group, and let $\Lambda$ be an $\ell$-adic coefficient ring (cf. \ref{chap30.0.0.1}). For each integer $j$, let $C^j(G,\Lambda^{\times})$ be the group of $\Lambda$-admissible maps from $G^j$ to $\Lambda^{\times}$. We define a complex
$$
C^1(G,\Lambda^{\times}) \xrightarrow[]{d^1} C^2(G,\Lambda^{\times}) \xrightarrow[]{d^2} C^3(G,\Lambda^{\times}),
$$
as follows: if $\lambda$ is an element of $C^1(G,\Lambda^{\times})$, we set
\begin{align}
\begin{split}
\label{chap3eq21}
d^1(\lambda) : G^2 &\rightarrow \Lambda^{\times} \\
(x,y) &\rightarrow \lambda(x)\lambda(y)\lambda(xy)^{-1},
\end{split}
\end{align}
which is indeed $\Lambda$-admissible, and if $\mu$ is an element of $C^2(G,\Lambda^{\times})$, we set
\begin{align}
\begin{split}
\label{chap3eq22}
d^2(\mu) : G^3 &\rightarrow \Lambda^{\times} \\
(x,y,z) &\rightarrow \mu(x,y)\mu(xy,z)\mu(x,yz)^{-1} \mu(y,z)^{-1},
\end{split}
\end{align}
which is $\Lambda$-admissible as well. If $\lambda$ is an element of $C^1(G,\Lambda^{\times})$, then we have
$$
d^1(\lambda)(x,y)d^1(\lambda)(xy,z) = \lambda(x) \lambda(y)\lambda(z)\lambda(xyz)^{-1} = d^1(\lambda)(x,yz) d^1(\lambda)(y,z),
$$
and thus $d^2 \circ d^1$ vanishes.

\begin{defi}\label{chap32cocyc} An admissible \textit{$2$-cocycle }(resp. \textit{$2$-boundary}) on $G$ with values in $\Lambda^{\times}$ is an element of the kernel of $d^2$ (resp. of the image of $d^1$). The \textit{second admissible cohomology group} of $G$ with coefficients in $\Lambda^{\times}$, denoted $H_{\adm}^2(G,\Lambda^{\times})$, is the quotient of the group of admissible $2$-cocycles on $G$ with values in $\Lambda^{\times}$, by the subgroup of admissible $2$-boundaries.
\end{defi}

We also have
$$
H_{\adm}^2(G,\Lambda^{\times}) = \mathrm{colim}_{\Lambda_0 } H_{\adm}^2(G,\Lambda_0^{\times}),
$$
where $\Lambda_0$ runs over the filtered set of admissible $\ell$-adic subrings of $\Lambda$ (cf. \ref{chap3filtered}), and the group $H_{\adm}^2(G,\Lambda_0^{\times})$ coincides with the second continuous cohomology group of $G$ with coefficients in $\Lambda_0^{\times}$. 

\begin{rema}\label{chap3finitecomp} If $G$ is finite, then any map from $G^j$ to $\Lambda^{\times}$ is $\Lambda$-admissible. Thus the group $H_{\adm}^2(G,\Lambda^{\times})$ coincides with the second cohomology group $H^2(G,\Lambda^{\times})$ of $G$ with coefficients in $\Lambda$.
\end{rema}

\subsection{\label{chap30.0}} Let $\Lambda$ be an $\ell$-adic coefficient ring (cf. \ref{chap30.0.0.1}). A \textit{$\Lambda$-admissible multiplier} on a topological group $G$ is an admissible $2$-cocycle $\mu$ with values in $\Lambda^{\times}$ (cf. \ref{chap32cocyc}), such that $\mu(1,1) = 1$. In particular, a multiplier $\mu$ satisfies the cocycle relation
\begin{align}
\label{chap3cocy} 
\mu(x,y)\mu(xy,z) = \mu(x,yz) \mu(y,z),
\end{align}
for all $x,y,z$ in $G$. By specializing this relation to $x = y = 1$, we obtain $\mu(1,z) = 1$ for any $z$ in $G$. Likewise, we have $\mu(x,1) = 1$ for any $x$ in $G$.

\begin{defi}\label{chap3unitary} A $\Lambda$-admissible multiplier $\mu$ on a topological group $G$ is said to be \textit{unitary} if there is an integer $r \geq 1$ such that $\mu^r = 1$.
\end{defi}

Since the group of $r$-th roots of unity in $\Lambda$ is a discrete subgroup of $\Lambda^{\times}$, any unitary multiplier on a topological group $G$ must be locally constant. If moreover $G$ is profinite, then any $\Lambda$-admissible multiplier on $G$ must be left and right $I$-invariant for some open normal subgroup $I$ of $G$ (cf. \ref{chap30.0.0.0.1}).

\begin{defi}\label{chap30.1} A \textit{$\Lambda$-twisted topological group} (resp. a \textit{$\Lambda$-twisted group}) is a pair $(G,\mu)$, where $G$ is a topological group (resp. a discrete group) and $\mu$ is a $\Lambda$-admissible multiplier on $G$. A morphism of $\Lambda$-twisted topological groups from $(G,\mu)$ to $(G',\mu')$ is a continuous group homomorphism $f : G \rightarrow G'$ such that $\mu(x,y) = \mu'(f(x),f(y))$ for all $x,y$ in $G$.
\end{defi}

\begin{rema}\label{chap3cobor0} Let $G$ be a topological group, and let $\lambda : G \rightarrow \Lambda^{\times}$ be a $\Lambda$-admissible map (cf. \ref{chap30.0.0.0}) such that $\lambda(1)=1$. Then the $2$-boundary $d^1(\lambda)$ (cf. \ref{chap3eq21}) is a multiplier on $G$. The quotient of the group of $\Lambda$-admissible multipliers on $G$ by the group of $2$-coboundaries $d^1(\lambda)$ such that $\lambda(1)=1$, is isomorphic to the second admissible cohomology group $H_{\adm}^2(G,\Lambda^{\times})$, since any $\Lambda$-admissible $2$-cocycle $\mu$ on $G$ factors as $\mu = c \mu'$, where $\mu'$ is a $\Lambda$-admissible multiplier on $G$ and $c = \mu(1,1)$ is a unit of $\Lambda$, and we have $c = d^1(c)$.

\end{rema}

%

\subsection{\label{chap30.2}} Let $(G,\mu)$ be a $\Lambda$-twisted topological group (cf. \ref{chap30.1}). A \textit{$\Lambda$-admissible representation of $(G,\mu)$} is a pair $(V, \rho)$, where $V$ is a free $\Lambda$-module of finite rank, and $\rho : G \rightarrow \Aut_{\Lambda}(V)$ is a $\Lambda$-admissible map (cf. \ref{chap30.0.0.0}) such that 
$$
\rho(x) \rho(y) = \mu(x,y) \rho(xy),
$$
for all $x,y$ in $G$. Since $\mu(1,1) = 1$, this relation implies $\rho(1) = \id_V$.

If $(V, \rho)$ and $(V',\rho')$ are both $\Lambda$-admissible representations of $(G,\mu)$, we define a morphism from $(V, \rho)$ to $(V',\rho')$ to be a homomorphism $f : V \rightarrow V'$ of $\Lambda$-modules such that $f \circ \rho(x) = \rho'(x) \circ f$ for all $x$ in $G$. 

We will denote by $\Rep_{\Lambda}(G,\mu)$ the category of $\Lambda$-admissible representations of $(G,\mu)$.

\begin{rema} If $(V, \rho)$ is a $\Lambda$-admissible representation of $(G, \mu)$, then the composition of $\rho$ with the projection $\Aut_{\Lambda}(V) \rightarrow \Aut_{\Lambda}(V)/{\Lambda}^{\times}$ to the projective general linear group of $V$ is a genuine group homomorphism, thereby defining a projective representation of $G$. Moreover, any $\Lambda$-admissible projective representation of a discrete group $G$ is obtained in this way from a $\Lambda$-admissible representation of $(G,\mu)$, for some multiplier $\mu$. However, the category $\Rep_{\Lambda}(G,\mu)$ is additive, unlike the category of projective representations of $G$. 
\end{rema}

\begin{rema}\label{chap3cobor} If $\lambda : G \rightarrow \Lambda^{\times}$ is a $\Lambda$-admissible map with $\lambda(1)=1$, then the functor $(V,\rho) \mapsto (V, \lambda \rho)$ is an equivalence of categories from $\Rep_{\Lambda}(G,\mu)$ to $\Rep_{\Lambda}(G,\mu d^1(\lambda))$ (cf. \ref{chap3cobor0}). These categories are thus equivalent, although non canonically, since the isomorphism just constructed depends on $\lambda$. In particular, $\Rep_{\Lambda}(G,\mu)$ depends only on the cohomology class of $\mu$, up to non unique equivalence.
\end{rema}

\begin{prop}\label{chap3restriction} Let $H$ be an open subgroup of finite index in a $\Lambda$-twisted topological group $(G,\mu)$. Let $V$ be a free $\Lambda$-module of finite rank, and let $\rho : G \rightarrow \Aut_{\Lambda}(V)$ be a map such that $\rho(x) \rho(y) = \mu(x,y) \rho(xy)$ for all $x,y$ in $G$. If $(V,\rho_{|H})$ is a $\Lambda$-admissible representation of $(H,\mu_{|H})$, then $(V,\rho)$ is a $\Lambda$-admissible representation of $(G,\mu)$.
\end{prop}

Indeed, there exists an admissible $\ell$-adic subring $\Lambda_0 \subseteq \Lambda$ and a $\Lambda_0$-form $V_0$ of $V$ such that $\rho_{|H}$ factors through $\Aut_{\Lambda_0}(V_0)$ and such that the induced map $\rho_{|H} : H \rightarrow \Aut_{\Lambda_0}(V_0)$ is continuous. Let $(g_i)_{i \in I}$ be a finite family of left $H$-cosets representatives. Up to replacing $\Lambda_0$ with a larger admissible $\ell$-adic subring of $\Lambda$, we can assume (and we do) that each $\rho(g_i)$ belongs to $\Aut_{\Lambda_0}(V_0)$, and that $\mu$ takes its values in $\Lambda_0^{\times}$. The restriction of $\rho$ to the open subset $Hg_i$ is then given by the formula
$$
\rho(hg_i) = \mu(h,g_i)^{-1} \rho(h) \rho(g_i).
$$ 
Thus the restriction $\rho_{|Hg_i}$ take its values in $\Aut_{\Lambda_0}(V_0)$ and is continuous. Therefore $\rho$ takes its values in $\Aut_{\Lambda_0}(V_0)$ and is continuous, and consequently $(V,\rho)$ is a $\Lambda$-admissible representation of $(G,\mu)$.

\subsection{\label{chap30.3.2}} Non zero $\Lambda$-admissible representations of $(G,\mu)$ may not exist for every $\mu$. Indeed, if $\Lambda$ is an algebraically closed field, then the cohomology class associated to $\mu$ (cf. \ref{chap3cobor0}) must have finite order for such a representation to exist, by the following proposition.

\begin{prop}\label{chap3unitary2} Assume that the $\ell$-adic coefficient ring $\Lambda$ is an algebraically closed field. If a $\Lambda$-twisted topological group $(G,\mu)$ admits a $\Lambda$-admissible representation of rank $r \geq 1$, then there exists a $\Lambda$-admissible map $\lambda : G \rightarrow \Lambda^{\times}$ (cf. \ref{chap30.0.0.0}) such that $\lambda (1)=1$, and such that $\mu d^1(\lambda)$ (cf. \ref{chap3cobor0}) is unitary (cf. \ref{chap3unitary}), with values in the group of $r$-th roots of unity in $\Lambda^{\times}$.
\end{prop}

Indeed, let $\Lambda_0 \subset \Lambda$ be an admissible $\ell$-adic subring of $\Lambda$ such that $\mu$ takes its values in $\Lambda_0^{\times}$, and such that there exists a $\Lambda_0$-admissible representation $(V,\rho)$ of $(G,\mu)$, of rank $r \geq 1$. We have
$$
\det(\rho(g)) \det(\rho(h)) = \mu(g,h)^r \det(\rho(gh)),
$$
for all $g,h$ in $G$. Thus, if $\lambda : G \rightarrow \Lambda^{\times}$ is a $\Lambda$-admissible map such that $\lambda (1)=1$ and $\lambda^r = \det \circ \rho$, then the multiplier $\mu d^1(\lambda)$ (cf. \ref{chap30.11}) is unitary, with values in $r$-th roots of unity. It remains to show the existence of such a continuous map $\lambda$. 

If $\Lambda_0$ is finite, then there exists a map $f$ from $\Lambda_0^{\times}$ to $\Lambda^{\times}$ such that $f(w)^r = w$ for any $w$ in $\Lambda_0^{\times}$. We can then choose $\lambda = f \circ \det \circ \rho$.

Otherwise, we can assume that $\Lambda_0$ is a finite extension $E$ of $\mathbb{Q}_{\ell}$, with ring of integers $\Ow_E$.
The subgroup $1 + r \ell^2 \Ow_E$ of $E^{\times}$ is open, and the map
\begin{align*}
1 + r \ell^2 \Ow_E &\rightarrow 1 + \ell \Ow_E \\
1 + x &\rightarrow (1+x)^{\frac{1}{r}} = \sum_{n \geq 0} \frac{r^{-1}(r^{-1}-1)\cdots (r^{-1}-n+1)}{n!} x^n,
\end{align*}
is continuous. Let $H$ be an open subgroup of $G$ such that $\det \circ \rho(H)$ is contained in $1 + r \ell^2 \Ow_E$. Let $(g_c)_{c \in G/H}$ be a set of representatives for the right cosets of $H$ in $G$, such that $g_1 = 1$, and let $(\lambda_c)_{c \in G/H}$ be elements of $\Lambda$ such that $\lambda_c^r = \det(\rho(g_c))$ and $\lambda_1 = 1$. Then, setting $\lambda( h g_c) = \det(\rho(h))^{\frac{1}{r}} \mu(h,g_c)^{-1} \lambda_c$ for $h$ and $c$ in $H$ and $G/H$ respectively yields a $\Lambda$-admissible map from $G$ to $\Lambda^{\times}$ such that $\lambda (1)=1$ and $\lambda^r = \det \circ \rho$, hence the result.

\subsection{\label{chap30.2.1}} Let $f : (G,\mu) \rightarrow (G',\mu')$ be a morphism of $\Lambda$-twisted topological groups (cf. \ref{chap30.1}), and let $(V,\rho)$ be a $\Lambda$-admissible representation of $(G',\mu')$ (cf. \ref{chap30.2}). Then $(V, \rho f)$ is a $\Lambda$-admissible representation of $(G,\mu)$, which we denote by $f^* V$. This yields a functor $f^*$ from $\Rep_{\Lambda}(G',\mu')$ to $\Rep_{\Lambda}(G,\mu)$.

\subsection{\label{chap30.2.2}} Let $(V_1,\rho_1)$ and $(V_2,\rho_2)$ be $\Lambda$-admissible representations of $\Lambda$-twisted topological groups $(G_1,\mu_1)$ and $(G_2,\mu_2)$ respectively. Then the formula
$$
(\mu_1 \otimes \mu_2)((x_1,x_2),(y_1,y_2)) = \mu_1(x_1,y_1) \mu_2(x_2,y_2)
$$
for $x_1,y_1$ in $G_1$ and $x_2,y_2$ in $G_2$, defines a multiplier on the topological group $G_1 \times G_2$. The free $\Lambda$-module of finite rank $V_1 \otimes_{\Lambda} V_2$ is then endowed with a structure of continuous linear representation of $(G_1 \times G_2,\mu_1 \otimes \mu_2 )$ over $\Lambda$, by defining $\rho(x_1,x_2) = \rho_1(x_1) \otimes \rho_2(x_2)$ for $(x_1,x_2)$ in $G_1 \times G_2$.

If $G_1 = G_2$, then the diagonal morphism $G \rightarrow G \times G$ is a morphism of $\Lambda$-twisted topological groups from $(G, \mu_1 \mu_2)$ to $(G_1 \times G_2,\mu_1 \otimes \mu_2 )$. The restriction of the $\Lambda$-admissible representation $V_1 \otimes_{\Lambda} V_2$ of $(G_1 \times G_2,\mu_1 \otimes \mu_2 )$ through this diagonal morphism then defines a $\Lambda$-admissible representation of $(G, \mu_1 \mu_2)$, still denoted by $V_1 \otimes_{\Lambda} V_2$.

\subsection{\label{chap30.3}} Let $(G,\mu)$ be a $\Lambda$-twisted group (cf. \ref{chap30.1}). The \textbf{twisted group algebra} $\Lambda[G,\mu]$ of $(G,\mu)$ over $\Lambda$ is given by a free $\Lambda$-module with basis $([x])_{x \in G}$, endowed with the $\Lambda$-bilinear product defined by
$$
[x] [y] = \mu(x,y) [xy],
$$
for all $x,y$ in $G$. The cocycle relation \ref{chap3cocy} is equivalent to the associativity of this product. Moreover, recall from \ref{chap30.0} that $\mu(x,1) = \mu(1,z) = 1$ for all $x,z$ in $G$, hence $[1]$ is a (left and right) neutral element, and thus $\Lambda[G,\mu]$ is a unital associative $\Lambda$-algebra.

\begin{prop}\label{chap3equiv} Let $(G,\mu)$ be a finite discrete $\Lambda$-twisted group. The functor which sends a left $\Lambda[G,\mu]$-module $V$ which is free of finite rank over $\Lambda$ to the $\Lambda$-admissible representation of $(G,\mu)$ on $V$ defined by the formula $\rho(x)(v) = [x]v$ for $x$ in $G$ and $v$ in $V$ is a $\Lambda$-linear equivalence of categories, from the category of left $\Lambda[G,\mu]$-modules which are free of finite rank over $\Lambda$ to the category $\Lambda$-admissible representations of $G$ over $\Lambda$ (cf. \ref{chap30.2}).
\end{prop} 

In particular, if $G$ is finite, then $\Lambda[G,\mu]$ is itself a non zero $\Lambda$-admissible representation of $(G,\mu)$. By combining this obervation with Proposition \ref{chap3unitary2}, we recover the well-known result that any multiplier on $G$ is cohomologous to a unitary multiplier with values in $|G|$-th roots of unity, and thus that the abelian group $H_{\adm}^2(G,\Lambda^{\times}) = H^2(G,\Lambda^{\times})$ (cf. \ref{chap3finitecomp}) is finite with exponent dividing $|G|$ whenever $\Lambda$ is an algebraically closed field. This observation admits the following reformulation in terms of central extensions:

\begin{prop}\label{chap3hty}Assume that the $\ell$-adic coefficient ring $\Lambda$ is an algebraically closed field. Let $G^*$ be a central extension of a finite group $G$ by $\Lambda^{\times}$. Then there exists a section $\sigma : G \rightarrow G^*$ with $\sigma(1) = 1$ such that the $2$-cocycle
\begin{align*}
\mu : G^2 &\rightarrow \Lambda^{\times} \\
(x,y) &\rightarrow \sigma(x) \sigma(y) \sigma(xy)^{-1},
\end{align*}
is unitary.
\end{prop}

Indeed, if $\sigma_0 : G \rightarrow G^*$ is an arbitrary section such that $\sigma_0(1) = 1$, with associated $2$-cocycle $\mu_0$, then we observed that there exists a map $\lambda : G \rightarrow \Lambda^{\times}$ such that $\lambda(1) = 1$ and such that $\mu_0 d^1(\lambda)$ is unitary. Since $\mu_0 d^1(\lambda)$ is the $2$-cocycle associated to the section $\lambda \sigma_0 : G \rightarrow G^*$, we can take $\sigma = \lambda \sigma_0$ in Proposition \ref{chap3hty}.

%

\subsection{\label{chap30.3.1}} Let $H$ be an open subgroup of finite index in a $\Lambda$-twisted topological group $(G,\mu)$. Then the restriction $\mu_{|H}$ of $\mu$ to $H \times H$ endows $H$ with a structure of $\Lambda$-twisted topological group, such that the inclusion $\iota : H \rightarrow G$ is a morphism of twisted topological groups. The functor $\iota^*$ (cf. \ref{chap30.2.1}) from $\Rep_{\Lambda}(G,\mu)$ to $\Rep_{\Lambda}(H,\mu_{|H})$ admits a left adjoint $\Ind_{H}^G$, given by
$$
\Ind_{H}^G(V) = \Lambda[G,\mu] \otimes_{\Lambda[H,\mu_{|H}]} V,
$$
cf. \ref{chap3equiv}. In order to verify that $\Ind_{H}^G(V)$ is indeed a $\Lambda$-admissible representation of $(G,\mu)$, it is sufficient by Proposition \ref{chap3restriction} to check that the restriction of $\Ind_{H}^G(V)$ to the finite index open subgroup $K = \cap_{g \in G/H} g H g^{-1}$ is a $\Lambda$-admissible representation of $(K,\mu_{|K})$. However, if $(g_c)_{c \in G/H}$ are left $H$-cosets representatives, then we have a decomposition
$$
\Lambda[G,\mu] = \bigoplus_{c \in G/H} [g_c] \Lambda[H,\mu_{|H}],
$$
as a right $\Lambda[H,\mu_{|H}]$-module, which yields in turn a decomposition
$$
\Ind_{H}^G(V) = \bigoplus_{c \in G/H} [g_c] V,
$$
where $[g_c] V$ is a $\Lambda$-admissible representation of $(K,\mu_{|K})$, since the action of $K$ on this $\Lambda$-module is given by 
\begin{align*}
[k] [g_c] v &= \mu(k,g_c) [k g_c] v\\ 
&= \mu(k,g_c) \mu(g_c,g_c^{-1}k g_c)^{-1} [g_c] [g_c^{-1}k g_c] v,
\end{align*}
for $k$ in $K$ and $v$ in $V$, so that $[g_c] V$ is isomorphic to the $\Lambda$-module $V$, endowed with the $\Lambda$-admissible map $\rho_c : K \rightarrow \Aut_{\Lambda}(V)$ given by $\rho_c(k) = \mu(k,g_c) \mu(g_c,g_c^{-1}k g_c)^{-1} \rho(g_c^{-1}k g_c)$.

\subsection{\label{chap30.2.0}} Let $(G,\mu)$ be a $\Lambda$-twisted topological group (cf. \ref{chap30.1}), and let $Z$ be a subgroup of $\Lambda^{\times}$ which contains the image of $\mu$, and which is contained in an admissible $\ell$-adic subring of $\Lambda$. Let us consider the central extension
$$
1 \rightarrow Z \xrightarrow[]{\iota} G^* \xrightarrow[]{\pi} G \rightarrow 1,
$$
associated to the $2$-cocycle $\mu$. The underlying topological space of $G^*$ is the product $Z \times G$, the group law is given by
$$
(\lambda_1,g_1) \cdot (\lambda_2,g_2) = (\lambda_1 \lambda_2 \mu(g_1,g_2), g_1 g_2),
$$
and we have $\iota(\lambda) = (\lambda,1)$ and $\pi(\lambda,g) = g$ for $(\lambda,g)$ in $Z \times G$.
The continuous map $\pi$ admits a distinguished continuous section, namely $\sigma : g \mapsto (1,g)$, which is a group homomorphism if and only if $\mu$ is trivial, i.e. $\mu = 1$.

If $(V,\rho)$ is a $\Lambda$-admissible representation of $(G,\mu)$ (cf. \ref{chap30.2}), then the continuous map
\begin{align*}
\rho^* : G^* &\rightarrow \Aut_{\Lambda}(V), \\
(\lambda,g) &\rightarrow \lambda \rho(g),
\end{align*}
is a group homomorphism, hence $(V,\rho^*)$ is a $\Lambda$-admissible representation of the topological group $G^*$. If $\zeta : Z \rightarrow \Lambda^{\times}$ is the character of $Z$ given by the inclusion, then the restriction of $(V,\rho^*)$ to $Z$ is $\zeta$-isotypical, i.e. $\rho^* \circ \iota(\lambda)(v) = \zeta(\lambda) v$ for all $(\lambda,v)$ in $Z \times V$. Conversely, any $\Lambda$-admissible representation of $G^*$ with $\zeta$-isotypical restriction to $Z$ yields a $\Lambda$-admissible representation of $(G,\mu)$ by composition with the section $\sigma$, and these two constructions are quasi-inverse to each other. We have obtained:

\begin{prop}\label{chap3projgroup} The functor $(V,\rho) \mapsto (V,\rho^*)$ is an equivalence of categories from $\Rep_{\Lambda}(G,\mu)$ to the category of $\Lambda$-admissible representations of $(G^*,1)$ whose restriction to $Z$ is $\zeta$-isotypical.
\end{prop}

If $H$ is a subgroup of $G$, endowed with the restriction of $\mu$ to $H$, then the corresponding group $H^*$ is the inverse image of $H$ by $\pi$. Let $\iota : H \rightarrow G$ be the inclusion, which is a morphism of $\Lambda$-twisted topological groups. Under the equivalence of Proposition \ref{chap3projgroup}, the functor $\iota^*$ (cf. \ref{chap30.2.1}) corresponds to the restriction functor from representations of $G^*$ to representations of $H^*$. By taking left adjoints when available (cf. \ref{chap30.3.1}), we obtain:

\begin{prop}\label{chap3projgroup2} Let $H$ be an open subgroup of finite index in $G$, with inverse image $H^*$ in $G^*$. Under the equivalence of Proposition \ref{chap3projgroup}, the functor $\Ind_{H}^G$ corresponds to $\Ind_{H^*}^{G^*}$.
\end{prop}


\subsection{\label{chap30.3.1.0}} Let us recall that if $(V,\rho)$ is a $\Lambda$-admissible representation of a $\Lambda$-twisted topological group $(G,\mu)$, then the composition
$$
G \xrightarrow[]{\rho} \Aut_{\Lambda}(V) \rightarrow \Aut_{\Lambda}(V)/ \Lambda^{\times},
$$
is a genuine group homomorphism. In particular, its image is a subgroup of $\Aut_{\Lambda}(V)/ \Lambda^{\times}$.

\begin{defi}\label{chap3fproj} A $\Lambda$-admissible representation $(V,\rho)$ of a twisted topological group $(G,\mu)$ is said to have \textit{finite projective image} if the composition of $\rho$ with the projection from $\Aut_{\Lambda}(V)$ to $\Aut_{\Lambda}(V)/{\Lambda}^{\times}$ has finite image. We denote by $\Rep_{\Lambda}^{\mathrm{fin}}(G,\mu)$ the full subcategory of $\Rep_{\Lambda}(G,\mu)$ whose objects are the $\Lambda$-admissible representations of $(G,\mu)$ with finite projective image.
\end{defi}

Under the equivalence of Proposition \ref{chap3projgroup}, the subcategory $\Rep_{\Lambda}^{\mathrm{fin}}(G,\mu)$ of $\Rep_{\Lambda}(G,\mu)$ is equivalent to the category of continuous linear representations of $G^*$ with finite projective image and with $\zeta$-isotypical restriction to $Z$ (with notation from \ref{chap3projgroup}).

\subsection{\label{chap30.3.1.1}} Let $(V,\rho)$ be an object of $\Rep_{\Lambda}^{\mathrm{fin}}(G,\mu)$, where $\Lambda$ is an $\ell$-adic coefficient ring \textit{in which $\ell$ is invertible}. Then the $\Lambda$-module $\End(V)$ is a $\Lambda$-admissible representation of $(G,1)$ under the action given by $g \cdot u = \rho(g) \circ u \circ \rho(g)^{-1}$. Moreover, this action factors through the image $G'$ of $G$ in $\Aut_{\Lambda}(V)/{\Lambda}^{\times}$, which is a finite group. Thus we can form the projector
\begin{align*}
P : \End(V) &\rightarrow \End(V) \\
u &\mapsto \frac{1}{|G'|} \sum_{g \in G'} g \cdot u,
\end{align*}
whose image is the space of endomorphisms of $(V,\rho)$. If an element $\pi$ of $\End(V)$ is a projector onto a $G$-stable subspace $W$, then $g \cdot \pi$ is a projector onto $W$ as well for each $g$ in $G'$, hence $P(\pi)$ is a projector onto $W$ which commutes with the action of $G$. Thus $W$ is a direct summand of $V$ in the additive category $\Rep_{\Lambda}(G,\mu)$. We thus obtain the following extension of Maschke's theorem:

\begin{prop}\label{chap3maschke2} Let $(G,\mu)$ be a $\Lambda$-twisted topological group, where $\Lambda$ is an $\ell$-adic coefficient ring in which $\ell$ is invertible. Then any object of $\Rep_{\Lambda}^{\mathrm{fin}}(G,\mu)$ is semisimple. 
\end{prop}

In particular, the indecomposable objects of $\Rep_{\Lambda}^{\mathrm{fin}}(G,\mu)$ are irreducible. Let us introduce the Grothendieck group of the category $\Rep_{\Lambda}^{\mathrm{fin}}(G,\mu)$:

\begin{defi}\label{chap3gr} Let $(G,\mu)$ be a $\Lambda$-twisted topological group. The Grothendieck group $K_0^{\mathrm{fin}}(G,\mu,\Lambda)$ is the quotient of the free abelian group with basis $( [V] )_{V}$ indexed by all $\Lambda$-admissible representations of $(G,\mu)$ with finite projective image and whose underlying $\Lambda$-module is $\Lambda^n$ for some integer $n$, by the relations
$$
[V'] + [V''] - [V],
$$ 
whenever $V$ is an extension of $V'$ by $V''$ in $\Rep_{\Lambda}^{\mathrm{fin}}(G,\mu)$. If $G$ is finite, the group $K_0^{\mathrm{fin}}(G,\mu,\Lambda)$ is simply denoted by $K_0(G,\mu,\Lambda)$.
\end{defi}

\begin{rema} The class of $\Lambda$-admissible representations of $(G,\mu)$ may not form a set, hence the unnatural restriction to $\Lambda$-modules which \textit{are} $\Lambda^n$ for some $n$, rather than being merely isomorphic to $\Lambda^n$ for some $n$.
\end{rema}

If $V$ is a $\Lambda$-admissible representation of $(G,\mu)$ with finite projective image, then $V$ is isomorphic to a representation $V'$ whose underlying $\Lambda$-module is $\Lambda^n$ for some integer $n$, and the class of $[V']$ in $K_0^{\mathrm{fin}}(G,\mu, \Lambda)$ depends only on $V$. This class will be simply denoted by $[V]$.
By Corollary \ref{chap3maschke2}, if $\ell$ is invertible in $\Lambda$ then the group $K_0^{\mathrm{fin}}(G,\mu, \Lambda)$ is a free abelian group, with basis given by the classes $[V]$ where $V$ is an (isomorphism class of) irreducible $\Lambda$-admissible representation of $(G,\mu)$ with finite projective image.

Before proceeding further, let us recall Brauer's induction theorem for finite groups:

\begin{teo}[\cite{S}, 10.5 Th. 20]\label{chap3brauer} If $G$ is a finite group, then the abelian group $K_0(G,1,C)$ is generated by the classes of representations of the form $\Ind_{H}^G V$, where $H$ is a subgroup of $G$ and $V$ is a one-dimensional $C$-linear representation of $H$.
\end{teo}

Brauer's induction theorem \ref{chap3brauer} is a consequence of the following two results:
\begin{enumerate}
\item If $G$ is a finite group, then $K_0(G,1,C)$ is generated by the classes of representations of the form $\Ind_{H}^G V$, where $H$ is a nilpotent subgroup of $G$, and where $V$ is an irreducible representation of $H$.
\item Any irreducible $C$-linear representation of a finite nilpotent group $G$ is isomorphic to $\Ind_{H}^G V$ for some subgroup $H$ and some one-dimensional representation $V$ of $H$.
\end{enumerate}
Moreover, it is sufficient for the first of these results to prove that the class of the trivial representation of $G$ is a linear combination in $K_0(G,1)$ of representations of the form $\Ind_{H}^G V$, where $H$ is a nilpotent subgroup of $G$. We refer to (\cite{S}, 10) for proofs of these results, and for a more complete discussion of Brauer's theorem. We will also need the following variant of Brauer's theorem:

\begin{teo}[\cite{De73}, Prop. 1.5]\label{chap3brauervirtuel} If $G$ is a finite group, then the abelian group $K_0(G,1,C)$ is generated by the class $[C]$ of the trivial representation of $G$ and by the classes of the form $[\Ind_{H}^G V] - [\Ind_{H}^G C]$, where $H$ is a subgroup of $G$ and $V$ is a one-dimensional $C$-linear representation of $H$.
\end{teo}

 We will need an extension of Brauer's theorem to continous representations of twisted groups with finite projective image:

\begin{teo}\label{chap3brauer2} Let $(G,\mu)$ be a $C$-twisted topological group. Then the group $K_0^{\mathrm{fin}}(G,\mu,C)$ (cf. \ref{chap3gr}) is generated by the classes of $C$-admissible representations of the form $\Ind_{H}^G V$, where $H$ is an open subgroup of finite index in $G$ and $V$ is a one-dimensional $C$-admissible representation of $(H, \mu_{|H})$.
\end{teo}

Indeed, let $(V,\rho)$ be an object of $\Rep_C^{\mathrm{fin}}(G,\mu)$, and let us consider the subgroup $G' = C^{\times} \rho(G)$ of $\Aut_C(V)$ generated by the image of $\rho$ and by homotheties. The topological group $G'$ is a central extension of the finite group $G'/C^{\times}$ by $C^{\times}$. By Proposition \ref{chap3hty}, there exists a set-theoretic section $\sigma$ from $G'/C^{\times}$ to $G'$ such that $\sigma(1)=1$, whose associated $2$-cocycle is unitary. There exists a unique continuous map $\lambda : G \rightarrow C^{\times}$ such that
$$
\sigma(\rho(g) C^{\times}) = \lambda(g) \rho(g),
$$
for all $g$ in $G$, and we have $\lambda(1)=1$. By replacing $(V,\rho)$ and $(G,\mu)$ by $(V,\lambda \rho)$ and $(G,d(\lambda) \mu)$ respectively (cf. \ref{chap3cobor}), we can therefore assume (and we do) that $\mu$ is unitary and that the set $\rho(G)$ is finite. 

Since $\mu$ is now assumed to be unitary, we can choose $Z$ to be a finite subroup of $C^{\times}$ in \ref{chap30.2.0}. Moreover, the representation $(V,\rho^*)$ of $G^*$ corresponding to $(V,\rho)$ by Proposition \ref{chap3projgroup} has finite image, namely $Z \rho(G)$. By applying Brauer's induction theorem \ref{chap3brauer} to the finite group $\rho^*(G^*)$, we obtain a decomposition
\begin{align}\label{chap3decbrauer}
[V] = \sum_{i \in I} n_i [\Ind_{H_i}^{G^*} \chi_i],
\end{align}
in $K_0^{\mathrm{fin}}(G^*,1)$, where $n_i$ is a (possibly negative) integer, $H_i$ is an open subgroup of finite index in $G^*$ and $\chi_i$ is a continuous character of $H_i$ with finite image. Since $Z$ is central in $G^*$, the decomposition of a representation into the isotypical components of its restriction to $Z$ yields a splitting
$$
K_0^{\mathrm{fin}}(G^*,1) = \bigoplus_{\zeta' : Z \rightarrow C^{\times}} K_0^{\mathrm{fin}}(G^*,1)[\zeta'],
$$
where $\zeta'$ runs through the set of characters of $Z$, and $K_0^{\mathrm{fin}}(G^*,1)[\zeta']$ is generated by continuous linear representations of $G^*$ with finite projective image and $\zeta'$-isotypical restriction to $Z$. By projecting the relation \ref{chap3decbrauer} onto the factor corresponding to $\zeta$, we obtain
$$
[V] = \sum_{i \in I} n_i [\Ind_{H_i}^{G^*} \chi_i[\zeta]],
$$
in the abelian group $K_0^{\mathrm{fin}}(G^*,1)[\zeta]$, where we denoted by $[\zeta]$ the $\zeta$-isotypical component.


For each $i$ in $I$, we can identify the representation $\Ind_{H_i}^{G^*} \chi_i$ with the $C$-vector space 
$$
\{ f : G^* \rightarrow C \ | \ \forall h \in H_i, g \in G^*, f(h g) = \chi_i(h) f(g) \},
$$
endowed with the action $(g \cdot f)(x) = f(xg)$ for $g,x$ in $G^*$. The $\zeta$-isotypical component is then the subspace of functions $f : G^* \rightarrow C$ such that $f(gh) = \zeta(h) f(g)$ for all $h$ in $Z$ and all $g$ in $G^*$. Thus the $\zeta$-isotypical component vanishes whenever $\chi_i$ and $\zeta$ do not coincide on $H_i \cap Z$. On the other hand, if $\chi_i$ and $\zeta$ do coincide on $H_i \cap Z$, then there exists a unique character $\psi_i$ of $Z H_i $ such that $\psi_{i|Z} = \zeta$ and $\psi_{i |H_i} = \chi_i$, in which case
\begin{align*}
\Ind_{H_i}^{G^*} \chi_i[\zeta] &= \{ f : G^* \rightarrow C \ | \ \forall h \in ZH_i, g \in G^*, f(h g) = \psi_i(h) f(g) \}, \\
&= \Ind_{ZH_i}^{G^*} \psi_i.
\end{align*}
If $G_i$ is the image in $G$ of $Z H_i$, then $G_i$ is an open subgroup of finite index in $G$, and $\psi_i$ yields a continuous one-dimensional linear representation of $(G_i,\mu_{|G_i})$, such that the relation
$$
V = \sum_{\substack{i \in I \\ \chi_{i|H_i \cap Z} = \zeta_{|H_i \cap Z} } } n_i [\Ind_{G_i}^G \psi_i]
$$
holds in $K_0^{\mathrm{fin}}(G,\mu)$ (cf. \ref{chap3projgroup}, \ref{chap3projgroup2}), hence the conclusion of Theorem \ref{chap3brauer2}.

%

\subsection{\label{chap30.4}} Let $(G,\mu)$ be a profinite $C$-twisted topological group (cf. \ref{chap30.1}), and let us consider an extension
$$
1 \rightarrow I \xrightarrow[]{\iota} Q \xrightarrow[]{\pi} G \rightarrow 1,
$$
of profinite topological groups. It is assumed that $G$ carries the quotient topology from $Q$. Then the formula $(\pi^*\mu)(x,y) = \mu(\pi(x),\pi(y))$ for $x,y$ in $Q$ defines a multiplier on $Q$, so that $(Q,\pi^* \mu)$ is a $C$-twisted topological group, and $\pi$ is a morphism of $C$-twisted topological groups (cf. \ref{chap30.1}). 

\begin{prop} \label{2.1.0} Assume that there exists a closed normal subgroup $I_{\ell}$ of $I$ of profinite order coprime to $\ell$ such that $I/I_{\ell}$ is a free pro-$\ell$-group on one generator. Let $(V,\rho)$ be an irreducible $C$-admissible representation of $(Q,\pi^*\mu)$. Then there exists a closed subgroup $I'$ of $I$, normal in $Q$, such that:
\begin{itemize}
\item[$(i)$] the image $\rho(I')$ of $I'$ by $\rho$ is finite;
\item[$(ii)$] the quotient group $I/I'$ is topologically of finite type, and its centralizer in $Q/I'$ is open.
\end{itemize}
\end{prop}

Indeed, the group $G = Q/I$ acts on the free pro-$\ell$-group $I/I_{\ell}$ through a continuous character $\chi : G \rightarrow \mathbb{Z}_{\ell}^{\times}$. If $\chi$ has finite image, then its kernel is open, hence the centralizer of $I/I_{\ell}$ is open in $Q/I_{\ell}$, and one may therefore take $I' = I_{\ell}$. One should note that $I_{\ell}$ must be normal in $Q$, since it is a characteristic subgroup of $I$.

We henceforth assume that $\chi$ has infinite image. We will show that in this case the image $\rho(I)$ is finite, so that one can take $I' = I$. Since $\rho(I_{\ell})$ is finite, there exists an open subgroup $Q'$ of $Q$ such that $\rho(Q' \cap I_{\ell})$ is trivial. By replacing $(Q,V)$ with $(Q'/Q' \cap I_{\ell},V')$, for irreducible subrepresentations $V'$ of the semisimple representation $V_{|Q'}$ of $(Q',\pi^*\mu_{|Q'})$, we can assume (and we henceforth do) that $I_{\ell}$ is trivial. Let $t$ be a topological generator of the free pro-$\ell$-group $I$. For any element $g$ of $G$, the automorphism $\rho(t)^{\chi(g)}$ is a conjugate of $\rho(t)$. Since $\chi$ has infinite image, this implies that the eigenvalues of $\rho(t)$ are roots of unity. In particular there exists an integer $m$ such that $\rho(t)^m$ is unipotent. The kernel of $\rho(t)^m - \mathrm{Id}_V$ is then a non zero sub-$C$-admissible representation of $(V,\rho)$, hence $\rho(t)^m = \mathrm{Id}_V$ by irreducibility of $(V,\rho)$. Thus the image of $I$ by $\rho$ is finite.

\begin{prop}\label{chap3irrfinite} \label{2.1.2} Let us assume that there exists a finite subgroup $I'$ of $I$, normal in $Q$, such that the quotient group $I/I'$ is topologically finitely generated, with open centralizer in $Q/I'$. Then the centralizer of $I$ in $Q$ is open. Moreover, any $C$-admissible representation of $(Q,\pi^*\mu)$ over $C$, which is irreducible as a representation of $I$, has finite projective image.
\end{prop}

By assumption, there exists an open subgroup $Q_0$ of $Q$ whose image in $Q/I'$ centralizes $I/I'$. Let $\Sigma$ be a finite subset of $I$ whose image in $I/I'$ generates a dense subgroup. By replacing $\Sigma$ with $\Sigma \cup I'$, we can assume (and we do) that $\Sigma$ generates a dense subgroup of $I$. For any element $t$ of $\Sigma$, the set of commutators $[Q_0,t]$ is contained in the finite group $I'$. By continuity of each of the maps $q \rightarrow [q,t] = qt q^{-1} t^{-1}$ from the profinite topological space $Q_0$ to the discrete topological space $I'$, there exists an open subgroup $Q_1$ of $Q_0$ such that we have $[q,t] = 1$ for any elements $q$ and $t$ of $Q_1$ and $\Sigma$ respectively. In particular, the open subgroup $Q_1$ of $Q$ centralizes $\Sigma$, as well as the closed subgroup generated by the latter finite set, namely $I$. The centralizer of $I$ in $Q$ thus contains $Q_1$, and is consequently open in $Q$.

The last assertion of Proposition \ref{2.1.2} is proved as follows: if $(V,\rho)$ is a $C$-admissible representation of $(Q,\pi^*\mu)$, which is irreducible as a representation of $I$, then Schur's lemma ensures that any element of $Q$ centralizing $I$ acts on $V$ as a homothety. Since the centralizer of $Q$ in $I$ is open, hence of finite index in $Q$, the projective image of $\rho$ is finite.


\begin{prop}\label{chap3dec}  \label{2.1.1}Let $(G,\mu)$ be a profinite $C$-twisted topological group, and let
$$
1 \rightarrow I \xrightarrow[]{\iota} Q \xrightarrow[]{\pi} G \rightarrow 1,
$$
be an extension of profinite topological groups, where $G$ (resp. $I$) carries the quotient topology (resp. induced topology) from $Q$. 

 Let $(V,\rho)$ be a $C$-admissible representation of $(Q,\pi^*\mu)$, such that $(V,\rho \iota)$ is a semisimple $C$-admissible representation of $I$ and let us assume that there exists a closed subgroup $I'$ of $I$, normal in $Q$, such that:
\begin{itemize}
\item[$(i)$] the image $\rho(I')$ of $I'$ by $\rho$ is finite;
\item[$(ii)$] the quotient group $I/I'$ is topologically of finite type, and its centralizer in $Q/I'$ is open.
\end{itemize}

Then there exists a finite family $(G_j)_{j \in J}$ of open subgroups of finite index in $G$, and for each $j$, a $C$-admissible multiplier $\nu_j$ on $G_j$, a $C$-admissible representation $W_j$ of $(G_j,\mu \nu_j^{-1})$, and a $C$-admissible representation $E_j$ of $(Q_j,\pi^* \nu_j)$, where $Q_j = \pi^{-1}(G_j)$, such that:
\begin{enumerate}
\item for each $j$ in $J$, the twisted representation $E_j$ has finite projective image, its restriction to $I$ is irreducible, and the action of $I$ on $E_j$ factors through $\rho(I)$;
\item the $C$-admissible representation $(V,\rho)$ of $Q$ is isomorphic to the direct sum
$$
\bigoplus_{j \in J} \Ind_{Q_j}^Q( E_j \otimes_C \pi^* W_j).
$$
\end{enumerate}
\end{prop}


\begin{rema} Even if we start with a genuine representation of $Q$, namely if $\mu = 1$, then the cocycles $\nu_j$ appearing in Proposition \ref{chap3dec} are usually not trivial, nor cohomologically trivial. 
\end{rema}

\begin{rema}\label{irreduciblerema} Let $(V,\rho)$ be a semisimple $C$-admissible representation of $(Q,\pi^*\mu)$. Then its restriction to $I$ is semisimple as well. If moreover the conditions $(i)$ and $(ii)$ of Proposition \ref{2.1.1} are fulfilled by some closed subgroup $I'$ of $I$, then we can apply Proposition \ref{2.1.1} to $(V,\rho)$. The $C$-admissible representations $W_j$ of $(G_j,\mu \nu_j^{-1})$ can then be taken to be semisimple. Indeed, we can assume $(V,\rho)$ to be irreducible, in which case each factor $W_j$ either vanishes or is irreducible.
\end{rema}

\begin{rema}\label{chap3unitary10} If $\mu$ is unitary, then the multipliers $\nu_j$ which appear in Proposition \ref{chap3dec} can be taken to be unitary as well. Indeed, one can assume that each $W_j$ is non zero, and thus by Proposition \ref{chap3unitary2}, there exists for each $j$ a $C$-admissible map $\lambda_j : G_j \rightarrow C^{\times}$ such that $\lambda_j(1)=1$ and such that $\nu_j d(\lambda_j)$ is unitary. Twisting $E_j$ and $W_j$ by $\lambda_j \circ \pi$ and $\lambda_j^{-1}$ respectively (cf. \ref{chap3cobor}) then leaves the tensor product $E_j \otimes_C \pi^* W_j$ unaltered.
\end{rema}

Let us prove Proposition \ref{chap3dec}. Let $(V,\rho)$ and $I'$ as in Proposition \ref{chap3dec}. The subgroup $K = I \cap \ker(\rho)$ of $Q$ is normal and closed. By replacing $(Q,I)$ with $(Q/K,I/K)$, we can assume (and we do) that the restriction of $(V,\rho)$ to $I$ is faithful, so that $I'$ is a finite subgroup of $I$ by the assumption $(i)$. Since the restriction of $\pi^*\mu$ to $I$ is trivial, the restriction $\iota^* V$ is a honest linear representation of $I$, which is moreover semisimple by assumption. 

For each isomorphism class $\chi$ of irreducible representations of $I$, let $V[\chi]$ be the $\chi$-isotypical component of $(V,\rho \iota)$. Let $X$ be the set of isomorphism classes of irreducible representations $\chi$ such that $V[\chi]$ is non zero, so that we have a decomposition
$$
V = \bigoplus_{\chi \in X} V[\chi].
$$
Let us denote by $\varphi : Q \rightarrow \Aut(I)$ the action of $Q$ on $I$ by automorphisms given by conjugation by an element of $Q$, i.e. sending an element $q$ of $Q$ to $(t \mapsto q t q^{-1})$. The homomorphism $\varphi$ is continuous with finite image, when $ \Aut(I)$ is endowed with the discrete topology, since its kernel is open in $Q$ by Proposition \ref{2.1.2}.

The group $Q$ acts continuously on the left on the finite discrete set $X$ by $q \cdot \chi = \chi \circ \varphi(q)^{-1}$, and $\rho(q)$ sends $V[\chi]$ onto $V[q \cdot \chi]$. This action factors through $\pi$, so that we obtain a continuous action of $G$ on $X$. Let $J$ be the set of $G$-orbits in $X$ and for each $j$ in $J$, let $\chi_j$ be a member of the orbit $j$. Let $G_j$ be the stabilizer of $\chi_j$ in $G$, and let $Q_j$ be its inverse image in $Q$. Then the action of $Q$ on $V$ restricts to an action of $Q_j$ on $V[\chi_j]$, and the homomorphism
\begin{align*}
\Ind_{Q_j}^Q(V[\chi_j]) &\rightarrow \bigoplus_{q \in Q/Q_j }V[q \chi_j] = \bigoplus_{\chi \in j}V[ \chi]\\
 [q] \otimes v &\mapsto \rho(q) v
\end{align*}
is an isomorphism. Thus $V$ is isomorphic to the direct sum 
$$
\bigoplus_{j \in J} \Ind_{Q_j}^Q( V[\chi_j]).
$$

Let $(E_j,\rho_{j})$ be an irreducible representation of $I$ in the class $\chi_j$, and let 
$$
W_{j} = \Hom_{C[I]}( E_{j}, V) = \Hom_{C[I]}( E_{j}, V[\chi_j]),
$$
considered as a trivial representation of $I$, so that the $C$-linear homomorphism
\begin{align}
\begin{split}
\label{chap3eq2}
E_{j} \otimes_C W_{j} &\rightarrow V[\chi_j]\\
e \otimes w &\mapsto w(e),
\end{split}
\end{align}
is an isomorphism of representations of $I$. 

For each $\psi$ in $\Aut(I)$ such that $\chi_j \circ \psi^{-1} = \chi_j$, the representations $(E_j,\rho_j)$ and $(E_j,\rho_j \circ \psi^{-1})$ are isomorphic, hence there exists a $C$-linear automorphism $\lambda_j
(\psi)$ of $E_j$ such that
$$
\rho_j \circ \psi^{-1}(t) = \lambda_j(\psi)^{-1} \rho_j(t) \lambda_j(\psi),
$$
for all $t$ in $I$. We can take $\lambda_j(\mathrm{id}_I) = \mathrm{id}_{E_j}$. Note that Schur's lemma implies that $\lambda_j(\psi)$ is uniquely determined up to multiplication by an invertible scalar, and that if we denote by $\Aut(I)_{\chi_j}$ the subgroup of $\Aut(I)$ formed by the elements $\psi$ of $\Aut(I)$ such that $\chi_j \circ \psi^{-1} = \chi_j$, then the composition
$$
\Aut(I)_{\chi_j} \xrightarrow[]{\lambda_j} \Aut_C(E_j) \rightarrow \Aut_C(E_j)/ C^{\times},
$$
is a group homomorphism.

Let $\sigma$ be a continuous section of the continuous map $\pi$, such that $\sigma(1)=1$. Such a section always exists by (\cite{Scg}, I.1.2, Prop. 1). For $q = \sigma(g) t$ in $Q_j$, with $t$ in $I$ and $g$ in $G_j$ (so that $g = \pi(q)$), we set
$$
\widetilde{\rho}_j(q) = (\lambda_j \circ \varphi \circ \sigma)(g) \rho_j(t).
$$
For $t$ in $I$, the automorphism $\lambda_j(\varphi(t))$ differs from $\rho_j(t)$ only by an invertible scalar, hence for each $q$ in $Q_j$, the automorphism $\widetilde{\rho}_j(q)$ differs from $\lambda_j \circ \varphi(q)$ by an invertible scalar. In particular, the projective image of $\widetilde{\rho}_j$ is finite, since $\varphi$ has finite image.

\begin{lem}\label{chap3lem1} For each $j$ in $J$, there exists a unique multiplier $\nu_j$ on $G_j$ such that $(E_j,\widetilde{\rho}_j)$ is a continuous linear representation of $(Q_j,\pi^* \nu_j)$ with finite projective image.
\end{lem}

Assuming the conclusion of Lemma \ref{chap3lem1}, we can conclude the proof of Proposition \ref{chap3dec} as follows. We first notice that the continuous map
\begin{align*}
Q_j &\rightarrow \Aut_C(W_j)\\
q &\mapsto (w \mapsto \rho(q) \circ w \circ \widetilde{\rho}_j(q)^{-1}).
\end{align*}
is right $I$-invariant, hence uniquely factors as $\tau_j \circ \pi_{|Q_j}$, where $\tau_j$ is a continuous map from $G_j$ to $\Aut_C(W_j)$. One then checks that $(W_j, \tau_j)$ is a $C$-admissible representation of $(G_j, \mu \nu_j^{-1})$, so that (\ref{chap3eq2}) is an isomorphism of $C$-admissible representations of $(Q_j, \pi^* \mu)$. Thus the decomposition
$$
V \simeq \bigoplus_{j \in J} \Ind_{Q_j}^Q( E_{j} \otimes_C \pi^* W_{j}).
$$
provides the wanted conclusion.

\subsection{\label{chap30.7}} Let us now prove Lemma \ref{chap3lem1}. The maps $\lambda_j$ and $\varphi_{|Q_j}$ are both continuous when $\Aut(I)_{\chi_j}$ is endowed with the discrete topology, hence $\widetilde{\rho}_j$ is continuous. Let us write
$$
c(g_1,g_2) = \sigma(g_2)^{-1} \sigma(g_1)^{-1} \sigma(g_1g_2) \in I.
$$
If $q_1 = \sigma(g_1) t_1$ and $q_2 = \sigma(g_2) t_2$ are elements of $G_j$, with $t_1,t_2$ in $I$, then we have
$$
q_1 q_2 = \sigma(g_1g_2) (c(g_1,g_2)^{-1} (\sigma(g_2)^{-1} t_1 \sigma(g_2)) t_2),
$$
so that we have
\begin{align*}
\widetilde{\rho}_j(q_1) \widetilde{\rho}_j(q_2) &= (\lambda_j \circ \varphi \circ \sigma)(g_1) \rho_j(t_1) (\lambda_j \circ \varphi \circ \sigma)(g_2) \rho_j(t_2) \\
&= (\lambda_j \circ \varphi \circ \sigma)(g_1) (\lambda_j \circ \varphi \circ \sigma)(g_2) \rho_j(\sigma(g_2)^{-1} t_1 \sigma(g_2) t_2) \\
&= (\lambda_j \circ \varphi \circ \sigma)(g_1) (\lambda_j \circ \varphi \circ \sigma)(g_2) (\rho_j \circ c)(g_1,g_2) (\lambda_j \circ \varphi \circ \sigma)(g_1 g_2)^{-1} \widetilde{\rho}_j(q_1 q_2).
\end{align*}
Let us define
$$
\nu_j(g_1,g_2) = (\lambda_j \circ \varphi \circ \sigma)(g_1) (\lambda_j \circ \varphi \circ \sigma)(g_2) (\rho_j \circ c)(g_1,g_2) (\lambda_j \circ \varphi \circ \sigma)(g_1 g_2)^{-1},
$$
so that we have 
\begin{align}\label{chap3calcul}
\widetilde{\rho}_j(q_1) \widetilde{\rho}_j(q_2) = \nu_j(g_1,g_2) \widetilde{\rho}_j(q_1q_2).
\end{align}
It remains to show that $\nu_j$ defines a multiplier on $G_j$, i.e. that $\nu_j(g_1,g_2)$ takes values in $C^{\times}$ and that it satisfies the cocycle formula (\ref{chap3cocy}). We noticed that for each $q$ in $Q_j$, the automorphism $\widetilde{\rho}_j(q)$ differs from $\lambda_j \circ \varphi(q)$ by an invertible scalar (cf. the discussion before Lemma \ref{chap3lem1}). Thus, we have $r \circ \widetilde{\rho}_j = r \circ \lambda_j \circ \varphi$, where $r$ is the projection from $\Aut_C(E_j)$ to $\Aut_C(E_j)/ C^{\times}$. In particular, since $r \circ \lambda_j$ is a group homomorphism, so is $r \circ \widetilde{\rho}_j$. This implies that $r \circ \nu_j$ is identically equal to $1$. Thus $\nu_j$ takes values in $C^{\times}$, and the formula (\ref{chap3calcul}) then implies that $\nu_j$ satisfies the cocycle condition (\ref{chap3cocy}). Consequently, $(E_j,\widetilde{\rho}_j)$ is a $C$-admissible representation of $(Q_j,\pi^* \nu_j)$.

\subsection{\label{chap30.6}} Let $(G,\mu)$ be a $C$-twisted profinite topological group (cf. \ref{chap30.1}). Let us consider an extension
$$
1 \rightarrow I \xrightarrow[]{\iota} Q \xrightarrow[]{\pi} G \rightarrow 1,
$$
of profinite topological groups, where $G$ (resp. $I$) carries the quotient topology (resp. induced topology) from $Q$. Let $K_0(Q,\pi^*\mu,C)$ be the Grothendieck group of $C$-admissible representations of $(Q,\pi^*\mu)$. Thus any $C$-admissible representation $V$ of $(Q,\pi^*\mu)$ has a well defined class $[V]$ in $K_0(Q,\pi^*\mu,C)$, and the latter is generated by such classes with relations $[V] = [V'] + [V'']$ for each short exact sequence
$$
0 \rightarrow V' \rightarrow V \rightarrow V'' \rightarrow 0,
$$
of $C$-admissible representations of $(Q,\pi^*\mu)$.

\begin{prop}\label{chap3generators5} Let us consider an extension
$$
1 \rightarrow I \xrightarrow[]{\iota} Q \xrightarrow[]{\pi} G \rightarrow 1,
$$
of profinite topological groups, where $G$ (resp. $I$) carries the quotient topology (resp. induced topology) from $Q$. Assume that there exists a closed normal subgroup $I_{\ell}$ of $I$ of profinite order coprime to $\ell$ such that $I/I_{\ell}$ is a free pro-$\ell$-group on one generator. 

 Then the abelian group
$$
\bigoplus_{\mu} K_0(Q,\pi^*\mu,C),
$$
where the sum runs over all unitary $C$-admissible multipliers $\mu$ on $G$, and $K_0(Q,\pi^*\mu,C)$ is as in \ref{chap30.6}, is generated by its subset of elements of the following two types:
\begin{enumerate}
\item the class $[C]$ in $K_0(Q,1,C)$,
\item for any unitary $C$-admissible multiplier $\mu$ on $G$, any open subgroup $Q'$ of $Q$, with image $G'$ in $G$, any unitary $C$-admissible multipliers $\mu_1$ and $\mu_2$ on $G'$ such that $\mu_1 \mu_2 = \mu_{| G'}$, any $C$-admissible representation $E$ of $(Q',\pi^* \mu_1)$ of rank $1$ and any irreducible $C$-admissible representation $W$ of $(G',\mu_2)$, the class 
$$
[\Ind_{Q'}^Q( E \otimes W)] - \rk(W) [\Ind_{Q'}^Q( C)],
$$
in the sum of $K_0(Q,\pi^* \mu,C)$ and $K_0(Q,1,C)$.
\end{enumerate}
\end{prop}

Let $V$ be a $C$-admissible representation of $(Q,\pi^*\mu)$, for some unitary $C$-admissible multiplier $\mu$ on $G$. We must prove that the class $[V]$ belongs to the group generated by the classes described in Proposition \ref{chap3generators5}. We can assume (and we do) that $V$ is irreducible. By Proposition \ref{2.1.0}, the representation $V$ then satisfies the assumptions of Proposition \ref{2.1.1}. Thus Proposition \ref{chap3dec}, together with Remarks \ref{irreduciblerema} and \ref{chap3unitary10}, yield that $V$ is of the form $\Ind_{Q'}^Q ( E \otimes \pi^* W)$, where $Q'$ is an open subgroup of $Q$ containing $I$, where $E$ is a $C$-admissible representation of $(Q',\pi^* \mu_1)$ with finite projective image, whose restriction to $I$ is irreducible, and $W$ is an irreducible $C$-admissible representation of $(G',\mu_2)$, for some unitary $C$-admissible multipliers $\mu_{1}$ and $\mu_{2}$ on the image $G'$ of $Q'$ in $G$, such that $\mu_{1} \mu_{2} = \mu_{| G'}$.

%


 By Proposition \ref{chap3irrfinite} and Theorem \ref{chap3brauer2}, we can assume (and we do) that $E$ is of the form $\Ind_{Q''}^{Q'} V_1$, where $Q''$ is an open subgroup of $Q'$, and where $V_1$ is a $C$-admissible representation of rank $1$ of $(Q'',\pi^* \mu_1)$. We then have an isomorphism 
 $$
 V = \Ind_{Q'}^Q( \Ind_{Q''}^{Q'} V_1 \otimes \pi^* W) \cong \Ind_{Q''}^{Q}( V_1 \otimes \pi^* W),
 $$
 hence a decomposition
 $$
 [V] = \left([\Ind_{Q''}^{Q}( V_1 \otimes \pi^* W)] - \rk(W)[\Ind_{Q''}^{Q} C] \right) + \rk(W) [\Ind_{Q''}^{Q} C].
 $$
The first term in this decomposition is of the required type $(2)$, while the last term $[\Ind_{Q''}^{Q} C]$ belongs to the subgroup of $K_0(Q,1,C)$ generated by elements of type $(1)$ or of type $(2)$ (with trivial multipliers $\mu,\mu_1,\mu_2$ and with trivial factor $W$): indeed, the group $Q$ acts on $\Ind_{Q''}^{Q} C$ through a finite quotient, hence the result follows from Theorem \ref{chap3brauervirtuel}.

\subsection{\label{CGTRES}} Let $G$ be a profinite group, and let $K_0(G,\Lambda)$ be the Grothendieck group $K_0(G, 1, \Lambda)$ of virtual $\Lambda$-admissible (untwisted) representations of $G$, cf. \ref{chap30.6}. We have a decomposition homomorphism
$$
d : K_0(G, \overline{\mathbb{Q}_{\ell}}) \rightarrow K_0(G, \overline{\mathbb{F}_{\ell}}),
$$
cf. (\cite{S}, 15.2), such that for any $\overline{\mathbb{Z}_{\ell}}$-admissible representation of $G$ on a free $\overline{\mathbb{Z}_{\ell}}$-module of finite rank, we have
$$
d( [V[\ell^{-1}]] ) = [V/ \m V],
$$
where $\m$ is the maximal ideal of $\overline{\mathbb{Z}_{\ell}}$. The homomorphism $d$ is surjective by (\cite{S}, 16.1 Th.33).

\subsection{\label{CGTRES2}} Let $G$ be a profinite group and let $\mu$ be a unitary $\overline{\mathbb{F}_{\ell}}$-admissible $2$-cocycle on $G$. We also denote by $\mu$ the Teichm\"{u}ller lift of $\mu$ to $\overline{\mathbb{Z}_{\ell}}$. Let $Z$ be a finite subgroup of $\overline{\mathbb{Z}_{\ell}}^{\times}$, of order coprime to $\ell$, containing the image of $\mu$, and let $G^*$ be the central extension of $G$ by $Z$ induced by $\mu$, cf. \ref{chap30.2.0}. For any $\Lambda$-admissible representation $V$ of $G^*$, we have a splitting
$$
V = \bigoplus_{\zeta : Z \rightarrow \Lambda^{\times}} V[\zeta],
$$
where the sum runs over homomorphisms from $Z$ to $\Lambda^{\times}$, and $V[\zeta]$ is the $\zeta$-isotypic component of the restriction of $V$ to $Z$. Let us denote by $K_0(G^*,\Lambda)[\zeta]$ the Grothendieck group of virtual $\Lambda$-admissible representations of $G^*$ on which $Z$ acts through $\zeta$. We thus have a splitting
$$
K_0(G^*,\Lambda) = \bigoplus_{\zeta : Z \rightarrow \Lambda^{\times}} K_0(G^*,\Lambda)[\zeta].
$$
The decomposition homomorphism 
$$
d : K_0(G^*, \overline{\mathbb{Q}_{\ell}}) \rightarrow K_0(G^*, \overline{\mathbb{F}_{\ell}}),
$$
cf. \ref{CGTRES}, respects this splitting. Since $d$ is surjective, it induces a surjective homomorphism
$$
d : K_0(G^*, \overline{\mathbb{Q}_{\ell}})[\zeta] \rightarrow K_0(G^*, \overline{\mathbb{F}_{\ell}})[\zeta],
$$
for any $\zeta$. If we take $\zeta$ to be the canonical inclusion from $Z$ into $\overline{\mathbb{Z}_{\ell}}^{\times}$, we obtain a surjective decomposition homomorphism
$$
K_0(G,\mu,\overline{\mathbb{Q}_{\ell}}) \rightarrow K_0(G,\mu,\overline{\mathbb{F}_{\ell}}),
$$
which sends the class $V[\ell^{-1}]$, for any $\mu$-twisted $\overline{\mathbb{Z}_{\ell}}$-admissible representation of $G$ on a free $\overline{\mathbb{Z}_{\ell}}$-module of finite rank, to the class of $V/ \m V$, where $\m$ is the maximal ideal of $\overline{\mathbb{Z}_{\ell}}$.

\section{Twisted \texorpdfstring{$\ell$}{l}-adic sheaves\label{chap3twistedsection}}

Let $k$ be a perfect field, and let $\bk$ be an algebraic closure of $k$. We denote by $G_k = \Gal(\bk/k)$ the Galois group of the extension $\bk/k$, and For any $k$-scheme $X$, and for any $k$-algebra $k'$, we denote by $X_{k'}$ the fiber product of $X$ and $\Spec(k')$ over $\Spec(k)$. The group $G_k$ acts on the left on $\bk$, and thus acts on the right on $X_{\bk}$. 
Let $\Lambda$ be an $\ell$-adic coefficient ring (cf. \ref{chap30.0.0.1}). We fix a unitary $\Lambda$-admissible mutiplier $\mu$ on the topological group $G_k$ (cf. \ref{chap3conv}, \ref{chap30.0}, \ref{chap3unitary}). 

\begin{defi}\label{chap3neutra} A finite Galois extension $k'/k$ contained in $\bk$ is said to neutralize $\mu$ if $\mu$ is the pullback of a multiplier on the finite quotient $\Gal(k'/k)$ of $G_k$. 
\end{defi}

By \ref{chap30.0.0.0.1}, the unitary multiplier $\mu$ is neutralized by some finite Galois extension of $k$.

%
%

\subsection{\label{chap31.0.1}} Assume that $\Lambda$ is a finite $\ell$-adic coefficient ring, and let $X$ be a $k$-scheme. We denote by $\Loc(X,\Lambda)$ the category of locally constant constructible $\Lambda$-modules on the small \'etale site of $X$. Moreover, we denote by $\Sh(X,\Lambda)$ the abelian category of constructible sheaves of $\Lambda$-modules on the small \'etale site of $X$. Thus an object $\F$ of $\Sh(X,\Lambda)$ is a $\Lambda$-module in the \'etale topos of $X$, such that for any affine open subset $U$ of $X$, there exists a finite partition $U = \sqcup_{j \in J} U_j$ into constructible locally closed subschemes, such that $\F_{|U_j}$ belongs to $\Loc(X,\Lambda)$.

\subsection{\label{chap31.0.2}} Assume that $\Lambda$ is the ring of integers in a finite subextension of $\mathbb{Q}_{\ell}$ in $C$, and let $X$ be a locally noetherian $k$-scheme. We denote by $\Sh(X,\Lambda)$ (resp. $\Loc(X,\Lambda)$) the inverse $2$-limit of the categories $\Sh(X,\Lambda/\ell^n)$ (resp. $\Loc(X,\Lambda/\ell^n)$) where $n$ ranges over all integers (cf. \ref{chap31.0.1}). Thus the objects of $\mathrm{Loc}^{\otimes}(G,\Lambda)$ are projective systems $(\F_n)_{n \geq 1}$ where $\F_n$ is a multiplicative $\Lambda/ \ell^n$-local system on $G$, such that the transition maps $\F_{n+1} \rightarrow \F_n$ induce isomorphisms
$\F_{n+1} \otimes \Lambda/ \ell^n \rightarrow \F_n,$
for each integer $n$.

 The category $\Sh(X,\Lambda)$ is abelian by (\cite{SGA5}, VI 1.1.3).

\subsection{\label{chap31.0.3}} Assume that $\Lambda$ is a finite subextension of $\mathbb{Q}_{\ell}$ in $C$, and let $X$ be a locally noetherian $k$-scheme. Let $\Lambda_0$ be the ring of integers in $\Lambda$. We denote by $\Sh(X,\Lambda)$ the quotient of the abelian category $\Sh(X,\Lambda_0)$ by the thick subcategory of torsion $\Lambda_0$-sheaves, cf. (\cite{De}, 1.1.1(c)), and by $\Loc(X,\Lambda)$ the essential image of $\Loc(X,\Lambda_0)$ in $\Sh(X,\Lambda)$. In particular, the category $\Sh(X,\Lambda)$ is abelian.

The natural functor $\Sh(X,\Lambda_0) \rightarrow \Sh(X,\Lambda)$ is essentially surjective, and will be denoted by $\otimes_{\Lambda_0} \Lambda$. If $X$ is noetherian then the natural homomorphism
$$
\Hom_{\Lambda_0}(\F,\G) \otimes_{\Lambda_0} \Lambda \rightarrow \Hom_{\Lambda}(\F \otimes_{\Lambda_0} \Lambda,\G \otimes_{\Lambda_0} \Lambda),
$$
is an isomorphism for any objects $\F$ and $\G$ of $\Sh(X,\Lambda_0)$.

\subsection{\label{chap31.0.4}} Let $X$ be a locally noetherian $k$-scheme. We denote by $\Sh(X,\Lambda)$ (resp. $\Loc(X,\Lambda)$) the $2$-colimit of the categories $\Sh(X,\Lambda_0)$, where $\Lambda_0$ ranges over admissible $\ell$-adic subrings of $\Lambda$ (cf. \ref{chap30.0.0.1}). For $\Lambda = C$, this coincides with (\cite{De}, 1.1.1(d)). We will simply refer to objects of $\Sh(X,\Lambda)$ (resp. $\Loc(X,\Lambda)$) as \textit{$\Lambda$-sheaves on $X$} (resp. \textit{$\Lambda$-local systems on $X$}). The category $\Sh(X,\Lambda)$ is abelian, since it is a filtered $2$-colimit of abelian categories.

\subsection{\label{chap31.0.5}} Let $X$ be a locally noetherian $k$-scheme. Let $\Lambda \rightarrow \Lambda'$ be a continuous homomorphism of $\ell$-adic coefficient rings. If $\Lambda$ and $\Lambda'$ are admissible, we define a functor
$$
\otimes_{\Lambda} \Lambda' : \Sh(X,\Lambda) \rightarrow \Sh(X,\Lambda'),
$$
as follows:
\begin{enumerate}
\item if $\Lambda$ is finite then so is $\Lambda'$, and $\otimes_{\Lambda} \Lambda'$ is the functor which sends a $\Lambda$-sheaf $\F$ to the tensor product $\F \otimes_{\Lambda} \Lambda'$,
\item if $\Lambda$ and $\Lambda'$ are rings of integers in finite extensions of $\mathbb{Q}_{\ell}$, then $\otimes_{\Lambda} \Lambda'$ is the functor which sends a projective system $(\F_n)_{n \geq 1}$ as in \ref{chap31.0.2} to $(\F_n \otimes_{\Lambda/\ell^n} \Lambda'/\ell^n)_{n \geq 1}$
\item if $\Lambda$ is a ring of integers in a finite extension of $\mathbb{Q}_{\ell}$ and if $\Lambda'$ is finite, then $\otimes_{\Lambda} \Lambda'$ is the functor which sends a projective system $(\F_n)_{n \geq 1}$ as in \ref{chap31.0.2} to $\F_n \otimes_{\Lambda/\ell^n} \Lambda'$, where $n$ is an integer such that $\ell^n$ vanishes in $\Lambda'$.
\item if $\Lambda$ is the ring of integers in a finite extension of $\mathbb{Q}_{\ell}$, and if $\Lambda'$ is a finite extension of $\mathbb{Q}_{\ell}$ with ring of integers $\Lambda_0'$, then $\otimes_{\Lambda} \Lambda'$ is the functor which sends a $\Lambda$-sheaf $\F$ to $(\F \otimes_{\Lambda} \Lambda_0') \otimes_{\Lambda_0'} \Lambda'$, cf. \ref{chap31.0.3}.
\item if $\Lambda$ and $\Lambda'$ are finite extensions of $\mathbb{Q}_{\ell}$, and if $\Lambda_0$ is the ring of integers in $\Lambda$, then $\otimes_{\Lambda} \Lambda'$ is a functor which sends $\F \otimes_{\Lambda_0} \Lambda$ to $\F\otimes_{\Lambda_0} \Lambda'$ for any $\Lambda_0$-sheaf $\F$.
\end{enumerate}

In general we let $\otimes_{\Lambda} \Lambda'$ be the $2$-colimit of the functors
$$
\otimes_{\Lambda_0} \Lambda_0' : \Sh(X,\Lambda_0) \rightarrow \Sh(X,\Lambda_0'),
$$
where $\Lambda_0$ and $\Lambda_0'$ are admissible $\ell$-adic subrings of $\Lambda$ and $\Lambda'$ respectively, such that $\Lambda_0'$ contains the image of $\Lambda_0$ in $\Lambda'$.


%

\subsection{\label{chap31.3}} Let $X$ be a locally noetherian $k$-scheme. Let $k'/k$ be a finite Galois subextension of $\bk$ neutralizing $\mu$ (cf. \ref{chap3neutra}). A \textbf{$\mu$-twisted $\Lambda$-sheaf on $X$}, is a pair $\F = (\F_{k'}, (\rho_{\F}(g))_{g \in \Gal(k'/k)})$, where $\F_{k'}$ is $\Lambda$-sheaf on $X_{k'}$ (cf. \ref{chap31.0.4}), and $\rho_{\F}(g) : g^{-1} \F_{k'} \rightarrow \F_{k'}$ is an isomorphism for each $g$ in $\Gal(k'/k)$, such that the diagram
\begin{center}
 \begin{tikzpicture}[scale=1]

\node (A) at (0,0) {$g^{-1} h^{-1} \F_{k'} $};
\node (B) at (3,0) {$g^{-1} \F_{k'}$};
\node (C) at (6,0) {$\F_{k'}$};

\path[->,font=\scriptsize]
(A) edge node[above]{$g^{-1} \rho_{\F}(h)$} (B)
(B) edge node[above]{$\rho_{\F}(g)$} (C)
(A) edge[bend right] node[above]{$ \mu(g,h) \rho_{\F}(gh)$} (C);
\end{tikzpicture} 
\end{center}
is commutative for any $g,h$ in $\Gal(k'/k)$. In particular, the endomorphism $\rho_{\F}(1)$ is the identity of $\F$. If $\F$ and $\G$ are $\mu$-twisted $\Lambda$-sheaves on $X$, a morphism from $\F$ to $\G$ is a morphism $f : \F_{k'} \rightarrow \G_{k'}$ in $\Sh(X_{k'},\Lambda)$ such that $f \circ \rho_{\F}(g) = \rho_{\G}(g) \circ (g^{-1}f)$ for any $g$ in $\Gal(k'/k)$.

\begin{rema} Since the action of $\Gal(k'/k)$ on $X_{k'}$ is a right action, we have $(gh)^{-1} = g^{-1} h^{-1}$ on $\Sh(X_{k'},\Lambda)$.
\end{rema}

If $k''/k'$ is a finite extension contained in $\bk$ such that $k''$ is a Galois extension of $k$, then $k''/k$ neutralizes $\mu$ as well. If $\pi$ is the projection from $\Gal(k''/k)$ to $\Gal(k'/k)$, then, by descent along the $\Gal(k''/k')$-torsor $X_{k''} \rightarrow X_{k'}$, the functor
$$
(\F_{k'}, (\rho_{\F}(g))_{g \in \Gal(k'/k)}) \mapsto (\F_{k' | X_{k''}}, (\rho_{\F}(\pi(g))_{g \in \Gal(k''/k)}),
$$
is an equivalence between the corresponding categories of $\mu$-twisted $\Lambda$-sheaves on $X$. We denote by $\Sh(X,\mu,\Lambda)$ the $2$-limit of these categories along the filtered set of Galois extension of $k$ contained in $\bk$ which neutralizes $\mu$. The category $\Sh(X,\mu,\Lambda)$ is abelian.

If $\Lambda \rightarrow \Lambda'$ is a continuous homomorphism of $\ell$-adic coefficient rings, the natural functor
$$
\Sh(X, \mu, \Lambda) \rightarrow \Sh(X, \mu, \Lambda'),
$$
will be denoted $\otimes_{\Lambda} \Lambda'$.

\subsection{\label{chap31.3.0}} Let $k''$ be a finite extension of $k$, and let $X$ be a locally noetherian $k''$-scheme. Let $\Sigma$ be the set of morphisms of $k$-algebras from $k''$ to $\bk$. Let $k'/k$ be a finite Galois subextension of $\bk$ neutralizing $\mu$ (cf. \ref{chap3neutra}) and containing the image of any element of $\Sigma$. We thus have a decomposition
$$
X \otimes_{k} k' = \coprod_{\iota \in \Sigma} X_{\iota},
$$
where $X_{\iota}$ is the $k'$-scheme $X \otimes_{k'',\iota} k'$. In particular, a $\Lambda$-sheaf on $X_{k'}$ can be considered as a collection $(\F_{\iota})_{\iota \in \Sigma}$, where $\F_{\iota}$ is a $\Lambda$-sheaf on $X_{\iota}$ for each $\iota$. Moreover, if $\F = (\F_{k'}, (\rho_{\F}(g))_{g \in \Gal(k'/k)})$ is a $\mu$-twisted $\Lambda$-sheaf on the $k$-scheme $X$, then $\F_{k'}$ is a $\Lambda$-sheaf on $X_{k'}$, and its component $(\F_{k'})_{\iota}$ is a $\mu_{|\Gal(\bk/\iota(k''))}$-twisted $\Lambda$-sheaf on the $k''$-scheme $X$.

\begin{prop}\label{chap3sanity} Let $X$ and $k''$ be as in \ref{chap31.3.0}. Let $\iota$ be an element of $\Sigma$, let $\Gal(\iota)$ be the Galois group of the extension $\iota : k'' \rightarrow \bk$, and let $\mu_{k''}$ be the restriction of $\mu$ to $G_{k''}$. The functor which sends a $\mu$-twisted $\Lambda$-sheaf $\F = (\F_{k'}, (\rho_{\F}(g))_{g \in \Gal(k'/k)})$ on the $k$-scheme $X$ to the $\mu_{k''}$-twisted $\Lambda$-sheaf $((\F_{k'})_{\iota}, (\rho_{\F}(g)_{\iota})_{g \in \Gal(\iota)})$ on the $k''$-scheme $X$, is an equivalence of categories.
\end{prop}

Indeed, a quasi-inverse to the functor from Proposition \ref{chap3sanity} can be described as follows. If $(\F_{\iota},\rho_{\F}(g))_{g \in \Gal(\iota)}$ is a $\mu_{k''}$-twisted $\Lambda$-sheaf on the $k''$-scheme $X$, then for each $\iota'$ in $\Sigma$, we define $\F_{\iota'}$ to be the sub-$\Lambda$-sheaf of 
$$
\prod_{\substack{g \in \Gal(k'/k) \\ g \iota = \iota' }} g^{-1} \F_{\iota},
$$
on $X_{\iota'}$, consisting of sections $(s_g)_g$ such that for any $g$ in $\Gal(k'/k)$ with $g\iota = \iota'$ and any element $h$ of $\Gal(\iota)$, we have
$$
s_{gh} = \mu(g,h) \rho_{\F}(h)^{-1} s_g.
$$
The collection $(\F_{\iota'})_{\iota' \in \Sigma}$ yields a $\Lambda$-sheaf on $X_{k'}$ which is naturally endowed with a structure $(\rho_{\F}(g))_{g \in \Gal(k'/k)}$ of $\mu$-twisted $\Lambda$-sheaf on the $k$-scheme $X$. For each $h$ in $\Gal(k'/k)$ and for each $\iota'$ in $\Sigma$, the morphism $\rho_{\F}(h)_{\iota'}$ sends a section $(s_g)_{g \iota = h^{-1} \iota'}$ of $h^{-1} \F_{h^{-1} \iota'}$ to the section $(\mu(h,h^{-1}g) s_{h^{-1}g})_{g \iota = \iota'}$ of $\F_{\iota'}$.

\subsection{\label{chap31.3.1}} Assume that $\Lambda$ is the ring of integers in a finite subextension of $\mathbb{Q}_{\ell}$ in $C$, let $X$ be a locally noetherian $k$-scheme. Then the natural functor
$$
\Sh(X,\mu,\Lambda) \rightarrow 2\text{-}\mathrm{lim}_{n} \ \Sh(X,\mu,\Lambda/ \ell^n),
$$
where $n$ runs over all positive integers, is an equivalence of categories. Indeed, if $k'/k$ is a finite Galois extension contained in $\bk$ which neutralizes $\mu$, then the natural functor
$$
\Sh(X_{k'},\Lambda) \rightarrow 2\text{-}\mathrm{lim}_{n} \ \Sh(X_{k'},\Lambda/ \ell^n),
$$
is an equivalence of categories by definition, cf. \ref{chap31.0.2}. Moreover, if $\F = (\F_n)_{n \geq 1}$ is an object of $\Sh(X_{k'},\Lambda) $, then we have
$$
\Hom_{\Lambda}(g^{-1} \F,\F) = \lim \Hom_{\Lambda/\ell^n}(g^{-1} \F_n,\F_n),
$$
for any $g$ in $\Gal(k'/k)$, hence a structure of $\mu$-twisted $\Lambda$-sheaf on $\F$ amounts to a compatible system of structures of $\mu$-twisted $\Lambda$-sheaves on each $\F_n$.
%

\subsection{\label{chap31.3.2}} Assume that $\Lambda$ is a finite subextension of $\mathbb{Q}_{\ell}$ in $C$, with ring of integers $\Lambda_0 \subseteq \Lambda$, let $X$ be a locally noetherian $k$-scheme. Then $\mu$ takes its values in $\Lambda_0$, and the natural functor
\begin{align}\label{chap3functor}
\Sh(X,\mu,\Lambda_0) \rightarrow \Sh(X,\mu,\Lambda)
\end{align}
induces an equivalence from the quotient of $\Sh(X,\mu,\Lambda_0)$ by its subcategory of torsion objects, to the category of $\mu$-twisted $\Lambda$-sheaves on $X$. Unlike \ref{chap31.3.1}, this statement is not entirely tautological.

By glueing, we can assume (and we do) that $X$ is noetherian. Let $k'/k$ be a finite Galois subextension of $\bk$, with Galois group $G = \Gal(k'/k)$, which neutralizes $\mu$ (cf. \ref{chap3neutra}). For any objects $\F$ and $\G$ of $\Sh(X,\mu,\Lambda_0)$, the natural homomorphism
$$
\Hom_{\Lambda_0}(\F_{k'},\G_{k'}) \otimes_{\Lambda_0} \Lambda \rightarrow \Hom_{ \Lambda}( \F_{k'} \otimes_{\Lambda_0} \Lambda, \G_{k'} \otimes_{\Lambda_0} \Lambda),
$$
is an isomorphism, cf. \ref{chap31.0.3}. Moreover, we have an exact sequence
$$
0 \rightarrow \Hom_{\Lambda_0}(\F,\G) \rightarrow \Hom_{\Lambda_0}(\F_{k'},\G_{k'}) \xrightarrow[]{f \mapsto (f \rho_{\F}(g) - \rho_{\G}(g) (g^{-1}f))_g} \prod_{g \in \Gal(k'/k)} \Hom_{\Lambda_0}(g^{-1}\F_{k'},\G_{k'}).
$$
By flatness of the ring homomorphism $\Lambda_0 \rightarrow \Lambda$, we deduce that the natural homomorphism
$$
\Hom_{\Lambda_0}(\F,\G) \otimes_{\Lambda_0} \Lambda \rightarrow \Hom_{ \Lambda}( \F \otimes_{\Lambda_0} \Lambda, \G \otimes_{\Lambda_0} \Lambda),
$$
is an isomorphism as well, hence the full faithfullness of the functor (\ref{chap3functor}).

 It remains to prove that the functor (\ref{chap3functor}) is essentially surjective. Let $(\F_{k'}, (\rho_{\F}(g))_{g \in G})$ be a pair representing an object $\F$ of $\Sh(X,\mu,\Lambda)$. Let $\Hl$ be a torsion free $\Lambda_0$-sheaf on $X_{k'}$ such that $\F_{k'}$ is isomorphic to $\Hl \otimes_{\Lambda_0} \Lambda$. For a sufficiently large integer $n$, we have for each element $g$ of $G$ a homomorphism $\widetilde{\rho}_{\F}(g)$ from $g^{-1} \Hl$ to $\Hl$, which induces the homomorphism $\ell^n \rho_{\F}(g)$ from $g^{-1} \F_{k'}$ to $\F_{k'}$. Let us consider the homomorphism
 $$
\theta : \bigoplus_{h \in G} h^{-1} \Hl \rightarrow \Hl,
$$	
given by $\widetilde{\rho}_{\F}(h)$ on the component $h^{-1} \Hl$. We endow the source $\Hl' = \bigoplus_{h \in G} h^{-1} \Hl$ of $\theta$ with a structure of $\mu$-twisted $\Lambda_0$-sheaf by equipping it for each $g$ in $G$ with the isomorphism
$$
 \rho_{\Hl'}(g) : g^{-1} \Hl' \rightarrow \Hl',
$$
which sends a tuple $(x_h)_{h \in G}$, where $x_h$ belongs to $g^{-1} h^{-1} \Hl$, to the tuple $(\mu(g,g^{-1}h)x_{g^{-1}h})_{h \in G}$. For each $g$ in $G$, the diagram
\begin{center}
 \begin{tikzpicture}[scale=1]

\node (A) at (0,2) {$g^{-1} \Hl' $};
\node (B) at (3,2) {$g^{-1} \Hl$};
\node (C) at (3,0) {$\Hl$};
\node (D) at (0,0) {$\Hl'$};

\path[->,font=\scriptsize]
(A) edge node[above]{$g^{-1} \theta$} (B)
(B) edge node[right]{$\widetilde{\rho}_{\F}(g)$} (C)
(A) edge node[right]{$\ell^n \rho_{\Hl'}(g) $} (D)
(D) edge node[above]{$\theta$} (C);
\end{tikzpicture} 
\end{center}
is commutative. Since $\mathcal{H}$ is torsion free, this implies that $ \rho_{\Hl'}(g)$ induces an isomorphism from the subobject $g^{-1} \ker(\theta)$ of $ g^{-1} \Hl'$ to the subobject $\ker(\theta)$ of $\Hl'$. Thus the kernel of $\theta$ is a $\mu$-twisted $\Lambda_0$-subsheaf of $\Hl'$. Consequently, the image of $\theta$ is a $\mu$-twisted $\Lambda_0$-sheaf as well. Since $\mathrm{Im}( \theta) \otimes_{\Lambda_0} \Lambda$ is isomorphic to $\F$, this proves the essential surjectivity of the functor (\ref{chap3functor}).

\subsection{\label{chap31.3.3}} Let $X$ be a locally noetherian $k$-scheme. Then the natural functor
$$
2 \text{-} \mathrm{colim}_{\Lambda_0} \ \Sh(X,\mu,\Lambda_0) \rightarrow \Sh(X,\mu,\Lambda),
$$
where $\Lambda_0$ runs over all admissible $\ell$-adic subrings of $\Lambda$ containing the image of $\mu$, is an equivalence of categories. 
 Indeed, if $k'/k$ is a finite Galois extension contained in $\bk$ which neutralizes $\mu$, then the natural functor
 $$
2 \text{-} \mathrm{colim}_{\Lambda_0} \ \Sh(X_{k'},\Lambda_0) \rightarrow \Sh(X_{k'},\Lambda),
$$
is an equivalence of categories by definition, cf. \ref{chap31.0.4}. Moreover, if $\F$ is an object of $\Sh(X_{k'},\Lambda_0) $, then the natural homomorphism
$$
\colim_{\Lambda_1} \Hom_{\Lambda_1}(g^{-1} \F \otimes_{\Lambda_0} \Lambda_1,\F \otimes_{\Lambda_0} \Lambda_1) \rightarrow \Hom_{\Lambda}(g^{-1} \F \otimes_{\Lambda_0} \Lambda,\F \otimes_{\Lambda_0} \Lambda),
$$
is an isomorphism for any $g$ in $\Gal(k'/k)$, where $\Lambda_1$ runs over all admissible $\ell$-adic subrings of $\Lambda$ containing $\Lambda_0$ and the image of $\mu$. Since $\Gal(k'/k)$ is finite, we obtain that a structure of $\mu$-twisted $\Lambda$-sheaf on $\F \otimes_{\Lambda_0} \Lambda$ amounts to a structure of $\mu$-twisted $\Lambda_1$-sheaf on $\F \otimes_{\Lambda_0} \Lambda_1$, for a large enough admissible $\ell$-adic subring $\Lambda_1$ of $\Lambda$.

\subsection{\label{chap31.3.4}} Let $f : X \rightarrow Y$ be a separated morphism of $k$-schemes of finite type. If $\F$ is a $\mu$-twisted $\Lambda$-sheaf on $X$, given by a tuple $(\F_{k'}, (\rho_{\F}(g))_{g \in \Gal(k'/k)})$ (cf. \ref{chap31.3}), then $R^\nu f_* \F_{k'}$ and $R^\nu f_! \F_{k'}$ are $\Lambda$-sheaves on $Y_{k'}$ for each integer $\nu$ by (\cite{SGA412}, 7.1.1), and the isomorphisms
$$
\rho_{\F}(g) : \F_{k'} \rightarrow g_* \F_{k'},
$$
on $X_{k'}$ yield by functoriality isomorphisms
\begin{align*}
R^\nu f_* \F_{k'} &\rightarrow R^\nu f_* g_*\F_{k'} = g_* R^\nu f_* \F_{k'} \\
R^\nu f_! \F_{k'} &\rightarrow R^\nu f_! g_*\F_{k'} = g_* R^\nu f_! \F_{k'},
\end{align*}
so that we obtain structures of $\mu$-twisted $\Lambda$-sheaves on $R^\nu f_* \F_{k'}$ and $R^\nu f_! \F_{k'}$. We denote by $R^\nu f_* \F$ and $R^\nu f_! \F$ the resulting $\mu$-twisted $\Lambda$-sheaves. In particular, by taking $Y = \Spec(k)$, we obtain a structure of $\mu$-twisted representation of $G_k$ on the cohomology groups $H_c^{\nu}(X_{\bk}, \F_{\bk})$.

\subsection{\label{chap31.4.0}} Let $X$ be a $k$-scheme. For any geometric point $\bar{x}$ of $X$, we denote by $\ev_{\bar{x}}$ the functor which to a finite \'etale $X$-scheme $Y$ associates the set $Y(\bar{x})$ of $\bar{x}$-points of $Y$. For a couple $(\overline{x}_0,\overline{x}_1)$ of geometric points of $X$, we denote by $\pi_1(X,\overline{x}_0,\overline{x}_1)$ the set of isomorphisms from $\ev_{\overline{x}_0}$ to $\ev_{\overline{x}_1}$, endowed with the coarsest topology such that for any finite \'etale $X$-scheme $Y$, the map
$$
\pi_1(X,\overline{x}_0,\overline{x}_1) \rightarrow \Hom(Y(\overline{x}_0),Y(\overline{x}_1) ),
$$
is continuous, where the target is endowed with the discrete topology. Thus $\pi_1(X,\overline{x}_0,\overline{x}_1)$ is a profinite topological space. For a triple $(\overline{x}_0,\overline{x}_1, \overline{x}_2)$ of geometric points of $X$, the composition induces a continuous map
$$
\pi_1(X,\overline{x}_1,\overline{x}_2) \times \pi_1(X,\overline{x}_0,\overline{x}_1) \rightarrow \pi_1(X,\overline{x}_0,\overline{x}_2),
$$
which satisfies the natural associativity condition. In particular, $\pi_1(X,\overline{x}_0,\overline{x}_0)$ is a profinite group, which will be simply denoted by $\pi_1(X,\overline{x}_0)$.

\subsection{\label{chap31.4.3}} Let $f : X \rightarrow Y$ be a morphism of $k$-schemes, and let $(\overline{x}_0,\overline{x}_1)$ be a couple of geometric points of $X$. We have natural isomorphisms $\ev_{f(\overline{x}_0)} \cong \ev_{\overline{x}_0} \circ f^{-1}$ and $\ev_{f(\overline{x}_1)} \cong \ev_{\overline{x}_1} \circ f^{-1}$, whence precomposition with $f^{-1}$ induces a continuous map
$$
f_* : \pi_1(X,\overline{x}_0,\overline{x}_1) \rightarrow \pi_1(X,f(\overline{x}_0),f(\overline{x}_1)),
$$
which is compatible with the composition laws. In particular, the formation of $\pi_1(X,\overline{x}_0,\overline{x}_1)$ (resp. $\pi_1(X,\overline{x}_0)$) is functorial in the triple $(X,\overline{x}_0,\overline{x}_1)$ (resp. in the couple $(X,\overline{x}_0)$).

\subsection{\label{chap31.4.1}} Let $X$ be a locally noetherian connected $k$-scheme, and let $\bar{x}$ be a geometric point of $X$. Then the fiber functor $\F \mapsto \F_{\bar{x}}$ is an equivalence of categories from $\Loc(X,\Lambda)$ to the category $\Rep_{\Lambda}(\pi_1(X,\overline{x}),1)$ of $\Lambda$-admissible representations of the profinite group $\pi_1(X,\overline{x})$ (cf. \ref{chap30.2}). This statement reduces to the case where $\Lambda$ is finite, which in turn follows from (\cite{SGA1}, V.7).

\subsection{\label{chap31.4.2}} Let $X$ be a locally noetherian geometrically connected $k$-scheme, let $k'$ be a finite Galois extension of $k$ and let $\bar{x}$ be a geometric point of $X_{k'}$. Let $u : X_{k'} \rightarrow X$ be the natural projection. By (\cite{SGA1}, IX 6.1), we have an exact sequence
\begin{align}\label{chap3exactseqs}
1 \rightarrow \pi_1(X_{k'},\bar{x}) \xrightarrow[]{u_*} \pi_1(X,u(\bar{x})) \xrightarrow[]{r} \Gal(k' / k) \rightarrow 1,
\end{align}
which we now proceed to describe in greater details. For each $g$ in $\Gal(k'/k)$, the geometric points $\bar{x}$ and $\bar{x} g$ have the same image by $u$, so that we obtain a continuous map
$$
u_* : \pi_1(X_{k'},\bar{x}, \bar{x}g) \rightarrow \pi_1(X,u(\bar{x}), u(\bar{x}g)) = \pi_1(X,u(\bar{x})).
$$
The collection of these continuous maps yields a homeomorphism
\begin{align}\label{chap3coproduit}
\coprod_{g \in \Gal(k'/k)} \pi_1(X_{k'},\bar{x}, \bar{x}g) \xrightarrow[]{u_*} \pi_1(X,u(\bar{x})),
\end{align}
whose composition with the homomorphism $r : \pi_1(X,u(\bar{x})) \rightarrow \Gal(k' / k)$ maps $\pi_1(X_{k'},\bar{x}, \bar{x}g)$ to $g$ for each $g$ in $\Gal(k'/k)$. The homeomorphism above can be promoted to an isomorphism of profinite groups if we endow its source with the following group law: for any elements $\alpha$ and $\beta$ of $\pi_1(X_{k'},\bar{x}, \bar{x}g)$ and $\pi_1(X_{k'},\bar{x}, \bar{x}h)$ respectively, we set $\alpha \cdot \beta = h_* \alpha \circ \beta$.

\subsection{\label{chap31.4.4}} Let $X$ be a locally noetherian geometrically connected $k$-scheme, let $k'/k$ be a finite Galois subextension of $\bk$ neutralizing $\mu$ (cf. \ref{chap3neutra}) and let $\bar{x}$ be a geometric point of $X_{k'}$. We abusively denote by $\mu$ the pullback of $\mu$ by the homomorphism $r : \pi_1(X,u(\bar{x})) \rightarrow \Gal(k' / k)$ from \ref{chap31.4.2}.

Let $\Loc(X,\mu, \Lambda)$ be the full subcategory of $\Sh(X,\mu, \Lambda)$ consisting of objects $\F$ such that $\F_{k'}$ is a $\Lambda$-local system on $X_{k'}$. Then the fiber functor $\F \mapsto \F_{\bar{x}}$ is an equivalence of categories from $\Loc(X,\mu, \Lambda)$ to the category $\Rep_{\Lambda}(\pi_1(X,u(\overline{x})),\mu)$ of $\Lambda$-admissible representations of the $\Lambda$-twisted profinite group $(\pi_1(X,u(\overline{x})),\mu)$ (cf. \ref{chap30.1}, \ref{chap30.2}). 

Indeed, by \ref{chap31.4.1} the functor $\F_{k'} \mapsto \F_{k',\bar{x}}$ realizes an equivalence of categories from $\Loc(X,\Lambda)$ to the category of $\Lambda$-admissible representations of the profinite group $(\pi_1(X_{k'},\overline{x}),1)$. Moreover, for each $g$ in $\Gal(k'/k)$, the morphism
\begin{align*}
\Hom_{\Lambda}(g^{-1} \F_{k'},\F_{k'}) &\rightarrow \Hom_{\pi_1(X_{k'}, \bar{x})}(r^{-1}(g),\Aut_{\Lambda}(\F_{k',\bar{x}})) \\
\rho &\mapsto (\alpha \mapsto \rho_{\bar{x}} \circ u_*^{-1}(\alpha)),
\end{align*}
where $u_*^{-1}$ is the inverse of the isomorphism from (\ref{chap3coproduit}), realizes an isomorphism onto the set of left and right $\pi_1(X_{k'}, \bar{x})$-equivariant maps from $r^{-1}(g)$ to $\Aut_{\Lambda}(\F_{k',\bar{x}})$. Thus, structures of $\mu$-twisted $\Lambda$-sheaf on $\F_{k'}$ correspond bijectively to structures of $\Lambda$-admissible representations of $(\pi_1(X,u(\overline{x})),\mu)$ on $\F_{k',\bar{x}}$, hence the result.
\begin{exemple}\label{chap3exampletwist} For $X = \Spec(k)$, the category $\Sh(X,\mu,\Lambda)$ is equivalent to the category of $\Lambda$-admissible representations of $(G_k,\mu)$ (cf. \ref{chap30.2}). 
\end{exemple}

\subsection{\label{chap31.4.5}} Let $Y$ be a locally noetherian geometrically connected $k$-scheme, and let $f : X \rightarrow Y$ be a finite \'etale morphism, where $X$ is a geometrically connected $k''$-scheme for some finite extension $k''$ of $k$ contained in $\bk$. Let $k'/k$ be a finite Galois subextension of $\bk$ neutralizing $\mu$ and containing $k''$ (cf. \ref{chap3neutra}), let $\bar{x}$ be a geometric point of $X_{k'}$ and let $u : X_{k'} \rightarrow X$ be the natural projection. The homomorphism
$$
f_* : \pi_1(X,u(\overline{x})) \rightarrow \pi_1(Y,fu(\overline{x})),
$$
is injective and its image is an open subgroup of finite index in $\pi_1(Y,fu(\overline{x}))$. Then the following diagram is commutative.
\begin{center}
 \begin{tikzpicture}[scale=1]

\node (A) at (0,2) {$\Loc(X,\mu,\Lambda)$};
\node (B) at (5,2) {$\Rep_{\Lambda}(\pi_1(X,u(\overline{x})),\mu)$};
\node (C) at (5,0) {$\Rep_{\Lambda}( \pi_1(Y,fu(\overline{x})),\mu)$};
\node (D) at (0,0) {$\Loc(Y,\mu,\Lambda)$};

\path[->,font=\scriptsize]
(A) edge node[above]{$\F \mapsto \F_{u(\overline{x})}$} (B)
(B) edge node[right]{$\Ind_{\pi_1(X,u(\overline{x}))}^{ \pi_1(Y,fu(\overline{x}))}$} (C)
(A) edge node[right]{$f_*$} (D)
(D) edge node[above]{$\F \mapsto \F_{fu(\overline{x})}$} (C);
\end{tikzpicture} 
\end{center}
In this diagram, the left vertical arrow is defined in \ref{chap31.4.0}, while the right vertical arrow is defined in \ref{chap30.3.1}. The horizontal arrows are equivalences of categories by \ref{chap31.4.4}.

\subsection{\label{chap31.4.8}} Assume $\Lambda$ is a finite subextension of $\mathbb{Q}_{\ell}$ in $C$, with ring of integers $\Lambda_0 \subseteq \Lambda$. We denote by $\mathfrak{m}_0$ the maximal ideal of $\Lambda_0$. Let $X$ be a noetherian $k$-scheme, and let $K_0(X,\mu,\Lambda)$ be the Grothendieck group of virtual classes of $\mu$-twisted $\Lambda$ sheaves on $X$. 

\begin{prop}\label{1.4.8.1} Let $\F$ be a $\mu$-twisted $\Lambda$-sheaf on $X$. If $\F_0$ and $\G_0$ are torsion free $\mu$-twisted $\Lambda_0$-sheaves on $X$ with image in $\mathrm{Sh}(X, \mu, \Lambda)$ isomorphic to $\F$, then the classes of $\F_0/\mathfrak{m}_0 \F_0 $ and $\G_0/\mathfrak{m}_0 \G_0 $ in $K_0(X,\mu,\Lambda_0/\mathfrak{m}_0)$ coincide.
\end{prop}

This is similar to (\cite{S}, 15.2): one can assume that we have an injective homomorphism $\varphi : \F_0 \rightarrow \G_0$ with torsion cokernel, and one then considers the $\mu$-twisted $\Lambda_0$-sheaves $\mathcal{H}_j = \varphi(\F_0) + \mathfrak{m}_0^j \G_0$. We have $\G_0 = \mathcal{H}_0$ and the morphism $\varphi$ induces an isomorphism from $\F_0$ to $\mathcal{H}_j$ for $j$ large enough. It is therefore enough to prove that $\mathcal{H}_j$ and $\mathcal{H}_{j+1}$ yield the same class in $K_0(X,\mu,\Lambda_0/\mathfrak{m}_0)$ for any integer $j$, which in turn follows from the exact sequence
$$
0 \rightarrow \mathcal{H}_{j}/ \mathcal{H}_{j+1} \xrightarrow[]{\pi} \mathcal{H}_{j+1} / \mathfrak{m}_0 \mathcal{H}_{j+1} \rightarrow \mathcal{H}_j / \mathfrak{m}_0 \mathcal{H}_j \rightarrow \mathcal{H}_j/ \mathcal{H}_{j+1} \rightarrow 0
$$
where $\pi$ is a generator of the ideal $\mathfrak{m}_0$ of $\Lambda_0$.

\begin{defi}\label{1.4.8.2} The decomposition homomorphism
$$
d : K_0(X,\mu,\Lambda) \rightarrow K_0(X,\mu,\Lambda_0/\mathfrak{m}_0),
$$
is the unique homomorphism sending the class of $\F \otimes_{\Lambda_0} \Lambda$ to the class of $\F \otimes_{\Lambda_0} \Lambda_0/\mathfrak{m}_0$ for any torsion free $\mu$-twisted $\Lambda_0$-sheaf on $X$. This is well defined by Proposition \ref{1.4.8.1}.
\end{defi}

\begin{prop}\label{1.4.8.3} Assume that $\mu$ is the Teichm\"{u}ller lift of unitary $\Lambda_0/\mathfrak{m}_0$-admissible mutiplier $\mu_0$ on $G_k$. Then for any noetherian $k$-scheme $X$ the decomposition homomorphism
$$
d : K_0(X,\mu,\Lambda) \rightarrow K_0(X,\mu_0,\Lambda_0/\mathfrak{m}_0),
$$
cf. \ref{1.4.8.2}, is surjective.
\end{prop}

We proceed by noetherian induction on $X$. Let $\F$ be a $\mu$-twisted $\Lambda_0/\mathfrak{m}_0$-sheaf on $X$. There exists a non empty open subscheme $j : U \rightarrow X$ of $X$ such that $j^{-1} \F$ belongs to $\Loc(U,\mu,\Lambda_0/\mathfrak{m}_0)$. By shrinking $U$ if necessary we can assume (and we do) that $U$ is a geometrically connected $k'$-scheme for some finite extension $k'$ of $k$. Let $Z$ be the complement of $U$ in $X$, endowed with its reduced scheme structure, and let $i : Z \rightarrow X$ be corresponding closed immersion. We then have
$$
[\F] = [i_* i^{-1} \F] + [j_! j^{-1} \F],
$$
in $K_0(X,\mu,\Lambda_0/\mathfrak{m}_0)$. By \ref{chap3sanity}, \ref{chap31.4.4}, and \ref{CGTRES2}, there exists an object $\G$ of $\Loc(U,\mu,\Lambda)$ such that $d([\G]) =  j^{-1} \F$. By noetherian induction, there exists a $\mu$-twisted $\Lambda$-sheaf $\mathcal{H}$ on $Z$ such that $d([\mathcal{H}]) =  i^{-1} \F$. We then have
$$
 [\F] = d( [i_* \mathcal{H}] + [j_! \G] ),
$$
hence the surjectivity of $d$.

\subsection{\label{chap31.4.6}} Let $s \rightarrow \Spec(k)$ be a finite extension of $k$, and let us fix a $k$-morphism $\overline{s} : \Spec(\bk) \rightarrow s$, corresponding to a $k$-linear embedding of $k(s)$ in $\bk$, so that the Galois group $G_{s} = \Gal(\bk/k(s))$ can be considered as an open subgroup of finite index in $G_k$. We still denote by $\mu$ the restriction of $\mu$ to $G_{s}$. We denote by $\delta_{s/k}$ the $\Lambda$-admissible character of rank $1$ of $G_k$ given by 
$$
\delta_{s/k} = \det \left( \Ind_{G_{s}}^{G_k} \Lambda \right).
$$
Thus $\delta_{s/k}$ is the signature character associated the left action of $G_k$ on the finite set $G_k / G_{s}$.

\begin{defi}\label{chap3verlagerung} Let $V$ be a $\Lambda$-admissible representation of $(G_{s},\mu)$. The \emph{transfer} or \emph{verlagerung} of $V$ with respect to the extension $k(s)/k$ is the $\Lambda$-admissible map $\Ver_{s/k}(V)$ from $G_k$ to $\Lambda^{\times}$ defined by
$$
\Ver_{s/k}(V) = \delta_{s/k}^{- \rk(V)} \det \left( \Ind_{G_{s}}^{G_k} V \right),
$$
where the induction is defined in \ref{chap30.3.1}.
\end{defi}

\begin{prop}\label{chap3transferhomo} Let $V$ be a $\Lambda$-admissible representation of $(G_{s},\mu)$. Then the determinant $\det(V)$ is a $\Lambda$-admissible representation of $(G_{s},\mu^{\rk(V)})$ and we have 
$$
\Ver_{s/k}(V) = \Ver_{s/k}(\det(V)).
$$
Moreover, if $\chi_1$ and $\chi_2$ are $\Lambda$-admissible representations of rank $1$ of $(G_{s},\mu_1)$ and $(G_s, \mu_2)$ respectively, for some multipliers $\nu_1,\nu_2$ on $G_k$, then
$$
\Ver_{s/k}(\chi_1 \chi_2) = \Ver_{s/k}(\chi_1) \Ver_{s/k}(\chi_2).
$$
\end{prop}

Indeed, let $(t_i)_{i \in I}$ be a family of $G_{s}$-left cosets representatives in $G_k$. Then we have a splitting (cf. \ref{chap30.3.1})
$$
\Ind_{G_{s}}^{G_k} V = \bigoplus_{i \in I} [t_i] V.
$$
Let $g$ be an element of $G_k$ and let us write $g t_i = t_{\sigma(g)(i)} h_{g,i}$ for some bijection $\sigma(g)$ of $I$ onto itself, and some element $h_{g,i}$ of $G_{s}$. For any element $v$ of $V$, we have
\begin{align*}
[g] [t_i] v &= \mu(g,t_i) [g t_i] v \\
&= \mu(g,t_i) \mu(t_{\sigma(g)(i)}, h_{g,i})^{-1} [t_{\sigma(g)(i)}] [h_{g,i}] v.
\end{align*}
Consequently, we have
$$
\det( g \ | \ \Ind_{G_{s}}^{G_k} V) = \mathrm{sgn}(\sigma(g))^{\rk(V)} \prod_{i \in I} \mu(g,t_i)^{\rk(V)} \mu(t_{\sigma(g)(i)}, h_{g,i})^{-\rk(V)} \det( h_{g,i} \ | \ V),
$$
where $\mathrm{sgn}$ is the signature homomorphism. The sign $\mathrm{sgn}(\sigma(g))$ is equal to $\delta_{s/k}(g)$ and thus
$$
\Ver_{s/k}(V)(g) = \prod_{i \in I} \mu(g,t_i)^{\rk(V)} \mu(t_{\sigma(g)(i)}, h_{g,i})^{-\rk(V)} \det( h_{g,i} \ | \ V),
$$
hence the conclusion of Proposition \ref{chap3transferhomo}.

\begin{rema} If $\mu = 1$ then $\det(V)$ is a group homomorphism from $G_s$ to $\Lambda^{\times}$ and the computation above yields
$$
\Ver_{s/k}(V) = \det(V) \circ \mathrm{ver}_{s/k},
$$
where $\mathrm{ver}_{s/k} : G_k^{\mathrm{ab}} \rightarrow G_s^{\mathrm{ab}}$ is the usual transfer homomorphism, cf. (\cite{JPS2}, VII.8 p.122) or (\cite{De73}, Prop. 1.2).
\end{rema}

\begin{cor}\label{chap3compositionverlag} Let $s' \rightarrow s$ be a finite extension. Then we have
\begin{align*}
\delta_{s'/k} &= \delta_{s/k}^{[k(s'):k(s)]} \Ver_{s/k}(\delta_{s'/s}), \\
\Ver_{s'/k} &= \Ver_{s/k} \circ \Ver_{s'/s}.
\end{align*}
\end{cor}

This follows from Proposition \ref{chap3transferhomo} and from the existence of a natural isomorphism
$$
 \Ind_{G_{s'}}^{G_k} V \cong \Ind_{G_{s}}^{G_k} \Ind_{G_{s'}}^{G_s} V.
$$

\subsection{\label{chap31.4.7}} Let $s \rightarrow \Spec(k)$ be a finite extension of $k$, and let us fix a $k$-morphism $\overline{s} : \Spec(\bk) \rightarrow s$, corresponding to a $k$-linear embedding of $k(s)$ in $\bk$, so that the Galois group $G_{s} = \Gal(\bk/k(s))$ can be considered as an open subgroup of finite index in $G_k$. We still denote by $\mu$ the restriction of $\mu$ to $G_{s}$. 

Let $X$ be a separated $s$-scheme of finite type. For any element $t$ of $G_k/G_s$, we denote by $t(\overline{s})$ the composition of $\overline{s}$ with the $k$-automorphism of $\Spec(\bk)$ induced by some element of $G_k$ lifting $t$, and by $^{t}X$ the $\bk$-scheme $X \times_{s, t(\overline{s})} \Spec(\bk )$. Let us consider the product
$$
\overline{X} = \prod_{t \in G_k/ G_s} {}^{t}X.
$$
 For each $t$ in $G_k/ G_s$ and each $g$ in $G_k$, the $k$-automorphism of $\Spec(\bk)$ induced by $g$ yields an isomorphism
$$
{}^{gt}X \rightarrow {}^{t}X.
$$
By taking the product over $t$, we obtain an automorphism of $\overline{X}$, and this produces a right action of $G_k$ on $\overline{X}$. We denote by $p_1$ be projection from $\overline{X}$ onto $X$ through its factor $X_1 = X_{\overline{s}}$.
Let $\F$ be a $\mu$-twisted $\Lambda$-sheaf on $X$. We construct a $\mu^{[k(s):k]}$-twisted $\Lambda$-sheaf $\overline{\F}$ on $\overline{X}$ as follows. For each choice $g_{\bullet} =(g_t)_{t \in G_k/G_s}$ of left $G_s$-cosets representatives in $G_k$, we set 
$$
\overline{\F}_{g_{\bullet}} = \bigotimes_{t \in G_k/G_s} g_t^{-1} p_{1}^{-1} \F.
$$
If $\tilde{g}_{\bullet} =(\tilde{g}_t)_{t \in G_k/G_s}$ is another set of left $G_s$-cosets representatives in $G_k$, then we can write $\tilde{g}_t = g_t h_t$ for some $h_t$ in $G_s$, hence an isomorphism
$$
 g_t^{-1} p_{1}^{-1} \rho_{\F}(h_t) : \tilde{g}_t^{-1} p_{1}^{-1} \F \cong g_t^{-1} h_t^{-1} p_{1}^{-1} \F \rightarrow g_t^{-1} p_{1}^{-1} \F,
$$
which yields an isomorphism $\overline{\F}_{\tilde{g}_{\bullet}} \rightarrow \overline{\F}_{g_{\bullet}}$ by taking the tensor product over $t$. We then define $\overline{\F}$ to be the inverse limit of $\overline{\F}_{g_{\bullet}}$ over all choices $g_{\bullet}$ of left $G_s$-cosets representatives in $G_k$, with respect to these transition maps.

If $g$ is an element of $G_k$, then $g g_{\bullet}$ is also a set of left $G_s$-cosets representatives in $G_k$, hence the transition maps
$$
g^{-1} \overline{\F}_{g_{\bullet}} \cong \overline{\F}_{g g_{\bullet}} \rightarrow \overline{\F}_{g_{\bullet}},
$$
yield an isomorphism $\rho_{\overline{\F}}(g) : g^{-1} \overline{\F} \rightarrow \overline{\F}$, which endows $\overline{\F}$ with a structure of $\mu^{[k(s):k]}$-twisted $\Lambda$-sheaf on $\overline{X}$.

\begin{prop}\label{chap3gradedvector} Let $s,\overline{s},X,\overline{X}$ be as in \ref{chap31.4.7}. Let $\F$ be a $\mu$-twisted $\Lambda$-sheaf on $X$, and let us assume that $R\Gamma_c(X_{\overline{s}},\F)$ is concentrated in a single cohomological degree $\nu$, and that $H_c^{\nu}(X_{\overline{s}},\F)$ is a free $\Lambda$-module of rank $1$. Let $\overline{\F}$ be the $\mu^{[k(s):k]}$-twisted $\Lambda$-sheaf on $\overline{X}$ constructed in \ref{chap31.4.7}. Then $R\Gamma_c(\overline{X},\overline{\F})$ is concentrated in degree $[k(s):k]\nu$, and the $\Lambda$-admissible representation
$$
H_c^{[k(s):k]\nu}(\overline{X},\overline{\F}),
$$
is of rank $1$, isomorphic to $\delta_{s/k}^{\nu} \Ver_{s/k} H_c^{\nu}(X_{\overline{s}},\F)$ (cf. \ref{chap31.4.6}).
\end{prop}

Indeed, for any choice $g_{\bullet} =(g_t)_{t \in G_k/G_s}$ of left $G_s$-cosets representatives in $G_k$, the canonical isomorphism $\overline{\F} \rightarrow \overline{\F}_{g_{\bullet}}$, together with K\"{u}nneth's formula (\cite{SGA4}, XVII 5.4.3), yields that $R\Gamma_c(\overline{X},\overline{\F})$ is concentrated in degree $[k(s):k]\nu$ and that we have a canonical isomorphism
\begin{align}\label{chap3kunneth}
\bigotimes_{t \in G_k/ G_s} H_c^{\nu}({}^t X,g_t^{-1} \F) \rightarrow H_c^{[k(s):k]\nu}(\overline{X},\overline{\F}).
\end{align}
It remains to understand how the action of $G_k$ on the target of this isomorphism translates on its source. Let $g$ be an element of $G_k$. Then we can write $g g_t = g_{\sigma(g)(t)} h_{g,t}$ for some bijection $\sigma(g)$ of $G_k/G_s$ onto itself and $h_{g,t}$ is some element of $G_s$. Since (\ref{chap3kunneth}) is an isomorphism of graded vector spaces (cf. (\cite{SGA4}, XVII 1.1.4), the element $g$ of $G_k$ acts on the source of \ref{chap3kunneth} by $\mathrm{sgn}(\sigma(g))^{\nu} = \delta_{s/k}^{\nu}(g)$ multiplied by 
$$
\prod_{t \in G_k/ G_s} \mu(g,g_t) \mu(g_{\sigma(g)(t)}, h_{g,t})^{-1} \mathrm{Tr}( h_{g,t} \ | \ H_c^{\nu}({}^t X,g_t^{-1} p_{1}^{-1} \F)).
$$
The latter is exactly $\Ver_{s/k} H_c^{\nu}(X_{\overline{s}},\F)$ evaluated at $g$, hence the conclusion of Proposition \ref{chap3gradedvector}.
\section{Gabber-Katz extensions \label{chap3GKsection}}

In this section, we review the theory of Gabber-Katz extensions from \cite{katz}, which will constitute our main tool in order to define geometric local $\varepsilon$-factors in Section \ref{chap3gfar}. 

Let $T$ be the spectrum of a $k$-algebra, which is a henselian discrete valuation ring $\Ow_T$, with maximal ideal $\m$, and whose residue field $\Ow_T/ \m$ is a finite extension of $k$ of degree $\deg(s)$. Let $j : \eta \rightarrow T$ be the generic point of $T$, and let $i : s \rightarrow T$ be its closed point, so that $T$ is canonically an $s$-scheme. We fix a uniformizer $\pi$ of $\Ow_T$, and we abusively denote by $\pi$ as well the morphism
$$
\pi : T \rightarrow \mathbb{A}^1_s,
$$
corresponding to the unique morphism $k(s)[t] \rightarrow \Ow_T$ of $k(s)$-algebras which sends $t$ to $\pi$. We fix a $k$-morphism $\overline{s}$ from $\Spec(\bk)$ to $s$, and a separable closure $\overline{\eta}$ of $\eta_{\bar{s}} = \eta \times_s \overline{s}$. We consider $\overline{\eta}$ as a geometric point of $\mathbb{A}^1_s$ through the morphism $\pi$. We denote by $G_s$ the Galois group of the extension $\bk/k(s)$, and by $G_{\eta}$ the Galois group of the extension $k(\overline{\eta})/k(\eta)$.


\subsection{\label{chap31.6.0}} Let $\widehat{T}$ be the spectrum of the $\m$-adic completion of $\Ow_T$, and let $\left(\mathbb{A}^1_{s}\right)_{(0)} $ (resp. $\widehat{\mathbb{A}^1_{s}}_{(0)}$) be the henselization (resp. spectrum of the completion) of $\mathbb{A}^1_s$ at $0$. Then we have the following commutative diagram.
\begin{center}
 \begin{tikzpicture}[scale=1]

\node (A) at (0,2) {$\widehat{T}$};
\node (B) at (3,2) {$T$};
\node (C) at (0,0) {$\widehat{\mathbb{A}^1_{s}}_{(0)}$};
\node (D) at (3,0) {$\left(\mathbb{A}^1_{s} \right)_{(0)}$};
\path[->,font=\scriptsize]
(A) edge (B)
(B) edge node[right]{$\pi$} (D)
(C) edge (D)
(A) edge node[left]{$\pi$} (C) ;
\end{tikzpicture} 
\end{center}
The left vertical morphism is an isomorphism, while the two horizontal morphisms induce isomorphisms on the corresponding \'etale sites, as it follows for example from Krasner's lemma, cf. (\cite{Stacks}, 09EJ). Thus the right vertical morphism induces an isomorphism on the corresponding \'etale sites.

\subsection{\label{chap31.6}} Following (\cite{katz}, 1.3.1), we define a category of special coverings of $\mathbb{G}_{m,s}$. For any scheme $X$, let $\Fet(X)$ be the category of finite \'etale $X$-schemes. By (\cite{SGA1}, V.7), the functor which sends an object $U \rightarrow \mathbb{G}_{m,s}$ of $\Fet(\mathbb{G}_{m,s})$ to the fiber $U_{\overline{\eta}}$, endowed with the natural action of $\pi_1(\mathbb{G}_{m,s}, \overline{\eta})$, realizes an equivalence of categories from $\Fet(\mathbb{G}_{m,s})$ to the category of finite sets endowed with a continuous left action of $\pi_1(\mathbb{G}_{m,s}, \overline{\eta})$.

\begin{defi}\label{chap3specialcover}(\cite{katz}, 1.3.1) The category $\Fet^{\diamondsuit}(\mathbb{G}_{m,s})$ of \textit{special} finite \'etale $\mathbb{G}_{m,s}$-schemes is the full subcategory of $\Fet(\mathbb{G}_{m,s})$ (cf. \ref{chap31.6}) whose objects are the finite \'etale morphisms $f : U \rightarrow \mathbb{G}_{m,s}$ such that:
\begin{enumerate}
\item the morphism $f$ is tamely ramified above $\infty$, i.e. the fiber product
$$
U \times_{\mathbb{G}_{m,s}} \Spec(k(s)((t^{-1}))) \rightarrow \Spec(k(s)((t^{-1}))),
$$
is a finite disjoint union of spectra of tamely ramified extensions of $k(s)((t^{-1}))$.
\item the geometric monodromy group of $f$, namely the image of the composition
$$
\pi_1(\mathbb{G}_{m,\overline{s}}, \overline{\eta}) \rightarrow \pi_1(\mathbb{G}_{m,s}, \overline{\eta}) \rightarrow \Aut(U_{\overline{\eta}}),
$$
has a unique $p$-Sylow subgroup.
\end{enumerate}
\end{defi}

\begin{teo}[\cite{katz}, 1.4.1]\label{chap3specialcover2} Let $\Fet^{\diamondsuit}(\mathbb{G}_{m,s})$ be as in \ref{chap3specialcover}. The functor
$$ 
\pi_{|\eta}^{-1} : \Fet^{\diamondsuit}(\mathbb{G}_{m,s}) \rightarrow \Fet(\eta),
$$
induced by the morphism $\pi_{|\eta} : \eta \rightarrow \mathbb{G}_{m,s}$, is an equivalence of categories.
\end{teo}

\subsection{\label{chap3gkexample}} For example, if $h$ is an element of $k(\eta)$ then Theorem \ref{chap3specialcover2} implies that the $\mathbb{F}_p$-torsor over $\eta$ defined by the equation
\begin{align}\label{chap3asex}
x - x^p = h,
\end{align}
must extend to an $\mathbb{F}_p$-torsor over $\mathbb{G}_{m,s}$, which is tamely ramified above $\infty$, hence unramified since $\mathbb{F}_p$ is a $p$-group. Let us write $h$ as $h_{\pi}(\pi) + r_{\pi}$, for some polynomial $h_{\pi}$ in $k(s)[t^{-1}]$ and some element $r_{\pi}$ of $\m$. Since $T$ is henselian, there exists a unique element $u_{\pi}$ of $\m$ such that $r_{\pi} = u_{\pi} - u_{\pi}^p$.
Then the $\mathbb{F}_p$-torsor over $\mathbb{G}_{m,s}$ defined by the equation
\begin{align}\label{chap3asex2}
y - y^p = h_{\pi}(t),
\end{align}
is unramified at $\infty$ and its pullback to $\eta$ by $\pi_{|\eta}$ is isomorphic to (\ref{chap3asex}). Therefore (\ref{chap3asex2}) is (up to isomorphism) the $\mathbb{F}_p$-torsor over $\mathbb{G}_{m,s}$ associated to (\ref{chap3asex}) by the equivalence in Theorem \ref{chap3specialcover2}. 

The Laurent polynomial $h_{\pi}(t)$ admits the following alternative description. Let $\nu \geq 1$ be an integer such that $h$ is of valuation strictly larger than $-\nu$, and let $u$ be the generator $1 \otimes 1 -t \otimes \pi^{-1}$ of the $\Ow_T/\m^{\nu}[t,t^{-1}]$-module $k(s)[t,t^{-1}] \otimes_{k(s)} \m^{-1}/\m^{-1 + \nu} $. Then we have
\begin{align}\label{chap3gkas}
h_{\pi}(t) = - \Res \left( h \frac{d u}{u} \right),
\end{align}
in $k(s)[t,t^{-1}]$. This formula should be understood as follows: we first take a lift $\tilde{u}$ of $u$ in $\pi^{-1} A[[\pi]]$, where $A = k(s)[t,t^{-1}]$, so that $\tilde{u}$ is invertible in $A((\pi))$, and we then set
$$
\Res \left( h \frac{d u}{u} \right) = \Res \left( h \frac{d \tilde{u}}{\tilde{u}} \right),
$$
where the right hand side is the specialization to $A = k(s)[t,t^{-1}]$ and $r=1$ of the following definition:

\begin{defi}\label{chap3residue} For any $k(s)$-algebra $A$, any non negative integer $r$ and any element $w = \sum_{n < r} w_n \otimes \pi^n$ of $A \otimes k(\eta)/\m^{r}$, we define
$$
d w = \sum_{n < r} w_n \otimes n \pi^{n-1} d\pi,
$$
in $A \otimes( k(\eta)/\m^{r-1}) d \pi$, and
$$
\Res(w d \pi) = w_{-1}.
$$
\end{defi}

Let us prove (\ref{chap3gkas}). We consider the lift $\tilde{u} = 1 - t \pi^{-1}$ of $u$ in $\pi^{-1} A[[\pi]]$. We have
\begin{align*}
\frac{d\tilde{u}}{\tilde{u}} &= \frac{t  d  \pi}{ \pi(\pi -t )} \\
&= - \frac{d \pi}{\pi(1- t^{-1}  \pi )} \\
&= - \sum_{n \geq 0} t^{-n} \pi^n \frac{d \pi}{\pi},
\end{align*}
in $A((\pi))d \pi$. If we write the image of $h$ in $\m^{1-\nu}/\m$ as $\sum_{n=0}^{\nu - 1} h_n \pi^{-n}$ for some elements $(h_n)_{0 \leq n < \nu}$ of $k(s)$, then this yields
\begin{align*}
\Res \left( h \frac{d u}{u} \right) &= - \sum_{0 \leq n < \nu} t^{-n} \otimes \Res \left( h \pi^n \frac{d \pi}{\pi} \right) \\
&= - \sum_{0 \leq n < \nu} h_n t^{-n},
\end{align*}
and the latter Laurent polynomial is exactly $-h_{\pi}(t)$, hence (\ref{chap3gkas}).

\subsection{\label{chap31.8}} Let $f: U \rightarrow \mathbb{G}_{m,s}$ be a connected special finite \'etale cover of $\mathbb{G}_{m,s}$ (cf. \ref{chap3specialcover}). Let $\eta' \rightarrow \eta$ be the pullback of $U$ to $\eta$ by $\pi$, and let us fix an $\eta$-morphism $\overline{\eta} \rightarrow \eta'$, so that we can consider $\overline{\eta}$ as a geometric point of $U$, henceforth denoted $\overline{\eta}'$. Similarly, let $s'$ be a finite \'etale extension of $s$ such that $U$ is a geometrically connected $s'$-scheme, and let us fix an $s$-morphism $\overline{s} \rightarrow s'$, so that we can consider $\overline{s}$ as a geometric point, denoted $\overline{s}'$, of the normalization $T'$ of $T$ in $\eta'$, or as a geometric point, also denoted $\overline{s}'$, of the normalization $X$ of $\mathbb{A}^1_s$ in $U$. We have a natural morphism
$$
\pi' : T' \rightarrow X,
$$
whose restriction to $\eta'$ is the natural morphism from $\eta'$ to $U$.

\begin{defi}\label{chap3specialcover5} A finite \'etale morphism $V \rightarrow U$ is \emph{$f$-special} if the composition
$$
V \rightarrow U \xrightarrow[]{f} \mathbb{G}_{m,s}
$$
is special (cf. \ref{chap3specialcover}). We denote by $\Fet^{\diamondsuit}(U,f)$ the category of $f$-special finite \'etale $U$-schemes, or equivalently the category of $U$-objects in $\Fet^{\diamondsuit}(\mathbb{G}_{m,s})$.
\end{defi}


Let us consider the fiber functor
\begin{align}
\begin{split}
\label{chap3fiberfunct}
\Fet^{\diamondsuit}(U,f) &\rightarrow \mathrm{Sets}\\
(V \rightarrow U) &\mapsto V_{\overline{\eta}'} = \Hom_{\eta'}(\overline{\eta}',V).
\end{split}
\end{align}
Let $\pi_1(U,f, \overline{\eta}')^{\diamondsuit}$ be the group of automorphisms of this functor, endowed with the coarsest topology such that for any object $V \rightarrow U$ of $\Fet^{\diamondsuit}(U,f)$, the natural group homomorphism
$$
\pi_1(U,f, \overline{\eta}')^{\diamondsuit} \rightarrow \Aut(V_{\overline{\eta}'}),
$$
is continuous, when the finite set $\Aut(V_{\overline{\eta}'})$ is endowed with the discrete topology. 

The topological group $\pi_1(U,f, \overline{\eta}')^{\diamondsuit}$ is profinite, and we have a natural surjective homomorphism
$$
\pi_1(U, \overline{\eta}') \rightarrow \pi_1(U,f, \overline{\eta}')^{\diamondsuit}.
$$
Moreover, the fiber functor (\ref{chap3fiberfunct}) realizes an equivalence from $\Fet^{\diamondsuit}(U,f)$ to the category of finite sets endowed with a continuous left action of $\pi_1(U,f, \overline{\eta}')^{\diamondsuit} $, cf. (\cite{katz}, 1.3.3). The natural surjective homomorphism
$$
\pi_1(U, \overline{\eta}') \rightarrow G_{s'},
$$
factors through $\pi_1(U,f, \overline{\eta}')^{\diamondsuit}$, and we denote by $\pi_1(U_{\overline{s}'},f, \overline{\eta}')^{\diamondsuit} $ the kernel of the resulting surjective homorphism from $\pi_1(U,f, \overline{\eta}')^{\diamondsuit}$ to $G_{s'}$.

By Theorem \ref{chap3specialcover2}, the homomorphism
\begin{align}\label{chap3profiniteisom}
(\pi'_{|\eta'})_* : G_{\eta'} \rightarrow \pi_1(U,f, \overline{\eta}')^{\diamondsuit},
\end{align}
is an isomorphism of profinite topological groups.

\subsection{\label{chap31.7}} Let $U,f,X$ be as in \ref{chap31.8}. Let $\Lambda$ be a finite admissible $\ell$-adic ring (cf. \ref{chap30.0.0.1}), and let $\mu$ be a $\Lambda$-admissible unitary multiplier on $G_k$ (cf. \ref{chap30.0}). The category $\Sh^{\diamondsuit}(X,f,\mu,\Lambda)$ of \textit{$f$-special} $\mu$-twisted $\Lambda$-sheaves on $X$ is the full subcategory of $\Sh(X,\mu,\Lambda)$ (cf. \ref{chap31.3}) whose objects are the $\mu$-twisted $\Lambda$-sheaves $\F = (\F_{k'}, (\rho_{\F}(g))_{g \in \Gal(k'/k)})$ on $X$ such that:
\begin{enumerate}
\item The restriction of $\F_{k'}$ to $U_{k'}$ is a $\Lambda$-local system (cf. \ref{chap31.0.1}).
\item The restriction of $\F_{k'}$ to any (equivalently, some) connected component of $U_{k'}$ is trivialized on a finite \'etale cover which is $f$-special, cf. \ref{chap3specialcover5}. Equivalently, the $\Lambda$-admissible representation $\F_{\overline{\eta}'}$ of $(\pi_1(U, \overline{\eta}'),\mu)$ (cf. \ref{chap31.4.4}) factors through a $\Lambda$-admissible representation of $(\pi_1(U,f, \overline{\eta}')^{\diamondsuit},\mu)$ (cf. \ref{chap31.8}).
\end{enumerate}

\subsection{\label{chap31.9}} Let $U,f,X$ be as in \ref{chap31.8}. Let $\Lambda$ be the ring of integers in a finite subextension of $\mathbb{Q}_{\ell}$ in $C$, and let $\mu$ be a $\Lambda$-admissible unitary multiplier on $G_k$ (cf. \ref{chap30.0}). The category $\Sh^{\diamondsuit}(X,f,\mu,\Lambda)$ of \textit{$f$-special} $\mu$-twisted $\Lambda$-sheaves on $X$ is the full subcategory of $\Sh(X,\mu,\Lambda)$ (cf. \ref{chap31.3.1}) consisting of its objects $(\F_n)_n$, where $\F_n$ is a $\mu$-twisted $\Lambda/\ell^n$-sheaf on $X$, such that $\F_n$ is special for each $n$ (cf. \ref{chap31.7}).

\subsection{\label{chap31.10}} Let $U,f,X$ be as in \ref{chap31.8}. Let $\Lambda$ be a finite subextension of $\mathbb{Q}_{\ell}$ in $C$, with ring of integers $\Lambda_0$, and let $\mu$ be a $\Lambda$-admissible unitary multiplier on $G_k$ (cf. \ref{chap30.0}). Then $\mu$ takes its values in $\Lambda_0^{\times}$. We define the category $\Sh^{\diamondsuit}(X,\mu,\Lambda)$ of \textit{$f$-special} $\mu$-twisted $\Lambda$-sheaves on $X$ to be the essential image of $\Sh^{\diamondsuit}(X,f,\mu,\Lambda_0)$ (cf. \ref{chap31.9}) in $\Sh(X,\mu,\Lambda)$ (cf. \ref{chap31.3.2}).

\subsection{\label{chap31.11}} Let $U,f,X$ be as in \ref{chap31.8}. Let $\Lambda$ be an $\ell$-adic coefficient ring (cf. \ref{chap30.0.0.1}), and let $\mu$ be a $\Lambda$-admissible unitary multiplier on $G_k$. We define the category $\Sh^{\diamondsuit}(X,f,\mu,\Lambda)$ of \textit{$f$-special} $\mu$-twisted $\Lambda$-sheaves on $X$ to be the full subcategory of $\Sh(X,f,\mu,\Lambda)$ (cf. \ref{chap31.3.3}) whose objects belong to the essential image of $\Sh(X,f,\mu,\Lambda_0)^{\diamondsuit}$ (cf. \ref{chap31.7}, \ref{chap31.9}, \ref{chap31.10}), for some admissible $\ell$-adic subring $\Lambda_0$ of $\Lambda$ containing the image of $\mu$.

We define the category $\Sh^{\diamondsuit}(U,f,\mu,\Lambda)$ of \textit{$f$-special} $\mu$-twisted $\Lambda$-sheaves on $U$ to be the full subcategory of $\Sh^{\diamondsuit}(X,f,\mu,\Lambda)$ consisting $f$-special $\mu$-twisted $\Lambda$-sheaves supported on the open subscheme $U$ of $X$.

\subsection{\label{chap31.12}} Let $U,f,X, \eta',s'$ be as in \ref{chap31.8}, and let $\Lambda$ and $\mu$ be as in \ref{chap31.11}. By \ref{chap31.4.4} and \ref{chap31.8}, the fiber functor $\F \mapsto \F_{\overline{\eta}'}$ realizes an equivalence from the category of $f$-special $\mu$-twisted $\Lambda$-sheaves on $U$ (cf. \ref{chap31.11}) to the category of $\Lambda$-admissible representations of $(\pi_1(U,f, \overline{\eta}')^{\diamondsuit},\mu)$.

More generally, for any object $\F$ of the category $\Sh^{\diamondsuit}(X,f,\mu,\Lambda)$, we have a cospecialization homomorphism
$$
c_{\F,\overline{s}', \overline{\eta}} : \F_{\overline{s}'} \rightarrow \F_{\overline{\eta}'},
$$
and the functor $\F \mapsto (\F_{\overline{s}'}, \F_{\overline{\eta}'}, c_{\F,\overline{s}', \overline{\eta}'})$ realizes an equivalence of categories from $\Sh^{\diamondsuit}(X,f,\mu,\Lambda)$ to the category of triples $(V_{s'}, V_{\eta'}, c)$, where $V_{\eta'}$ is a $\Lambda$-admissible representation of $(\pi_1(U,f, \overline{\eta}')^{\diamondsuit},\mu)$, where $V_{s'}$ is a $\Lambda$-admissible representation of $(G_{s'},\mu)$, and where 
$$
c : V_{s'} \rightarrow V_{\eta'}^{\pi_1(U_{\overline{s}'},f, \overline{\eta}')^{\diamondsuit}},
$$
is a homomorphism of $\Lambda$-admissible representations of $(G_{s'},\mu)$ from $V_{s'}$ to the subrepresentation of $V_{\eta'}$ consisting of its $\pi_1(U_{\overline{s}'}, f, \overline{\eta}')^{\diamondsuit}$-invariant elements (cf. \ref{chap31.8}).

\subsection{\label{chap31.5}} Let $U,f,X, T', \eta',s'$ be as in \ref{chap31.8}, and let $\Lambda$ and $\mu$ be as in \ref{chap31.11}. By \ref{chap31.4.4}, the fiber functor $\F \mapsto \F_{\overline{\eta}'}$ realizes an equivalence from the category of $\mu$-twisted $\Lambda$-sheaves on $\eta'$ (cf. \ref{chap31.3}) to the category of $\Lambda$-admissible representations of $(G_{\eta'},\mu)$.

More generally, the functor $\F \mapsto (\F_{\overline{s}'}, \F_{\overline{\eta}'}, c_{\F,\overline{s}', \overline{\eta}'})$, where $c_{\F,\overline{s}', \overline{\eta}'} : \F_{\overline{s}'} \rightarrow \F_{\overline{\eta}'} $ is the cospecialization homomorphism, realizes an isomorphism from $\Sh(T',\mu,\Lambda)$ to the category of triples $(V_{s'}, V_{\eta'}, c)$, where $V_{\eta'}$ is a $\Lambda$-admissible representation of $(G_{\eta'},\mu)$, where $V_{s'}$ is a $\Lambda$-admissible representation of $(G_{s'},\mu)$, and where 
$$
c : V_{s'} \rightarrow V_{\eta'}^{G_{\eta_{\overline{s}'}}},
$$
is a homomorphism of $\Lambda$-admissible representations of $(G_{s'},\mu)$ from $V_{s'}$ to the subrepresentation of $V_{\eta'}$ consisting of its $G_{\eta_{\overline{s}'}}$-invariant elements.

This remark, combined with the Galoisian description of $f$-special $\mu$-twisted $\Lambda$-sheaves on $X$ (cf. \ref{chap31.12}) and with the isomorphism (\ref{chap3profiniteisom}), yields:

\begin{teo}\label{chap3GK} Let $U,f,X,T', \pi'$ be as in \ref{chap31.8}, let $\Lambda$ be an $\ell$-adic coefficient ring (cf. \ref{chap30.0.0.1}), and let $\mu$ be a $\Lambda$-admissible unitary multiplier on $G_k$ (cf. \ref{chap30.0}). The pullback functor
$$
(\pi')^{-1} : \Sh^{\diamondsuit}(X,f,\mu,\Lambda) \rightarrow \Sh(T',\mu,\Lambda),
$$
is an equivalence from the category of $f$-special $\mu$-twisted $\Lambda$-sheaves on $X$ (cf. \ref{chap31.11}) to the category of $\mu$-twisted $\Lambda$-sheaves on $T'$.
\end{teo}

By restricting to $\mu$-twisted $\Lambda$-sheaves with vanishing fiber at $s'$, we similarly obtain:

\begin{teo}\label{chap3GK2} Let $U,f,\eta'$ be as in \ref{chap31.8}, let $\Lambda$ be an $\ell$-adic coefficient ring (cf. \ref{chap30.0.0.1}), and let $\mu$ be a $\Lambda$-admissible unitary multiplier on $G_k$ (cf. \ref{chap30.0}). The pullback functor
$$
(\pi_{|\eta}')^{-1} : \Sh^{\diamondsuit}(U,f,\mu,\Lambda) \rightarrow \Sh(\eta',\mu,\Lambda),
$$
is an equivalence from the category of $f$-special $\mu$-twisted $\Lambda$-sheaves on $U$ (cf. \ref{chap31.11}) to the category of $\mu$-twisted $\Lambda$-sheaves on $\eta'$.
\end{teo}

\subsection{\label{chap31.13}} For $U = \mathbb{G}_{m,s}$ and $f = \mathrm{id}$, let us simply refer to $f$-special sheaves as special sheaves. This agrees with the terminology in (\cite{katz}, 1.5). As a particular case of Theorem \ref{chap3GK}, we have the following extension result:

\begin{teo}\label{chap3GK3} Let $\Lambda$ be an $\ell$-adic coefficient ring (cf. \ref{chap30.0.0.1}), and let $\mu$ be a $\Lambda$-admissible unitary multiplier on $G_k$ (cf. \ref{chap30.0}). The pullback functor
$$
\pi^{-1} : \Sh^{\diamondsuit}(\mathbb{A}^1_s,\mu,\Lambda) \rightarrow \Sh(T,\mu,\Lambda),
$$
is an equivalence from the category of special $\mu$-twisted $\Lambda$-sheaves on $\mathbb{A}^1_s$ (cf. \ref{chap31.11}) to the category of $\mu$-twisted $\Lambda$-sheaves on $T$.
\end{teo}

By restricting to $\mu$-twisted $\Lambda$-sheaves with vanishing special fiber we similarly obtain:

\begin{teo}\label{chap3GK4} Let $\Lambda$ be an $\ell$-adic coefficient ring (cf. \ref{chap30.0.0.1}), and let $\mu$ be a $\Lambda$-admissible unitary multiplier on $G_k$ (cf. \ref{chap30.0}). The pullback functor
$$
\pi_{|\eta}^{-1} : \Sh^{\diamondsuit}(\mathbb{G}_{m,s},\mu,\Lambda) \rightarrow \Sh(\eta,\mu,\Lambda),
$$
is an equivalence from the category of special $\mu$-twisted $\Lambda$-sheaves on $\mathbb{G}_{m,s}$ (cf. \ref{chap31.11}) to the category of $\mu$-twisted $\Lambda$-sheaves on $\eta$.
\end{teo}

When $\mu =1$, the latter theorem matches (\cite{katz}, Th. 1.5.6).

\subsection{\label{chap31.14}} Let $\Lambda$ be an $\ell$-adic coefficient ring (cf. \ref{chap30.0.0.1}), let $\mu$ be a $\Lambda$-admissible unitary multiplier on $G_k$ (cf. \ref{chap30.0}), and let $U,f,\eta'$ be as in \ref{chap31.8}. 

\begin{lem}\label{chap3pushlemma} If $\F$ is an object of $ \Sh^{\diamondsuit}(U,f,\mu,\Lambda)$, then its pushforward $f_* \F$ (cf. \ref{chap31.3.4}) belongs to $\Sh^{\diamondsuit}(\mathbb{G}_{m,s},\mu,\Lambda)$.
\end{lem}

We can assume (and we do) that $\Lambda$ is finite, in which case the $\Lambda$-local system $\F$ is represented by a finite \'etale morphism $g : V \rightarrow U$ which is $f$-special (cf. \ref{chap3specialcover5}). The $\Lambda$-local system $f_* \F$ is then represented by the finite \'etale morphism $fg : V \rightarrow \mathbb{G}_{m,s}$, which is special, hence the conclusion of Lemma \ref{chap3pushlemma}.

It follows from Lemma \ref{chap3pushlemma} that we have the following commutative diagram (up to natural isomorphisms).

\begin{center}
 \begin{tikzpicture}[scale=1]

\node (A) at (0,2) {$\Sh(\eta',\mu,\Lambda)$};
\node (B) at (5,2) {$ \Sh^{\diamondsuit}(U,f,\mu,\Lambda)$};
\node (C) at (5,0) {$\Sh^{\diamondsuit}(\mathbb{G}_{m,s},\mu,\Lambda)$};
\node (D) at (0,0) {$\Sh(\eta,\mu,\Lambda)$};

\path[->,font=\scriptsize]
(B) edge node[above]{$(\pi_{|\eta}')^{-1}$} (A)
(B) edge node[right]{$f_*$} (C)
(A) edge node[right]{$(f_{| \eta'})_*$} (D)
(C) edge node[above]{$\pi_{|\eta}^{-1}$} (D);
\end{tikzpicture} 
\end{center}
The rows of this diagram are equivalences of categories by Theorems \ref{chap3GK2} and \ref{chap3GK4}.

\section{Geometric class field theory\label{chap3gcftsection}}

We review in this section global and local geometric class field theory. Let $\Lambda$ be an $\ell$-adic coefficient ring (cf. \ref{chap3conv}, \ref{chap30.0.0.1}). The purpose of geometric class field theory is to establish equivalences between groupoids of $\Lambda$-local systems of rank $1$ on curves over $k$, or over germs of curves, and groupoids of multiplicative local systems over certain group schemes. The notion of multiplicative local system, which geometrizes the notion of character, is reviewed in \ref{chap32.0}, \ref{chap32.0.1}, \ref{chap32.0.2} and \ref{chap32.0.3} below. 

\subsection{\label{chap32.0}} Let us assume that $\Lambda$ is a finite $\ell$-adic coefficient ring (cf. \ref{chap30.0.0.1}). Let $S$ be a $k$-scheme and let $G$ be a commutative $S$-group scheme, with multiplication $m : G \times_S G \rightarrow G$. A \textit{multiplicative $\Lambda$-local system on $G$} is a $\Lambda$-local system $\Lc$ on $G$, of rank $1$, together with an isomorphism $\theta : p_1^{-1} \Lc \otimes p_2^{-1} \Lc \rightarrow m^{-1} \Lc$ of $\Lambda$-local systems on $G \times G$ where $p_1$ and $p_2$ are the canonical projections, which satisfy the following two properties, cp. (\cite{G18}, Def. 2.5).
\begin{enumerate}
\item Symmetry: if $\sigma$ is the involution of $G \times G$ which switches the two factors, then the isomorphism
$$
 p_2^{-1} \Lc \otimes p_1^{-1} \Lc \rightarrow \sigma^{-1}(p_1^{-1} \Lc \otimes p_2^{-1} \Lc) \xrightarrow[]{\sigma^{-1} \theta} \sigma^{-1} m^{-1} \Lc \rightarrow m^{-1} \Lc
$$
is the composition of $\theta$ with the canonical isomorphism $ p_2^{-1} \Lc \otimes p_1^{-1} \Lc \rightarrow p_1^{-1} \Lc \otimes p_2^{-1} \Lc$.
\item Associativity: if $q_i : G \times G \times G \rightarrow G$ (resp. $q_{ij} : G \times G \times G \rightarrow G \times G$) is the projection on the $i$-th factor for $i \in \llbracket1,3\rrbracket$ (resp. on the $i$-th and $j$-th factors for $(i,j) \in \llbracket1,3\rrbracket^2$ such that $i<j$) and if $m_3 : G \times G \times G \rightarrow G$ is the multiplication morphism, then the diagram of $\Lambda$-local systems on $G \times G \times G$
\begin{center}
 \begin{tikzpicture}[scale=1]

\node (A) at (0,2) {$q_1^{-1} \Lc \otimes q_2^{-1} \Lc \otimes q_3^{-1} \Lc$};
\node (B) at (3,4) {$q_1^{-1} \Lc \otimes (m q_{23})^{-1} \Lc$};
\node (C) at (3,0) {$(m q_{12})^{-1} \Lc \otimes q_3^{-1} \Lc$.};
\node (D) at (6,2) {$m_3^{-1} \Lc$};
\path[->,font=\scriptsize]
(A) edge node[above left]{$id \otimes q_{23}^{-1} \theta $} (B)
(A) edge node[below left]{$ q_{12}^{-1} \theta \otimes id$} (C)
(B) edge node[above right]{$ (q_1 \times mq_{23})^{-1} \theta $} (D)
(C) edge node[below right]{$ (mq_{12} \times q_3)^{-1} \theta$} (D);

\end{tikzpicture} 
\end{center}
is commutative.
\end{enumerate}
A morphism between multiplicative $\Lambda$-local systems $(\Lc_1,\theta_1)$ and $(\Lc_2,\theta_2)$ on $G$ is an isomorphism $\alpha : \Lc_1 \rightarrow \Lc_2$ of $\Lambda$-local systems such that the diagram
\begin{center}
 \begin{tikzpicture}[scale=1]

\node (A) at (0,0) {$p_1^{-1} \Lc_1 \otimes p_2^{-1} \Lc_1$};
\node (B) at (0,2) {$p_1^{-1} \Lc_2 \otimes p_2^{-1} \Lc_2$};
\node (C) at (3,0) {$m^{-1} \Lc_1$.};
\node (D) at (3,2) {$m^{-1} \Lc_2$};
\path[->,font=\scriptsize]
(A) edge node[above]{$\theta_1$} (C)
(B) edge node[above]{$\theta_2$} (D)
(A) edge node[left]{$ p_1^{-1} \alpha \otimes p_2^{-1} \alpha$} (B)
(C) edge node[right]{$ m^{-1} \alpha$} (D);

\end{tikzpicture} 
\end{center}
is commutative.

We denote by $\mathrm{Loc}^{\otimes}(G,\Lambda)$ the groupoid of multiplicative $\Lambda$-local systems on $G$. The group of automorphisms of an object of $\mathrm{Loc}^{\otimes}(G,\Lambda)$ is given by
$$
\Aut_{\mathrm{Loc}^{\otimes}(G,\Lambda)}(\Lambda) = \mathrm{Hom}_{\mathrm{Grp}/S}(G,\Lambda^{\times}_S),
$$
where $\Lambda^{\times}_S$ is the constant $S$-group scheme associated to $\Lambda^{\times}$. If $S$ is connected and if $G$ is an extension of a constant $S$-group scheme, associated to a discrete group $\pi_0(G)$, by an $S$-group scheme with connected geometric fibers, then we also have
$$
\Aut_{\mathrm{Loc}^{\otimes}(G,\Lambda)}(\Lambda) = \mathrm{Hom}_{\mathrm{Grp}}(\pi_0(G),\Lambda^{\times}).
$$

\begin{rema}\label{chap3dico} The functor $\Lc \mapsto \mathcal{I}som(\Lambda, \Lc)$ which sends a multiplicative $\Lambda$-local system $\Lc$ on $G$ to the $\Lambda^{\times}$-torsor of its local trivializations realizes an equivalence from the category $\mathrm{Loc}^{\otimes}(G,\Lambda)$, to the groupoid of multiplicative $\Lambda^{\times}$-torsors on $G$ in the sense of (\cite{G18}, Def. 2.5). By (\cite{G18}, Def. 2.9), the latter groupoid is equivalent to the groupoid of extensions of $G$ by $\Lambda^{\times}$ in the category of commutative $S$-group schemes. Namely, to an extension
$$
0 \rightarrow \Lambda^{\times} \rightarrow E \rightarrow G \rightarrow 0,
$$
of commutative $S$-group schemes, one associate the $\Lambda$-local system of rank $1$ on $G$ corresponding to the $\Lambda^{\times}$-torsor $E$ over $G$, where $\Lambda^{\times}$ acts by left multiplication on $E$; this $\Lambda$-local system of rank $1$ on $G$ is naturally endowed with a structure of multiplicative $\Lambda$-local system on $G$, cf. (\cite{G18}, Def. 2.4).
\end{rema}

\begin{exemple} Assume that $k$ is finite of cardinality $q$, and that $G$ is a connected commutative $k$-group scheme. The $q$-Frobenius morphism $F : G \rightarrow G$ is then an homomorphism of $k$-group schemes. Moreover, the sequence
$$
0 \rightarrow G(k) \rightarrow G \xrightarrow[]{1-F} G \rightarrow 0,
$$
is exact, cf. (\cite{Laumon}, 1.1.3). For any homomorphism $\chi : G(k) \rightarrow \Lambda^{\times}$, the pushout of this exact sequence provides an extension of $G$ by $\Lambda^{\times}$, which yields in turn a multiplicative $\Lambda$-local system $\Lc_{\chi}$ on $G$ by Remark \ref{chap3dico}.
\end{exemple}

%

\subsection{\label{chap32.0.1}} Let us assume that $\Lambda$ is the ring of integers in a finite extension of $\mathbb{Q}_{\ell}$. Let $S$ be a connected $k$-scheme and let $G$ be a commutative $S$-group scheme, which is an extension of a constant $S$-group scheme, associated to a discrete group $\pi_0(G)$, by an $S$-group scheme with connected geometric fibers. We define the groupoid $\mathrm{Loc}^{\otimes}(G,\Lambda)$ of multiplicative $\Lambda$-local systems on $G$ to be the $2$-limit of the categories $\mathrm{Loc}^{\otimes}(G,\Lambda/ \ell^n)$ (cf. \ref{chap32.0}), where $n$ runs over over all positive integers. Thus the objects of $\mathrm{Loc}^{\otimes}(G,\Lambda)$ are given by projective systems $(\Lc_n)_{n \geq 1}$ where $\Lc_n$ is a multiplicative $\Lambda/ \ell^n$-local system on $G$, such that the transition maps $\Lc_{n+1} \rightarrow \Lc_n$ induce isomorphisms
$\Lc_{n+1} \otimes \Lambda/ \ell^n \rightarrow \Lc_n,$
for each integer $n$.

The group of automorphisms of an object of $\mathrm{Loc}^{\otimes}(G,\Lambda)$ is given by
$$
\Aut_{\mathrm{Loc}^{\otimes}(G,\Lambda)}(\Lambda) = \lim_n \mathrm{Hom}_{\mathrm{Grp}}(\pi_0(G), (\Lambda/ \ell^n)^{\times}) = \mathrm{Hom}_{\mathrm{Grp}}(\pi_0(G),\Lambda^{\times}).
$$

\subsection{\label{chap32.0.2}} Let us assume that $\Lambda$ is a finite extension of $\mathbb{Q}_{\ell}$, with ring of integers $\Lambda_0 \subseteq \Lambda$. Let $S$ be a connected $k$-scheme and let $G$ be a commutative $S$-group scheme, which is an extension of a constant $S$-group scheme, associated to a discrete group $\pi_0(G)$, by an $S$-group scheme with connected geometric fibers. We define the groupoid $\mathrm{Loc}^{\otimes}(G,\Lambda)$ of multiplicative $\Lambda$-local systems on $G$ to be the groupoid whose objects are those of $\mathrm{Loc}^{\otimes}(G,\Lambda_0)$ (cf. \ref{chap32.0.1}), and whose morphisms are given by
$$
\mathrm{Isom}_{\Lambda}(\Lc_1,\Lc_2) = \mathrm{Hom}_{\mathrm{Grp}}(\pi_0(G),\Lambda^{\times}) \wedge \mathrm{Isom}_{\Lambda_0}(\Lc_1,\Lc_2),
$$
namely the quotient of $\mathrm{Hom}_{\mathrm{Grp}}(\pi_0(G),\Lambda^{\times}) \times \mathrm{Isom}_{\Lambda_0}(\Lc_1,\Lc_2)$ by the action of $\mathrm{Hom}_{\mathrm{Grp}}(\pi_0(G),\Lambda_0^{\times})$ given by $u(\lambda, \varphi) = (u^{-1} \lambda, u\varphi)$ for $u$ in $\mathrm{Hom}_{\mathrm{Grp}}(\pi_0(G),\Lambda_0^{\times})$ and $(\lambda, \varphi)$ in $\mathrm{Hom}_{\mathrm{Grp}}(\pi_0(G),\Lambda^{\times}) \times \mathrm{Isom}_{\Lambda_0}(\Lc_1,\Lc_2)$. In particular, the group of automorphisms of an object of $\mathrm{Loc}^{\otimes}(G,\Lambda)$ is given by
$$
\Aut_{\mathrm{Loc}^{\otimes}(G,\Lambda)}(\Lambda) = \mathrm{Hom}_{\mathrm{Grp}}(\pi_0(G),\Lambda^{\times}).
$$

Isomorphisms between multiplicative $\Lambda$-local systems $\Lc_1$ and $\Lc_2$ on $G$ can be alternatively described as isomorphisms $\alpha : \Lc_1 \rightarrow \Lc_2$ of $\Lambda$-local systems (cf. \ref{chap31.0.3}) such that the diagram
\begin{center}
 \begin{tikzpicture}[scale=1]

\node (A) at (0,0) {$p_1^{-1} \Lc_1 \otimes p_2^{-1} \Lc_1$};
\node (B) at (0,2) {$p_1^{-1} \Lc_2 \otimes p_2^{-1} \Lc_2$};
\node (C) at (3,0) {$m^{-1} \Lc_1$.};
\node (D) at (3,2) {$m^{-1} \Lc_2$};
\path[->,font=\scriptsize]
(A) edge node[above]{$\theta_1$} (C)
(B) edge node[above]{$\theta_2$} (D)
(A) edge node[left]{$ p_1^{-1} \alpha \otimes p_2^{-1} \alpha$} (B)
(C) edge node[right]{$ m^{-1} \alpha$} (D);

\end{tikzpicture} 
\end{center}
is commutative. Indeed, if $\pi$ is a uniformizer of $\Lambda$, then such an isomorphism $\alpha$ must be of the form $\alpha = \pi^{v} \varphi$, where $\varphi : \Lc_1 \rightarrow \Lc_2$ is an isomorphism of $\Lambda_0$-local systems, and where $v$ is a map from $\pi_0(G)$ to $\Z$. The commutativity of the diagram above implies that $d^1(\pi^v) = \pi^{d^1(v)}$ (cf. \ref{chap3eq21}) takes its values in $\Lambda_0^{\times}$, hence is trivial. Thus $\lambda = \pi^v$ is a group homomorphism from $\pi_0(G)$ to $\Lambda^{\times}$ and $\varphi$ is an isomorphism of multiplicative $\Lambda_0$-local systems.

\subsection{\label{chap32.0.3}} We now consider an arbitrary $\ell$-adic coefficient ring $\Lambda$. Let $S$ be a connected $k$-scheme and let $G$ be a commutative $S$-group scheme, which is an extension of a constant $S$-group scheme, associated to a discrete group $\pi_0(G)$, by an $S$-group scheme with connected geometric fibers. We define the groupoid $\mathrm{Loc}^{\otimes}(G,\Lambda)$ of multiplicative $\Lambda$-local systems on $G$ to be the $2$-colimit of the groupoids $\mathrm{Loc}^{\otimes}(G,\Lambda_0)$, where $\Lambda_0$ runs over all admissible $\ell$-adic subrings of $\Lambda$ (cf. \ref{chap32.0}, \ref{chap32.0.1} and \ref{chap32.0.2}). The group of automorphisms of an object of $\mathrm{Loc}^{\otimes}(G,\Lambda)$ is given by the group of \textit{$\Lambda$-admissible} group homomorphisms from $\pi_0(G)$ to $\Lambda^{\times}$ (cf. \ref{chap30.0.0.0}).

\begin{rema} If the discrete group $\pi_0(G)$ is finitely generated, then any group homomorphism from $\pi_0(G)$ to $\Lambda^{\times}$ is $\Lambda$-admissible.
\end{rema}

\begin{exemple}\label{chap3zmultloc} Let us consider the discrete group scheme $G = \Z_S$. For each integer $n$, let $n : S \rightarrow \Z_S$ be the section corresponding to the element $n$ of $\Z$. Then the pullback by the section $1$ realizes an equivalence from $\mathrm{Loc}^{\otimes}( \Z_S,\Lambda)$ to the groupoid of $\Lambda$-local systems of rank $1$ on $S$. A quasi-inverse to this functor is given by sending a $\Lambda$-local system $\F$ of rank $1$ on $S$ to the multiplicative $\Lambda$-local system on $\Z_S$ whose pullback by the section $n$ is $\F^{\otimes n}$, for any integer $n$, together with the isomorphism
$\theta : p_1^{-1} \Lc \otimes p_2^{-1} \Lc \rightarrow m^{-1} \Lc,$
whose pullback by a section $(n,m)$ of $\Z_S \times_S \Z_S$ is the canonical isomorphism
$$
\F^{\otimes n} \otimes \F^{\otimes m} \rightarrow \F^{\otimes (n+m)}.
$$
\end{exemple}

\subsection{\label{chap32.0.4}} Let $S$ be a connected $k$-scheme and let $G_1,G_2$ be commutative $S$-group schemes as in \ref{chap32.0.3}. Then $G_1 \times_S G_2$ is an extension of the discrete group $\pi_0(G_1) \times \pi_0(G_2)$ by an $S$-group scheme with connected geometric fibers. Let $p_1,p_2$ be the natural projections from $G_1 \times_S G_2$ to $G_1$ and $G_2$ respectively. Then the functor
\begin{align*}
\mathrm{Loc}^{\otimes}(G_1,\Lambda) \times \mathrm{Loc}^{\otimes}(G_2,\Lambda) &\rightarrow \mathrm{Loc}^{\otimes}(G_1 \times_S G_2,\Lambda) \\
(\Lc_1, \Lc_2) &\mapsto p_1^{-1} \Lc_1 \otimes p_2^{-1} \Lc_2,
\end{align*}
is an equivalence of categories. Indeed, if $\iota_1,\iota_2$ are the inclusions of the factors $G_1$ and $G_2$ respectively in $G_1 \times_S G_2$, then for any multiplicative $\Lambda$-local system $(\Lc,\theta)$ on $G_1 \times_S G_2$, the isomorphism $\theta$ produces an isomorphism
$$
\Lc \rightarrow p_1^{-1} \iota_1^{-1} \Lc \otimes p_2^{-1} \iota_2^{-1} \Lc.
$$

\subsection{\label{chap32.0.5}} Let $G$ be a commutative $k$-group scheme as in \ref{chap32.0.3}. Let $i : H \rightarrow G$ be a closed connected sub-$k$-group scheme of $G$, so that $G/H$ is also an extension of the constant $k$-group scheme $\pi_0(G)$ by a quasi-compact connected $k$-group scheme. Let $r : G \rightarrow G/H$ be the canonical projection. Then the pullback functor $r^{-1}$ induces an equivalence from the groupoid
$$
\mathrm{Loc}^{\otimes}(G/H,\Lambda)
$$
of multiplicative $\Lambda$-local systems on $G/H$ (cf. \ref{chap32.0.3}) to the groupoid of triples $(\Lc,\theta, \zeta)$, where $(\Lc,\theta)$ is a multiplicative $\Lambda$-local systems on $G$ and $\zeta : \Lambda_{H} \rightarrow i^{-1} \Lc$ is an isomorphism of multiplicative $\Lambda$-local systems on $H$. Indeed, if $(\Lc,\theta, \zeta)$ is such a triple, if $p_1,p_2,m$ are as in \ref{chap32.0}, and if $\varphi$ is the isomorphism
\begin{align*}
\varphi : G \times_k H &\rightarrow G \times_{G/H} G \\
(g,h) &\rightarrow (g,gh)
\end{align*}
then we have a sequence of isomorphisms
$$
\varphi^{-1} p_1^{-1} \Lc \rightarrow p_1^{-1} \Lc_{|G \times_k H} \xrightarrow[]{\mathrm{id} \otimes \zeta} p_1^{-1} \Lc \otimes p_2^{-1} \Lc_{|G \times_k H} \xrightarrow[]{\theta_{|G\times_k H}} m^{-1} \Lc_{|G \times_k H} \xrightarrow[]{} \varphi^{-1} p_2^{-1} \Lc,
$$
which yields a descent datum $p_1^{-1} \Lc \rightarrow p_2^{-1} \Lc$ on $G \times_{G/H} G$, with respect $r$, which is a morphism of effective descent for the fibered category of $\Lambda$-local systems, as well as for the fibered category of multiplicative $\Lambda$-local systems.

\subsection{\label{chap32.1}} Let $X$ be a smooth geometrically connected projective curve of genus $g$ over $k$, let $i : D \rightarrow X$ be an effective divisor of degree $d \geq 1$ on $X$, and let $U$ be the open complement of $D$ in $X$. Our aim is to describe $\Lambda$-local systems of rank $1$ on $U$. One first introduces a measure of the ramification at infinity of such a local system:

\begin{defi}\label{chap3ramdefi} A $\Lambda$-local system $\F$ of rank $1$ on $U$ has \textit{ramification bounded by $D$} if for any point $x$ of $D$, the Swan conductor of the restriction of $F$ to the spectrum of the fraction field of the completed local ring of $X$ at $x$ is strictly less than the multiplicity of $D$ at $x$.
\end{defi}

The main theorem of geometric class field theory, namely Theorem \ref{chap3ggcft} below, states an equivalence between the groupoid of $\Lambda$-local systems of rank $1$ on $U$ with ramification bounded by $D$, and the groupoid of multiplicative $\Lambda$-local systems on a $k$-group scheme, the generalized Picard scheme, which we now introduce:

\begin{defi}\label{chap32.1.0.0} The generalized Picard functor $\Pic_k(X,D)$ associated to $(X,D)$ is the functor which to a $k$-scheme $S$ associates the group of isomorphism classes of pairs $(\Lc, \alpha)$ where $\Lc$ is an invertible $\Ow_{X_S}$-module and $\alpha : \Ow_{D_S} \rightarrow i_S^* \Lc$ is an isomorphism of $\Ow_{D_S}$-modules. Here, by $X_S,D_S,i_S$ we denote the base change of $X,D,i$ along $S \rightarrow \Spec(k)$.
\end{defi}

If the effective divisor $D$ is given by a single $k$-point $x$ of $X$ with multiplicity $1$, then the morphism
\begin{align*}
\Pic_k(X,x) &\rightarrow \Pic_k(X) \\
(\Lc, \alpha) &\rightarrow \Lc,
\end{align*}
is an isomorphism, where $\Pic_k(X)$ is the Picard functor of $X$. The latter is well-known to be representable by a $k$-group scheme, namely an extension of $\mathbb{Z}_k$ by the Jacobian scheme of $X$, which is an abelian $k$-scheme of dimension $g$. In general, we have:

\begin{prop}[\cite{G18}, Prop. 4.8] The generalized Picard functor $\Pic_k(X,D)$ is representable by a smooth separated $k$-group scheme of dimension $d + g - 1$.
\end{prop}


Let us consider the Abel-Jacobi morphism
\begin{align}\label{chap3abeljacob}
\Phi : U \rightarrow \Pic_k(X,D),
\end{align}
which sends a section $x$ of $U$ to the pair $(\Ow(x),1)$, where $1 : \Ow_D \rightarrow \Ow(x) \otimes_{\Ow_X} \Ow_D$ is the trivialization of $\Ow(x)$ on $D$ induced by the canonical section $1 : \Ow_X \hookrightarrow \Ow(x)$. Global geometric class field theory can then be stated as follows:

\begin{teo}[Global geometric class field theory]\label{chap3ggcft} Let $\F$ be a $\Lambda$-local system of rank $1$ on $U$, with ramification bounded by $D$ (cf. \ref{chap3ramdefi}). Then, there exists a unique (up to unique isomorphism) pair $(\chi_{\F},\beta)$, where $\chi_{\F}$ is a multiplicative $\Lambda$-local system on $\Pic_k(X,D)$ (cf. \ref{chap32.0.3}), and $\beta : \Phi^{-1} \chi_{\F} \rightarrow \F$ is an isomorphism. The functor $\F \mapsto \chi_{\F}$ is an equivalence from the groupoid of $\Lambda$-local systems of rank $1$ on $U$, with ramification bounded by $D$, to the groupoid of multiplicative $\Lambda$-local systems on $\Pic_k(X,D)$.
\end{teo}

This theorem reduces to the case where $\Lambda$ is finite, which was originally proved by Serre and Lang, cf. (\cite{Lang}, 6) and \cite{JPS3}, by using the Albanese property of Rosenlicht's generalized Picard schemes \cite{Ros}. Deligne gave another proof in the tamely ramified case. An exposition of Deligne's proof in the unramified case over a finite field can be found in \cite{Laumon}. Deligne's approach was later extended to allow arbitrary ramification simultaneously by the author (\cite{G18}, Th. 1.1) and by Takeuchi (\cite{T18}, Th. 1.1).

%

\subsection{\label{chap32.3}} Let $T$ be the spectrum of a $k$-algebra, which is a henselian discrete valuation ring whose residue field is a finite extension of $k$. Let $\eta$ be the generic point of $T$, and let $s$ be its closed point, so that $k(\eta)$ is a henselian discrete valuation field, with valuation subring $\Ow_{T,s}$, and with residue field $k(s)$ which is a finite extension of $k$. By Hensel's lemma, there exists a unique morphism $T \rightarrow s$ of $k$-schemes whose composition with the immersion $s \rightarrow T$ is the identity. We can thus consider $T$ as an $s$-scheme. 

\begin{defi}\label{chap3ramdefiloc} Let $D$ be a closed subscheme of $T$ supported on $s$. A $\Lambda$-local system $\F$ of rank $1$ on $\eta$ has \textit{ramification bounded by $D$} if its Swan conductor is strictly less than the multiplicity of $D$ at $s$, namely the length of $\Ow_{D,s}$ as an $\Ow_{T,s}$-module.
\end{defi}

Before proceeding further, we need the following result:

\begin{prop}\label{chap3lemmadiag} Let $D$ be a closed subscheme of $T$ supported on $s$. The kernel $\mathcal{I}$ of the homomorphism
\begin{align*}
\Ow_T \otimes_{k(s)} \Ow_D &\rightarrow \Ow_D \\
f_1 \otimes f_2 &\rightarrow f_1 f_2,
\end{align*}
is an invertible ideal of $\Ow_T \otimes_{k(s)} \Ow_D$, which generates the unit ideal of $k(\eta) \otimes_{k(s)} \Ow_D$.
\end{prop}

Indeed, if $\pi$ is a uniformizer of the discrete valuation field $k(\eta)$, then the kernel of the multiplication from $\Ow_T \otimes_{k(s)} \Ow_D$ to $\Ow_D$ is generated by $\pi \otimes 1 - 1 \otimes \pi$. Since $1 \otimes \pi$ is nilpotent in $\Ow_T \otimes_{k(s)} \Ow_D$, the section $\pi \otimes 1 - 1 \otimes \pi$ becomes a unit in $k(\eta) \otimes_{k(s)} \Ow_D$. Since the natural homomorphism from $\Ow_T \otimes_{k(s)} \Ow_D$ to $k(\eta) \otimes_{k(s)} \Ow_D$ is injective, this proves that $\pi \otimes 1 - 1 \otimes \pi$ generates an invertible ideal of $\Ow_T \otimes_{k(s)} \Ow_D$, and this concludes the proof of Proposition \ref{chap3lemmadiag}.

\subsection{\label{chap32.4}} We aim at describing the $\Lambda$-local systems $\F$ of rank $1$ on $\eta$ with ramification bounded by $D$ in terms of multiplicative $\Lambda$-local systems on a certain group scheme (cf. \ref{chap32.0.3}), the local Picard group, which we now introduce. 

\begin{defi}\label{chap3localpic} Let $D$ be a closed subscheme of $T$ supported on $s$, and let $\mathcal{I}$ be the invertible ideal of $\Ow_T \otimes_{k(s)} \Ow_D$ from Proposition \ref{chap3lemmadiag}. The local Picard group $\Pic(T,D)$ associated to the pair $(T,D)$ is the functor which sends a $T$-scheme $S$ to the group of pairs $(d,u)$, where $d$ is a locally constant $\mathbb{Z}$-valued map on $S$, and 
$$u : \Ow_S \otimes_{k(s)} \Ow_D \rightarrow \Ow_S \otimes_{\Ow_T} \mathcal{I}^{-d} $$
is an isomorphism of $\Ow_S \otimes_{k(s)} \Ow_D$-modules.
\end{defi}

One can informally think of a section $(d,u)$ of $\Pic(T,D)$ over a $T$-scheme $S$ as a trivialization of the line bundle $\Ow(d \Delta)$ along the effective Cartier divisor $S \times_{s} D$ on the germ of $S$-curve $S \times_{s} T$, where $\Delta : S \rightarrow S \times_{s} T$ is the diagonal embedding. 

Sending a section $(d,u)$ of $\Pic(T,D)$ to $d$ defines a homomorphism 
$$
\Pic(T,D) \rightarrow \Z_T,
$$
of group valued functors. We denote by $\Pic^0(T,D)$ the kernel of this homomorphism. The special fiber $\Pic^0(T,D)_s$ is the functor which sends an $s$-scheme $S$ to the group of units in $\Ow_S \otimes_{k(s)} \Ow_D$. The natural homomorphism
$$
\Pic^0(T,D)_s \times_s T \rightarrow \Pic^0(T,D),
$$
which sends a unit $u$ of $\Ow_S \otimes_{k(s)} \Ow_D$ to the pair $(0,u)$ is an isomorphism.

\begin{prop}\label{chap3exactseq} Let $D$ be a closed subscheme of $T$ supported on $s$. The functor $\Pic(T,D)$ is representable by a $T$-group scheme, which fits into a (split) exact sequence
$$
1 \rightarrow \Pic^0(T,D) \rightarrow \Pic(T,D) \rightarrow \Z_T \rightarrow 0,
$$
where $\Pic^0(T,D)$ is representable by a smooth separated affine $T$-group scheme, with geometrically connected fibers of dimension equal to the multiplicity of $D$ at $s$.
\end{prop}

Indeed, if $\pi$ is a uniformizer of the discrete valuation field $k(\eta)$, then the ideal $\mathcal{I}$ of $\Ow_T \otimes_{k(s)} \Ow_D$ is generated $1 \otimes \pi - \pi \otimes 1$, hence for any $S$-point $(d,u)$ of $\Pic(T,D)$ the isomorphism $u$ can be uniquely written as a sum
$$
u = \left(\sum_{0 \leq n < \nu} u_n \pi^{n} \right) (1 \otimes \pi - \pi \otimes 1)^{-d},
$$
where $\nu$ is the multiplicity of $D$ at $s$, and $(u_n)_{0 \leq n < \nu}$ are sections of $\Ow_S$, such that $u_{0}$ is invertible. Thus $\Pic(T,D)$ is representable by a product $\mathbb{Z}_{T} \times_{T} \mathbb{G}_{m,T} \times_{T} \mathbb{A}^{\nu-1}_{T}$, and this concludes the proof of Proposition \ref{chap3exactseq}.

\begin{rema} The exact sequence
$$
1 \rightarrow \Pic^0(T,D) \rightarrow \Pic(T,D) \rightarrow \Z_T \rightarrow 0,
$$
is split, but the splitting constructed in the proof of Proposition \ref{chap3exactseq} depends on a choice of uniformizer, and is therefore non canonical.
\end{rema}

\subsection{\label{chap32.10}} By Proposition \ref{chap3lemmadiag}, the natural homomorphism $\Ow_T \otimes_{k(s)} \Ow_D \rightarrow \mathcal{I}^{-1}$ induces an isomorphism
$$
u_{\mathrm{can}} : k(\eta)\otimes_{k(s)} \Ow_D \rightarrow k(\eta) \otimes_{\Ow_T} \mathcal{I}^{-1}.
$$
The pair $(1,u_{\mathrm{can}})$ yields a $k(\eta)$-point of $\Pic(T,D)$, corresponding to a morphism
$$
\Phi_{\eta} : \eta \rightarrow \Pic(T,D),
$$
of $T$-schemes, which plays the role of a local Abel-Jacobi morphism in Theorem \ref{chap3lgcft1}. Recall that we have an isomorphism
 $$
\alpha : \Pic^0(T,D)_s \times_s T \rightarrow \Pic^0(T,D),
$$
which sends a unit $u$ of $\Ow_S \otimes_{k(s)} \Ow_D$ to the pair $(0,u)$. Let us also denote by $p_1$ the projection of $\Pic^0(T,D)_s \times_s T$ onto the first factor. 

\begin{defi}\label{chap3lgcft0} We define the groupoid $\Trip(T,D,\Lambda)$ to be the category of triples $(\chi,\widetilde{\chi}, \beta)$, consisting of
\begin{enumerate}

\item a multiplicative $\Lambda$-local system $\chi$ on the $s$-group scheme $\Pic^0(T,D)_s$ (cf. \ref{chap32.0.3}),
\item a multiplicative $\Lambda$-local system $\widetilde{\chi}$ on the $T$-group scheme $\Pic(T,D)$ (cf. \ref{chap32.0.3}),
\item an isomorphism $\beta : \alpha^{-1} \widetilde{\chi} \rightarrow p_1^{-1} \chi$ on $\Pic^0(T,D)_s \times_s T $.

\end{enumerate}
\end{defi}

With this definition at hand, the main theorem of local geometric class field theory can be stated as follows:

\begin{teo}[Local geometric class field theory]\label{chap3lgcft1} Let $D$ be a closed subscheme of $T$ supported on $s$. Then, the functor $\Phi_{\eta}^{-1}$, which sends an object $(\chi,\widetilde{\chi}, \beta)$ of $\Trip(T,D)$ to the pullback $\Phi_{\eta}^{-1} \widetilde{\chi}$, is an equivalence from the groupoid $\Trip(T,D,\Lambda)$ to the groupoid of $\Lambda$-local systems of rank $1$ on $\eta$, with ramification bounded by $D$.
\end{teo}

We postpone the proof of Theorem \ref{chap3lgcft1} to the paragraph \ref{chap32.2} below. We now provide an equivalent version of Theorem \ref{chap3lgcft1}, whose formulation is somewhat simpler, although non canonical, as it depends on a choice of uniformizer. Let $\pi$ be a uniformizer of $k(\eta)$. Then $1 \otimes \pi - \pi \otimes 1$ is a generator of $\mathcal{I}$, and $1 \otimes \pi$ is a generator of $k(s) \otimes_{\Ow_T} \mathcal{I}$ as an $\Ow_D$-module. Thus we obtain an isomorphism
\begin{align*}
\alpha_{\pi} : \Pic(T,D)_s \times_s T &\rightarrow \Pic(T,D) \\
(d,u) &\mapsto (d, u (1 \otimes \pi)^d (1 \otimes \pi - \pi \otimes 1)^{-d}),
\end{align*}
whose restriction to $\Pic^0(T,D)_s \times_s T $ coincides with $\alpha$. Here, a section $u$ over an $s$-scheme $S$ is considered as an isomorphism 
$$
\Ow_S  \otimes_{k(s)} \Ow_D \rightarrow \Ow_S  \otimes_{k(s)} \Ow_D \otimes_{\Ow_T} \m^{-d},
$$
where $\m$ denotes the defining ideal of $s$ in $T$. If $(\chi,\widetilde{\chi}, \beta)$ is a triple as in Theorem \ref{chap3lgcft1}, then we have a splitting
$$
\Pic(T,D)_s \times_s T \cong (\Pic^0(T,D)_s \times_s T) \times_T \Z_T,
$$
of $T$-group schemes (cf. \ref{chap3exactseq}), and a corresponding decomposition of $\alpha_{\pi}^{-1} \widetilde{\chi}$ as $\alpha^{-1} \widetilde{\chi} \boxtimes \gamma$, where $\gamma$ is a multiplicative $\Lambda$-local system on $\Z_T$ (cf. \ref{chap32.0.4}). The first factor $\alpha^{-1} \widetilde{\chi}$ is isomorphic to $p_1^{-1} \chi$ by $\beta$. Since $T$ is henselian, the pullback by the morphism $T \rightarrow s$ is an equivalence from the groupoid of multiplicative $\Lambda$-local systems on $\Z_s$ to the groupoid of multiplicative $\Lambda$-local systems on $\Z_T$, hence $\gamma$ descends to $\Z_s$ (cf. \ref{chap3zmultloc}). Thus $\alpha_{\pi}^{-1} \widetilde{\chi}$ canonically descends to a multiplicative $\Lambda$-local system on $\Pic(T,D)_s$. If we further note that $\alpha_{\pi}^{-1} \circ \Phi_{\eta}$ is given by the $k(\eta)$-point $(1,u)$ of $\Pic(T,D)_s \times_s T$, where
$$
u = (1 \otimes \pi)^{-1} (1 \otimes \pi - \pi \otimes 1) = 1 - \pi \otimes \pi^{-1},
$$
then we obtain that Theorem \ref{chap3lgcft1} is equivalent to the following:

\begin{teo}[Local geometric class field theory, second version]\label{chap3lgcft2} Let $D$ be a closed subscheme of $T$ supported on $s$, and let $\pi$ be a uniformizer of $k(\eta)$. Let $\Phi_{\eta,\pi} : \eta \rightarrow \Pic(T,D)_s$ be the morphism corresponding to the $k(\eta)$-point $(1,1 - \pi \otimes \pi^{-1})$ of $\Pic(T,D)_s$. Then, the functor $\Phi_{\eta, \pi}^{-1}$ is an equivalence from the groupoid of multiplicative $\Lambda$-local systems on $\Pic(T,D)_s$ to the groupoid of $\Lambda$-local systems of rank $1$ on $\eta$, with ramification bounded by $D$.
\end{teo}

\subsection{\label{chap32.10.1}} Our deduction of Theorem \ref{chap3lgcft2} from Theorem \ref{chap3lgcft1} also shows that the functor
$$
\Trip(T,D,\Lambda) \rightarrow \mathrm{Loc}^{\otimes}(\Pic(T,D)_s,\Lambda),
$$
which sends a triple $(\chi,\widetilde{\chi}, \beta)$ to the restriction of $\widetilde{\chi}$ to the special fiber $\Pic(T,D)_s$ of $\Pic(T,D)$, is an equivalence of groupoids. Moreover, the composition of $\Phi_{\eta,\pi}^{-1}$ with this restriction functor coincides with 
the functor $\Phi_{\eta}^{-1}$ from Theorem \ref{chap3lgcft1}. In particular, the equivalence $\Phi_{\eta,\pi}^{-1}$ from Theorem \ref{chap3lgcft2} does not depend on $\pi$, up to natural isomorphism. We denote by $\F \mapsto \chi_{\F}$ a quasi-inverse to $\Phi_{\eta,\pi}^{-1}$, which is well-defined up to natural isomorphism.

\subsection{\label{chap32.10.2}} If $(\chi_1,\widetilde{\chi}_1, \beta_1)$ and $(\chi_2,\widetilde{\chi}_2, \beta_2)$ are objects of $\Trip(T,D,\Lambda)$ such that $\chi_1 = \chi_2$, then $\Phi_{\eta}^{-1} \widetilde{\chi}_2$ is isomorphic to $\Phi_{\eta}^{-1} \widetilde{\chi}_1 \otimes \G$, where $\G$ is the pullback to $\eta$ of a $\Lambda$-local systems of rank $1$ on $s$. We thus obtain a simpler (although weaker) version of Theorem \ref{chap3lgcft1} by ignoring twists by unramified $\Lambda$-local systems of rank $1$ on $\eta$:

\begin{teo}[Local geometric class field theory, third version]\label{chap3lgcft3} Let $D$ be a closed subscheme of $T$ supported on $s$. If $\F$ is a $\Lambda$-local system of rank $1$ on $\eta$, with ramification bounded by $D$, then there exists a unique (up to isomorphism) multiplicative $\Lambda$-local system $\chi$ on the $s$-group scheme $\Pic^0(T,D)_s$, such that the $\Lambda$-local system $\chi \boxtimes \F$ on the product
$$
\Pic^0(T,D)_s \times_s \eta \xrightarrow[\sim]{\alpha \cdot \Phi_{\eta}} \Pic^1(T,D)_{\eta},
$$
extends to a $\Lambda$-local system on $\Pic^1(T,D)$. This provides a bijection from the group of isomorphism classes of $\Lambda$-local systems of rank $1$ on $\eta$ with ramification bounded by $D$, up to twist by unramified $\Lambda$-local systems of rank $1$ on $\eta$, to the group of isomorphism classes of multiplicative $\Lambda$-local systems on the $s$-group scheme $\Pic^0(T,D)_s$.
\end{teo}

One should note that the restriction functor from $\Pic^1(T,D)$ to $\Pic^1(T,D)_{\eta}$ realizes an equivalence between the groupoid of $\Lambda$-local systems on $\Pic^1(T,D)$ and a full subcategory of the groupoid of $\Lambda$-local systems on $\Pic^1(T,D)_{\eta}$. The formulation of Theorem \ref{chap3lgcft3} is due to Gaitsgory, and can also be found in Bhatt's Oberwolfach report (\cite{Bhatt}, Th. 11).

\subsection{\label{chap32.5}} We now describe the relation between our version of local geometric class field theory, namely Theorem \ref{chap3lgcft1}, and Contou-Carrere's theory of the local Jacobian. Let $\m$ be the defining ideal of $s$, so that $D$ is defined by $\m^{\nu}$ for some nonnegative integer $\nu$. Then $k(s) \otimes_{\Ow_T} \mathcal{I}^{-d}$ is naturally isomorphic to $\m^{-d}/\m^{-d+\nu}$ as a module over $\Ow_D = \Ow_T/\m^\nu$, and we can identify $\Pic(T,D)_s$ with the functor which sends an $s$-scheme $S$ to the group of pairs $(d,u)$, where $d$ is a locally constant $\mathbb{Z}$-valued map on $S$, and 
$$u : \Ow_S \otimes_{k(s)} \Ow_T/ \m^\nu \rightarrow \Ow_S \otimes_{k(s)} \m^{-d} / \m^{-d + \nu} $$
is an isomorphism of $\Ow_S \otimes_{k(s)} \Ow_T/ \m^{\nu}$-modules.
The inverse limit over $\nu$ of these $s$-group schemes can then be identified with the functor $J(\eta)$ which sends a $k(s)$-algebra $A$ to the group of units in $A \widehat{\otimes}_{k(s)} k(\eta) = (A \widehat{\otimes}_{k(s)} \Ow_T) \otimes_{\Ow_T} k(\eta)$, where $A \widehat{\otimes}_{k(s)} \Ow_T$ is the $\m$-adic completion of $A \otimes_{k(s)} \Ow_T$, which generate a sub-$(A \widehat{\otimes}_{k(s)} \Ow_T)$-module of the form $\m^{-d} (A \widehat{\otimes}_{k(s)} \Ow_T)$ for some locally constant function $d$ on $\Spec(A)$. 

\begin{rema} A unit of $A \widehat{\otimes}_{k(s)} k(\eta)$ may not necessarily belong to $J(\eta)(A)$. For example, if $a$ is a nilpotent element of $A$ and if $u$ is an element of $k(\eta)$ which does not belong to $\Ow_T$, then $1-a \otimes u$ is a unit of $A \widehat{\otimes}_{k(s)} k(\eta)$ which does not belong to $J(\eta)(A)$. However, if $A$ is reduced, then $J(\eta)(A)$ is simply the group of units in $A \widehat{\otimes}_{k(s)} k(\eta)$ by (\cite{CC}, 0.8) or (\cite{G18}, Prop. 3.4).
\end{rema}

\begin{prop} The functor $J(\eta)$ is representable by an $s$-group scheme, which fits into an exact sequence
$$
1 \rightarrow J(\eta)^0 \rightarrow J(\eta) \rightarrow \Z_s \rightarrow 0,
$$
where $J(\eta)^0$ is representable by a geometrically connected affine $s$-group scheme.
\end{prop}

Indeed, $J(\eta)^0$ is the limit of the inverse system $(\Pic^0(T,\m^{\nu})_s)_{\nu \geq 0}$ of affine $s$-group schemes.

\begin{teo}[Local geometric class field theory, fourth version]\label{chap3lgcft4} Let $\pi$ be a uniformizer of $k(\eta)$. Let $\Psi_{\eta,\pi} : \eta \rightarrow J(\eta)$ be the morphism corresponding to the $k(\eta)$-point $(1,1 - \pi \otimes \pi^{-1})$ of $J(\eta)$. Then, the functor $\Psi_{\eta, \pi}^{-1}$ is an equivalence from the groupoid of multiplicative $\Lambda$-local systems on $J(\eta)$ to the groupoid of $\Lambda$-local systems of rank $1$ on $\eta$.
\end{teo}

If $\Lambda$ is finite, then Theorem \ref{chap3lgcft4} follows from Theorem \ref{chap3lgcft2}, since the category of finite etale $J(\eta)$-schemes is the $2$-colimit of the categories of finite etale $\Pic(T,\m^{\nu})_s$-schemes when $\nu$ ranges over all integers, and thus the groupoid $\mathrm{Loc}^{\otimes}(J(\eta),\Lambda)$ is the $2$-limit of the groupoids $\mathrm{Loc}^{\otimes}(\Pic(T,\m^{\nu})_s,\Lambda)$. If $\Lambda$ is the ring of integers in a finite extension of $\mathbb{Q}_{\ell}$, then the conclusion of Theorem \ref{chap3lgcft4} holds for the finite $\ell$-adic coefficient rings $\Lambda/\ell^n$ for each $n$, and thus for $\Lambda$ as well by taking $2$-limits. This implies the validity of Theorem \ref{chap3lgcft4} when $\Lambda$ is a finite extension of $\mathbb{Q}_{\ell}$, and by taking $2$-colimits this yields the result when $\Lambda$ is an arbitrary $\ell$-adic coefficient ring.


The $s$-group scheme $J(\eta)$ coincides Contou-Carrere's local Jacobian, and the morphism $\Psi_{\eta,\pi}$ in Theorem \ref{chap3lgcft3} is the morphism studied by Contou-Carrere or considered by Deligne in his 1974 letter to Serre (\cite{bloch}, p.74). Contou-Carrere established an Albanese property for the morphism $\Psi_{\eta,\pi}$, which was used by Suzuki (\cite{TS}, Th. A (1)) in order to give a different proof of Theorem \ref{chap3lgcft4}. Moreover, Suzuki (op. cit.) showed that the equivalence constructed by Serre in \cite{JPS} when $k$ is algebraically closed, is a quasi-inverse to the equivalence in Theorem \ref{chap3lgcft4}.


\subsection{\label{chap32.2}} We now prove Theorem \ref{chap3lgcft1}, by combining the Gabber-Katz extension theorem \ref{chap31.5} with global geometric class field theory, namely Theorem \ref{chap3ggcft}. More precisely, we prove its equivalent version \ref{chap3lgcft2}.
Let $D$ be a closed subscheme of $T$ supported on $s$, and let $\pi$ be a uniformizer of $k(\eta)$. The uniformizer $\pi$ provides a morphism $k(s)[t,t^{-1}] \rightarrow k(\eta)$ sending $t$ to $\pi$, corresponding to a morphism
$$
\pi : \eta \rightarrow \mathbb{G}_{m,s},
$$
of $s$-schemes.

By Theorem \ref{chap3GK4}, the restriction of the pullback functor $\pi^{-1}$ to the category of special $\Lambda$-sheaves on $\mathbb{A}^1_{s}$ vanishing at $0$ is an equivalence with the category of $\Lambda$-sheaves on $\eta$ (cf. \ref{chap31.5}). Let $\pi_{\diamondsuit}$ be a quasi-inverse to this equivalence. Let $\F$ be a $\Lambda$-sheaf on $\eta$ with ramification bounded by $D$. Then $\pi_{\diamondsuit} \F$ is a $\Lambda$-local system on the open subscheme $\mathbb{G}_{m,s}$ of $\mathbb{P}^1_{s}$, extended by zero at $0$ and $\infty$, with ramification bounded by the divisor $D' = D + [\infty]$. Let us consider the Abel-Jacobi morphism
$$
\Phi : \mathbb{G}_{m,s} \rightarrow \Pic_{s}(\mathbb{P}^1_{s},D'),
$$
which sends a section $t$ of $\mathbb{G}_{m,s}$ to the pair $(\Ow(t),1)$, cf. (\ref{chap3abeljacob}). Since the Picard scheme of $\mathbb{P}^1_{s}$ is the constant group scheme $\mathbb{Z}_{s}$, we can identify the connected component of degree $d$ of $\Pic_{s}(\mathbb{P}^1_{s},D')$ with the functor which to an $s$-scheme $S$ associates the quotient by $ \mathbb{G}_{m,s}(S)$ of the group of isomorphisms $\alpha : \Ow_{D'_S} \rightarrow i_S^{'*} \Ow(d[0])$, where $i' : D' \rightarrow \mathbb{P}^1_{s}$ is the inclusion. The latter functor can be further identified with the functor which to an $s$-scheme $S$ associates the group of isomorphisms $\alpha : \Ow_{D_S} \rightarrow i_S^{*} \Ow(d[0])$, where $i : D \rightarrow \mathbb{P}^1_{s}$ is the inclusion. Using the uniformizer $\pi$, we obtain an isomorphism
$$
\theta : \Pic_{s}(\mathbb{P}^1_{s},D') \rightarrow \Pic(T,D)_s.
$$
Consequently, if $\F$ has rank $1$, then Theorem \ref{chap3ggcft} implies that $\pi_{\diamondsuit} \F$ is isomorphic to the pullback by $\theta \circ \Phi$ of a multiplicative $\Lambda$-local system on $\Pic(T,D)_s$. We therefore deduce from Theorems \ref{chap3GK} and \ref{chap3ggcft} that the pullback by $\theta \circ \Phi \circ \pi$ induces an equivalence from the groupoid of multiplicative $\Lambda$-local systems on $\Pic(T,D)_s$ to the groupoid of $\Lambda$-local systems of rank $1$ on $\eta$, with ramification bounded by $D$.

It remains to check that the composition $\theta \circ \Phi \circ \pi$ coincides with the morphism $\Phi_{\eta,\pi}$ in Theorem \ref{chap3lgcft2}. If $t$ is a section of $\mathbb{G}_{m,s}$ over an $s$-scheme $S$, then the isomorphism $\Ow([t]) \rightarrow \Ow([0])$ given by multiplication by $1 - t x^{-1}$, where $x$ is the coordinate on $\mathbb{G}_{m,s}$, sends the canonical trivialization $1 : \Ow_{D'_S} \rightarrow i_S^{\prime *} \Ow([t])$ to the trivialization $\alpha : \Ow_{D'_S} \rightarrow i_S^{\prime  *} \Ow([0])$ corresponding to $1-t x^{-1}$. Thus $\theta \circ \Phi$ sends $t$ to the $S$-point $\Pic(T,D)_s$ corresponding to $1 - t \otimes \pi^{-1}$. By taking $S= \eta$ and $t = \pi$, we obtain that $\theta \circ \Phi \circ \pi$ coincides with $\Phi_{\eta,\pi}$. This concludes our proof of Theorem \ref{chap3lgcft2}, which in turn implies Theorems \ref{chap3lgcft1}, \ref{chap3lgcft3} and \ref{chap3lgcft4}.

\begin{rema} This proof of the main theorem of local geometric class field theory \ref{chap3lgcft1} uses global geometric class field theory. The latter admits geometric proofs which do not use the local theory, cf. for example \cite{T18}, hence the argument is not circular. Moreover, the use of local geometric class field theory in (\cite{G18}, Prop. 3.14) can be avoided by resorting to a computation with Artin-Schreier-Witt theory, as in \cite{T18}.
\end{rema}

\subsection{\label{chap32.6}} In this paragraph, we describe the compatibility between local and global geometric class field theory, namely Theorems \ref{chap3ggcft} and \ref{chap3lgcft1}. Let $X$ be a smooth geometrically connected projective curve of genus $g$ over $k$, let $i : D \rightarrow X$ be an effective Cartier divisor on $X$, and let $U$ be the open complement of $D$ in $X$. We introduced in \ref{chap32.1.0.0} the generalized Picard group scheme $\Pic_k(X,D)$, and in (\ref{chap3abeljacob}) the Abel-Jacobi morphism
$$
\Phi : U \rightarrow \Pic_k(X,D),
$$
which sends a section $x$ of $U$ to the pair $(\Ow(x),1)$, where $1 : \Ow_D \rightarrow \Ow(x) \otimes_{\Ow_X} \Ow_D$ is the canonical trivialization of $\Ow(x)$ on $D$, cf. (\ref{chap3abeljacob}).

Let $x$ be a point of $D$, and let $X_{(x)}$ be the henselisation of $X$ at $x$, with generic point $\eta_x$. We identify the closed point of $X_{(x)}$ with $x$, and we denote by $D_x$ the pullback of $D$ to $X_{(x)}$, which is a closed subscheme of $X_{(x)}$ supported on $x$. Let $\widetilde{x}$ be the $X_{(x)}$-point of $X \times_k X_{(x)}$ given by the diagonal embedding. The restriction of $\widetilde{x}$ to $\eta_x$ factors through $U \times_k X_{(x)}$. We now define a morphism
$$
\tau : \Pic(X_{(x)},D_x) \rightarrow \Pic_k(X,D)\times_k X_{(x)}
$$
of $X_{(x)}$-group schemes as follows. Let $(d,u)$ be a point of $\Pic(X_{(x)},D_x)$ over an $X_{(x)}$-scheme $S$ (cf. \ref{chap3localpic}). The image in $\Ow_D \otimes_k \Ow_S$ of the kernel of the natural multiplication homomorphism 
$$
\Ow_X \otimes_k \Ow_S \rightarrow \Ow_S,
$$
is the ideal $\mathcal{I}$ as in Proposition \ref{chap3lemmadiag}. Thus the pullback to $D \times_k S$ of the line bundle $\Ow(d\widetilde{x})$ on the $S$-curve $X \times_k S$ is given by the invertible module $\Ow_S \otimes_{\Ow_{X_{(x)}}} \mathcal{I}^{-d}$. Consequently, $u$ provides a trivialization of $\Ow(d\widetilde{x})$ on the divisor $D_x \times_k S$ of $X \times_k S$. Moreover, the canonical section $1 : \Ow_{X \times_k S} \rightarrow \Ow(d\widetilde{x})$ provides a trivialization of $\Ow(d\widetilde{x})$ on the divisor $(D \setminus D_x) \times_k S$. We thus obtain a trivialization $\beta : \Ow_{D \times_k S} \rightarrow (i \times \mathrm{id}_S)^* \Ow(d\widetilde{x})$, and the pair $(\Ow(d\widetilde{x}),\beta)$ defines an $S$-point of $\Pic_k(X,D)\times_k X_{(x)}$. This construction is functorial in $S$, and thus defines a morphism $\tau$ as above.

We also let $\tau_x : \Pic(X_{(x)},D_x)_x \rightarrow \Pic_k(X,D) $ be the restriction of $\tau$ to the special fiber. We then have the following commutative diagram.
\begin{center}
 \begin{tikzpicture}[scale=1]

\node (A) at (6,0) {$\Pic^0(X_{(x)},D_x)_x \times_x X_{(x)}$};
\node (B) at (8,2) {$\Pic^0(X_{(x)},D_x)_x$};
\node (C) at (0,2) {$\eta_x$};
\node (D) at (4,2) {$\Pic(X_{(x)},D_x)$};
\node (E) at (0,4) {$U \times_k X_{(x)}$};
\node (F) at (4,4) {$\Pic_k(X,D)\times_k X_{(x)}$};
\node (G) at (8,4) {$\Pic_k(X,D)$};
\path[->,font=\scriptsize]
(A) edge node[above left]{$p_1$} (B)
(B) edge node[left]{$\tau_x$} (G)
(A) edge node[above right]{$ \alpha$} (D)
(C) edge node[above]{$\Phi_{\eta}$} (D)
(C) edge node[left]{$\widetilde{x}_{|\eta_x}$} (E)
(D) edge node[left]{$\tau$} (F)
(E) edge node[above]{$\Phi \times \mathrm{id}$} (F)
(F) edge node[above]{$p_1$} (G);
\end{tikzpicture} 
\end{center}

The morphisms $\alpha$ and $\Phi_{\eta}$ in this diagram are defined in \ref{chap32.10}, while $p_1$ always denotes the projection onto the first factor. Let $\mathrm{Loc}_1(U,D,\Lambda)$ (resp. $\mathrm{Loc}_1(\eta_x,D_x,\Lambda)$) be the groupoid of $\Lambda$-local systems of rank $1$ on $U$ with ramification bounded by $D$ (resp. on $\eta_x$ with ramification bounded by $D_x$). If $\Lc$ is an object of $\mathrm{Loc}^{\otimes}(\Pic_k(X,D),\Lambda)$, we obtain an object $(\widetilde{\chi},\chi,\theta)$ of $\Trip(X_{(x)},D_x,\Lambda)$ (cf. \ref{chap3lgcft0}) as follows: we set $\widetilde{\chi} = (p_1 \circ \tau )^{-1} \Lc $ and $\chi = \tau_x^{-1} \Lc$, while $\theta : \alpha^{-1} \widetilde{\chi} \rightarrow p_1^{-1} \chi$ is the natural isomorphism resulting from the commutativity of diagram above. We thus obtain a functor from $\mathrm{Loc}^{\otimes}(\Pic_k(X,D),\Lambda)$ to $\Trip(X_{(x)},D_x,\Lambda)$, which we abusively denote by $\tau^{-1}$ for simplicity. This functor $\tau^{-1}$ fits into the following diagram, which is commutative up to natural isomorphism.

\begin{center}
 \begin{tikzpicture}[scale=1]

\node (A) at (0,0) {$\mathrm{Loc}_1(\eta_x,D_x,\Lambda)$};
\node (B) at (4,0) {$\Trip(X_{(x)},D_x,\Lambda)$};
\node (C) at (0,2) {$\mathrm{Loc}_1(X,D,\Lambda)$};
\node (D) at (4,2) {$\mathrm{Loc}^{\otimes}(\Pic_k(X,D),\Lambda)$};

\path[->,font=\scriptsize]
(B) edge node[above]{$\Phi_{\eta}^{-1}$} (A)
(C) edge node[left]{} (A)
(D) edge node[left]{$\tau^{-1}$} (B)
(D) edge node[above]{$\Phi^{-1}$} (C);
\end{tikzpicture} 
\end{center}

The rows of this diagram are equivalences of groupoids by Theorem \ref{chap3ggcft} and \ref{chap3lgcft1}. Thus the restriction functor, which is the left vertical arrow in this diagram, corresponds to the functor $\tau^{-1}$ in terms of multiplicative $\Lambda$-local systems.

\subsection{\label{chap32.19}} We now describe the functoriality property of geometric local class field theory. Let $T$ (resp. $T'$) be the spectrum of a $k$-algebra, which is a henselian discrete valuation ring $\Ow_T$ (resp. $\Ow_{T'}$), whose residue field is a finite extension of $k$. Let $\eta$ (resp. $\eta'$) be the generic point of $T$ (resp. $T'$), and let $s$ (resp. $s'$) be its closed point, so that $k(\eta)$ (resp. $k(\eta')$) is a henselian discrete valuation field, with valuation subring $\Ow_{T,s}$ (resp. $\Ow_{T',s'}$), and with residue field $k(s)$ (resp. $k(s')$) which is a finite extension of $k$.

Let $f : T' \rightarrow T$ be a finite morphism, such that the restriction $f_{|\eta'} : \eta' \rightarrow \eta$ is \'etale, namely such that the finite extension $k(\eta) \rightarrow k(\eta')$ induced by $f$ is separable. Let $D$ be a closed subscheme of $T$ supported on $s$, and let $D'$ be its pullback to $T'$. Let $\mathcal{I}$ (resp. $\mathcal{I}'$) be the kernel of the homomorphism
\begin{align*}
\Ow_T \otimes_{k(s)} \Ow_D &\rightarrow \Ow_D \\
f_1 \otimes f_2 &\rightarrow f_1 f_2,
\end{align*}
(resp. of the homomorphism $\Ow_{T'} \otimes_{k(s')} \Ow_{D'} \rightarrow \Ow_{D'}$), which, by Proposition \ref{chap3lemmadiag}, is a principal invertible ideal of $\Ow_T \otimes_{k(s)} \Ow_D$ (resp. $\Ow_{T'} \otimes_{k(s')} \Ow_{D'}$) generating the unit ideal of $k(\eta) \otimes_{k(s)} \Ow_D$ (resp. $k(\eta') \otimes_{k(s')} \Ow_{D'}$).

For any $T'$-scheme $S$, the $\Ow_{S} \otimes_{k(s)} \Ow_D$-algebra $\Ow_S \otimes_{k(s')} \Ow_{D'}$ is free of finite rank equal to the ramification index $e_f$ of the extension $k(\eta')/k(\eta)$, hence we can consider the norm map
\begin{align}\label{chap3norm}
N_{f} : \Ow_S \otimes_{k(s')} \Ow_{D'} \rightarrow \Ow_{S} \otimes_{k(s)} \Ow_D,
\end{align}
which sends a section $u$ of $ \Ow_S \otimes_{k(s')} \Ow_{D'} $ to the determinant of the $\Ow_{S} \otimes_{k(s)} \Ow_D$-linear endomorphism $x \mapsto ux$ of $ \Ow_S \otimes_{k(s')} \Ow_{D'} $.
The norm map $N_f$ is homogeneous of degree $e_f$, and therefore the image by $N_f$ of a principal ideal of $\Ow_S \otimes_{k(s')} \Ow_{D'}$ generates a principal ideal of $\Ow_{S} \otimes_{k(s)} \Ow_D$. 

\begin{lem}\label{chap3lemanorm} The ideal of $\Ow_{S} \otimes_{k(s)} \Ow_D$ generated by the image by $N_f$ of the ideal $\Ow_S \otimes_{\Ow_{T'}} \mathcal{I}'$ of $\Ow_S \otimes_{k(s')} \Ow_{D'}$ is $\Ow_S \otimes_{\Ow_{T}} \mathcal{I}$.
\end{lem}

Indeed, if $\pi'$ is a uniformizer of $k(\eta')$ then the ideal $\Ow_S \otimes_{\Ow_{T'}} \mathcal{I}'$ of $\Ow_S \otimes_{k(s')} \Ow_{D'}$ is generated by $\pi' \otimes 1 - 1 \otimes \pi' $. If $P_{\pi'}(X) = X^{e_f} + a_{1} X^{e_f -1 } + \cdots + a_{e_f}$ is the characteristic polynomial of $\pi'$ in the totally ramified extension $k(\eta')/k(\eta_{s'})$, where $\eta_{s'} = \eta \times_s s'$, then the ideal generated by $N_f(\Ow_S \otimes_{\Ow_{T'}} \mathcal{I}')$ is generated by the element
$$
(1 \otimes P_{\pi'})( \pi' \otimes 1 ) = \pi'^{e_f} \otimes 1 + \pi'^{e_f-1} \otimes a_1 + \dots + 1 \otimes a_{e_f},
$$
of $\Ow_{S} \otimes_{k(s')} \Ow_{D_{s'}}$. Since $P_{\pi'}$ is an Eisenstein polynomial, the elements $(a_i)_{i \leq e_f}$ of $\Ow_{T_{s'}}$ belong to the maximal ideal, and $a_{e_f}$ is a uniformizer of $k(\eta_{s'})$. In particular, we can write $a_i = b_i a_{e_f}$, for some elements $(b_i)_{i \leq e_f}$ of $\Ow_{T_{s'}}$. We obtain a decomposition
$$
(1 \otimes P_{\pi'})( \pi' \otimes 1 ) = P_{\pi'}(\pi') \otimes 1 + (u \otimes 1)(1 \otimes a_{e_f} - a_{e_f} \otimes 1) + (1 \otimes a_{e_f}) v,
$$
where we have set
\begin{align*}
u &= 1 + \sum_{j=1}^{e_f - 1} b_{e_f -j} \pi'^j \\
v &= \sum_{j = 1}^{e_f} (\pi'^{e_f - j} \otimes 1)( 1 \otimes b_{j} - b_j \otimes 1).
\end{align*}
The term $P_{\pi'}(\pi')$ vanishes by the Cayley-Hamilton theorem. The elements $1 \otimes a_{e_f} - a_{e_f} \otimes 1$ and $v$ belong to the ideal $\Ow_S \otimes_{\Ow_{T}} \mathcal{I}$, while $1 \otimes a_{e_f} - a_{e_f} \otimes 1$ generates the latter. Moreover, the elements $u \otimes 1$ and $1 \otimes a_{e_f}$ are respectively invertible and nilpotent in $\Ow_{S} \otimes_{k(s)} \Ow_D$. Thus $(1 \otimes P_{\pi'})( \pi' \otimes 1 )$ is a generator of $\Ow_S \otimes_{\Ow_{T}} \mathcal{I}$, hence the result.

\begin{defi} \label{chap3norm2} The \textit{norm morphism} associated to $f$ is the homomorphism
$$
N_f : \Pic(T',D') \rightarrow \Pic(T,D) \times_T T',
$$
of $T'$-group schemes (cf. \ref{chap3localpic}) which sends a section $(d,u)$ of $\Pic(T',D')$ over a $T'$-scheme $S$ to the $S$-point $(d, N_f(u))$ of $\Pic(T,D)$, where $N_f(u)$ is the trivialization of the invertible $\Ow_{S} \otimes_{k(s)} \Ow_D$-module generated by $N_f(\Ow_S \otimes_{\Ow_{T'}} \mathcal{I}'^{-d})$, which coincides with $\Ow_S \otimes_{\Ow_{T}} \mathcal{I}^{-d}$ by Lemma \ref{chap3lemanorm}, obtained by applying the norm map (\ref{chap3norm}) to $u$.
\end{defi}

Let $\Phi_{\eta}$ (resp. $\Phi_{\eta'}$) be the local Abel-Jacobi morphism for $(T,D)$ (resp. $(T,D')$), cf. \ref{chap32.10}. The following diagram is clearly commutative.

\begin{center}
 \begin{tikzpicture}[scale=1]

\node (A) at (0,0) {$\eta$};
\node (B) at (3,0) {$\Pic(T,D)$};
\node (C) at (0,1.854) {$\eta'$};
\node (D) at (3,1.854) {$\Pic(T',D')$};

\path[->,font=\scriptsize]
(A) edge node[above]{$\Phi_{\eta}$} (B)
(C) edge node[left]{$f_{|\eta'}$} (A)
(D) edge node[left]{$N_f$} (B)
(C) edge node[above]{$\Phi_{\eta'}$} (D);
\end{tikzpicture} 
\end{center}

Correspondingly, the following diagram is commutative up to natural isomorphism.

\begin{center}
 \begin{tikzpicture}[scale=1]

\node (A) at (0,0) {$\mathrm{Loc}_1(\eta,D,\Lambda)$};
\node (B) at (3.3,0) {$\Trip(T,D,\Lambda)$};
\node (C) at (0,1.854) {$\mathrm{Loc}_1(\eta',D',\Lambda)$};
\node (D) at (3.3,1.854) {$\Trip(T',D',\Lambda)$};

\path[->,font=\scriptsize]
(B) edge node[above]{$\Phi_{\eta}^{-1}$} (A)
(A) edge node[left]{$f_{|\eta'}^{-1}$} (C)
(B) edge node[left]{$N_f^{-1}$} (D)
(D) edge node[above]{$\Phi_{\eta'}^{-1}$} (C);
\end{tikzpicture} 
\end{center}

The rows of this diagram are equivalences of groupoids by Theorem \ref{chap3lgcft0}.

 \subsection{\label{chap32.20}} Let $\bk$ be an algebraic closure of $k$. We denote by $G_k = \Gal(\bk/k)$ the Galois group of the extension $\bk/k$ and by $\mu$ a unitary $\Lambda$-admissible mutiplier on the topological group $G_k$ (cf. \ref{chap30.0}, \ref{chap3unitary}). Let $k'/k$ be a neutralizing extension of $k$ contained in $\bk$, cf. \ref{chap3neutra}.

\begin{defi}\label{chap3twistedmult} Let $S$ be a connected $k$-scheme and let $G$ be a commutative $S$-group scheme, which fits into an exact sequence
$$
1 \rightarrow G^0 \rightarrow G \xrightarrow[]{d} \Z_S \rightarrow 0,
$$
 where $G^0$ is an $S$-group scheme with connected geometric fibers. A \textbf{$\mu$-twisted multiplicative $\Lambda$-local system on $G$}, is a pair $\Lc = (\Lc_{k'}, (\rho_{\Lc}(g))_{g \in \Gal(k'/k)})$, where $\Lc_{k'}$ is a multiplicative $\Lambda$-local system on the $S_{k'}$-group scheme $G_{k'}$ (cf. \ref{chap31.0.4}), and $\rho_{\Lc}(g) : g^{-1} \Lc_{k'} \rightarrow \Lc_{k'}$ is an isomorphism of multiplicative $\Lambda$-local systems for each $g$ in $\Gal(k'/k)$, such that the diagram
\begin{center}
 \begin{tikzpicture}[scale=1]

\node (A) at (0,0) {$g^{-1} h^{-1} \Lc_{k'} $};
\node (B) at (3,0) {$g^{-1} \Lc_{k'}$};
\node (C) at (6,0) {$\Lc_{k'}$};

\path[->,font=\scriptsize]
(A) edge node[above]{$g^{-1} \rho_{\Lc}(h)$} (B)
(B) edge node[above]{$\rho_{\Lc}(g)$} (C)
(A) edge[bend right] node[above]{$ \mu(g,h)^d \rho_{\Lc}(gh)$} (C);
\end{tikzpicture} 
\end{center}
is commutative for any $g,h$ in $\Gal(k'/k)$. Here $\mu(g,h)^d$ is the section of $\Lambda^{\times}$ on $G$ which is constant equal to $\mu(g,h)^r$ on the inverse image of $r$ by the given homomorphism $d :G \rightarrow \Z_S$, for each integer $r$.

 If $\Lc$ and $\mathcal{M}$ are $\mu$-twisted multiplicative $\Lambda$-local systems on $G$, a morphism from $\Lc$ to $\mathcal{M}$ is an isomorphism $f : \Lc_{k'} \rightarrow \mathcal{M}_{k'}$ of multiplicative $\Lambda$-local systems such that $f \circ \rho_{\Lc}(g) = \rho_{\mathcal{M}}(g) \circ (g^{-1}f)$ for any $g$ in $\Gal(k'/k)$.
\end{defi}

\begin{rema}\label{chap3sanity2} If $S$ is a also a $k''$-scheme, for some finite extension $k''$ of $k$, and if $\iota : k'' \rightarrow \bk$ is a $k$-linear embedding, with Galois group $\Gal(\iota)$, then as in \ref{chap3sanity}, the groupoid of $\mu$-twisted multiplicative $\Lambda$-local systems on $G$, where $S$ is considered as a $k$-scheme, is equivalent to the groupoid of $\mu_{|\Gal(\iota)}$-twisted multiplicative $\Lambda$-local systems on the $S$-group scheme $G$, where $S$ is now considered as a $k''$-scheme.
\end{rema}

\subsection{\label{chap32.21}} Let $T, \eta$ and $s$ be as in \ref{chap32.3}, and let $\Sigma$ be the set of $k$-linear embeddings of $k(s)$ in $\bk$, and assume that $k'$ contains the image of any element of $\Sigma$. We then have a decomposition
$$
T_{k'} = \coprod_{\iota \in \Sigma} T_{\iota},
$$
where $T_{\iota} = T \otimes_{k(s),\iota} k'$ is the spectrum of an henselian discrete valuation $k'$-algebra with residue field $k'$. If $\F = (\F_{k'}, (\rho_{\F}(g))_{g \in \Gal(k'/k)})$ is a $\mu$-twisted $\Lambda$-local system of rank $1$ on $\eta$, then the Swan conductor of the restriction $\F_{\iota}$ of $\F_{k'}$ to $T_{\iota}$ is independent of $\iota$.

\begin{defi}\label{chap3ramdefiloc2} Let $D$ be a closed subscheme of $T$ supported on $s$. A $\mu$-twisted $\Lambda$-local system $\F = (\F_{k'}, (\rho_{\F}(g))_{g \in \Gal(k'/k)})$ of rank $1$ on $\eta$ has \textit{ramification bounded by $D$} if for each $\iota$ in $\Sigma$ (or equivalently, some $\iota$), the Swan conductor of $\F_{\iota}$ is strictly less than the multiplicity of $D$ at $s$.
\end{defi}

\begin{teo}[Local geometric class field theory, twisted version]\label{chap3lgcft2twisted} Let $D$ be a closed subscheme of $T$ supported on $s$, and let $\pi$ be a uniformizer of $k(\eta)$. Let $\Phi_{\eta,\pi} : \eta \rightarrow \Pic(T,D)_s$ be the morphism corresponding to the $k(\eta)$-point $(1,1 - \pi \otimes \pi^{-1})$ of $\Pic(T,D)_s$. Then, the functor $\Phi_{\eta, \pi}^{-1}$ is an equivalence from the groupoid of $\mu$-twisted multiplicative $\Lambda$-local systems on $\Pic(T,D)_s$ (cf. \ref{chap3twistedmult}) to the groupoid of $\mu$-twisted $\Lambda$-local systems of rank $1$ on $\eta$, with ramification bounded by $D$.
\end{teo}

This follows immediately from Theorem \ref{chap3lgcft2}, and from the functoriality of local geometric class field theory, cf. \ref{chap32.19}.

\subsection{\label{chap32.22}} Let $X$ be a smooth geometrically connected projective curve over $k$, let $i : D \rightarrow X$ be an effective Cartier divisor on $X$, and let $U$ be the open complement of $D$ in $X$. 

\begin{defi}\label{chap3ramdefi3} A $\mu$-twisted $\Lambda$-local system $\F$ of rank $1$ on $U$ has \textit{ramification bounded by $D$} if for any point $x$ of $D$, the restriction of $\F$ to the generic point of the henselization $X_{(x)}$ of $X$ at $x$ has ramification bounded by the restriction of $D$ to $X_{(x)}$ (cf. \ref{chap3ramdefiloc2}).
\end{defi}

Let us consider the Abel-Jacobi morphism
\begin{align*}
\Phi : U \rightarrow \Pic_k(X,D),
\end{align*}
which sends a section $x$ of $U$ to the pair $(\Ow(x),1)$, cf. \ref{chap3abeljacob}. As in \ref{chap32.21}, we have the following twisted version of the main theorem of global geometric class field theory:

\begin{teo}[Global geometric class field theory, twisted version]\label{chap3ggcfttwisted} Let $X,U,D$ be as in \ref{chap32.22}. Then the pullback $\Phi^{-1}$ by the Abel-Jacobi morphism of $(X,D)$ realizes an equivalence from the groupoid of $\mu$-twisted multiplicative $\Lambda$-local systems on $\Pic_k(X,D)$ (cf. \ref{chap3twistedmult}) to the category of groupoid of $\mu$-twisted $\Lambda$-local systems of rank $1$ on $U$, with ramification bounded by $D$ (cf. \ref{chap3ramdefi3}).
\end{teo}

\section{Extensions of additive groups\label{chap3extensions}}

Let $A$ be a perfect $\mathbb{F}_p$-algebra, and let $S$ be its spectrum. In this section, we denote by $S_{\fppf}$ the topos of sheaves of sets on the small $\fppf$ site of $S$, cf. (\cite{SGA4}, VII.4.2), and by $\abx$ the category of commutative group objects in $S_{\fppf}$. The purpose of this section is to study extensions of the group scheme $\Ga_{a,S}$ by a finite abelian group $\Gamma$, or equivalently to study short exact sequences
$$
0 \rightarrow \Gamma \rightarrow E \rightarrow \Ga_{a,S} \rightarrow 0,
$$
of abelian groups in $\abx$, where $\Gamma$ is considered as a constant $S$-group scheme. In such an exact sequence, the action of $\Gamma$ on $E$ by left multiplication turns $E$ into a left $\Gamma$-torsor over $\Ga_{a,S}$, hence $E$ is representable by a finite \'etale $\Ga_{a,S}$-scheme.

\subsection{\label{chap34.1}} Let $\Gamma$ be a finite abelian group, and let $G$ be a finitely presented $S$-group scheme with geometrically connected fibers. The extensions of $G$ by $\Gamma$ in $\abx$ are classified by the elements of the abelian group $\Ext_{\abx}(G, \Gamma)$.

\begin{prop}\label{chap3inj} Let $i : \Gamma' \rightarrow \Gamma$ be an injective homomorphism of finite abelian groups and let $G$ be a finitely presented $S$-group scheme with geometrically connected fibers. Then the natural homomorphism
$$
\Ext_{\abx}(G, \Gamma') \xrightarrow[]{i} \Ext_{\abx}(G, \Gamma),
$$
is injective.
\end{prop}

Indeed, we have an exact sequence
$$
\Hom_{\abx}(G, \Gamma/\Gamma') \rightarrow \Ext_{\abx}(G, \Gamma') \rightarrow \Ext_{\abx}(G, \Gamma),
$$
whose first term vanishes since $G$ has geometrically connected fibers over $S$.

\begin{prop}\label{chap3ppart} Let $G$ be a finitely presented $S$-group scheme with geometrically connected fibers, annihilated by an integer $n \geq 1$. Let $\Gamma$ be a finite abelian group and let $\Gamma[n]$ be its subgroup of $n$-torsion elements. Then the natural homomorphism
$$
\Ext_{\abx}(G, \Gamma[n]) \rightarrow \Ext_{\abx}(G, \Gamma),
$$
is an isomorphism.
\end{prop}

Indeed, if $\Gamma' \subseteq \Gamma$ is the image of the multiplication by $n$ in $\Gamma$, then we have a short exact sequence
$$
0 \rightarrow \Gamma[n] \rightarrow \Gamma \xrightarrow[]{n} \Gamma' \rightarrow 0,
$$
which yields an exact sequence
$$
\Hom_{\abx}(G, \Gamma') \rightarrow \Ext_{\abx}(G, \Gamma[n]) \rightarrow \Ext_{\abx}(G, \Gamma) \xrightarrow[]{n} \Ext_{\abx}(G, \Gamma'),
$$
whose first term vanishes since $G$ has geometrically connected fibers over $S$, and whose last homomorphism vanishes as well, since its composition with the injective homomorphism
$$
\Ext_{\abx}(G, \Gamma') \rightarrow \Ext_{\abx}(G, \Gamma),
$$
cf. \ref{chap3inj}, is the multiplication by $n$ on $\Ext_{\abx}(G, \Gamma)$, which is zero since $n$ annihilates $G$. 

\subsection{\label{chap34.8}} Let $G$ be a finitely presented affine commutative $S$-group scheme, and let $H$ be an object of $\abx$ such that any $H$-torsor over a finitely presented affine $S$-scheme is trivial. Let us consider an extension 
\begin{align}\label{chap3ext10}
0 \rightarrow H \xrightarrow[]{\iota} E \xrightarrow[]{r} G \rightarrow 0,
\end{align}
of $G$ by $H$ in $\abx$. The action of $H$ by left multiplication on $E$ turns the latter into a left $H$-torsor over $G$. Since $G$ is affine and finitely presented over $S$, this torsor is trivial, and thus we can assume (and we do) that $E$ is $H \times_S G$ as an object of $S_{\fppf}$, that $r(h,g) = g$ for any local sections $h$ and $g$ of $H$ and $G$, and that the left action of $H$ on $E$ is given by $(h',0) + (h,g) = (h'+h,g)$. The addition on $E$ must then take the form
\begin{align}\label{chap3ext11}
(h_1,g_1) + (h_2,g_2) = (h_1 + h_2 + c(g_1,g_2), g_1 + g_2),
\end{align}
for some morphism of sheaves $c$ from $G \times_S G$ to $H$. Since $E$ is a commutative group under the law (\ref{chap3ext11}), we have the relations
\begin{align}\label{chap3ext12}
&c(g_1,g_2) = c(g_2,g_1), \\
\label{chap3ext14} &c(g_1,g_2+g_3) + c(g_2,g_3) = c(g_1,g_2) + c(g_1+g_2,g_3),
\end{align}
corresponding respectively to the commutativity and the associativity of the law (\ref{chap3ext11}).
 
 Conversely, if $c$ is a morphism from $G \times_S G$ to $H$ which satisfies the relations (\ref{chap3ext12}) and (\ref{chap3ext14}), then the formula (\ref{chap3ext11}) defines an extension $E_c$ of $G$ by $H$ in $\abx$, whose underlying $S$-scheme is $H \times_S G$.

\begin{prop}\label{chap3ext5} Let $G$ be a finitely presented affine commutative $S$-group scheme, and let $H$ be an object of $\abx$ such that any $H$-torsor over a finitely presented affine $S$-scheme is trivial. Let $\mathcal{C}(G,H)$ be the group of morphisms from $G \times_S G$ to $H$ in $S_{\fppf}$ which satisfy the relations (\ref{chap3ext12}) and (\ref{chap3ext14}). We then have an exact sequence
$$
0 \rightarrow \Hom_{\abx}(G,H) \rightarrow \Hom_{S_{\fppf}}(G,H) \xrightarrow[]{d} \mathcal{C}(G,H) \xrightarrow[]{c \mapsto E_c} \Ext_{\abx}(G,H) \rightarrow 0,
$$
where $d : \Hom_{S_{\fppf}}(G,H) \rightarrow \mathcal{C}(G,H)$ is the homomorphism given on sections by
$$
d(f)(g_1,g_2) = f(g_1 + g_2) - f(g_1) - f(g_2).
$$
\end{prop}
 
Indeed, we already know that the map $c \mapsto E_c$ from $\mathcal{C}(G,H)$ to $\Ext_{\abx}(G,H)$ is surjective, and its kernel consists of the elements $c$ of $\mathcal{C}(G,H)$ for which $E_c$ is a trivial extension of $G$ by $H$. If $c$ is such a morphism, then the surjective homomorphism $E_c \rightarrow G$ has a section, which must take the form $g \mapsto (f(g),g)$ for some morphism $f$ from $G$ to $H$. We have
$$
(f(g_1 + g_2),g_1 + g_2) = (f(g_1),g_1) + (f(g_2),g_2) = (f(g_1) +f(g_2) + c(g_1,g_2),g_1 + g_2),
$$
for any local sections $g_1,g_2$ of $G$, hence $d(f) = c$. Conversely, any element $f$ of $\Hom_{S_{\fppf}}(G,H)$ such that $d(f) = c$ provides a section $g \mapsto (f(g),g)$ of the extension $E_c$. Thus the sequence in Proposition \ref{chap3ext5} is exact at $\mathcal{C}(G,H)$. The result then follows from the description of homomorphisms from $G$ to $H$ as elements $f$ of $\Hom_{S_{\fppf}}(G,H)$ such that $d(f)$ vanishes.

\begin{exemple}\label{chap3witt} Let $c$ be the image in $\mathbb{F}_p[U_1,U_2]$ of the polynomial
$$
\frac{ U_1^p + U_2^p - (U_1+U_2)^p}{p} = -\sum_{i=1}^{p-1} \frac{(p-1)!}{i! (p-i)!} U_1^i U_2^{p-i} \in \mathbb{Z}[U_1,U_2].
$$
Then the morphism from $\Ga_{a,\mathbb{F}_p} \times_{\mathbb{F}_p}\Ga_{a,\mathbb{F}_p}$ to $\Ga_{a,\mathbb{F}_p}$ corresponding to $c$ belongs to $\mathcal{C}(\Ga_{a,\mathbb{F}_p},\Ga_{a,\mathbb{F}_p})$. The corresponding extension $E_c$ of $\Ga_{a,\mathbb{F}_p}$ by itself is isomorphic to the $\mathbb{F}_p$-group scheme of Witt vectors of length $2$.
\end{exemple} 

\subsection{\label{chap34.4}} The group $G = H = \mathbb{G}_{a,S}$ satisfy the assumptions of Proposition \ref{chap3ext5}, and thus satisfy its conclusion. We therefore have a homomorphism
$$
 \Hom_{S_{\fppf}}(\mathbb{G}_{a,S},\mathbb{G}_{a,S}) \xrightarrow[]{d} \mathcal{C}(\mathbb{G}_{a,S},\mathbb{G}_{a,S}),
$$
with kernel $\Hom_{\abx}(\mathbb{G}_{a,S},\mathbb{G}_{a,S})$ and with cokernel $\Ext_{\abx}(\mathbb{G}_{a,S},\mathbb{G}_{a,S})$. Moreover, the group $\Hom_{S_{\fppf}}(\mathbb{G}_{a,S},\mathbb{G}_{a,S})$ can be identified with the group $A[U]$ of polynomials in one variable over $A$, while $\mathcal{C}(\mathbb{G}_{a,S},\mathbb{G}_{a,S}) $ can be identified with the group $A[U_1,U_2]$ of polynomials $c$ in two variables over $A$, which satisfy the relations
 \begin{align*}
&c(U_1,U_2) = c(U_2,U_1), \\
&c(U_1,U_2+U_3) + c(U_2,U_3) = c(U_1,U_2) + c(U_1+U_2,U_3).
\end{align*}

\begin{prop}[\cite{JPS2}, V.5]\label{chap3newton} A polynomial $f$ in $A[U]$ belongs to $\Hom_{\abx}(\mathbb{G}_{a,S},\mathbb{G}_{a,S})$ if and only if it is of the form
$$
f(U) = \sum_{r \geq 0} a_r U^{p^r},
$$
for some elements $(a_r)_{r \geq 0}$ of $A$.
\end{prop}

Indeed, a polynomial $f(U) = \sum_{m \geq 0} b_m U^m$ of $A[U]$ belongs to $\Hom_{\abx}(\mathbb{G}_{a,S},\mathbb{G}_{a,S})$ if and only if the element
$$
d(f) = \sum_{m \geq 0} b_m d(U^m),
$$
of $A[U_1,U_2]$ vanishes. Since each $b_m d(U^m)$ is a homogeneous polynomial of degree $m$, it is the homogeneous part of degree $m$ of $d(f)$. Consequently, $d(f)$ vanishes if and only if so does $b_m d(U^m)$ for each $m$. The conclusion of Proposition \ref{chap3newton} then follows from the fact that for each integer $m \geq 1$, we have $d(U^m) = 0$ if and only if $m$ is a power of $p$. Indeed, if $m = p^v n$ with $n$ prime to $p$, then $d(U^m) = d(U^n)^{p^v}$ vanishes if and only if $d(U^n)$ does. We have $d(U) = 0$, and if $n > 1$, then the coefficient of $U_1 U_2^{n-1}$ in $d(U^n)$ is $n$, whence $d(U^n)$ is non zero since $n$ is prime to $p$.
 
\begin{prop}\label{chap3ext6} Let $E$ be an extension of $\Ga_{a,S}$ by itself in $\abx$, whose pushout by the homomorphism $t \mapsto t - t^p$ from $\Ga_{a,S}$ to itself is a trivial extension of $\Ga_{a,S}$ by itself. Then $E$ is a trivial extension of $\Ga_{a,S}$ by itself.
\end{prop}

Indeed, let $c$ be an element of $\mathcal{C}(\mathbb{G}_{a,S},\mathbb{G}_{a,S}) $ such that $E_{c - c^p}$ a trivial extension of $\Ga_{a,S}$ by itself. By Proposition \ref{chap3ext5}, there exists an element $f = \sum_{n \geq 0} b_n U^n$ of $A[U]$ such that we have $d(f) = c - c^p$. We thus have
$$
c(U_1,U_2) = c(U_1,U_2)^p + \sum_{n \geq 0} b_n d(U^n).
$$
In particular, we have $c(0,0) = c(0,0)^p - b_0$ since $d(1) = -1$, and thus
$$
c(U_1,U_2) - c(0,0) = (c(U_1,U_2) - c(0,0) )^p + \sum_{n \geq 1} b_n d(U^n).
$$
By iterating this identity, we obtain a relation
$$
c(U_1,U_2) - c(0,0) = \sum_{v \geq 0} \left( \sum_{n \geq 1} b_n d(U^n) \right)^{p^v} = \sum_{m \geq 1} c_m d(U^m),
$$
in the power series ring $A[[U_1,U_2]]$, where $c_m = \sum_{m = p^vn} b_n^{p^v}$. For each integer $m$, the polynomial $d(U^m)$ is homogeneous of degree $m$. In particular, the polynomial $c_m d(U^m)$ is the homogeneous part of degree $m$ of $c(U_1,U_2) - c(0,0)$. Since $c(U_1,U_2) -c(0,0)$ is a polynomial, we must have $c_m d(U^m)= 0$ for $m$ large enough. Thus the power series
$$
-c(0,0) + \sum_{\substack{m \geq 1 \\ d(U^m) \neq 0}} c_m U^m,
$$
is a polynomial, whose image by $d$ is $c$. Consequently, $E_c$ is a trivial extension of $\Ga_{a,S}$ by itself.

\subsection{\label{chap34.2}} Let $F : \Ga_{a,S} \rightarrow \Ga_{a,S}$ be the Frobenius homomorphism, given on sections by $t \mapsto t^p$. We then have the so-called Artin-Schreier exact sequence
\begin{align}\label{chap3as}
0 \rightarrow \mathbb{F}_{p,S} \rightarrow \Ga_{a,S} \xrightarrow[]{1-F} \Ga_{a,S} \rightarrow 0.
\end{align}
By applying the functor $\Hom(\Ga_{a,S},-)$ to this short exact sequence, we obtain a long exact sequence

\begin{center}
\begin{tikzpicture}[descr/.style={fill=white,inner sep=1.5pt}]
 \matrix (m) [
 matrix of math nodes,
 row sep=3em,
 column sep=2.5em,
 text height=1.5ex, text depth=0.25ex
 ]
 { \Hom_{\abx}(\Ga_{a,S},\mathbb{F}_{p}) & \Hom_{\abx}(\Ga_{a,S},\Ga_{a,S}) & \Hom_{\abx}(\Ga_{a,S},\Ga_{a,S}) \\
 \Ext_{\abx}(\Ga_{a,S},\mathbb{F}_{p}) & \Ext_{\abx}(\Ga_{a,S},\Ga_{a,S}) & \Ext_{\abx}(\Ga_{a,S},\Ga_{a,S}) \\
 };

 \path[overlay,->, font=\scriptsize,>=latex]
 (m-1-1) edge (m-1-2)
 (m-1-2) edge node[above]{$1-F$} (m-1-3)
 (m-1-3) edge[out=355,in=175] node[descr,yshift=0.3ex] {$\delta$} (m-2-1)
 (m-2-1) edge (m-2-2)
 (m-2-2) edge node[above]{$1-F$} (m-2-3);
\end{tikzpicture}

\end{center}
whose first term vanishes since $\Ga_{a,S}$ has geometrically connected fibers over $S$, and whose last homomorphism is injective by Proposition \ref{chap3ext6}. We thus have a short exact sequence
 \begin{align}\label{chap3exactseq10}
0 \rightarrow \Hom_{\abx}(\Ga_{a,S},\Ga_{a,S}) \xrightarrow[]{1-F} \Hom_{\abx}(\Ga_{a,S},\Ga_{a,S}) \xrightarrow[]{\delta} \Ext_{\abx}(\Ga_{a,S},\mathbb{F}_{p}) \rightarrow 0.
 \end{align}

\begin{prop}\label{chap3desc} If we denote, for each element $a$ of the ring $A = \Gamma(S,\Ow_S)$, by $m_a$ the endomorphism of $\Ga_{a,S}$ which sends a section $t$ to $at$, then the homomorphism
\begin{align*}
A \oplus \Hom_{\abx}(\Ga_{a,S},\Ga_{a,S}) &\rightarrow \Hom_{\abx}(\Ga_{a,S},\Ga_{a,S}) \\
(a,u) &\mapsto m_a + F(u)-u,
\end{align*}
is an isomorphism.
 \end{prop}
 
Indeed, the group $\Hom_{\abx}(\Ga_{a,S},\Ga_{a,S})$ is the group of polynomials $v$ in $A[T]$ which are additive, in the sense that
$$
v(T+S) = v(T) + v(S).
$$
These are exactly the polynomials of the form $v(T) = \sum_{j \geq 0} a_j T^{p^j}$, where $(a_j)_{j \geq 0}$ is a finite family of elements of $A$, cf. Proposition \ref{chap3newton}. By writing successively monomials of the form $a T^{p^j}$ with $j \geq 1$ as $a^{p^{-1}} T^{p^{j-1}} + (F-1)(u)$, where $u = a^{p^{-1}} T^{p^{j-1}}$, we obtain that any such polynomial $v$ can be decomposed as
$$
v(T) = a T + (F-1)(u),
$$
for some additive polynomial $u$, and with $a = \sum_{j \geq 0} a_j^{p^{-j}}$. The homomorphism in Proposition \ref{chap3desc} is thus surjective.

 On the other hand, if an additive polynomial $u = \sum_{j \geq 0} u_j T^{p^j}$ is such that $(F-1)(u)$ is of the form $aT$ for some element $a$ of $A$, then we have
 $$
 aT = \sum_{j \geq 0} u_j^p T^{p^{j+1}} - \sum_{j \geq 0} u_j T^{p^j} = -u_0 T + \sum_{j \geq 1} (u_{j-1}^p - u_j) T^{p^j},
 $$
 so that $a = -u_0$ and $u_{j-1}^p = u_j$ for each $j \geq 1$. This implies that for each $j \geq 0$, there exists an integer $N$ such that $u_j^{p^N} = 0$, and thus $u_j = 0$ since $A$ is reduced. Thus such an additive polynomial $u$ must be zero, which proves that the homomorphism in Proposition \ref{chap3desc} is injective.
 
 \begin{cor}\label{chap3as2} The homomorphism of abelian groups
\begin{align*}
 A &\rightarrow \Ext_{\abx}(\Ga_{a,S},\mathbb{F}_{p}) \\
 a &\mapsto \delta(m_a)
\end{align*}
is an isomorphism.
 \end{cor}

This follows immediately from Proposition \ref{chap3desc} and from the exact sequence (\ref{chap3exactseq10}). For each element $a$ of $A$, the extension $\delta(m_a)$ can be explicitly described as the pullback of the extension (\ref{chap3as}) by the homomorphism $m_a$.

\subsection{\label{chap34.3}} Let $\Lambda$ be an $\ell$-adic coefficient ring, and let $\psi : \mathbb{F}_{p} \rightarrow \Lambda^{\times}$ be a non trivial (hence injective) homomorphism. The pushout by $\psi$ of the Artin-Schreier $\mathbb{F}_p$-torsor (cf. \ref{chap3as}) yields a multiplicative $\Lambda$-local system on $\Ga_{a,S}$ (cf. \ref{chap32.0.3}), which we denote by $\Lc_{\psi}$. More generally, if $f : X \rightarrow \Ga_{a,S}$ is a morphism of $S$-schemes, we denote by $\Lc_{\psi} \lbrace f \rbrace$ the $\Lambda$-local system of rank $1$ on $X$ given by the pullback of $\Lc_{\psi}$ by $f$. If moreover $X$ is an $S$-group scheme and if $f$ is a homomorphism of $S$-group schemes, then $\Lc_{\psi} \lbrace f \rbrace$ is a multiplicative $\Lambda$-local system on $X$. 

\begin{prop}\label{chap3charsheafas} Let $V$ be a finitely generated projective $A$-module, and let $\mathcal{V}$ be the corresponding $S$-group scheme. Then any multiplicative $\Lambda$-local system on $\mathcal{V}$ is isomorphic to $\Lc_{\psi}\lbrace v^* \rbrace$, for a unique $A$-linear homomorphism $v^* : V \rightarrow A$, considered as a homomorphism from $\mathcal{V}$ to $\Ga_{a,S}$.
\end{prop}

By \ref{chap32.0.4}, we can assume (and we do) that $V$ is the $A$-module $A$, and thus that $\mathcal{V}$ is $\Ga_{a,S}$. Moreover, we can assume (and we do as well) that $\Lambda$ is finite, in which case multiplicative $\Lambda$-local systems on $\Ga_{a,S}$ correspond to extensions in $S_{\fppf}$ of $\Ga_{a,S}$ by the finite abelian group $\Lambda^{\times}$. Since $\psi$ realizes an isomorphism from $\mathbb{F}_{p}$ to the $p$-torsion subgroup of $\Lambda^{\times}$, we deduce from Proposition \ref{chap3ppart} that the homomorphism
$$
\Ext_{\abx}(\Ga_{a,S}, \mathbb{F}_{p}) \xrightarrow[]{\psi} \Ext_{\abx}(\Ga_{a,S}, \Lambda^{\times}),
$$
induced by $\psi$ is an isomorphism. The conclusion then follows from Corollary \ref{chap3as2}.

\begin{cor}\label{chap3vanish13} Assume that $A = k$ is an algebraically closed field. Let $V$ be a $k$-vector space of finite dimension $r \geq 1$, let $\mathcal{V}$ be the corresponding $k$-group scheme and let $\mathcal{M}$ be a multiplicative $\Lambda$-local system on $\mathcal{V}$. Then the cohomology group
$$
H_c^{\nu}(\mathcal{V}, \mathcal{M})
$$
vanishes for each integer $\nu$, unless $\nu = 2r$ and $\mathcal{M}$ is trivial, in which case it is a free $\Lambda$-module of rank $1$.
\end{cor}

By Proposition \ref{chap3charsheafas}, we can assume that $V$ is $k^r$ and that $\mathcal{M}$ is $\Lc_{\psi}\lbrace x_1 \rbrace$, where $x_1 : k^r \rightarrow k$ is the first coordinate. If $\overline{\mathbb{F}}_p$ is the algebraic closure of $\mathbb{F}_p$ in $k$, then for each integer $\nu$ the group 
$$
H_c^{\nu}(\mathcal{V}, \mathcal{M}) = H_c^{\nu}(\Ga_{a,k}^r, \Lc_{\psi}\lbrace x_1 \rbrace)
$$
is isomorphic to $H_c^{\nu}(\Ga_{a,\overline{\mathbb{F}}_p}^r, \Lc_{\psi}\lbrace x_1 \rbrace)$, which vanishes by (\cite{SGA412}, [Sommes trig.] Th. 2.7).

\begin{cor}\label{chap3vanish14} Assume that $A = k$ is an algebraically closed field of characteristic $p \neq 2$. Let $V$ be a $k$-vector space of finite dimension $r \geq 1$, let $\mathcal{V}$ be the corresponding $k$-group scheme, let $\gamma : V \rightarrow k$ be a non zero linear form and let $\mathcal{M}$ be a multiplicative $\Lambda$-local system on $\mathcal{V}$. Then the cohomology group
$$
H_c^{\nu}(\mathcal{V}, \mathcal{M} \otimes \Lc_{\psi}\lbrace \gamma^2 \rbrace)
$$
vanishes for each integer $\nu$, unless $\nu = 2r-1$ and $\mathcal{M}$ is isomorphic to $ \Lc_{\psi}\lbrace \alpha \gamma \rbrace$ for some element $\alpha$ of $k$, in which case it is a free $\Lambda$-module of rank $1$.
\end{cor}

Indeed, if $t$ is the coordinate on $\Ga_{a,S}$, then the projection formula yields an isomorphism
$$
R\gamma_! (\mathcal{M} \otimes \Lc_{\psi}\lbrace \gamma^2 \rbrace) \cong R\gamma_! (\mathcal{M} )\otimes \Lc_{\psi}\lbrace t^2 \rbrace.
$$
By Corollary \ref{chap3vanish13} and by multiplicativity of $\mathcal{M}$, the geometric fibers of this complex vanish unless the restriction of $\mathcal{M}$ to the kernel of $\gamma$ is trivial. By Proposition \ref{chap3charsheafas}, the multiplicative $\Lambda$-local system $\mathcal{M}$ is isomorphic to $ \Lc_{\psi}\lbrace v^* \rbrace$ for a unique linear form $v^*$ on $V$. We thus obtain that $R\gamma_! (\mathcal{M} \otimes \Lc_{\psi}\lbrace \gamma^2 \rbrace)$ vanishes unless the restriction of $v^*$ to the kernel of $\gamma$ vanishes, in which case we have $v^* = \alpha \gamma$ for some element $\alpha$ of $k$.

When $\mathcal{M}$ is isomorphic to $ \Lc_{\psi}\lbrace \alpha \gamma \rbrace$ for some element $\alpha$ of $k$, we have 
\begin{align*}
R\gamma_! (\mathcal{M} \otimes \Lc_{\psi}\lbrace \gamma^2 \rbrace) &\cong R\gamma_! (\Lambda) \otimes \Lc_{\psi}\lbrace \alpha t + t^2 \rbrace \\
&\cong \Lc_{\psi}\lbrace \alpha t + t^2 \rbrace (-r)[-2r],
\end{align*}
and the conclusion then follows from the fact that the cohomology group
$$
H_c^{\nu}(\Ga_{a,k}, \Lc_{\psi}\lbrace \alpha t + t^2 \rbrace)
$$
is of rank $1$ for $\nu = 1$, and vanishes otherwise. Indeed, this group vanishes for $\nu = 0$ since $\Lc_{\psi}\lbrace \alpha t + t^2 \rbrace$ has no punctual sections, it vanishes for $\nu =2$ by Poincar\'e duality, and its Euler characteristic is $-1$ by the Grothendieck-Ogg-Shafarevich formula, since the Swan conductor of $\Lc_{\psi}\lbrace \alpha t + t^2 \rbrace$ at $\infty$ is $2$.

\subsection{\label{chap34.5}} We now assume that $S$ is of characteristic $p=2$. The element $c = U_1 U_2$ of $A[U_1,U_2]$ belongs to $\mathcal{C}(\Ga_{a,S},\Ga_{a,S} )$ (cf. \ref{chap3ext5}), and thus defines an extension $G=E_c$ of $\Ga_{a,S}$ by itself, cf. \ref{chap34.4}. Thus $G$ is $\Ga_{a,S} \times_S \Ga_{a,S}$ as an $S$-scheme, endowed with the multiplication
$$
(t_1,u_1) + (t_2,u_2) = (t_1 + t_2 + u_1u_2, u_1 + u_2),
$$
for sections $t_1,u_1,t_2,u_2$ of $\Ga_{a,S}$. Equivalently, the $S$-group scheme $G$ is the pullback to $S$ of the group of Witt vectors of length $2$ over $\mathbb{F}_2$, cf. \ref{chap3witt}.

The group $G$ satisfies the assumptions of Proposition \ref{chap3ext5}, and thus satisfies its conclusion. We therefore have a homomorphism
$$
 \Hom_{S_{\fppf}}(G,G) \xrightarrow[]{d} \mathcal{C}(G,G),
$$
with kernel $\Hom_{\abx}(G,G)$ and with cokernel $\Ext_{\abx}(G,G)$. The group $\Hom_{S_{\fppf}}(G,G)$ can be identified with the group of couples $f = (f_0,f_1)$ of elements of $A[T,U]$.

\begin{prop}\label{chap3newton2}For $p=2$, let $G$ be the extension of $\Ga_{a,S}$ by itself defined by $c(U_1,U_2) = U_1 U_2$ (cf. \ref{chap34.4}). A couple $f = (f_0,f_1)$ of elements of $A[T,U]$ belongs to $\Hom_{\abx}(G,G)$ if and only if it is of the form
$$
\left( a \left( T^{\frac{1}{2}} \right)^2 + \widetilde{a}(U) + b(U), a(U) \right),
$$
where $a(U) = \sum_{r \geq 0} a_r U^{2^r}$ and $b(U)$ are additive polynomials, and where
\begin{align*}
 a \left( T^{\frac{1}{2}} \right)^2 &= \sum_{r \geq 0} a_r^2 T^{2^r},\\
\widetilde{a}(U) &= \sum_{r_1 > r_2 \geq 0} a_{r_1} a_{r_2} U^{2^{r_1}+2^{r_2}}.
\end{align*} 
\end{prop}

Indeed, such a couple $f = (f_0,f_1)$ belongs to $\Hom_{\abx}(G,G)$ if and only if $d(f)$ vanishes, namely if and only if the relations
\begin{align}
\label{chap3mor1}f_0(T_1+T_2 + U_1U_2, U_1 + U_2) &= f_0(T_1,U_1) + f_0(T_2,U_2) + f_1(T_1,U_1) f_1(T_2,U_2) \\
\label{chap3mor2}f_1(T_1+T_2 + U_1U_2, U_1 + U_2) &= f_1(T_1,U_1) + f_1(T_2,U_2),
\end{align}
hold in $A[T_1,U_1,T_2,U_2]$. Setting $U_1 = T_2 = 0$ in (\ref{chap3mor2}), we obtain
$$
f_1(T_1,U_2) = f_1(T_1,0) + f_1(0,U_2).
$$
Setting $U_1 = U_2 = 0$ in (\ref{chap3mor2}), we obtain that $f_1(T,0)$ is an additive polynomial, so that we can write
$$
f_1(T,0) = \sum_{r \geq 0} x_r T^{2^r},
$$ 
for some elements $(x_r)_r$ of $A$ (cf. Proposition \ref{chap3newton}), while setting $T_1 = T_2 = 0$ in (\ref{chap3mor2}) yields
$$
f_1(0,U_1+U_2) + f_1(U_1U_2,0) = f_1(0,U_1) + f_1(0,U_2).
$$
Writing $f_1(0,U)$ as $\sum_{n \geq 0} y_n U^n$ for some elements $(y_n)_n$ of $A$, we obtain
$$
\sum_{n \geq 0} y_n \left( (U_1 + U_2)^{n} - U_1^n - U_2 ^n \right) + \sum_{r \geq 0} x_r U_1^{2^r} U_2^{2^r} = 0.
$$
For each integer $r \geq 0$, the homogeneous part of degree $2^{r+1}$ in this relation yields the vanishing of $x_r U_1^{2^r} U_2^{2^r}$, hence $x_r =0$. We thus have
$$
f_1(T,U) = a(U),
$$
and $a(U) = f_1(0,U)$ is an additive polynomial.

Setting $U_1 = T_2 = 0$ in (\ref{chap3mor1}), we obtain
$$
f_0(T_1,U_2) = f_0(T_1,0) + f_0(0,U_2) + f_1(T_1,0) f_1(0,U_2) = f_0(T_1,0) + f_0(0,U_2),
$$
since $f_1(T_1,0) = a(0)$ vanishes. Setting $U_1 = U_2 = 0$ in (\ref{chap3mor1}), we obtain that $f_0(T,0)$ is an additive polynomial, so that we can write
$$
f_0(T,0) = \sum_{r \geq 0} c_r T^{2^r},
$$ 
for some elements $(c_r)_r$ of $A$ (cf. Proposition \ref{chap3newton}), while setting $T_1 = T_2 = 0$ in (\ref{chap3mor1}) yields
\begin{align}
\label{chap3mor3} f_0(0,U_1+U_2) + f_0(U_1U_2,0) = f_0(0,U_1) + f_0(0,U_2) + a(U_1)a(U_2).
\end{align}
Let us write the additive polynomial $a(U)$ as $\sum_{r \geq 0} a_r U^{2^r}$, cf. Proposition \ref{chap3newton}, and let us write the polynomial $f_0(0,U)$ as $\sum_{n \geq 0} b_n U^n$ for some elements $(b_n)_n$ of $A$, so that (\ref{chap3mor3}) can be written as
$$
\sum_{n \geq 0} b_n \left( (U_1 + U_2)^{n} - U_1^n - U_2 ^n \right) = \sum_{r_1,r_2 \geq 0} a_{r_1} a_{r_2} U_1^{2^{r_1}} U_2^{2^{r_2}} - \sum_{r \geq 0} c_r U_1^{2^r} U_2^{2^r}.
$$
For each integer $r \geq 0$, the homogeneous part of degree $2^{r+1}$ in this relation yields the vanishing of $(a_r^2 - c_r) U_1^{2^r} U_2^{2^r}$, hence $c_r =a_r^2$. Furthermore, for each integer $n$ which is not a power of $2$, the homogeneous part of degree $n$ in this relation yields $b_n = a_{r_1} a_{r_2}$ if $n = 2^{r_1} + 2^{r_2}$ for some pair of distinct integers $(r_1,r_2)$, and $b_n = 0$ otherwise. We thus have
$$
f_0(0,U) = \widetilde{a}(U) + b(U),
$$
where $b(U) = \sum_{r \geq 0} b_{2^r} U^{2^r}$ is an additive polynomial and where $\widetilde{a}(U) = \sum_{r_1 > r_2 \geq 0} a_{r_1} a_{r_2} U^{2^{r_1}+2^{r_2}} $. Moreover, the relation $c_r = a_r^2$ yields $f_0(T,0) = a ( T^{\frac{1}{2}})^2$, so that
$$
f_0(T,U) = f_0(T,0) + f_0(0,U) = a \left( T^{\frac{1}{2}} \right)^2 + \widetilde{a}(U) + b(U),
$$
hence the conclusion of Proposition \ref{chap3newton2}.

\begin{prop}\label{chap3ext17}For $p=2$, let $G$ be the extension of $\Ga_{a,S}$ by itself defined by $c(U_1,U_2) = U_1 U_2$ (cf. \ref{chap34.4}), and let $F$ be the endomorphism $(t,u) \mapsto (t^2,u^2)$ of $G$. Let $E$ be an extension of $G$ by itself in $\abx$, whose pushout by the endomorphism $1-F$ of $G$ is a trivial extension of $G$ by itself. Then $E$ is a trivial extension of $G$ by itself.
\end{prop}

We prove Proposition \ref{chap3ext17} by an argument similar to the one we used to prove Proposition \ref{chap3ext6}. We endow $A[T,U]$ (resp. $A[T_1,U_1,T_2,U_2]$) with a structure of $\mathbb{N}$-graded $A$-algebra by assigning weight $2$ to the variable $T$ (resp. $T_1,T_2$), and weight $1$ to the variable $U$ (resp. $U_1,U_2$). If $B$ is an $\mathbb{N}$-graded $A$-algebra, an element $f = (b_0,b_1)$ of $G(B)$ is said to be \textit{homogeneous of degree $n$} if $b_0$ and $b_1$ are homogeneous elements of degrees $n$ and $\frac{n}{2}$ respectively in $B$. In particular, we have $b_1 = 0$ if $n$ is odd. One should note that for each integer $n$, the subset of $G(B)$ consisting of homogeneous elements of degree $n$ is a subgroup of $G(B)$. Any element $f$ of $G(B)$ can be uniquely written as a finite sum
$$
f = \sum_{n \geq 0} f_n,
$$
where $f_n$ is a homogeneous element of degree $n$ in $G(B)$. The element $f_n$ of $G(B)$ will be referred to as the \textit{homogeneous part of degree $n$} of $f$. 

Let $\gamma$ be an element of $\mathcal{C}(G,G) $ such that $E_{\gamma - F(\gamma)}$ is a trivial extension of $G$ by itself. By Proposition \ref{chap3ext5}, there exists an element $f$ of $G(A[T,U])$ such that $d(f) = \gamma - F(\gamma)$. Let us write $f = \sum_{n \geq 0} f_n$ as the sum of its homogeneous parts, as above. We have
$$
\gamma(T_1,U_1,T_2,U_2) = F(\gamma(T_1,U_1,T_2,U_2)) + \sum_{n \geq 0} d(f_n),
$$
in $A[T_1,U_1,T_2,U_2]$. In particular, we have $\gamma(0) = F(\gamma(0)) + d(f_0)$, and thus
$$
\gamma(T_1,U_1,T_2,U_2) - \gamma(0) = F(\gamma(T_1,U_1,T_2,U_2) - \gamma(0) ) + \sum_{n \geq 1} d(f_n).
$$
By iterating this identity, we obtain a relation
$$
\gamma(T_1,U_1,T_2,U_2) - \gamma(0) = \sum_{v \geq 0} F^v\left( \sum_{n \geq 1} d(f_n) \right) = \sum_{m \geq 1} d(g_m),
$$
in the group $G(A[[T_1,U_1,T_2,U_2]])$, where $g_m = \sum_{m = p^vn} F^v(f_n)$ is homogeneous of degree $m$. For each integer $m$, the element $d(g_m)$ of $G(A[T_1,U_1,T_2,U_2])$ is homogeneous of degree $m$. In particular, the homogeneous part of degree $m$ of $\gamma(T_1,U_1,T_2,U_2) - \gamma(0) $ is $d(g_m)$. Since the element $\gamma(T_1,U_1,T_2,U_2) - \gamma(0) $ of $G(A[T_1,U_1,T_2,U_2])$ has a non zero homogeneous part of degree $m$ for only finitely many integers $m$, we must have $ d(g_m)= 0$ for $m$ large enough. Thus the element
$$
-\gamma(0) + \sum_{\substack{m \geq 1 \\ d(g_m) \neq 0}} g_m,
$$
of $G(A[[T_1,U_1,T_2,U_2]])$ belongs to $G(A[T_1,U_1,T_2,U_2])$, and its image by $d$ is $\gamma$. Consequently, $E_{\gamma}$ is a trivial extension of $G$ by itself.

\subsection{\label{chap34.6}}For $p=2$, let $G$ be the extension of $\Ga_{a,S}$ by itself defined by $c(U_1,U_2) = U_1 U_2$ (cf. \ref{chap34.4}), and let $F$ be the endomorphism $(t,u) \mapsto (t^2,u^2)$ of $G$. We then have the Lang-Artin-Schreier exact sequence
\begin{align}\label{chap3as3}
0 \rightarrow G(\mathbb{F}_{2}) \rightarrow G \xrightarrow[]{1-F} G \rightarrow 0.
\end{align}
By applying the functor $\Hom(G,-)$ to this short exact sequence, we obtain a long exact sequence

\begin{center}
\begin{tikzpicture}[descr/.style={fill=white,inner sep=1.5pt}]
 \matrix (m) [
 matrix of math nodes,
 row sep=3em,
 column sep=2.5em,
 text height=1.5ex, text depth=0.25ex
 ]
 { \Hom_{\abx}(G,G(\mathbb{F}_{2}) ) & \Hom_{\abx}(G,G) & \Hom_{\abx}(G,G) \\
 \Ext_{\abx}(G,G(\mathbb{F}_{2}) ) & \Ext_{\abx}(G,G) & \Ext_{\abx}(G,G) \\
 };

 \path[overlay,->, font=\scriptsize,>=latex]
 (m-1-1) edge (m-1-2)
 (m-1-2) edge node[above]{$1-F$} (m-1-3)
 (m-1-3) edge[out=355,in=175] node[descr,yshift=0.3ex] {$\delta$} (m-2-1)
 (m-2-1) edge (m-2-2)
 (m-2-2) edge node[above]{$1-F$} (m-2-3);
\end{tikzpicture}

\end{center}
whose first term vanishes since $G$ has geometrically connected fibers over $S$, and whose last homomorphism is injective by Proposition \ref{chap3ext17}. We thus have a short exact sequence
 \begin{align}\label{chap3exactseq12}
0 \rightarrow \Hom_{\abx}(G,G) \xrightarrow[]{1-F} \Hom_{\abx}(G,G) \xrightarrow[]{\delta} \Ext_{\abx}(G,G(\mathbb{F}_{2})) \rightarrow 0.
 \end{align}
 
\begin{rema}\label{chap3trivia} There is a unique isomorphism of abelian groups
$$
\Z / 4 \Z \rightarrow G(\mathbb{F}_2)
$$
which sends $1$ to $(0,1)$.
\end{rema}

 For any additive polynomials $a(U) = \sum_{r \geq 0} a_r U^{2^r}$ and $b(U)$ with coefficients in $A$, let us denote by $\langle a,b \rangle$ the element of $\Hom_{\abx}(\Ga_{a,S},\Ga_{a,S})$ given by
$$
\langle b,a \rangle =\left( a \left( T^{\frac{1}{2}} \right)^2 + \widetilde{a}(U) + b(U), a(U) \right)
$$
where we have set $\widetilde{a}(U) = \sum_{r_1 > r_2 \geq 0} a_{r_1} a_{r_2} U^{2^{r_1}+2^{r_2}}$.
By Proposition \ref{chap3newton2}, any endomorphism of $G$ is of the form $\langle b,a \rangle$ for a (necessarily unique) couple $(a,b)$ of additive polynomials with coefficients in $A$.

\begin{prop}\label{chap3desc2} The map
\begin{align*}
A^2 \times \Hom_{\abx}(G,G) &\rightarrow \Hom_{\abx}(G,G) \\
((\beta,\alpha),f) &\mapsto \langle \beta U,\alpha U \rangle + f - F \circ f,
\end{align*}
is bijective.
 \end{prop}
 
Let $\langle b,a \rangle$ be an endomorphism of $G$, where $a(U) = \sum_{r \geq 0} a_r U^{2^r}$ and $b(U)$ are additive polynomials with coefficients in $A$. By Proposition \ref{chap3desc}, there exists a unique element $\alpha$ of $A$ and a unique additive polynomial $g(U)$ with coefficients in $A$ such that $a(U)$ is equal to $\alpha U + g(U) - g(U)^2$. We thus have
$$
\langle b,a \rangle = \langle b',\alpha U \rangle + (1-F)(\langle 0,g \rangle),
$$
for a uniquely determined additive polynomial $b'$. By Proposition \ref{chap3desc} again, there exists a unique element $\beta$ of $A$ and a unique additive polynomial $h(U)$ with coefficients in $A$ such that $b'(U)$ is equal to $\beta U + h(U) - h(U)^2$. We then have
\begin{align*}
\langle b,a \rangle &= \langle \beta U, \alpha U \rangle + (1-F)(\langle h,0 \rangle) + (1-F)(\langle 0,g\rangle) \\
&= \langle \beta U, \alpha U \rangle + (1-F)(\langle h,g \rangle),
\end{align*}
hence the conclusion of Proposition \ref{chap3desc2}.

 \begin{cor}\label{chap3as4}For $p=2$, let $G$ be the extension of $\Ga_{a,S}$ by itself defined by $c(U_1,U_2) = U_1 U_2$ (cf. \ref{chap34.4}). Then the map
\begin{align*}
 G(A) &\rightarrow \Ext_{\abx}(G,G(\mathbb{F}_{2}) ) \\
(\beta,\alpha) &\mapsto \delta(\langle \beta^{\frac{1}{2}} U, \alpha U \rangle )
\end{align*}
is an isomorphism of abelian groups.
 \end{cor}

The bijectivity of the map in Corollary \ref{chap3as4} follows immediately from Proposition \ref{chap3desc2} and from the exact sequence (\ref{chap3exactseq12}). The fact that this map is a group homomorphism follows from the following computation: if $(\beta_1,\alpha_1)$ and $(\beta_2,\alpha_2)$ are elements of $G(A)$, then we have
$$
\langle \beta_1^{\frac{1}{2}} U, \alpha_1 U \rangle + \langle \beta_2^{\frac{1}{2}} U, \alpha_2 U \rangle = \langle (\beta_1 + \beta_2)^{\frac{1}{2}} U + \alpha_1 \alpha_2 U^2, (\alpha_1 +\alpha_2) U \rangle,
$$
and the right hand side can be decomposed as 
$$
\langle (\beta_1 + \beta_2)^{\frac{1}{2}} U + \alpha_1 \alpha_2 U^2, (\alpha_1 +\alpha_2) U \rangle = \langle (\beta_1 + \beta_2 + \alpha_1 \alpha_2)^{\frac{1}{2}} U, (\alpha_1 +\alpha_2) U \rangle - (1-F)( \langle \alpha_1^{\frac{1}{2}} \alpha_2^{\frac{1}{2}} U,0 \rangle),
$$
hence the conclusion.

\subsection{\label{chap34.7}}For $p=2$, let $G$ be the extension of $\Ga_{a,S}$ by itself defined by $c(U_1,U_2) = U_1 U_2$ (cf. \ref{chap34.4}). Let $\Lambda$ be an $\ell$-adic coefficient ring, and let $\xi :G(\mathbb{F}_2) \rightarrow \Lambda^{\times}$ be an injective homomorphism of abelian groups; this amounts to a choice of primitive fourth root of unity in $\Lambda$, cf. Remark \ref{chap3trivia}.

The pushout by $\xi$ of the Lang-Artin-Schreier $G(\mathbb{F}_2)$-torsor (cf. \ref{chap3as}) yields a multiplicative $\Lambda$-local system on $G$ (cf. \ref{chap32.0.3}), which we denote by $\Lc_{\xi}$. More generally, if $f = (f_0,f_1) : X \rightarrow G$ is a morphism of $S$-schemes, we denote by $\Lc_{\xi} \lbrace f_0, f_1 \rbrace$ the $\Lambda$-local system of rank $1$ on $X$ given by the pullback of $\Lc_{\xi}$ by $f$. If moreover $X$ is an $S$-group scheme and if $f = (f_0,f_1)$ is a homomorphism of $S$-group schemes, then $\Lc_{\xi} \lbrace f_0,f_1 \rbrace$ is a multiplicative $\Lambda$-local system on $X$. 

If $f = (f_0,f_1)$ and $f' = (f'_0,f'_1)$ are $S$-morphisms from an $S$-scheme $X$ to $G$, then multiplicativity of $\Lc_{\xi}$ on $G$ yields an isomorphism
\begin{align} \label{chap3multipli57}
\Lc_{\xi} \lbrace f_0, f_1 \rbrace \otimes \Lc_{\xi} \lbrace f'_0, f'_1 \rbrace \cong \Lc_{\xi} \lbrace f_0 + f'_0 + f_1 f_1', f_1 +f_1' \rbrace
\end{align}
of $\Lambda$-local systems on $X$.

\begin{rema}\label{chap3as89} For any morphism $f : X \rightarrow \Ga_{a,S}$ of $S$-schemes, the $\Lambda$-local system $\Lc_{\xi} \lbrace f,0 \rbrace$ is isomorphic to the Artin-Schreier local system $\Lc_{\psi} \lbrace f \rbrace$ from \ref{chap34.3}, where $\psi$ is the restriction of $\xi$ to the subgroup $\mathbb{F}_2$ of $G(\mathbb{F}_2)$. Moreover, the composition of $\psi$ with the surjective homomorphism $G(\mathbb{F}_{2}) \rightarrow \mathbb{F}_2$ is equal to $\xi^2$, and we have isomorphisms by \ref{chap3multipli57}
$$
\Lc_{\xi} \lbrace f_0,f_1 \rbrace^{\otimes 2} \cong \Lc_{\xi} \lbrace f_1^2, 0 \rbrace \cong \Lc_{\psi} \lbrace f_1^2 \rbrace \cong \Lc_{\psi} \lbrace f_1 \rbrace,
$$
for any morphism $f = (f_0,f_1) : X \rightarrow G$ of $S$-schemes.
\end{rema}

\begin{prop}\label{chap3charsheafas2} Let $V$ be a finitely generated projective $A$-module, and let $\mathcal{V}$ be the corresponding $S$-group scheme. Let $\gamma : V \rightarrow A$ be a surjective $A$-linear homomorphism, and let $\widetilde{\mathcal{V}}$ be the extension of $\mathcal{V}$ by $\Ga_{a,S}$ defined by the element $c : (v_1,v_2) \mapsto \gamma(v_1) \gamma(v_2)$ of $\mathcal{C}(\mathcal{V},\Ga_{a,S})$, cf. \ref{chap3ext5}. We denote by $t : \widetilde{\mathcal{V}} \rightarrow \Ga_{a,S}$ and by $v : \widetilde{\mathcal{V}} \rightarrow \mathcal{V}$ the canonical projections.

Then any multiplicative $\Lambda$-local system on $\widetilde{\mathcal{V}}$ is isomorphic to $\Lc_{\xi}\lbrace \alpha^2 t + v^*(v), \alpha \gamma(v) \rbrace$, for a unique $A$-linear homomorphism $v^* : V \rightarrow A$, considered as a homomorphism from $\mathcal{V}$ to $\Ga_{a,S}$, and a unique element $\alpha$ of $A$.
\end{prop}

By \ref{chap32.0.4} and by Proposition \ref{chap3charsheafas}, we can assume (and we do) that $V$ is the $A$-module $A$, that $\gamma$ is the identity, and thus that $\widetilde{\mathcal{V}}$ is $G$. Moreover, we can assume (and we do as well) that $\Lambda$ is finite, in which case multiplicative $\Lambda$-local systems on $G$ correspond to extensions in $S_{\fppf}$ of $G$ by the finite abelian group $\Lambda^{\times}$. Since $\xi$ realizes an isomorphism from $G(\mathbb{F}_{2})$ to the $4$-torsion subgroup of $\Lambda^{\times}$, we deduce from Proposition \ref{chap3ppart} that the homomorphism
$$
\Ext_{\abx}(G, G(\mathbb{F}_{2})) \xrightarrow[]{\xi} \Ext_{\abx}(G, \Lambda^{\times}),
$$
induced by $\xi$ is an isomorphism. The conclusion then follows from Corollary \ref{chap3as4}.

\begin{cor}\label{chap3vanish18} Assume that $A = k$ is an algebraically closed field of characteristic $p=2$. Let $V$ a $k$-vector space of finite dimension $r \geq 1$, let $\mathcal{V}$ be the corresponding $k$-group scheme, let $\gamma : V \rightarrow k$ be a non zero linear form, and let $\widetilde{\mathcal{V}}$ be the extension of $\mathcal{V}$ by $\Ga_{a,S}$ defined by the element $c : (v_1,v_2) \mapsto \gamma(v_1) \gamma(v_2)$ of $\mathcal{C}(\mathcal{V},\Ga_{a,S})$, cf. \ref{chap3ext5}. We denote by $t : \widetilde{\mathcal{V}} \rightarrow \Ga_{a,S}$ and by $v : \widetilde{\mathcal{V}} \rightarrow \mathcal{V}$ the canonical projections.

 Let $\mathcal{M}$ be a multiplicative $\Lambda$-local system on $\mathcal{V}$, and let $\alpha$ be an element of $k$. Then the cohomology group
$$
H_c^{\nu}(\widetilde{\mathcal{V}}, \mathcal{M} \otimes \Lc_{\xi}\lbrace \alpha^2 t,0 \rbrace)
$$
vanishes for each integer $\nu$, unless $\nu =2r+2$, $\alpha = 0$ and $\mathcal{M}$ is trivial, or $\nu = 2r+1$, $\alpha$ is non zero and $\mathcal{M}$ is isomorphic to $\Lc_{\xi}\lbrace \alpha^2 t + \delta \gamma(v), \alpha \gamma(v) \rbrace$ for some element $\delta$ of $k$, in which case it is a free $\Lambda$-module of rank $1$.
\end{cor}

By Proposition \ref{chap3charsheafas2}, the multiplicative $\Lambda$-local system $\mathcal{M}$ isomorphic to $\Lc_{\xi}\lbrace \beta^2 t + v^*(v), \beta \gamma(v) \rbrace$ for some $k$-linear form $v^* : V \rightarrow k$, and some element $\beta$ of $k$. The projection formula yields an isomorphism
$$
Rv_! (\mathcal{M} \otimes \Lc_{\xi}\lbrace \alpha^2 t,0 \rbrace) \cong Rv_!(\Lc_{\xi}\lbrace (\beta^2 + \alpha^2) t, 0 \rbrace) \otimes \Lc_{\xi}\lbrace v^*(v), \beta \gamma(v) \rbrace,
$$
and this complex is quasi-isomorphic to $0$ by Corollary \ref{chap3vanish13} and Remark \ref{chap3as89}, unless $\beta = \alpha $, in which case it is quasi-isomorphic to $\Lc_{\xi}\lbrace v^*(v), \beta \gamma(v) \rbrace[-2]$. 
We can thus assume (and we do) that $ \beta = \alpha $, and we must prove that
$$
H_c^{\nu}(\mathcal{V},\Lc_{\xi}\lbrace v^*(v), \alpha \gamma(v) \rbrace)
$$
vanishes for each integer $\nu$, unless $\nu =2r$ and $\mathcal{M}$ is trivial, or $\nu = 2r-1$ and $v^* = \delta \gamma$ for some element $\delta$ of $k$, in which case it is a free $\Lambda$-module of rank $1$.

If $\alpha = 0$, this follows from Corollary \ref{chap3vanish13} and Remark \ref{chap3as89}. We now assume that $\alpha$ is non zero. In this case, by Corollary \ref{chap3vanish13} and Remark \ref{chap3as89} again, the complex
$$
R\gamma_!(\Lc_{\xi}\lbrace v^*(v), \alpha \gamma(v) \rbrace) \cong R\gamma_!(\Lc_{\xi}\lbrace v^*(v),0 \rbrace) \otimes \Lc_{\xi}\lbrace 0, \alpha x \rbrace,
$$
where $x$ is the coordinate on $\Ga_{a,k}$, vanishes unless the restriction of $v^*$ to the kernel of $\gamma$ vanishes, namely if and only $v^* = \delta \gamma$ for some element $\delta$ of $k$, in which case it is isomorphic by the projection formula to
$$
R\gamma_!(\Lambda) \otimes \Lc_{\xi}\lbrace \delta x, 0 \rbrace \otimes \Lc_{\xi}\lbrace 0, \alpha x \rbrace \cong \Lc_{\xi}\lbrace \delta x, \alpha x \rbrace[2-2r].
$$
It remains to prove that 
$$
H_c^{\nu}(\Ga_{a,k},\Lc_{\xi}\lbrace \delta x, \alpha x \rbrace)
$$
vanishes for each integer $\nu$, unless $\nu = 1$, in which case it is of rank $1$. Since $\Lc_{\xi}\lbrace \delta x, \alpha x \rbrace$ has no punctual sections, this group vanishes for $\nu =0$, and by Poincar\'e duality it vanishes as well for $\nu = 2$. In order to conclude, it remains to compute the Euler characteristic with compact supports of $\Lc_{\xi}\lbrace \delta x, \alpha x \rbrace$ on $\Ga_{a,k}$.

The Swan conductor of $\Lc_{\xi}\lbrace \delta x, \alpha x \rbrace$ at infinity is equal to the highest ramification jump of the extension $k((x^{-1}))[t,u]$ of $k((x^{-1}))$ where
\begin{align*}
u - u^2 &= \alpha x, \\
t - t^2 &= u^3 + \delta x,
\end{align*}
corresponding to the equation $(t,u) - (t^2,u^2) = (\delta x,\alpha x)$ in $G$. The only ramification jump of the extension $k((x^{-1}))[t,u]/ k((x^{-1}))[u]$ (resp. $k((x^{-1}))[u]/k((x^{-1}))$), which is a degree $2$ Galois extension, is $3$ (resp. $1$). Thus the extension $k((x^{-1}))[t,u]$ of $k((x^{-1}))$ has two ramification jumps, namely $1$ and some rational number $j > 1$. The slope between $1$ and $3$ of the Herbrand function of the extension $k((x^{-1}))[t,u]/k((x^{-1}))$, cf. (\cite{JPS2}, IV.3), is equal $\frac{1}{2}$, namely the inverse of the degree of the subextension $k((x^{-1}))[u]/k((x^{-1}))$. This slope is also equal to $\frac{j-1}{3-1} = \frac{j-1}{2}$ hence $j-1 = 1$, and thus the second ramification jump of the extension $k((x^{-1}))[t,u]/k((x^{-1}))$ is equal to $j=2$. Consequently, the Swan conductor of $\Lc_{\xi}\lbrace \delta x, \alpha x \rbrace$ at infinity is equal to $2$, and the Grothendieck-Ogg-Shafarevich formula implies that the Euler characteristic with compact supports of $\Lc_{\xi}\lbrace \delta x, \alpha x \rbrace$ on $\Ga_{a,k}$ is equal to $-1$, which concludes the proof of Proposition \ref{chap3vanish18}.

\section{Geometric local \texorpdfstring{$\varepsilon$}{epsilon}-factors for sheaves of generic rank at most 1 \label{chap3lgf1}}


Let $\Lambda$ be an $\ell$-adic coefficient ring (cf. \ref{chap3conv}, \ref{chap30.0.0.1}) which is a field, and let $\psi : \mathbb{F}_{p} \rightarrow \Lambda^{\times}$ be a non trivial homomorphism. We fix a unitary $\Lambda$-admissible mutiplier $\mu$ on the topological group $G_k$ (cf. \ref{chap30.0}, \ref{chap3unitary}).

Let $T$ be the spectrum of a $k$-algebra, which is a henselian discrete valuation ring $\Ow_T$, with maximal ideal $\m$, and whose residue field $\Ow_T/ \m$ is a finite extension of $k$. Let $j : \eta \rightarrow T$ be the generic point of $T$, and let $i : s \rightarrow T$ be its closed point, so that $T$ is canonically an $s$-scheme, as in \ref{chap32.3}. We fix a $\bk$-point $\overline{s} : \Spec(\bk) \rightarrow T$ of $T$ above $s$, so that the Galois group $G_s = \Gal(\bk / k(s))$ can be considered as a subgroup of $G_k$. We still denote by $\mu$ the restriction of $\mu$ to $G_s$. We also fix a geometric point $\overline{\eta}$ of $T$ above $\eta$. 


\subsection{\label{chap32.12}} We denote by $\Omega^1_{\eta} = \Omega^1_{\eta/k}$ the one-dimensional $k(\eta)$-vector space of $1$-forms over $\eta$; it is endowed with a differential $d : k(\eta) \rightarrow \Omega^1_{\eta}$, which is continuous for the valuation topology, and such that $d \pi$ is non zero for any uniformizer $\pi$ of $k(\eta)$. We also denote by $\Omega^{1,\times}_{\eta}$ the $k(\eta)^{\times}$-torsor of non zero elements of $\Omega^1_{\eta}$. If $\omega$ is an element of $\Omega^{1,\times}_{\eta}$, we denote by $v(\omega)$ the unique integer such that for any uniformizer $\pi$ of $k(\eta)$, the element $\frac{\omega}{d\pi}$ of $k(\eta)^{\times}$ has valuation $v(\omega)$. 

\subsection{\label{chap32.11}} Let $\F$ be a $\Lambda$-sheaf on $T$. The \textit{conductor} of the couple $(T,\F)$ is the integer
$$
a(T,\F) = \rk(\F_{\overline{\eta}}) + \sw(\F_{\overline{\eta}}) - \rk(\F_{\overline{s}}).
$$
Following Laumon (\cite{La87}, 3.1.5.1), if $\omega$ is an element of $\Omega^{1,\times}_{\eta}$ (cf. \ref{chap32.12}), we define the conductor of the triple $(T,\F,\omega)$ to be the integer
$$
a(T,\F,\omega) = a(T,\F) + \rk(\F_{\overline{\eta}}) v(\omega).
$$

\subsection{\label{chap32.11998}} Let $f : T' \rightarrow T$ be a finite generically \'etale morphism, where $T'$ is the spectrum of a henselian discrete valuation ring, with generic extension $\eta' \rightarrow \eta$ and residual extension $s' \rightarrow s$. Let $\F$ be a $\Lambda$-sheaf on $T'$, and let $\omega$ be an element of $\Omega^{1,\times}_{\eta'}$. We then have
$$
a(T,f_* \F) = [s':s] a(T,\F)  + v(\partial_{\eta'/\eta}) \rk(\F_{\overline{\eta}}),
$$
by (\cite{JPS2}, IV.2 Prop. $4$), where $\partial_{\eta'/\eta}$ denotes the discriminant of the separable extension $\eta' \rightarrow \eta$, cf. (\cite{JPS2}, III.3). We have $v(\partial_{\eta'/\eta}) = [s':s] v(\partial_{\eta'/\eta_{s'}})$, and
$$
[\eta':\eta_{s'}]v(\omega) + v(\partial_{\eta'/\eta_{s'}}) = v(f^* \omega),
$$
hence the formula
$$
a(T,f_* \F,\omega)  = [s':s] a(T,\F,f^*\omega).
$$

\subsection{\label{chap32.7}} Let $\F$ be a $\mu$-twisted $\Lambda$-sheaf on $T$ \emph{supported on $s$} (cf. \ref{chap31.3}), where $T$ is considered either as a $k$-scheme or as an $s$-scheme (cf. \ref{chap3sanity}), and let $\omega$ be an element of $\Omega^{1,\times}_{\eta}$. For any element $\omega$ of $\Omega^{1,\times}_{\eta}$, we define the \textit{$\varepsilon$-factor} of the triple $(T,\F,\omega)$ to be the $\Lambda$-admissible map
\begin{align*}
\varepsilon_{ \overline{s}}(T,\F,\omega) : G_s &\rightarrow \Lambda^{\times} \\
g &\mapsto \det \left( g \ | \ \F_{ \overline{s}} \right)^{-1}
\end{align*}

 The map $\varepsilon_{ \overline{s}}(T,\F,\omega)$ defines a $\Lambda$-admissible representation of rank $1$ of $(G_s,\mu^{-\rk \F_{ \overline{s}}})$ which is isomorphic to $\det \left(\F_{ \overline{s}} \right)^{-1}$. In particular, we have
$$
d^1(\varepsilon_{ \overline{s}}(T,\F,\omega)) = \mu^{-\rk \F_{ \overline{s}}} = \mu^{a(T,\F,\omega)},
$$
cf. \ref{chap30.11} and \ref{chap32.11} for the notation.

\subsection{\label{chap32.13}} Let $\F$ be a $\mu$-twisted $\Lambda$-sheaf on $T$ (cf. \ref{chap31.3}, \ref{chap3sanity}), \emph{supported on $\eta$}, such that $j^{-1}\F$ is of rank $1$, and let $\omega$ be an element of $\Omega^{1,\times}_{\eta}$. Then $j^{-1} \F$ is a $\mu$-twisted $\Lambda$-local system of rank $1$ on $\eta$. Let $D = \nu s$ be an effective Cartier divisor on $T$ such that $j^{-1}\F$ has ramification bounded by $D$ (cf. \ref{chap3ramdefiloc2}). Theorem \ref{chap3lgcft2twisted} then produces a $\mu$-twisted multiplicative $\Lambda$-local system $\chi_{j^{-1}\F}$ (cf. \ref{chap3twistedmult}, \ref{chap3sanity2}) on the $s$-group scheme $\Pic(T,D)_s$ (cf. \ref{chap3localpic}).

Recall from \ref{chap32.5} that $\Pic(T,D)_s$ is naturally isomorphic to the functor which sends an $s$-scheme $S$ to the group of pairs $(d,u)$, where $d$ is a locally constant $\mathbb{Z}$-valued map on $S$, and 
$$u : \Ow_S \otimes_{k(s)} \Ow_T/ \m^{\nu} \rightarrow \Ow_S \otimes_{k(s)} \m^{-d} / \m^{-d + \nu} $$
is an isomorphism of $\Ow_S \otimes_{k(s)} \Ow_T/ \m^{\nu}$-modules. Denoting by $\Pic^d(T,D)_s$ the component of degree $d$ of $\Pic(T,D)_s$, we consider the morphism
\begin{align*}
\Res_{\omega} : \Pic^{a(T,\F,\omega)}(T,D)_s &\rightarrow \mathbb{G}_{a,s} \\
u &\mapsto \Res(u \omega),
\end{align*}
cf. \ref{chap3residue}, which is well defined since $\nu - a(T,\F,\omega)$ is greater than or equal to $-v(\omega)$.

\begin{prop}\label{chap3rankcohom} The cohomology group 
$$
H_c^j \left(\Pic^{a(T,\F,\omega)}(T,D)_{\overline{s}},\chi_{j^{-1}\F} \otimes \Lc_{\psi} \lbrace \Res_{\omega} \rbrace \right),
$$
vanishes for $j \neq 2\nu - a(T,\F) $, and is a $\Lambda$-module of rank $1$ if $j = 2\nu - a(T,\F)$.
\end{prop}

The Artin-Schreier sheaf $\Lc_{\psi} \lbrace \Res_{\omega} \rbrace$ is defined in the paragraph \ref{chap34.3}. An equivalent version of Proposition \ref{chap3rankcohom} appears with a lacunary proof in Section $g.(B)$ of Deligne's $1974$ letter to Serre, published as an appendix in \cite{bloch}. We postpone the proof of Proposition \ref{chap3rankcohom} to the paragraphs \ref{chap32.14}, \ref{chap32.14.1}, \ref{chap32.14.2}, \ref{chap32.14.3} and \ref{chap32.14.4} below.

\begin{defi}\label{chap3localfact1}Let $\F$ be a $\mu$-twisted $\Lambda$-sheaf on $T$, \emph{supported on $\eta$}, such that $j^{-1}\F$ is of rank $1$, and let $D$ be an effective Cartier divisor on $T$ such that $j^{-1}\F$ has ramification bounded by $D$ (cf. \ref{chap3ramdefiloc}), namely $\sw(\F_{\overline{\eta}})$ is strictly less than the multiplicity $\nu$ of $D$ at $s$. Let $\omega$ be an element of $\Omega^{1,\times}_{\eta}$. The \textit{$\varepsilon$-factor} of the triple $(T,\F,\omega)$ is the $\Lambda$-admissible map $\varepsilon_{ \overline{s}}(T,\F,\omega) : G_s \rightarrow \Lambda^{\times}$ such that
\begin{align*}
\varepsilon_{ \overline{s}}(T,\F,\omega)(g) = \mathrm{Tr} \left( g \ | \ H_c^{2\nu - a(T,\F)} \left(\Pic^{a(T,\F,\omega)}(T,D)_{ \overline{s}},\chi_{j^{-1}\F} \otimes \Lc_{\psi} \lbrace \Res_{\omega} \rbrace (\nu - a(T,\F,\omega)) \right) \right),
\end{align*}
for any $g$ in $G_s$.
\end{defi}

Since $\chi_{j^{-1}\F} \otimes \Lc_{\psi} \lbrace \Res_{\omega} \rbrace$ is a $\mu^{a(T,\F,\omega)}$-twisted $\Lambda$-sheaf on $\Pic^{a(T,\F,\omega)}(T,D)_{s}$ (cf. \ref{chap3twistedmult}), and since the cohomology group in \ref{chap3localfact1} is of rank $1$ by Proposition \ref{chap3rankcohom} we obtain
$$
d^1(\varepsilon_{ \overline{s}}(T,\F,\omega)) = \mu^{a(T,\F,\omega)},
$$
cf. \ref{chap31.3.4}.

The notation $\varepsilon_{ \overline{s}}(T,\F,\omega)$ suggests that the choice of $D$ is irrelevant. We have indeed:

\begin{prop}\label{chap3rankcohom2}Let $\F,D, \omega$ be as in \ref{chap3localfact1}. Then the map $\varepsilon_{ \overline{s}}(T,\F,\omega)$ is independent of the choice of $D$.
\end{prop}

Indeed, if $\delta$ is a positive integer, let us consider the projection morphism
$$
\tau : \Pic^{a(T,\F,\omega)}(T,D + \delta s)_s \rightarrow \Pic^{a(T,\F,\omega)}(T,D)_s,
$$
which sends for any $k$-scheme $S$ a trivialization $u$ of $\Ow_S \otimes_{k(s)} \m^{-a(T,\F,\omega)}/ \m^{-a(T,\F,\omega) + \nu + \delta}$ to its image in $\Ow_S \otimes_{k(s)} \m^{-a(T,\F,\omega)}/ \m^{-a(T,\F,\omega) + \nu}$. It is a (trivial) fibration in affine spaces of dimension $\delta$. Let $\G$ be the $\Lambda$-sheaf $\chi_{j^{-1}\F} \otimes \Lc_{\psi} \lbrace \Res_{\omega} \rbrace$ on $\Pic^{a(T,\F,\omega)}(T,D)_s$. The trace homomorphism
$$
R \tau_! \tau^{-1} \G (\delta)[2 \delta]\rightarrow \G,
$$
is an isomorphism, and thus the Leray spectral sequence for $(\tau, \G)$ yields that the $\Lambda$-admissible representation given by
$$
H_c^{2\nu + 2\delta - a(T,\F)} \left(\Pic^{a(T,\F,\omega)}(T,D + \delta)_{ \overline{s}},\tau^{-1} \G(\nu + \delta - a(T,\F,\omega)) \right),
$$
is isomorphic to
$$
H_c^{2\nu - a(T,\F)} \left(\Pic^{a(T,\F,\omega)}(T,D )_{ \overline{s}}, \G(\nu - a(T,\F,\omega)) \right),
$$
hence the result.


\subsection{\label{chap32.15}} Let $\F$ be a $\mu$-twisted $\Lambda$-sheaf on $T$, such that $\F_{\overline{\eta}}$ is of rank at most $1$ over $\Lambda$, and let $\omega$ be an element of $\Omega^{1,\times}_{\eta}$. If $\F$ is supported on a single point of $T$, then we defined the $\varepsilon$-factor of $(T,\F,\omega)$ in \ref{chap32.7} and \ref{chap3localfact1}. We combine these two definitions as follows:

\begin{defi}\label{chap3localfact2} Let $\F$ be a $\mu$-twisted $\Lambda$-sheaf on $T$, such that $j^{-1}\F$ is of rank at most $1$ over $\Lambda$. The \textit{$\varepsilon$-factor} of the triple $(T,\F,\omega)$ is the $\Lambda$-admissible map from $G_s$ to $\Lambda^{\times}$ given by
$$
\varepsilon_{ \overline{s}}(T,\F,\omega) = \varepsilon_{ \overline{s}}(T,j_{!} j^{-1} \F, \omega) \varepsilon_{ \overline{s}}(T,i_* i^{-1} \F, \omega).
$$
\end{defi}

It follows from \ref{chap32.7} and \ref{chap3localfact1} that we have
$$
d^1(\varepsilon_{ \overline{s}}(T,\F,\omega)) = \mu^{a(T,\F,\omega)}.
$$

\subsection{\label{chap32.30}} Let $\F$ be a $\mu$-twisted $\Lambda$-sheaf of rank $1$ on $\eta$, with ramification bounded by $D = \nu s$ (cf. \ref{chap3ramdefiloc2}). If $z$ is an element of valuation $d$ in $k(\eta)^{\times}$, then the image of $z^{-1}$ in $\m^{-d} / \m^{-d + \nu}$ yields an $s$-point of $\Pic^d(T,D)_s$. The restriction $\chi_{\F | z^{-1}}$ of $\chi_{\F}$ to this $s$-point is a $\Lambda$-local system of rank $1$ on $s$, and we define 
\begin{align*}
\langle \chi_{\F} \rangle (z) : G_s &\rightarrow \Lambda^{\times} \\
g &\mapsto \det \left( g \ | \ (\chi_{\F | z^{-1}})_{ \overline{s}} \right),
\end{align*}
so that $d^1(\langle \chi_{\F} \rangle (z) ) = \mu^{ d}$. The map $\langle \chi_{\F} \rangle (z) $ depends only on $\F$ and $z$, and not on the choice of $D$. By multiplicativity of the local system $\chi_{\F}$ (cf. \ref{chap32.0}), we have
$$
\langle \chi_{\F} \rangle (z_1 z_2) = \langle \chi_{\F} \rangle (z_1) \langle \chi_{\F} \rangle (z_2) 
$$
for any elements $z_1,z_2$ of $k(\eta)^{\times}$. 

\begin{prop}\label{chap3changeforme} Let $\F$ be a $\mu$-twisted $\Lambda$-sheaf on $T$, such that $j^{-1}\F$ is of rank $1$ over $\Lambda$. For any element $\alpha$ of $k(\eta)^{\times}$, of valuation $v(\alpha)$, we have
$$
\varepsilon_{ \overline{s}}(T,\F,\alpha \omega) = \langle \chi_{j^{-1}\F}\rangle(\alpha) \chi_{\cyc}^{- v(\alpha)} \varepsilon_{ \overline{s}}(T,\F,\omega),
$$
where $\langle \chi_{j^{-1}\F}\rangle$ is as in \ref{chap32.30}, and $\chi_{\cyc}$ is the $\ell$-adic cyclotomic character of $k$ (cf. \ref{chap3conv}).
\end{prop}

Indeed, we have an isomorphism
\begin{align*}
\theta : \Pic^{a(T,\F, \omega)}(T,D)_s &\rightarrow \Pic^{a(T,\F,\alpha \omega)}(T,D)_s \\
u &\mapsto \alpha^{-1} u,
\end{align*}
of $s$-schemes, such that the pullback of $\chi_{j^{-1}\F} \otimes \Lc_{\psi} \lbrace \Res_{\alpha \omega} \rbrace (\nu - a(T,\F, \alpha \omega))$ by $\theta$ is isomorphic to the twist of $\chi_{j^{-1}\F} \otimes \Lc_{\psi} \lbrace \Res_{\omega} \rbrace (\nu - a(T,\F,\omega))$ by $(\chi_{j^{-1}\F})_{|\alpha^{-1}}(- v(\alpha))$.

\begin{prop}\label{chap3computation} Let $n \geq 1$ be an integer prime to $p$, and let $h$ be an element of $k(\eta)$ of valuation $-n$, and let us consider the Artin-Schreier $\Lambda$-sheaf $\Lc_{\psi} \lbrace h \rbrace$ on $\eta$ (cf. \ref{chap34.3}). Let $\Mac$ be a $\mu$-twisted $\Lambda$-sheaf on $\eta$ of rank $1$ with ramification bounded by the divisor $\lceil \frac{n}{2} \rceil s$, and let us consider the $\mu$-twisted $\Lambda$-sheaf $\F = \Mac \otimes \Lc_{\psi} \lbrace -h \rbrace$. For any element $\omega$ of $\Omega^{1,\times}_{\eta}$, we have:
\begin{itemize}
\item if $n = 2n' - 1$ is odd then 
$$
\varepsilon_{ \overline{s}}(T,j_!\F, \omega) = \langle \chi_{\F}\rangle \left( \frac{\omega}{dh}\right) \chi_{\cyc}^{- v\left( \frac{\omega}{dh}\right) + n'},
$$
\item if $n = 2n'$ is even then $p$ is odd and we have
$$
\varepsilon_{ \overline{s}}(T,j_!\F, \omega) = \langle \chi_{\F}\rangle \left( \frac{\omega}{dh}\right) \chi_{\cyc}^{- v\left( \frac{\omega}{dh}\right) + n' + 1} \gamma_{\psi}(-n h_0),
$$
where $h_0$ is an element of $k(s)^{\times}$ such that $\frac{h}{h_0}$ is a square in $k(\eta)^{\times}$, and where, for any element $c$ of $k(s)$, we have set
\begin{align*}
\gamma_{\psi}(c) : G_s &\rightarrow \Lambda^{\times}\\
g &\mapsto \det \left( g \ | \ H_c^1 \left( \Ga_{a, \overline{s}}, \Lc_{\psi}\lbrace \frac{c t^2}{2} \rbrace) \right)\right).
\end{align*}
\end{itemize}
\end{prop}

When $k$ is a finite field of odd characteristic, the conclusion of Proposition \ref{chap3computation} follows from (\cite{abbes}, Prop. 8.7) and from Proposition \ref{chap3compafini} below.

We note that the Swan conductor of a $\mu$-twisted $\Lambda$-sheaf $\F$ as in Proposition \ref{chap3computation} is $n$. By Proposition \ref{chap3changeforme}, it is sufficient to prove Proposition \ref{chap3computation} when $\omega = \frac{d \pi}{\pi^{n+1}}$, for a fixed uniformizer $\pi$, in which case $\frac{\omega}{d h}$ has valuation $0$. The conclusion then follows from the proof of Proposition \ref{chap3rankcohom} below, cf. Porisms \ref{chap3por1} and \ref{chap3por2}.

\subsection{\label{chap32.14}} We now prove Proposition \ref{chap3rankcohom}. Let us first consider the projection morphism
$$
\tau : \Pic^{a(T,\F,\omega)}(T,D)_s \rightarrow \Pic^{a(T,\F,\omega)}(T,a(T,\F)s)_s,
$$
of relative dimension $\delta = \nu - a(T,\F)$, which sends for any $k$-scheme $S$ a trivialization $u$ of $\Ow_S \otimes_{k(s)} \m^{-a(T,\F,\omega)}/ \m^{-a(T,\F,\omega) + \nu}$ to its image in $\Ow_S \otimes_{k(s)} \m^{-a(T,\F,\omega)}/ \m^{-v(\omega)}$. Let $\G$ be the $\Lambda$-sheaf $\chi_{j^{-1}\F} \otimes \Lc_{\psi} \lbrace \Res_{\omega} \rbrace$ on the $s$-group scheme $\Pic^{a(T,\F,\omega)}(T,a(T,\F)s)_s$. As in the proof of \ref{chap3rankcohom2}, we have an isomorphism
$$
R \tau_! \tau^{-1} \G (\delta)[2 \delta]\rightarrow \G,
$$
hence the Leray spectral sequence for $(\tau, \G)$ implies that it is enough to prove Proposition \ref{chap3rankcohom} when the multiplicity $\nu$ of $D$ at $s$ is exactly $a(T,\F)$.

Let us now assume that $\nu$ is equal to $a(T,\F)$. Let $\pi$ be a uniformizer of $k(\eta)$. Let us write $\omega = \alpha^{-1} \frac{d \pi}{\pi^{\nu}}$ for some element $\alpha$ of $k(\eta)^{\times}$ of valuation $-a(T,\F,\omega) = - \nu - v(\omega)$, and let us consider the isomorphism
\begin{align*}
\theta : \Pic^{0}(T,D)_s &\rightarrow \Pic^{a(T,\F,\omega)}(T,D)_s \\
u &\mapsto \alpha u.
\end{align*}
The pullback of $\chi_{j^{-1}\F}$ by $\theta$ coincides with $\chi_{j^{-1}\F}$ on $\Pic^{0}(T,D)_s$, up to twist by the fiber $\alpha^{-1}\chi_{j^{-1}\F }$ of $\chi_{j^{-1}\F }$ at the $s$-point $\alpha$, while the pullback of $\Lc_{\psi} \lbrace \Res_{\omega} \rbrace$ is isomorphic to $\Lc_{\psi} \lbrace \Res_{\pi^{-\nu} d \pi} \rbrace$. Thus we can assume that $\omega$ is equal to $\frac{d \pi}{\pi^{\nu}}$, so that we have $a(T,\F,\omega) = 0$.

\subsection{\label{chap32.14.1}} Let us prove Proposition \ref{chap3rankcohom} when $\omega$ is equal to $\frac{d \pi}{\pi^{\nu}}$ and when $\nu = a(T,\F)$ is equal to $1$. In this case, the $s$-group scheme $\Pic^{0}(T,D)_s$ is simply the multiplicative group $\mathbb{G}_{m,s}$, and $\Lc_{\psi} \lbrace \Res_{\omega} \rbrace$ coincides with $\Lc_{\psi} \lbrace t \rbrace$. Any multiplicative $\Lambda$-local system on $\mathbb{G}_{m,s}$, such as $\chi_{j^{-1} \F}$, is tamely ramified at $0$ and $\infty$, cf. for example (\cite{G18}, 3.15). Thus the $\Lambda$-local system $\chi_{j^{-1}\F} \otimes \Lc_{\psi} \lbrace t \rbrace$ on $\mathbb{G}_{m,s}$ has Swan conductor $0$ at $0$, and $1$ at infinity. By the Grothendieck-Ogg-Shafarevich formula, we have
$$
\chi_c(\mathbb{G}_{m,\overline{s}},\chi_{j^{-1}\F} \otimes \Lc_{\psi} \lbrace t \rbrace ) = 2 - 1 - 2 = -1.
$$
Moreover, the cohomology groups $H_c^0(\mathbb{G}_{m,\overline{s}},\chi_{j^{-1}\F} \otimes \Lc_{\psi} \lbrace t \rbrace )$ and $H^0(\mathbb{G}_{m,\overline{s}},\chi_{j^{-1}\F}^{-1} \otimes \Lc_{\psi} \lbrace -t \rbrace )$ both vanish, and so does the group
$$H_c^2(\mathbb{G}_{m,\overline{s}},\chi_{j^{-1}\F} \otimes \Lc_{\psi} \lbrace t \rbrace ),$$
by Poincar\'e duality. This proves that the cohomology group
$$
H_c^j(\mathbb{G}_{m,\overline{s}},\chi_{j^{-1}\F} \otimes \Lc_{\psi} \lbrace t \rbrace ),
$$
vanishes when $j \neq 1$, and is of rank $1$ when $j$ is equal to $1$.

\subsection{\label{chap32.14.2}} Let us prove Proposition \ref{chap3rankcohom} when $\omega$ is equal to $\frac{d \pi}{\pi^{\nu}}$ and when $\nu = a(T,\F)$ is even, hence of the form $2 \nu'$ for some integer $\nu' \geq 1$. Let us consider the projection morphism
$$
\sigma : \Pic^{0}(T,D)_s \rightarrow \Pic^{0}(T,\nu' s)_s,
$$
which is a homomorphism of $s$-group schemes. Let $\V$ be the additive $s$-group scheme associated to the finite dimensional $s$-vector space $V = \m^{\nu'}/\m^{\nu}$. Then the morphism
\begin{align*}
r : \V &\rightarrow \Pic^{0}(T,D)_s \\
x &\mapsto 1 + x,
\end{align*}	
realizes an isomorphism of $s$-group schemes from $\V$ onto the kernel of $\sigma$. By Proposition \ref{chap3charsheafas}, the multiplicative $\Lambda$-local system $r^{-1}\chi_{j^{-1} \F}$ is isomorphic to an Artin-Schreier sheaf $\Lc_{\psi} \lbrace - v^* \rbrace$, for some linear form $v^*$ on $V$. The $k(s)$-linear map
\begin{align}
\begin{split}\label{chap3dual}
\Ow_T/\m^{\nu'} &\rightarrow \Hom_{k(s)}(V,k(s)) \\
y &\mapsto \Res_{y \omega} = (x \mapsto \Res(x y \omega)),
\end{split}
\end{align}
is an isomorphism onto the $k(s)$-linear dual of $V = \m^{\nu'}/\m^{\nu}$. Thus there exists a unique element $y$ of $\Ow_T/\m^{\nu'}$ such that $v^*$ is equal to $\Res_{y \omega}$.

Since the ramification of $\chi_{j^{-1} \F}$ is not bounded by the divisor $(\nu -1)s = \sw(\F_{\overline{\eta}}) s$, it follows from Theorem \ref{chap3lgcft2twisted} and from \ref{chap32.0.5} that the restriction of $r^{-1}\chi_{j^{-1} \F}$ to the sub-$s$-group scheme corresponding to $\m^{\nu-1}/\m^{\nu}$ is non trivial. Thus the restriction of $v^* = \Res_{y \omega}$ to $\m^{\nu-1}/\m^{\nu}$ is non trivial, and consequently $y$ is a unit of $\Ow_T/\m^{\nu'}$. In particular, $y$ defines an $s$-point of $\Pic^{0}(T,\nu' s)_s$.

Let $\overline{t}$ be the spectrum of an algebraically closed extension of $k(s)$, and let $u$ be a $\overline{t}$-point of $\Pic^{0}(T,D)_{\overline{t}}$. The morphism
\begin{align*}
ur : \V_{\overline{t}} &\rightarrow \Pic^{0}(T,D)_{\overline{t}} \\
x &\mapsto u(1 + x),
\end{align*}
realizes an isomorphism onto the fiber of $\sigma$ above $\sigma(u)$. Let $\G$ be the $\Lambda$-sheaf $\chi_{j^{-1}\F} \otimes \Lc_{\psi} \lbrace \Res_{\omega} \rbrace$ on the $s$-group scheme $\Pic^{0}(T,D)_s$. The pullback of $\G$ by $ur$ is isomorphic to $\Lc_{\psi} \lbrace \Res_{(u-y) \omega} \rbrace$, with notation as in (\ref{chap3dual}), up to twist by the stalk $\G_u $ of $\G$ at the geometric point $u$ of $\Pic^{0}(T,D)_s$. Together with Proposition \ref{chap3vanish13}, this implies that the complex
$$
R\Gamma_c( \sigma^{-1}(\sigma(u)), \G),
$$
vanishes unless $\sigma(u)$ is equal to $y$, in which case it is concentrated in degree $2 \nu' = \nu$, and the cohomology group $H_c^{\nu}( \sigma^{-1}(\sigma(u)), \G)$ is of rank $1$. We obtain that $R \sigma_{!} \G$ is of the form $y_* \Lc[-\nu]$, where $\Lc$ is a $\Lambda$-sheaf of rank $1$ on $s$, and thus that the complex
$$
R\Gamma_c(\Pic^{0}(T,D)_{\overline{s}}, \G) \cong R\Gamma_c(\Pic^{0}(T,\nu' s)_{\overline{s}}, R \sigma_{!} \G),
$$
is isomorphic to $\Lc_{\overline{s}}[-\nu]$, hence the conclusion of Proposition \ref{chap3rankcohom}.

\begin{por}\label{chap3por1} Let us assume that $\F = \Mac \otimes \Lc_{\psi} \lbrace - h \rbrace$ is as in Proposition \ref{chap3computation}, so that $\nu = n+1$ and $n = 2 \nu' - 1$. We keep the notation from \ref{chap32.14.2}. Since $\Mac$ has ramification bounded by $\nu'$, the restriction $r^{-1} \chi_{\Mac}$ is trivial. Moreover, by \ref{chap3gkexample} and \ref{chap3lgcft1}, we have $\chi_{\F} = \chi_{\Mac} \otimes \Lc_{\psi} \lbrace \Res(h\frac{du}{u}) \rbrace$, where $u$ is the universal unit parametrized by $\Pic(T,D)_s$. For any section $x$ of $\V$, we have 
\begin{align*}
\Res \left( h\frac{d(1+x)}{1+x} \right) &= \Res(h dx) - \Res \left(xh\frac{d x}{1+x} \right) \\
&= -\Res(x dh) - \Res \left(xh\frac{d x}{1+x} \right),
\end{align*}
which is equal to $ -\Res(x dh)$ since $\frac{xh dx}{1+x}$ has nonnegative valuation. Thus the element $y$ of $\Ow_T/\m^{\nu'}$ which appears in the proof above is equal to $\frac{dh}{\omega}$. Thus $R^{\nu} \sigma_{!} \G$ is concentrated on the $s$-point $\frac{dh}{\omega}$ of $\Pic^{0}(T,\nu' s)_s$, and its restriction to the $s$-point $\frac{dh}{\omega}$ is $\chi_{\F | \frac{dh}{\omega}}(-\nu')$, hence the conclusion of Proposition \ref{chap3computation} when $n$ is odd.
\end{por}

\subsection{\label{chap32.14.3}} Let us prove Proposition \ref{chap3rankcohom} when $\omega$ is equal to $\frac{d \pi}{\pi^{\nu}}$, when $\nu = a(T,\F)$ is of the form $2 \nu' + 1$ for some integer $\nu' \geq 1$, and when the characteristic $p$ of $k$ is odd. As in \ref{chap32.14.2}, let us consider the projection morphisms
%
%
\begin{center}
 \begin{tikzpicture}[scale=1]

\node (A) at (0,0) {$\Pic^{0}(T,D)_s$};
\node (B) at (4,0) {$\Pic^{0}(T,(\nu'+1) s)_s$};
\node (C) at (8,0) {$\Pic^{0}(T,\nu' s)_s,$};

\path[->,font=\scriptsize]
(A) edge node[above]{$\sigma_1$} (B)
(B) edge node[above]{$\sigma_2$} (C)
(A) edge[bend right] node[above]{$ \sigma$} (C);
\end{tikzpicture} 
\end{center}
so that $\sigma, \sigma_1$ and $\sigma_2$ are all homomorphisms of $s$-group schemes. Let $\V$ be the additive $s$-group scheme associated to the finite dimensional $s$-vector space $V = \m^{\nu'}/\m^{\nu}$. Then the morphism
\begin{align*}
r : \V &\rightarrow \Pic^{0}(T,D)_s \\
x &\mapsto 1 + x + \frac{x^2}{2},
\end{align*}	
realizes an isomorphism of $s$-group schemes from $\V$ to the kernel of $\sigma$. By Proposition \ref{chap3charsheafas}, the multiplicative $\Lambda$-local system $r^{-1}\chi_{j^{-1} \F}$ is isomorphic to an Artin-Schreier sheaf $\Lc_{\psi} \lbrace - v^* \rbrace$, for some linear form $v^*$ on $V$. Let $y$ be the unique element of $\Ow_T/\m^{\nu'+1}$ such that $v^*$ coincides with the linear form
\begin{align*}
\Res_{y\omega} : \m^{\nu'}/\m^{\nu} &\rightarrow k(s) \\
x &\mapsto \Res(x y \omega).
\end{align*}
As in \ref{chap32.14.2}, the element $y$ is a unit of $\Ow_T/\m^{\nu'+1}$, whence $y$ defines an $s$-point of $\Pic^{0}(T,(\nu'+1) s)_s$.

Let $\overline{t}$ be the spectrum of an algebraically closed extension of $k(s)$, and let $u$ be a $\overline{t}$-point of $\Pic^{0}(T,D)_{\overline{t}}$. The morphism
\begin{align*}
ur : \V_{\overline{t}} &\rightarrow \Pic^{0}(T,D)_{\overline{t}} \\
x &\mapsto u(1 + x + \frac{x^2}{2}),
\end{align*}
realizes an isomorphism onto the fiber of $\sigma$ above $\sigma(u)$. Let $\G$ be the $\Lambda$-sheaf $\chi_{j^{-1}\F} \otimes \Lc_{\psi} \lbrace \Res_{\omega} \rbrace$ on the $s$-group scheme $\Pic^{0}(T,D)_s$. The pullback of $\G$ by $ur$ is isomorphic to $\Lc_{\psi} \lbrace \Res_{(u-y) \omega} + \alpha \gamma^2 \rbrace$, up to twist by the stalk $\G_u $ of $\G$ at the geometric point $u$ of $\Pic^{0}(T,D)_s$, where $\alpha$ is the image of $\frac{u}{2}$ in $k(\overline{t})^{\times}$ and $\gamma : V \rightarrow k(s)$ is the linear form which sends an element $x$ to the image of $\pi^{-\nu'} x$ in $k(s)$. Together with Proposition \ref{chap3vanish14}, this implies that the complex
$$
R\Gamma_c( \sigma^{-1}(\sigma(u)), \G),
$$
vanishes unless the linear form $\Res_{(u-y) \omega}$ on $V$ is proportional to $\gamma$, namely unless $\sigma(u) = \sigma_2(y)$, in which case this complex is concentrated in degree $2(\nu'+1)-1 = \nu$, and the cohomology group $H_c^{\nu}( \sigma^{-1}(\sigma(u)), \G)$ is of rank $1$. We obtain that $R \sigma_{!} \G$ is of the form $\sigma_2(y)_* \Lc[-\nu]$, where $\Lc$ is a $\Lambda$-sheaf of rank $1$ on $s$, and thus that the complex
$$
R\Gamma_c(\Pic^{0}(T,D)_{\overline{s}}, \G) \cong R\Gamma_c(\Pic^{0}(T,\nu' s)_{\overline{s}}, R \sigma_{!} \G),
$$
is isomorphic to $\Lc_{\overline{s}}[-\nu]$, hence the conclusion of Proposition \ref{chap3rankcohom}.

\begin{por}\label{chap3por2} Let us assume that $\F = \Mac \otimes \Lc_{\psi} \lbrace - h \rbrace$ is as in Proposition \ref{chap3computation}, so that $\nu = n+1$ and $n = 2 \nu' $. We keep the notation from \ref{chap32.14.3}. Since $\Mac$ has ramification bounded by $\nu'$, the restriction $r^{-1} \chi_{\Mac}$ is trivial. Moreover, by \ref{chap3gkexample} and \ref{chap3lgcft1}, we have $\chi_{\F} = \chi_{\Mac} \otimes \Lc_{\psi} \lbrace \Res(h\frac{du}{u}) \rbrace$, where $u$ is the universal unit parametrized by $\Pic(T,D)_s$ (cf. \ref{chap32.5}). For any section $x$ of $\V$, we have 
\begin{align*}
\Res \left( h\frac{d(1+x + \frac{x^2}{2})}{1+x + \frac{x^2}{2}} \right) &= \Res(h dx) - \Res \left(x^2h\frac{d x}{2(1+x)} \right) \\
&= -\Res(x dh) - \Res \left(x^2h\frac{d x}{2(1+x)} \right),
\end{align*}
which is equal to $ -\Res(x dh)$ since $x^2h\frac{d x}{2(1+x)} $ has nonnegative valuation. Thus the element $y$ of $\Ow_T/\m^{\nu'+1}$ which appears in the proof above is equal to $\frac{dh}{\omega}$. Thus $R^{\nu} \sigma_{!} \G$ is concentrated on the $s$-point $\frac{dh}{\omega}$ of $\Pic^{0}(T,\nu' s)_s$. Moreover, choosing $u = y = \frac{dh}{\omega}$ in the computation above, the restriction $(ur)^{-1} \chi_{\F} \otimes \Lc_{\psi} \lbrace \Res_{\omega} \rbrace$ is isomorphic to $\chi_{\F | \frac{dh}{\omega}} \otimes \Lc_{\psi} \lbrace \alpha \gamma^2 \rbrace$, where $\alpha$ is the image of $\frac{dh}{2 \omega}$ in $k(s)^{\times}$. If $h_0$ is the image in $k(s)^{\times}$ of $\pi^{n} h$, then we have $\alpha = - \frac{n h_0}{2}$. We thus have
$$
\varepsilon_{ \overline{s}}(T,j_!\F, \omega)(g) = \langle \chi_{\F}\rangle \left( \frac{\omega}{dh}\right)(g) \chi_{\cyc}(g)^{\nu' + 1} \mathrm{Tr} \left( g \ | \ H_c^{\nu} \left( \V_{ \overline{s}}, \Lc_{\psi}\lbrace - \frac{n h_0 \gamma^2}{2} \rbrace) \right)\right),
$$
for any $g$ in $G_s$, hence the conclusion of Proposition \ref{chap3computation} when $n$ is even.
\end{por}

\subsection{\label{chap32.14.4}} Let us prove Proposition \ref{chap3rankcohom} when $\omega$ is equal to $\frac{d \pi}{\pi^{\nu}}$, when $\nu = a(T,\F)$ is of the form $2 \nu' + 1$ for some integer $\nu' \geq 1$, and when $k$ is of characteristic $p =2$. As in \ref{chap32.14.3}, let us consider the projection morphisms
%
%
\begin{center}
 \begin{tikzpicture}[scale=1]

\node (A) at (0,0) {$\Pic^{0}(T,D)_s$};
\node (B) at (4,0) {$\Pic^{0}(T,(\nu'+1) s)_s$};
\node (C) at (8,0) {$\Pic^{0}(T,\nu' s)_s,$};

\path[->,font=\scriptsize]
(A) edge node[above]{$\sigma_1$} (B)
(B) edge node[above]{$\sigma_2$} (C)
(A) edge[bend right] node[above]{$ \sigma$} (C);
\end{tikzpicture} 
\end{center}
so that $\sigma, \sigma_1$ and $\sigma_2$ are all homomorphisms of $s$-group schemes. Let $\V$ be the additive $s$-group scheme associated to the finite dimensional $k(s)$-vector space $V = \m^{\nu'}/\m^{\nu}$. Let $\gamma : V \rightarrow k(s)$ be the $k(s)$-linear form which sends an element $v$ of $V$ to the image in $k(s)$ of the element $\pi^{-\nu'} v$ of $\Ow_T/ \m^{\nu'+1}$, and let $\widetilde{\mathcal{V}}$ be the extension of $\mathcal{V}$ by $\Ga_{a,s}$ defined by the element $c : (v_1,v_2) \mapsto \gamma(v_1) \gamma(v_2)$ of $\mathcal{C}(\mathcal{V},\Ga_{a,S})$, cf. \ref{chap3ext5}. For any sections $(t_1,v_1)$ and $(t_2,v_2)$ of the $k(s)$-scheme $\widetilde{\mathcal{V}} = \Ga_{a,s} \times_s \V$, we have
$$
(t_1,v_1) + (t_2,v_2) = (t_1 + t_2 + \gamma(v_1) \gamma(v_2), v_1 + v_2),
$$
in $\widetilde{\mathcal{V}}$, and we have 
$$
(1 + v_1 + \pi^{2 \nu'} t_1)(1 + v_2 + \pi^{2 \nu'} t_2) = 1 + v_1 + v_2 + \pi^{2 \nu'} \left( t_1 + t_2 + (\pi^{-\nu'} v_1)(\pi^{-\nu'} v_2) \right),
$$
in $\Pic^{0}(T,D)_s$. Thus the morphism
\begin{align*}
r : \widetilde{\mathcal{V}} &\rightarrow \Pic^{0}(T,D)_s \\
(t,v) &\mapsto 1 + v + \pi^{2 \nu'} t,
\end{align*}
of $s$-schemes is a homomorphism of $s$-group schemes.	It surjects onto the kernel of $\sigma$, and its kernel is isomorphic to $\Ga_{a,s}$ (cf. \ref{chap3kernelisom} below).

Let $G$ be the extension of $\Ga_{a,S}$ by itself defined by $c(U_1,U_2) = U_1 U_2$ (cf. \ref{chap34.4}), and let $\xi : G(\mathbb{F}_2) \rightarrow \Lambda^{\times}$ be an injective character (cf. \ref{chap34.7}), whose restriction to the subgroup $\mathbb{F}_2$ of $ G(\mathbb{F}_2)$ is $\psi$. By Proposition \ref{chap3charsheafas2}, the multiplicative $\Lambda$-local system $r^{-1}\chi_{j^{-1} \F}$ is isomorphic to $\Lc_{\xi} \lbrace \alpha^2 t + v^*(v), \alpha \gamma(v) \rbrace$ (cf. \ref{chap34.7} for the notation), for some linear form $v^*$ on $V$ and some element $\alpha$ of $k$. Let $y$ be the unique element of $\Ow_T/\m^{\nu'+1}$ such that $v^*$ coincides with the linear form
\begin{align*}
\Res_{y\omega} : \m^{\nu'}/\m^{\nu} &\rightarrow k(s) \\
x &\mapsto \Res(x y \omega).
\end{align*}
Since the homomorphism of $s$-group schemes
\begin{align}
\begin{split} \label{chap3kernelisom}
\tau: \Ga_{a,s}&\rightarrow \widetilde{\mathcal{V}} \\
t &\mapsto (t,\pi^{2\nu'}t),
\end{split}
\end{align}	
realizes an isomorphism onto the kernel of $r$, the pullback $\Lc_{\psi} \lbrace (\alpha^2 + v^*(\pi^{2 \nu'}) ) t \rbrace$ of $r^{-1}\chi_{j^{-1} \F}$ by $\tau$ is trivial, hence $\alpha^2 = v^*(\pi^{2 \nu'})$ by Proposition \ref{chap3charsheafas}.

Since the ramification of $\chi_{j^{-1} \F}$ is not bounded by the divisor $(\nu -1)s = \sw(\F_{\overline{\eta}}) s$, it follows from Theorem \ref{chap3lgcft2twisted} and from \ref{chap32.0.5} that the restriction of $r^{-1}\chi_{j^{-1} \F}$ to the sub-$s$-group scheme $\Ga_{a,s}$ of $\widetilde{\mathcal{V}}$ is non trivial, hence $\alpha$ is non zero. This implies that the scalar $ v^*(\pi^{2 \nu'}) = \alpha^2$ is non zero, and consequently that $y$ is a unit of $\Ow_T/\m^{\nu'+1}$. In particular, $y$ defines an $s$-point of $\Pic^{0}(T,(\nu'+1) s)_s$.

Let $\overline{x}$ be the spectrum of an algebraically closed extension of $k(s)$, and let $u$ be an $\overline{x}$-point of $\Pic^{0}(T,D)_{\overline{x}}$. The morphism
\begin{align*}
ur : \widetilde{\mathcal{V}}_{\overline{x}} &\rightarrow \Pic^{0}(T,D)_{\overline{x}} \\
(t,v) &\mapsto u(1 + v + \pi^{2 \nu'} t),
\end{align*}
surjects onto the fiber of $\sigma$ above $\sigma(u)$. Let $\G$ be the $\Lambda$-sheaf $\chi_{j^{-1}\F} \otimes \Lc_{\psi} \lbrace \Res_{\omega} \rbrace$ on the $s$-group scheme $\Pic^{0}(T,D)_s$. The pullback of $\G$ by $ur$ is isomorphic to $\Lc_{\xi} \lbrace (\alpha^2 + \beta) t + \Res_{(u-y) \omega}(v), \alpha \gamma(v) \rbrace$, up to twist by the stalk $\G_u $ of $\G$ at the geometric point $u$ of $\Pic^{0}(T,D)_s$, where $\beta$ is the image in $k(\overline{x})$ of $u$. The trace morphism
$$
R(ur)_! (ur)^{-1} \G [2](1) \rightarrow \G_{|\sigma^{-1}(\sigma(u))},
$$
is an isomorphism, since $ur$ is a $\Ga_{a}$-torsor over $\sigma^{-1}(\sigma(u))$. Together with Proposition \ref{chap3vanish18}, this implies that the complex
$$
R\Gamma_c( \sigma^{-1}(\sigma(u)), \G),
$$
vanishes unless $\beta = \alpha^2$ and the linear form $\Res_{(u-y) \omega}$ on $V$ is proportional to $\gamma$, namely unless $\sigma(u) = \sigma_2(y)$, in which case it is concentrated in degree $2( \nu' + 1) -1 = \nu$, and the cohomology group $H_c^{\nu}( \sigma^{-1}(\sigma(u)), \G)$ is of rank $1$. We obtain that $R \sigma_{!} \G$ is of the form $ \sigma_2(y)_* \Lc[-\nu]$, where $\Lc$ is a $\Lambda$-sheaf of rank $1$ on $s$, and thus that the complex
$$
R\Gamma_c(\Pic^{0}(T,D)_{\overline{s}}, \G) \cong R\Gamma_c(\Pic^{0}(T,\nu' s)_{\overline{s}}, R \sigma_{!} \G),
$$
is isomorphic to $\Lc_{\overline{s}}[-\nu]$, hence the conclusion of Proposition \ref{chap3rankcohom}.

This, together with the paragraphs \ref{chap32.14}, \ref{chap32.14.1}, \ref{chap32.14.2}, and \ref{chap32.14.3} concludes the proof of Proposition \ref{chap3rankcohom}.

\subsection{\label{chap32.8}} We now assume that $k$ is a finite field of cardinality $q$, that $\Lambda$ is $C$, and that $\mu = 1$ is the trivial $C$-admissible mutiplier on $G_k$. The $k$-automorphism $x \mapsto x^q$ is a topological generator of $G_k$, and we denote by $\Frob_k$ its inverse. Similarly, $\Frob_{s} = \Frob_k^f$ is a topological generator of the subgroup $G_s$ of $G_k$, where $f$ is the degree of the extension $k(s)/k$.

If $\F$ is $C$-local system of rank $1$ on $\eta$, and let $D$ be an effective Cartier divisor on $T$ such that $\F$ has ramification bounded by $D$ (cf. \ref{chap3ramdefiloc}), namely $\sw(\F_{\overline{\eta}})$ is strictly less than the multiplicity $\nu$ of $D$ at $s$. Then Theorem \ref{chap3lgcft2} produces a multiplicative $C$-local system $\chi_{\F}$ on the $s$-group scheme $\Pic(T,D)_s$. Moreover, the map
\begin{align*}
\chi_{\F} : k(\eta)^{\times} &\rightarrow C^{\times} \\
z &\mapsto \langle \chi_{\F} \rangle(z)(\Frob_s) 
\end{align*}
is a group homomorphism, cf. \ref{chap32.30}.

\begin{prop}\label{chap3compafini} Let $\F$ be a $C$-local system of rank $1$ on $\eta$, and let $c$ be an arbitrary element of $k(\eta)$ of valuation $a(T,j_* \F,\omega)$. Then we have
$$
(-1)^{a(T,j_* \F)} \varepsilon_{ \overline{s}}(T,j_* \F,\omega)(\Frob_s) = \int_{c^{-1} \Ow_T^{\times}} \chi_{\F}^{-1}(z) \psi(\Tr_{k/\mathbb{F}_p}\Res ( z \omega)) d z,
$$
if $\F$ is ramified, where $dz$ is the Haar measure on $k(\eta)$ normalized so that $\int_{\Ow_T} dz = 1$. If $\F$ is unramified, then we have $a(T,j_* \F) = 0$ and
$$
\varepsilon_{ \overline{s}}(T,j_* \F,\omega)(\Frob_s) = \chi_{\F}(c) q^{f v(\omega)}.
$$
\end{prop}

Let $\Psi_{\omega} : k(\eta) \rightarrow \Lambda^{\times}$ be the additive character given by $z \mapsto \psi(\Tr_{k/\mathbb{F}_p}(z \omega))$. Proposition \ref{chap3compafini} can then be summarized as an equality
$$
(-1)^{a(T,j_* \F)} \varepsilon_{ \overline{s}}(T,j_* \F,\omega) (\Frob_s) = \varepsilon(\chi_{\F}, \Psi_{\omega}),
$$
where $ \varepsilon(\chi_{\F}, \Psi_{\omega})$ is the automorphic $\varepsilon$-factor of the pair $(\chi_{\F}, \Psi_{\omega})$, cf (\cite{La87}, 3.1.3.2).

\subsection{\label{chap32.17}} We now prove Proposition \ref{chap3compafini}. Let us first assume that $\F$ is unramified, so that $\F$ is the pullback to $\eta$ of a $C$-local system $\G$ of rank $1$ on $s$. We can take $D= s$ and $\nu = 1$ above. For each integer $d$, the multiplicative $C$-local system $\chi_{\F}$ is given on the component $\Pic^d(T,s)_s$ by the pullback of $\G^{\otimes d}$. The $C$-admissible representations of rank $1$ corresponding to $\varepsilon_k(T,i_* i^{-1} j_* \F,\omega)$ and $\varepsilon_k(T,j_! j^{-1} j_* \F,\omega)$ are then respectively isomorphic to $\G_{\overline{s}}^{-1}$ and to
\begin{align*}
&H_c^{2\nu - a(T,j_! \F)} \left(\Pic^{a(T,j_!\F,\omega)}(T,s)_{\overline{s}}, \Lc_{\psi} \lbrace \Res_{\omega} \rbrace(\nu - a(T,\F,\omega)) \right) \otimes \G_{\overline{s}}^{\otimes a(T,j_!\F,\omega)}\\
 =& H_c^{1} \left(\Pic^{1 + v(\omega)}(T,s)_{\overline{s}}, \Lc_{\psi} \lbrace \Res_{\omega} \rbrace (-v(\omega)) \right) \otimes \G_{\overline{s}}^{\otimes(1 + v(\omega))}.
\end{align*}
Let $\pi$ be a uniformizer of $k(\eta)$, and let us write $\omega$ as $\alpha \frac{d \pi}{ \pi}$ for some element $\alpha$ of $k(\eta)^{\times}$ of valuation $1 + v(\omega)$. We have an isomorphism	
\begin{align*}
\theta : \mathbb{G}_{m,s} &\rightarrow \Pic^{1 + v(\omega)}(T,s)_s \\
t &\mapsto t \alpha^{-1},
\end{align*}
so that $\theta^{-1} \Lc_{\psi} \lbrace \Res_{\omega} \rbrace$ is isomorphic to $\Lc_{\psi} \lbrace t\rbrace $. By Proposition \ref{chap3rankcohom} and by the Grothendieck-Lefschetz trace formula (\cite{Gr66}, \'eq. (25)), we have
\begin{align*}
\Tr \left( \Frob_s \ | \ H_c^{1} \left(\Pic^{1 + v(\omega)}(T,s)_{\overline{s}}, \Lc_{\psi} \lbrace \Res_{\omega} \rbrace \right) \right) &= \Tr \left( \Frob_s \ | \ H_c^{1} \left(\mathbb{G}_{m,\overline{s}}, \Lc_{\psi} \lbrace t\rbrace\right) \right) \\
&= - \sum_{t \in k(s)^{\times}} \psi(\Tr_{k/\mathbb{F}_p}(t))\\
&= 1.
\end{align*}
We thus obtain that the quantity $ \varepsilon_{\overline{s}}(T,j_! \F,\omega) ( \Frob_s )$ is given by $\Tr \left( \Frob_s \ | \ \G_{\overline{s}}\right)^{1 + v(\omega)} q^{f v(\omega)}$. This implies that the value of $\varepsilon_{\overline{s}}(T,j_* \F,\omega)$ at $\Frob_k$ is given by $\Tr \left( \Frob_s \ | \ \G_{\overline{s}}\right)^{v(\omega)} q^{f v(\omega)}$, hence the result since we have
$$
\chi_{\F}(c) = \Tr \left( \Frob_s \ | \ \G_{\overline{s}}\right)^{v(c)} = \Tr \left( \Frob_s \ | \ \G_{\overline{s}}\right)^{v(\omega)}.
$$

\subsection{\label{chap32.18}} We now prove Proposition \ref{chap3compafini} in the case where $\F$ is ramified. In this situation, the $C$-sheaf $j_* \F$ is supported on $\eta$, and the quantity $(-1)^{a(T,j_* \F)} \varepsilon_{\overline{s}}(T,j_* \F,\omega) (\Frob_s )$ is thus equal to
$$
 (-1)^{a(T,j_! \F)} q^{f a(T,j_! \F, \omega) - f \nu} \Tr \left( \Frob_s \ | \ H_c^{2\nu - a(T,j_!\F)} \left(\Pic^{a(T,j_!\F,\omega)}(T,D)_{\overline{s}},\chi_{\F} \otimes \Lc_{\psi} \lbrace \Res_{\omega} \rbrace \right) \right).
$$ 
By Proposition \ref{chap3rankcohom} and by the same Grothendieck-Lefschetz trace formula (\cite{Gr66}, \'eq. (25)), the latter quantity coincides with
$$
q^{f a(T,j_! \F, \omega) - f \nu} \sum_{u \in c^{-1}(\Ow_T/\m^{\nu})^{\times}} \chi_{\F}^{-1}(u) \psi(\Tr_{k/\mathbb{F}_p}(u \omega)).
$$
The factor $q^{f a(T,j_! \F, \omega) - f \nu}$ is equal to $\int_{c^{-1}(u + \m^{\nu} \Ow_T)} dz$ for any element $u$ of $(\Ow_T/\m^{\nu})^{\times}$, so that we obtain
\begin{align*}
(-1)^{a(T,j_* \F)} \varepsilon_{\overline{s}}(T,j_* \F,\omega) (\Frob_s ) &= \sum_{u \in (\Ow_T/\m^{\nu})^{\times}} \int_{c^{-1}(u + \m^{\nu} \Ow_T)} \chi_{\F}^{-1}(z) \psi(\Tr_{k/\mathbb{F}_p}(z \omega)) d z \\
&= \int_{c^{-1} \Ow_T^{\times}} \chi_{\F}^{-1}(z) \psi(\Tr_{k/\mathbb{F}_p}(z \omega)) d z, 
\end{align*}
hence the result.

\section{The product formula for sheaves of generic rank at most 1 (after Deligne)\label{chap3productsection1}}


We review in this section Deligne's computation of the determinant of the cohomology of rank $1$ local systems on curves, as exposed in his $1974$ letter to Serre, which is published as an appendix in \cite{bloch}. The material of this section is thus entirely due to Deligne, besides the terminology regarding twisted sheaves.

\subsection{\label{chap33.0}} Let us recall that the base field $k$ is assumed throughout to be a perfect field of characteristic $p$. Let $\Lambda$ be an $\ell$-adic coefficient ring (cf. \ref{chap3conv}, \ref{chap30.0.0.1}) which is a field, and let $\psi : \mathbb{F}_{p} \rightarrow \Lambda^{\times}$ be a non trivial homomorphism. We fix a unitary $\Lambda$-admissible mutiplier $\mu$ on the topological group $G_k$ (cf. \ref{chap30.0}, \ref{chap3unitary}).

\begin{defi}\label{chap3globaleps} Let $X$ be a connected smooth curve over $k$, and let $\F$ be a $\mu$-twisted $\Lambda$-sheaf on $X$. The global $\varepsilon$-factor of the pair $(X,\F)$ is the $\Lambda$-admissible map $\varepsilon_{\bk}(X,\F) : G_k \rightarrow \Lambda^{\times}$ defined by
\begin{align*}
\varepsilon_{\bk}(X,\F)(g) = \det( g, R\Gamma_c(X_{\bk}, \F) )^{-1},
\end{align*}
for any $g$ in $G_k$, cf. \ref{chap31.3.4}.

\end{defi}

For any smooth connected projective curve over $k$, we denote by $|X|$ the set of closed points of $X$ and by $X_{(x)}$ the henselization of $X$ at a closed point $x$.

\begin{teo}\label{chap3productform1} Let $X$ be a smooth connected projective curve of genus $g$ over $k$, let $\omega$ be a non zero global meromorphic differential $1$-form on $X$ and let $\F$ be a $\mu$-twisted $\Lambda$-sheaf on $X$ of generic rank $\rk(\F)$ at most $1$. Then, for all but finitely many closed points $x$ of $X$ the $\varepsilon$-factor of the triple $(X_{(x)},\F_{|X_{(x)}},\omega_{|X_{(x)}})$ is identically equal to $1$, and we have
$$
\varepsilon_{\bk}(X, \F) = \chi_{\cyc}^{N(g-1)\rk(\F)} \prod_{x \in |X|} \delta_{x/k}^{a(X_{(x)},\F_{|X_{(x)}})} \Ver_{x/k} \varepsilon_{\overline{x}}(X_{(x)},\F_{|X_{(x)}},\omega_{|X_{(x)}}),
$$
 where $N$ is the number of connected components of $X_{\bk}$, where $\overline{x}$ is a $k$-morphism from $\Spec(\bk)$ to $x$, and where $\delta_{x/k}$ and $\Ver_{x/k}$ are defined in \ref{chap31.4.6}. 
\end{teo}

In Theorem \ref{chap3productform1}, the image by $d^1$ of the left hand side (cf. \ref{chap30.11}) is 
$$
d^1( \varepsilon_{\bk}(X, \F) ) = \mu^{- \chi(X_{\bk}, \F)},
$$
cf. \ref{chap31.3.4}, while the image by $d^1$ of the right hand side is $\prod_{x \in |X|} \mu^{\deg(x) a(X_{(x)},\F_{|X_{(x)}},\omega_{|X_{(x)}}) }$ (cf. \ref{chap32.11}, \ref{chap32.15}), where the conductor $ a(X_{(x)},\F_{|X_{(x)}},\omega_{|X_{(x)}}) 	$ vanishes for all but finitely many closed points $x$ of $X$. Thus the conclusion of Theorem \ref{chap3productform1} is consistent with the identity
$$
- \chi(X_{\bk}, \F) = \sum_{x \in |X|} \deg(x) a(X_{(x)},\F_{|X_{(x)}},\omega_{|X_{(x)}}),
$$
which results from the Grothendieck-Ogg-Shafarevich formula.

\subsection{\label{chap33.5}} We now describe Deligne's proof of Theorem \ref{chap3productform1}. Let $X,g,N,\omega, \F$ be as in \ref{chap3productform1}. By replacing $k$ with a finite extension if necessary, we can assume (and we do) that $X$ is geometrically connected over $k$, so that $N=1$. Let $j : U \rightarrow X$ be a non empty open subscheme such that $j^{-1} \F$ is a $\mu$-twisted $\Lambda$-local system on $U$. We have an exact sequence
$$
0 \rightarrow j_! j^{-1} \F \rightarrow \F \rightarrow \bigoplus_{x \in X \setminus U}i_{x*} i_x^{-1} \F \rightarrow 0,
$$
where $i_x : x \rightarrow X$ is the inclusion of a closed point of $X$. The product formula \ref{chap3productform1} holds for $i_{x*} i_x^{-1} \F$ for each $x$ in $|X|$, hence we can assume (and we now do) that $\F$ vanishes outside $U$, i.e. $\F = j_! j^{-1} \F $. By replacing $U$ with a smaller non empty open subscheme of $X$ if necessary, we can further assume (and we do as well) that the complement $X \setminus U$ contains at least two closed points and that $j^{-1} \F$ is of rank $\rk(\F) = 1$.

Let us consider the effective Cartier divisor
$$
D = \sum_{x \in X \setminus U} (1 + \sw_x(\F)) x,
$$
where $\sw_x(\F)$ is the Swan conductor of $\F$ at $x$. Equivalently, we have
$$
D = \sum_{x \in |X|} a(X_{(x)},\F_{|X_{(x)}}) x,
$$
cf. \ref{chap32.11}. The Grothendieck-Ogg-Shafarevich formula then implies that the Euler characteristic of $\F$, namely
$$
\chi_c(U,j^{-1}\F) = \sum_{i \in \Z} (-1)^i \dim H_c^i(U,j^{-1}\F),
$$
is equal to $-d$, where $d = \deg(D) - 2 + 2g$ is a nonnegative integer. Since the canonical line bundle $\omega_X$ of $X$ has degree $2g-2$, the integer $d$ is also the degree of the line bundle $\omega_X(D)$.

The $\mu$-twisted $\Lambda$-local system $j^{-1} \F$ of rank $1$ on $U$ has ramification bounded by $D$, cf. \ref{chap3ramdefi3}, and consequently there exists a $\mu$-twisted multiplicative $\Lambda$-local system $\chi_{\F}$ (cf. \ref{chap3twistedmult}) on $\Pic_k(X,D)$ (cf. \ref{chap32.1.0.0}) whose pullback by the Abel-Jacobi morphism
\begin{align*}
\Phi : U \rightarrow \Pic_k(X,D),
\end{align*}
cf. \ref{chap3abeljacob}, is isomorphic to $j^{-1} \F$.

Let $\Sym_k^d(U)$ be the quotient of $U^d$ by the group of bijection of $\{1,\dots,d \}$ onto itself, acting by permuting the $d$ factors of $U^d$, cf. (\cite{G18}, Prop. 2.27). The $\Lambda$-local system $p_1^{-1} j^{-1} \F \otimes \cdots \otimes p_d^{-1} j^{-1} \F$ on $U^d$, where $(p_i)_{i=1}^d$ are the projections on each factor, descends to a $\Lambda$-local system $\F^{[d]}$ of rank $1$ on $\Sym_k^d(U)$, cf. (\cite{G18}, Prop. 2.32). The symmetric K\"unneth formula (\cite{SGA4}, XVII 5.5.21) implies that we have a natural isomorphism
\begin{align}\label{chap3isopre1}
R\Gamma_c(\Sym_k^d(U)_{\bk}, \F^{[d]})[d] \cong L \Gamma^d \left( R\Gamma_c(U_{\bk}, j^{-1}\F) \right)[d],
\end{align}
where the functor $L \Gamma^d$ is defined in (\cite{Illusie}, I.4.2.2.6). By (\cite{Illusie}, I.4.3.2.1), we have Quillen's shift formula
\begin{align}\label{chap3isopre2}
L \Gamma^d \left( R\Gamma_c(U_{\bk}, j^{-1}\F) \right)[d] \cong L \Lambda^d\left( R\Gamma_c(U_{\bk}, j^{-1}\F)[1] \right),
\end{align}
where the derived $d$-th exterior power $L \Lambda^d$ is also defined in (\cite{Illusie}, I.4.2.2.6). The shifted complex $R\Gamma_c(U_{\bk}, j^{-1}\F)[1]$ is of rank $d$, and the isomorphism
$$
\det( \Lambda^{\rk(V)} V ) \cong \det(V),
$$
valid for any finite dimensional $\Lambda$-vector space $V$, extends to an isomorphism
\begin{align}\label{chap3isopre3}
\det \left( L \Lambda^d\left( R\Gamma_c(U_{\bk}, j^{-1}\F)[1] \right) \right) \cong \det( R\Gamma_c(U_{\bk}, j^{-1}\F)[1] ).
\end{align}
By combining \ref{chap3isopre1}, \ref{chap3isopre2} and \ref{chap3isopre3}, we obtain a natural isomorphism
\begin{align}\label{chap3isomorphismproof0}
\det R\Gamma_c(U_{\bk}, j^{-1}\F)^{-1} \cong \det R\Gamma_c(\Sym_k^d(U)_{\bk}, \F^{[d]})^{(-1)^d}.
\end{align}
Moreover, if we denote by 
\begin{align}\label{chap3abeljacobi2}
\Phi_d : \Sym_k^d(U) \rightarrow \Pic_k(X,D),
\end{align}
the $d$-th Abel-Jacobi map, whose composition with the canonical projection $U^d \rightarrow \Sym_k^d(U)$ sends a section $(x_i)_{i=1}^d$ of $U^d$ to $\prod_{i=1}^d \Phi(x_i)$, then the multiplicativity of $\chi_{\F}$ implies that the pullback of $\Phi_d^{-1} \chi_{\F}$ to $U^d$ is isomorphic
$$
\left( (\Phi p_1)(\Phi p_2) \dots (\Phi p_d) \right)^{-1} \chi_{\F} \cong p_1^{-1} \Phi^{-1} \chi_{\F}\otimes \cdots \otimes p_d^{-1} \Phi^{-1} \chi_{\F} \cong p_1^{-1} j^{-1} \F \otimes \cdots \otimes p_d^{-1} j^{-1} \F,
$$
and thus the pullback $\Phi_d^{-1} \chi_{\F}$ is isomorphic to $\F^{[d]}$.

The Leray spectral sequence for $(\Phi_d,\F^{[d]})$ and the projection formula then yield
\begin{align}
\begin{split}
\label{chap3isomorphismproof}
\det R\Gamma_c(\Sym_k^d(U)_{\bk}, \F^{[d]}) &\cong \otimes_{q \in \Z} \det R\Gamma_c(\Pic_k(X,D)_{\bk}, R^q \Phi_{d!} \Phi_d^{-1} \chi_{\F})^{(-1)^q} \\
 &\cong \otimes_{q \in \Z} \det R\Gamma_c(\Pic_k(X,D)_{\bk}, \chi_{\F} \otimes R^q \Phi_{d!} \Lambda)^{(-1)^q}.
\end{split}
\end{align}

\subsection{\label{chap33.2}} Let $i :D \rightarrow X$ be the closed immersion of $D$ into $X$, and let $J$ (resp. $J^0$) be the functor which associates to a $k$-scheme $S$ the set of isomorphisms $\alpha : \Ow_{D_S} \rightarrow i_{S}^{*} \omega_X(D)$ of $\Ow_{D_S}$-modules (resp. of automorphisms of $\Ow_{D_S}$ as a module over itself). Then $J^0$ is representable by a smooth connected affine group scheme of dimension $\deg(D)$ over $k$, and $J$ is a $J^0$-torsor. In particular, since the action of $\mathbb{G}_{m,k}$ by multiplication on $\Ow_D$ turns $\mathbb{G}_{m,k}$ into a sub-$k$-group scheme of $J^0$, the $J^0$-torsor $J$ is naturally endowed with an action of $\mathbb{G}_{m,k}$ by left multiplication. Moreover, the morphism
\begin{align*}
f : J' &\rightarrow \Pic_k(X,D)\\
\alpha &\rightarrow (\omega_X(D),\alpha),
\end{align*}
where $J' = J/ \mathbb{G}_{m,k}$, is a closed immersion, its image being the fiber of the canonical projection $\Pic_k(X,D) \rightarrow \Pic_k(X)$ at $\omega_X( D)$. 

Let us denote by $\Res_D : i^* \omega_X(D) \rightarrow k$ the residue homomorphism, given by 
$$
\Res_D(\alpha) = \sum_{x \in |D|} \mathrm{Tr}_{k(x)/k} \Res_x(\alpha),
$$
where $\Res_x$ denotes the residue homomorphism at a closed point $x$ (cf. \ref{chap3residue}). We denote by $g : \Sigma \rightarrow J$ the closed subscheme consisting of isomorphisms $\alpha$ such that $\Res_D(\alpha) = 0$, and by $g' : \Sigma' \rightarrow J'$ its quotient by $\mathbb{G}_{m,k}$, which is a closed immersion as well.

\begin{lem}[\cite{bloch}, p.82]\label{chap3delignecalcul1} The $\Lambda$-sheaves $R^q \Phi_{d!} \Lambda$ (cf. \ref{chap3abeljacobi2}) on $\Pic^d_k(X,D)$ admit the following description:
\begin{enumerate}
\item for $q = 2g -2$, there exists a short exact sequence
$$
0 \rightarrow R^{2g-2} \Phi_{d!} \Lambda \rightarrow \Lambda(1-g) \rightarrow f_* \Lambda(1-g) \rightarrow 0,
$$
of $\Lambda$-sheaves on $\Pic^d_k(X,D)$.
\item for $q= 2g$, the $\Lambda$-sheaf $R^{2g} \Phi_{d!} \Lambda$ is isomorphic to $(fg')_* \Lambda(-g)$,
\item the $\Lambda$-sheaf $R^q \Phi_{d!} \Lambda$ vanishes if $q$ is not equal to $2g$ or $2g-2$.
\end{enumerate}
\end{lem}

Let $\overline{t}$ be the spectrum of an algebraically closed extension of $k$, and let $(\Lc,\alpha)$ be a $\overline{t}$-point of $\Pic^d_k(X,D)$ (cf. \ref{chap32.1.0.0}). The fiber of $\Phi_d$ above $(\Lc,\alpha)$ parametrizes sections $\sigma$ in $H^0(X_{\overline{t}}, \Lc)$ whose image in $H^0(D_{\overline{t}}, \Lc)$ is $\alpha$, cf. (\cite{G18}, Prop. 4.12). The degree of the line bundle $\Lc(-D)$ is $2g-2$, and the Riemann-Roch theorem yields:
\begin{enumerate}
\item if $\Lc$ is not isomorphic to $\omega_X(D)$, then the fiber of $\Phi_d$ above $(\Lc,\alpha)$ is a torsor under the additive $\overline{t}$-group scheme associated to the $(g-1)$-dimensional $k(\overline{t})$-vector space $H^0(X_{\overline{t}}, \Lc(-D))$.
\item if $\Lc = \omega_X(D)$ and if $\Res_D(\alpha)=0$ then the fiber of $\Phi_d$ above $(\Lc,\alpha)$ is a torsor under the additive $\overline{t}$-group scheme associated to the $g$-dimensional $k(\overline{t})$-vector space $H^0(X_{\overline{t}}, \Lc(-D))$.
\item if $\Lc = \omega_X(D)$ and if $\Res_D(\alpha) \neq 0$ then the fiber of $\Phi_d$ above $(\Lc,\alpha)$ is empty.
\end{enumerate} 

Let $w : W \rightarrow \Pic^d_k(X,D)$ be the open complement of the image of the closed immersion $f$. Then we have a distinguished triangle
$$
w_! w^{-1} R \Phi_{d!} \Lambda \rightarrow R \Phi_{d!} \Lambda \rightarrow f_* f^{-1} R \Phi_{d!} \Lambda \xrightarrow[]{[1]}.
$$
Above $W$, the morphism $\Phi_d$ is a fibration in affine spaces, of relative dimension $g-1$, hence $w^{-1} R \Phi_{d!} \Lambda$ is quasi-isomorphic to $\Lambda(1-g)[2-2g]$. Moreover, the description above of the fibers of $\Phi_d$ imply that $f^{-1} R \Phi_{d!} \Lambda$ is supported on $\Sigma'$. Above $\Sigma'$, the morphism $\Phi_d$ is a fibration in affine spaces, of relative dimension $g$, hence $f^{-1} R \Phi_{d!} \Lambda$ is quasi-isomorphic to $g'_*\Lambda(-g)[-2g]$. Thus $R^q \Phi_{d!} \Lambda$ vanishes if $q$ is not equal to $2g$ or $2g-2$, and we have isomorphisms
\begin{align*}
R^{2g-2} \Phi_{d!} \Lambda &\cong w_!\Lambda(1-g), \\
R^{2g} \Phi_{d!} \Lambda &\cong (fg')_*\Lambda(-g),
\end{align*}
hence the conclusion of Lemma \ref{chap3delignecalcul1}.

\subsection{\label{chap33.3}} By combining the formula \ref{chap3isomorphismproof} with Lemma \ref{chap3delignecalcul1}, we obtain that the determinant of the complex 
$$
R\Gamma_c(\Sym_k^d(U)_{\bk}, \F^{[d]}),
$$
is isomorphic to
\begin{align}\label{chap3isomorphismproof2}
 \frac{ \det R\Gamma_c(\Pic^d_k(X,D)_{\bk}, \chi_{\F}(1-g)) \det R\Gamma_c(\Sigma'_{\bk}, (fg')^{-1}\chi_{\F}(-g))}{\det R\Gamma_c(J'_{\bk}, f^{-1}\chi_{\F}(1-g)) }.
\end{align}

\begin{lem}[\cite{bloch}, p.82]\label{chap3delignecalcul2} The factor $\det R\Gamma_c(\Pic^d_k(X,D)_{\bk}, \chi_{\F}(1-g)) $ is isomorphic to $\Lambda$, as a $\Lambda$-admissible representation of $G_k$ of rank $1$.
\end{lem}

When $\F$ is everywhere tamely ramified, then $\Pic^d_k(X,D)_{\bk}$ is a torsor under the $\bk$-group scheme $\Pic^0_k(X,D)_{\bk}$, which is an extension of an abelian scheme of dimension $g$ by a torus of dimension $\deg(D)$. The result follows in this case from the fact that the determinant of the cohomology of a tame $\Lambda$-local system on an extension of an abelian $k$-scheme by a torus of dimension at least $2$ is canonically trivial. We refer to (\cite{bloch}, Constr. 1 p.70) for a proof of the latter result.

 When $\F$ is not everywhere tamely ramified, there exists a closed point $x$ on $D$ such that the multiplicity $\nu$ of $D$ at $x$ is at least $2$. Since $\Pic^d_k(X,D)$ is a torsor over the $k$-group scheme $\Pic^0_k(X,D)$ and since $\chi_{\F}$ is multiplicative, it is sufficient to prove that the complex
 $$
 R\Gamma_c(\Pic^0_k(X,D)_{\bk}, \chi_{\F}),
$$
is quasi-isomorphic to $0$. The projection
$$
\tau : \Pic^0_k(X,D) \rightarrow \Pic^0_k(X,D-x),
$$
is a homomorphism of $k$-group schemes, whose kernel is the additive $k$-group scheme associated to the $k$-vector space $\m_x^{\nu-1}/\m_x^{\nu}$, where $\m_x$ is the maximal ideal of $\Ow_{X,x}$. Since the ramification of $\F$ is not bounded by the divisor $D-x$, Theorem \ref{chap3ggcfttwisted} and \ref{chap32.0.5} imply that the restriction of $\chi_{\F}$ to the kernel of $\tau$ is non trivial. Together with Proposition \ref{chap3vanish13} and with the multiplicativity of $\F$, this implies the vanishing of $R\tau_! \chi_{\F}$, hence the result.

\subsection{\label{chap33.4}} By combining the formula \ref{chap3isomorphismproof2} with Lemma \ref{chap3delignecalcul2}, we obtain an isomorphism
\begin{align}\label{chap3isomorphismproof3}
\det R\Gamma_c(\Sym_k^d(U)_{\bk}, \F^{[d]}) \cong \det R\Gamma_c(\Sigma'_{\bk}, (fg')^{-1}\chi_{\F}(-g)) \det R\Gamma_c(J'_{\bk}, f^{-1}\chi_{\F}(1-g))^{-1}.
\end{align}
Let $\tau : J \rightarrow J'$ be the natural projection (cf. \ref{chap33.2}). 

\begin{lem}[\cite{bloch}, p.84]\label{chap3delignecalcul3} The $\Lambda$-sheaves $R^q \tau_{!} \Lc_{\psi}\{ \Res \}$ admits the following description:
\begin{enumerate}
\item the $\Lambda$-sheaf $R^1 \tau_{!} \Lc_{\psi}\{ \Res \}$ is isomorphic to the constant sheaf $\Lambda$ on $J'$,
\item the $\Lambda$-sheaf $R^2 \tau_{!} \Lc_{\psi}\{ \Res \}$ is isomorphic to $g'_* \Lambda(-1)$, with $g'$ as in \ref{chap33.2}.
\item the $\Lambda$-sheaf $R^q \tau_{!} \Lc_{\psi}\{ \Res \}$ vanishes when $q$ is not equal to $1$ or $2$.
\end{enumerate}
\end{lem}

Recall that $\tau$ is a $\mathbb{G}_m$-torsor. Let $\overline{J}$ be the quotient of $\Ga_{a,k} \times_k J$ by the action of $\mathbb{G}_{m,k}$ given by $t\cdot (y, \alpha) = (t^{-1} y, t \alpha)$, and let $\overline{\Res} : \overline{J} \rightarrow \Ga_{a,k}$ be the morphism which sends the class in $\overline{J}$ of a section $(y,\alpha)$ of $J$ to $y \Res(\alpha)$. Let $u : J \rightarrow \overline{J}$ be the open immersion which sends a section $\alpha$ to the class of $(1,\alpha)$, and let $\overline{\tau} : \overline{J} \rightarrow J'$ be the morphism which sends the class of a section $(y,\alpha)$ to $\tau(\alpha)$. Let $i : J' \rightarrow \overline{J}$ be the section of $\overline{\tau}$ which sends $\tau(\alpha)$ to the class of $(0,\alpha)$ in $\overline{J}$, for any section $\alpha$ of $J$. We then have an exact sequence
$$
0 \rightarrow u_! \Lc_{\psi}\{ \Res \} \rightarrow \Lc_{\psi}\{ \overline{\Res} \} \rightarrow i_* \Lambda \rightarrow 0.
$$
By applying the funtor $R\overline{\tau}_!$, we obtain a distinguished triangle
\begin{align}\label{chap3isomproof21}
R \tau_{!} \Lc_{\psi}\{ \Res \} \rightarrow R \overline{\tau}_{!} \Lc_{\psi}\{ \overline{\Res} \} \rightarrow \Lambda[0] \xrightarrow[]{[1]}.
\end{align}
Moreover, the fiber of $\overline{\tau}$ over a geometric point $\overline{\alpha}$ of $J'$ is isomorphic to $\Ga_a$ and the restriction of $\Lc_{\psi}\{ \overline{\Res} \}$ to this fiber is a multiplicative $\Lambda$-local system, which is trivial if and only if $\overline{\alpha}$ factors through $g'$ (cf. \ref{chap33.2}). Thus Proposition \ref{chap3vanish13} implies that $R \overline{\tau}_{!} \Lc_{\psi}\{ \overline{\Res} \}$ is supported on the image $\Sigma'$ of $g'$, and consequently
$$
R \overline{\tau}_{!} \Lc_{\psi}\{ \overline{\Res} \} \cong g'_* g^{'-1} R \overline{\tau}_{!} \Lc_{\psi}\{ \overline{\Res} \} \cong g'_* R (\overline{\tau}_{| \overline{\tau}^{-1}(\Sigma')})_! \Lambda.
$$
Above $\Sigma'$, the morphism $\overline{\tau}$ is a fibration in affine spaces, of relative dimension $1$, hence $R \overline{\tau}_{!} \Lc_{\psi}\{ \overline{\Res} \}$ is quasi-isomorphic to $g'_* \Lambda(-1)[-2]$. The conclusion of Lemma \ref{chap3delignecalcul3} then follows from \ref{chap3isomproof21}.

\subsection{\label{chap33.6}} By combining the formula \ref{chap3isomorphismproof3} with Lemma \ref{chap3delignecalcul3}, we obtain an isomorphism
\begin{align}
\begin{split}
\label{chap3isomorphismproof4}
\det R\Gamma_c(\Sym_k^d(U)_{\bk}, \F^{[d]}) &\cong \det R\Gamma_c(J'_{\bk}, f^{-1}\chi_{\F}(1-g) \otimes R\tau_{!} \Lc_{\psi}\{ \Res \}) \\
&\cong \det R\Gamma_c(J_{\bk}, (f \tau)^{-1}\chi_{\F}(1-g) \otimes \Lc_{\psi}\{ \Res \})
\end{split}
\end{align}

By Proposition \ref{chap3changeforme} and Proposition \ref{chap3transferhomo}, the product formula \ref{chap3productform1} for some non zero meromorphic $1$-form $\omega$ on $X$ implies the product formula for all such $1$-forms. In particular, we can assume (and we do) that $\omega$ is a global section of $\omega_X(D)$ such that $i^* \omega : \Ow_D \rightarrow i^* \omega_X(D)$ is an isomorphism. For any closed point $x$ of $X$, we denote by $i_x : D_x \rightarrow X_{(x)}$ the restriction of $D$ to the henselization $X_{(x)}$ of $X$ at $x$, by $\F_x$ and $\omega_x$ the restrictions of $\F$ and $\omega$ to $X_{(x)}$, and by $\Sigma_x$ the set of $k$-linear embeddings of $k(x)$ into $\bk$.

The sections of $J$ over a $\bk$-scheme $S$ consist of all isomorphisms $ \alpha : \Ow_{D_S} \rightarrow i_{S}^{*} \omega_X(D)$ of $\Ow_{D_S}$-modules. We have a decomposition
$$
\Ow_{D_S} \cong (\Ow_D \otimes_k \bk ) \otimes_{\bk} \Ow_S \cong \prod_{x \in D} \prod_{\iota \in \Sigma_x} (\Ow_{D_x} \otimes_{k(x),\iota} \Ow_S),
$$
hence $\alpha$ can be identified with a tuple $(\alpha_{x,\iota})_{x \in D, \iota \in \Sigma_x}$, where each $\alpha_{x,\iota}$ is a trivialization of the $(\Ow_{D_x} \otimes_{k(x),\iota} \Ow_S)$-module $i_x^* \omega_X(D_x) \otimes_{k(x),\iota} \Ow_S$. In particular, the morphism
\begin{align*}
\delta : \prod_{x \in D} \prod_{\iota \in \Sigma_x} \Pic^0(X_{(x)},D_x)_{\iota, \bk} &\rightarrow J_{\bk} \\
(u_x)_{x \in D, \iota \in \Sigma_x} &\rightarrow (u_x \omega_x)_{x \in D, \iota \in \Sigma_x},
\end{align*}
is an isomorphism of $\bk$-schemes, which fits into a commutative diagram
\begin{center}
 \begin{tikzpicture}[scale=1]

\node (A) at (0,2) {$\prod_{x \in D} \prod_{\iota \in \Sigma_x} \Pic^0(X_{(x)},D_x)_{\iota, \bk}$};
\node (B) at (4,2) {$J_{\bk}$};
\node (E) at (4,1) {$J_{\bk}'$};
\node (C) at (4,0) {$\Pic^d_k(X,D)_{\bk},$};
\node (D) at (0,0) {$\Pic^0_k(X,D)_{\bk}$};

\path[->,font=\scriptsize]
(A) edge node[above]{$\delta$} (B)
(B) edge node[right]{$\tau$} (E)
(E) edge node[right]{$f$} (C)
(A) edge  (D)
(D) edge  (C);
\end{tikzpicture} 
\end{center}
where the bottom horizontal arrow is the translation by the $\bk$-point of $\Pic^d_k(X,D)$ which is the image by $\Phi_d$ of the $\bk$-point of $\Sym^d(U)$ corresponding to $\sum_{x \notin D} \sum_{\iota \in \Sigma_x} v_x(\omega) \iota(x)$. If $(p_{x,\iota})_{x \in D, \iota \in \Sigma_x}$ are the natural projections from the source of $\delta$ onto each factor, then we have decompositions
\begin{align*}
\delta^{-1} \Lc_{\psi}\{ \Res \} &\cong \bigotimes_{x \in D} \bigotimes_{\iota \in \Sigma_x} p_{x,\iota}^{-1}\Lc_{\psi}\{ \Res_{\omega_x} \},\\
(f \tau \delta)^{-1}\chi_{\F} &\cong \bigotimes_{x \in D} \bigotimes_{\iota \in \Sigma_x} p_{x,\iota}^{-1}\chi_{\F_x} \otimes \bigotimes_{x \notin D} \bigotimes_{\iota \in \Sigma_x} \F_{x,\iota}^{\otimes v(\omega_x)},
\end{align*}
by the compatibility of local and global geometric class field theory, cf. \ref{chap32.6}. By Proposition \ref{chap3rankcohom}, each complex 
$$
R\Gamma_c(\Pic^0(X_{(x)},D_x)_{\iota, \bk},\chi_{\F_x}\otimes \Lc_{\psi}\{ \Res_{\omega_x} \})
$$
is a one-dimensional $\Lambda$-vector space concentrated in degree $a_x = a(X_{(x)},\F_{|X_{(x)}})$. By K\"unneth's formula (\cite{SGA4}, Thm. 5.4.3), and by Proposition \ref{chap3gradedvector}, we obtain that the complex 
$$
R\Gamma_c(J_{\bk}, (f \tau)^{-1}\chi_{\F} \otimes \Lc_{\psi}\{ \Res \}) 
$$
is concentrated in degree $\sum_{x \in D} a_x = \deg(D)$, and that $G_k$ acts on its cohomology group of degree $\deg(D)$ through the $\Lambda$-admissible map given by 
\begin{align*}
\prod_{x \in D} \delta_{x/k}^{a_x} \Ver_{x/k} H_c^{a_x}(\Pic^0(X_{(x)},D_x)_{\overline{x}}, \chi_{\F_x}\otimes \Lc_{\psi}\{ \Res_{\omega_x} \}) \prod_{x \notin D} \Ver_{x/k} \F_{\overline{x}}^{\otimes v(\omega_x)},
\end{align*}
with notation as in \ref{chap31.4.6}, and the latter is equal to
$$
 \prod_{x \in |X|} \delta_{x/k}^{a_x} \Ver_{x/k} \left( \chi_{\cyc}^{v(\omega_x)} \varepsilon_{\overline{x}}(X_{(x)},\F_{|X_{(x)}},\omega_{|X_{(x)}}) \right).
$$
Since $d$ has the same parity as $\deg(D)$, the latter map is equal by (\ref{chap3isomorphismproof4}) to the trace function on $G_k$ of the $\Lambda$-admissible representation $\det R\Gamma_c(\Sym_k^d(U)_{\bk}, \F^{[d]})^{(-1)^d}(g-1)$ of $(G_k,\mu^d)$, and thus to the trace function of $\det R\Gamma_c(U_{\bk},j^{-1}\F)^{-1}(g-1) $ as well by (\ref{chap3isomorphismproof0}), hence the conclusion of Theorem \ref{chap3productform1}, since we have
$$
 \prod_{x \in |X|} \Ver_{x/k} \left( \chi_{\cyc}^{v(\omega_x)} \right) =\prod_{x \in |X|} \chi_{\cyc}^{[k(x):k]v(\omega_x)} = \chi_{\cyc}^{2g-2}.
$$


\section{Geometric local \texorpdfstring{$\varepsilon$}{epsilon}-factors in arbitrary rank \label{chap3gfar}}

Let $\Lambda$ be an $\ell$-adic coefficient ring (cf. \ref{chap3conv}, \ref{chap30.0.0.1}) which is a field. Let $\psi : \mathbb{F}_{p} \rightarrow \Lambda^{\times}$ be a non trivial homomorphism, hence producing a multiplicative $\Lambda$-local system $\Lc_{\psi}$ on $\Ga_{a,k}$, cf. \ref{chap34.3}. 

As in \ref{chap3lgf1}, let $T$ be the spectrum of a $k$-algebra, which is a henselian discrete valuation ring $\Ow_T$, with maximal ideal $\m$, and whose residue field $\Ow_T/ \m$ is a finite extension of $k$ of degree $\deg(s)$. Let $j : \eta \rightarrow T$ be the generic point of $T$, and let $i : s \rightarrow T$ be its closed point, so that $T$ is canonically an $s$-scheme, as in \ref{chap32.3}. We fix a $\bk$-point $\overline{s} : \Spec(\bk) \rightarrow T$ of $T$ above $s$, so that the Galois group $G_s = \Gal(\bk / k(s))$ can be considered as a subgroup of $G_k$. We fix a unitary $\Lambda$-admissible mutiplier $\mu$ on the topological group $G_s$ (cf. \ref{chap30.0}, \ref{chap3unitary}).


\subsection{\label{chap35.0}} Let $\pi$ be a uniformizer of $\Ow_T$. We abusively denote by $\pi$ as well the morphism
$$
\pi : T \rightarrow \mathbb{A}^1_s,
$$
corresponding to the unique morphism $k(s)[t] \rightarrow \Ow_T$ of $k(s)$-algebras which sends $t$ to $\pi$. By Theorem \ref{chap3GK3}, the pullback functor $\pi^{-1}$ realizes an equivalence from the category of special $\mu$-twisted $\Lambda$-sheaves on $\mathbb{A}^1_s$ to the category of $\Lambda$-sheaves on $T$. We denote by $\pi_{\diamondsuit}$ a quasi-inverse to this equivalence.

\begin{defi}\label{chap3localepsfact} Let $\F$ be a $\mu$-twisted $\Lambda$-sheaf on $T$, and let $\omega$ be an element of $\Omega^{1,\times}_{\eta}$ (cf. \ref{chap32.12}). Then the \textit{$\varepsilon$-factor} of the triple $(T,\F,\omega)$ is the $\Lambda$-admissible map $\varepsilon_{\pi,\overline{s}}(T,\F,\omega)$ from $G_s$ to $\Lambda^{\times}$ defined by
\begin{align*}
G_s &\rightarrow \Lambda^{\times} \\
g &\mapsto \langle \chi_{\det(j^{-1}\F)}\rangle(\frac{\omega}{d\pi})(g) \chi_{\cyc}(g)^{- v(\omega) \rk(j^{-1} \F)} \det \left( g \ | \ R\Gamma_c\left( \mathbb{A}^1_{\overline{s}}, \pi_{\diamondsuit }\F \otimes \Lc_{\psi}^{-1} \right) \right)^{-1},
\end{align*}
cf. \ref{chap32.30} and \ref{chap3changeforme}, or equivalently by 
$$
\varepsilon_{\pi,\overline{s}}(T,\F,\omega) = \langle \chi_{\det(j^{-1}\F)}\rangle(\frac{\omega}{d\pi}) \chi_{\cyc}^{- v(\omega) \rk(j^{-1} \F)} \varepsilon_{\overline{s}}(\mathbb{A}^1_s, \pi_{\diamondsuit }\F \otimes \Lc_{\psi}^{-1}),$$
cf. \ref{chap3globaleps}.
\end{defi}

It is clear from Definition \ref{chap3localepsfact} that $\varepsilon$-factors are multiplicative in short exact sequences:

\begin{prop}\label{chap3multiplicativity} For any exact sequence
$$
0 \rightarrow \F' \rightarrow \F \rightarrow \F'' \rightarrow 0,
$$
of $\Lambda$-sheaves on $T$, and for any element $\omega$ of $\Omega^{1,\times}_{\eta}$, we have 
$$
\varepsilon_{\pi,\overline{s}}(T,\F,\omega) = \varepsilon_{\pi,\overline{s}}(T,\F',\omega) \varepsilon_{\pi,\overline{s}}(T,\F'',\omega).
$$
\end{prop}

\begin{prop}\label{chap3sanity5}Let $\F$ be a $\mu$-twisted $\Lambda$-sheaf on $T$, and let $\omega$ be an element of $\Omega^{1,\times}_{\eta}$. The coboundary of $\varepsilon_{\pi,\overline{s}}(T,\F,\omega)$ (cf. \ref{chap30.11}) is then given by
$$
d^1 \left(\varepsilon_{\pi,\overline{s}}(T,\F,\omega) \right) = \mu^{ a(T,\F,\omega)},
$$
cf. \ref{chap32.11} for the notation.
\end{prop}

Indeed, $\det(j^{-1}\F)$ is a $\mu^{\rk(j^{-1}\F)}$-twisted $\Lambda$-sheaf of rank $1$ on $\eta$, hence
$$
d^1 \left( \langle \chi_{\det(j^{-1}\F)}\rangle(\frac{\omega}{d\pi})\right) = \mu^{v(\omega)\rk(j^{-1}\F)}.
$$
Moreover, $\pi_{\diamondsuit }\F \otimes \Lc_{\psi}^{-1}$ is a $\mu$-twisted $\Lambda$-sheaf on $\mathbb{A}^1_s$, and consequently we have
$$
d^1 \left( \varepsilon_{\overline{s}}(\mathbb{A}^1_s, \pi_{\diamondsuit }\F \otimes \Lc_{\psi}^{-1}) \right) = \mu^{-\chi_c(\mathbb{A}^1_{\overline{s}}, \pi_{\diamondsuit }\F \otimes \Lc_{\psi}^{-1})}.
$$
The Swan conductor of $\pi_{\diamondsuit }\F \otimes \Lc_{\psi}^{-1}$ at infinity is equal to $\rk(j^{-1} \F)$, hence the conductor of the pair $(\pi_{\diamondsuit }\F \otimes \Lc_{\psi}^{-1}, dt)$ at infinity is equal to $\rk(j^{-1} \F)(2 + v_{\infty}(dt)) = 0$, so that the Grothendieck-Ogg-Shafarevich formula yields
$$
-\chi_c(\mathbb{A}^1_{\overline{s}}, \pi_{\diamondsuit }\F \otimes \Lc_{\psi}^{-1}) = a(T,\F,d\pi),
$$
hence the conclusion of Proposition \ref{chap3sanity5}.

\begin{prop}\label{chap3unramtwist} Let $\F$ be a $\mu$-twisted $\Lambda$-sheaf on $T$, and let $\mathcal{G}$ be a $\nu$-twisted $\Lambda$-local system on $T$, for some $\Lambda$-admissible multiplier $\nu$ on $G_s$. We then have
$$
\varepsilon_{\pi,\overline{s}}(T,\F \otimes \G,\omega) = \det(\G_{\overline{s}})^{a(T,\F,\omega)} \varepsilon_{\pi,\overline{s}}(T,\F,\omega)^{\rk(\G)},
$$
where $\det \left( \G_{\overline{s}} \right)$ is the $\Lambda$-admissible map on $G_{s}$ which sends an element $g$ of $G_s$ to the determinant of its action on the stalk $\G_{\overline{s}}$ of $\G$ at $\overline{s}$.
\end{prop}

We can assume (and we do) that $\omega$ is $d \pi$, in which case the conclusion of Proposition \ref{chap3unramtwist} follows from the definition \ref{chap3localepsfact} and from the fact that $\pi_{\diamondsuit}(\F \otimes \G) = \pi_{\diamondsuit} \F \otimes \pi_{\diamondsuit} \G$ is the twist of $\pi_{\diamondsuit} \F$ by $\pi_{\diamondsuit} \G$, the geometrically constant $\nu$-twisted $\Lambda$-sheaf on $\mathbb{A}^1_s$ associated to the fiber $i^{-1} \G$.

\begin{prop}\label{chap3onapoint} Let $\F$ be a $\mu$-twisted $\Lambda$-sheaf on $s$. Then we have
$$
\varepsilon_{\pi,\overline{s}}(T,i_* \F,\omega) = \det \left(\F_{\overline{s}} \right)^{-1},
$$
where $\det \left( \F_{\overline{s}} \right)$ is the $\Lambda$-admissible map on $G_{s}$ which sends an element $g$ to the determinant of its action on the stalk $\F_{\overline{s}}$ of $\F$ at $\overline{s}$. In particular, $\varepsilon_{\pi,\overline{s}}(T,i_* \F,\omega)$ agrees with the local $\varepsilon$-factor $\varepsilon_{\overline{s}}(T,i_* \F,\omega)$ defined in \ref{chap32.7}.
\end{prop}

Indeed, the $\mu$-twisted $\Lambda$-sheaf $\pi_{\diamondsuit} i_* \F$ is supported on the $s$-point $0$ of $\mathbb{A}^1_s$, with restriction $\F$ to $0$, and consequently we have
$$
R\Gamma_c\left( \mathbb{A}^1_{\overline{s}}, \pi_{\diamondsuit }\F \otimes \Lc_{\psi}^{-1} \right) \cong \F_{\overline{s}}[0],
$$
hence the conclusion of Proposition \ref{chap3onapoint}.

\subsection{\label{chap35.1}} Let $k(\overline{\eta})$ be a separable closure of $k(\eta_{\overline{s}})$ and let $\overline{\eta} \rightarrow \eta_{\overline{s}}$ be the corresponding morphism of $k$-schemes. We then have an exact sequence
$$
1 \rightarrow I_{\eta} \rightarrow G_{\eta} \rightarrow G_s \rightarrow 1,
$$
where $G_s$ is the Galois group of the extension $\bk/k(s)$ and $G_{\eta} = \pi_1(\eta, \overline{\eta})$ is the Galois group of the extension $k(\overline{\eta})/k(\eta)$, while $I_{\eta} = \pi_1(\eta_{\overline{s}},\overline{\eta})$ is the inertia group of this extension. The functor $\F \mapsto \F_{\overline{\eta}}$ is then an equivalence from the category of $\mu$-twisted $\Lambda$-sheaves on $\eta$ to the category of $\Lambda$-admissible representations of $(G_{\eta},\mu)$, cf. \ref{chap31.4.4}.

\begin{prop}\label{localstructure} There exists a (unique) closed normal subgroup $I_{\eta,\ell}$ of $I_{\eta}$, of profinite order coprime to $\ell$, such that $I_{\eta}/I_{\eta,\ell}$ is a free pro-$\ell$-group on one generator.
\end{prop}

It follows from standard results in the ramification theory of henselian discretely valued fields, see for example (\cite{JPS2}, IV), that $I_{\eta}$ admits a unique pro-$p$-Sylow subgroup $P$, the \textit{wild inertia group} of $\eta$, and that there exists an isomorphism
$$
I_{\eta}/P \rightarrow \prod_{\ell' \neq p} \mathbb{Z}_{\ell'}(1).
$$
The conclusion then follows by taking $I_{\eta,\ell}$ to be the kernel of the projection $I_{\eta}/P \rightarrow \mathbb{Z}_{\ell}(1)$.

\subsection{\label{chap35.5}} We henceforth assume that $\Lambda$ is an algebraically closed field.

\begin{prop}\label{chap3generators}  Let $K_0(T,\mu,\Lambda)$ be as in \ref{chap30.6}. The abelian group
$$
\bigoplus_{\nu} K_0(T,\nu,\Lambda),
$$
where the sum runs over all unitary $\Lambda$-admissible multipliers on $G_s$, is generated by its subset of elements of the following three types:
\begin{enumerate}
\item the class $[j_! \Lambda]$ in $K_0(T,1,\Lambda)$,
\item for any unitary $\Lambda$-admissible multiplier $\nu$ on $G_s$ and any $\nu$-twisted $\Lambda$-sheaf $\G$ on $s$, the class $[i_* \G]$ in $K_0(T,\nu,\Lambda)$,
\item for any unitary $\Lambda$-admissible multiplier $\nu$ on $G_s$, any connected finite \'etale cover $ \eta' \rightarrow \eta$ with normalization $f : T' \rightarrow T$, for any $s$-morphism $\overline{s}' : \Spec(\bk) \rightarrow T'$ over the closed point $s'$ of $T'$ such that $f(\overline{s}') = \overline{s}$, any unitary $\Lambda$-admissible multipliers $\nu_1$ and $\nu_2$ on $G_{s'}$ such that $\nu_1 \nu_2 = \nu_{| G_{s'}}$, any $\nu_1$-twisted $\Lambda$-sheaf $\F_1$ of rank $1$ over $\eta'$, and any unramified $\nu_2$-twisted $\Lambda$-sheaf $\F_2$ over $\eta'$, the class 
$$
[f_* j'_!(\F_1 \otimes \F_2)] - \rk(\F_2) [f_* j'_!\Lambda],
$$
in the sum of $K_0(T,\nu,\Lambda)$ and $K_0(T,1,\Lambda)$, where $j' : \eta' \rightarrow T'$ is the canonical open immersion.

\end{enumerate}
\end{prop}

We have an isomorphism
$$
K_0(T,\nu,\Lambda) \xrightarrow[]{(i^{-1},j^{-1})} K_0(s,\nu,\Lambda) \oplus K_0(\eta,\nu,\Lambda).
$$
For $\Lambda=C$, the conclusion follows from \ref{chap31.4.4}, \ref{chap31.4.5}, and from Proposition \ref{chap3generators5}, whose assumptions are satisfied by Proposition \ref{localstructure}. The case $\Lambda$ is an algebraic closure of $\mathbb{F}_{\ell}$ follows from the previous case by Proposition \ref{CGTRES2} and by \ref{chap31.4.4} (or by \ref{1.4.8.3}).

\begin{prop}\label{chap3sanity54} Let $ \eta' \rightarrow \eta$ be a connected finite \'etale cover of $\eta$, with normalization $f : T' \rightarrow T$. Let $\overline{s}' : \Spec(\bk) \rightarrow T'$ be an $s$-morphism over the closed point $s'$ of $T'$ such that $f(\overline{s}') = \overline{s}$, and let $j' : \eta' \rightarrow T'$ be the natural inclusion. Let $\mu_1$ and $\mu_2$ be unitary $\Lambda$-admissible multipliers on $G_{s'}$ such that $\mu_1 \mu_2 = \mu_{| G_{s'}}$, let $\F_1$ be a $\mu_1$-twisted $\Lambda$-sheaf of rank $1$ over $\eta'$ and let $\F_2$ be a geometrically constant $\mu_2$-twisted $\Lambda$-sheaf over $T'$. Then for any element $\omega$ of $\Omega^{1,\times}_{\eta}$, there exists a $\Lambda$-admissible homomorphism $\lambda_{f,\pi}(\omega)$ from $G_s$ to $\Lambda^{\times}$, depending only on $f,\pi$ and $\omega$, such that the map $\varepsilon_{\pi,\overline{s}}(T,f_* (j'_!\F_1 \otimes \F_2), \omega) $ is equal to
$$
 \lambda_{f,\pi}(\omega)^{\rk(\F_2)} \delta_{s'/s}^{a(T',j'_! \F_1, f^{*} \omega)\rk(\F_2)} \Ver_{s'/s} \left( \det \left( \F_{2,\overline{s'}} \right)^{a(T',j'_! \F_1, f^{*} \omega)} \varepsilon_{\overline{s}'}(T',j'_! \F_1, f^{*} \omega)^{\rk(\F_2)} \right),
$$
where $\varepsilon_{\overline{s}'}(T',j'_! \F_1, f^{*} \omega)$ is the $\varepsilon$-factor defined in \ref{chap3localfact1} and where $\delta_{s'/s},\Ver_{s'/s}$ are defined in \ref{chap31.4.6} and \ref{chap3verlagerung} respectively. When $\eta' = \eta$, the factor $\lambda_{f,\pi}(\omega)$ is identically equal to $1$.

\end{prop}


By Proposition \ref{chap3changeforme} and Proposition \ref{chap3transferhomo}, we can further assume (and we do) that $\omega = d \pi$. Let us now consider the special cover $f' : U \rightarrow \mathbb{G}_{m,s}$ (cf. \ref{chap3specialcover}) associated to $\pi$ and to the extension $\eta'$ of $\eta$, cf. Theorem \ref{chap3specialcover2}. Let $\G_2$ be a $\nu_2$-twisted $\Lambda$-sheaf on $s'$ whose pullback to $\eta'$ is isomorphic to $\F_2$. Let $\G_1$ be the $f'$-special $\nu_1$-twisted $\Lambda$-local system of rank $1$ on $U$ associated to $\F_1$ by Theorem \ref{chap3GK2}, so that the pullback of $\G_1$ to $\eta'$ is isomorphic to $\F_1$. The $\nu_1$-twisted $\Lambda$-local system $\pi_{\diamondsuit} j_! f_*(\F_1)$ vanishes at $0$, and its restriction to $\mathbb{G}_{m,s}$ is isomorphic to $f'_* \G$, cf. \ref{chap31.14}. We thus have
\begin{align*}
\varepsilon_{\pi,\overline{s}}(T,f_* j'_!(\F_1 \otimes \F_2), d \pi) &= \varepsilon_{\overline{s}}(\mathbb{G}_{m,s},f'_* (\G_1 \otimes \G_2) \otimes \Lc_{\psi} \lbrace -t \rbrace ) \\ 
&=\varepsilon_{\overline{s}}(U, \G_1 \otimes \G_2 \otimes \Lc_{\psi} \lbrace -f' \rbrace ) \\
&= \delta_{s'/s}^{-\chi_c(U_{\overline{s}'}, \G_1 \otimes \G_2 \otimes \Lc_{\psi} \lbrace -f' \rbrace)} \Ver_{s'/s} \left( \varepsilon_{\overline{s}'}(U, \G_1 \otimes \G_2 \otimes \Lc_{\psi} \lbrace -f' \rbrace ) \right).
\end{align*}

Let $u : U \rightarrow X$ be a smooth compactification of $U$, so that $X$ is a geometrically connected smooth projective curve over $s$. Let $D$ be the closed complement in $X$ of the union of $U$ and the closed point of $T'$. For any point $x$ of $D$, we have
\begin{align*}
a(X_{(x)},u_! \G_1 \otimes \G_2 \otimes \Lc_{\psi} \lbrace -f' \rbrace, d f' ) &= (1 + \sw_x(\Lc_{\psi} \lbrace -f' \rbrace) + v_x(df')) \rk(\F_2) \\
&= 0,
\end{align*}
since $\G_1 \otimes \G_2 $ is tamely ramified at $x$ (cf. \ref{chap3specialcover}) and since $\sw_x(\Lc_{\psi} \lbrace -f' \rbrace)$ is equal to the valuation $v_x(f')$ of $f'$ at $x$, the latter being an integer prime to $p$. Consequently, the Grothendieck-Ogg-Shafarevich formula (\cite{La87}, 3.1.5.3) yields
$$
-\chi_c(U_{\overline{s}'}, \G_1 \otimes \G_2 \otimes \Lc_{\psi} \lbrace -f' \rbrace) = \alpha \rk(\F_2),
$$
where $\alpha = a(T',j'_! \F_1, f^{-1} d \pi)$. We also have 
$$
\varepsilon_{\overline{s}'}(U, \G_1 \otimes \G_2 \otimes \Lc_{\psi} \lbrace -f' \rbrace ) = \det(\G_{2,\overline{s}'})^{-\chi_c(U_{\overline{s}'}, \G_1 \otimes \Lc_{\psi} \lbrace -f' \rbrace)} \varepsilon_{\overline{s}'}(U, \G_1 \otimes \Lc_{\psi} \lbrace -f' \rbrace )^{\rk(\F_2)},
$$
with $-\chi_c(U_{\overline{s}'}, \G_1 \otimes \Lc_{\psi} \lbrace -f' \rbrace) = \alpha$, hence
\begin{align}\label{chap3eq56989}
\varepsilon_{\pi,\overline{s}}(T,f_* j'_!(\F_1 \otimes \F_2), d \pi) =  \delta_{s'/s}^{\alpha \rk(\F_2)} \Ver_{s'/s} \left( \det \left( \F_{2,\overline{s'}} \right)^{\alpha} \varepsilon_{\overline{s}'}(U, \G_1 \otimes \Lc_{\psi} \lbrace -f' \rbrace )^{\rk(\F_2)} \right).
\end{align}
If we have a formula of the form
\begin{align}\label{chap3eq56990}
\varepsilon_{\pi,\overline{s'}}(T_{s'},f_* j'_! \F_1, d \pi) = \lambda  \varepsilon_{\overline{s}'}(T',j'_! \F_1, f^{*} \omega),
\end{align}
for some $\Lambda$-admissible homomorphism $\Lambda$ from $G_{s'}$ to $\Lambda^{\times}$, then (\ref{chap3eq56989}) applied to $s=s'$ and $\F_2 = \Lambda$ yields
$$
\lambda  \varepsilon_{\overline{s}'}(T',j'_! \F_1, f^{*} \omega)  = \varepsilon_{\overline{s}'}(U, \G_1 \otimes \Lc_{\psi} \lbrace -f' \rbrace ),
$$
and inserting the latter formula in (\ref{chap3eq56989}) yields the required result with $\lambda_{f,\pi}(d \pi) = \Ver_{s'/s}(\lambda)$. It is therefore sufficient to prove (\ref{chap3eq56990}), so that we can assume that $s = s'$ and $\F_2 = \Lambda$ (and we henceforth do), in which case the formula to be proved is
$$
\varepsilon_{\pi,\overline{s}}(T,f_* j'_! \F_1, d \pi) = \lambda_{f,\pi}(d \pi) \varepsilon_{\overline{s}}(T',j'_! \F_1, f^{*} d \pi),
$$
for some $\Lambda$-admissible map $\lambda_{f,\pi}(d \pi)$ from $G_s$ to $\Lambda^{\times}$. We now apply the product formula \ref{chap3productform1} with $\omega = f^{\prime *} d t = df'$, according to which the map $\varepsilon_{\overline{s}}(U, \G_1 \otimes \Lc_{\psi} \lbrace -f' \rbrace)$ is equal to
$$
\lambda_{f,\pi}(d \pi) \varepsilon_{\overline{s}}(T',j'_! \F_1, f^{*} d \pi),
$$ 
where we have set 
$$
\lambda_{f,\pi}(d \pi) = \chi_{\cyc}^{g-1} \prod_{x \in D} \delta_{x/s}^{1 - v_x(f')} \Ver_{x/s} \left( \varepsilon_{\overline{x}}(X_{(x)},u_! \G_1 \otimes \Lc_{\psi} \lbrace - f' \rbrace_{|X_{(x)}}, df' ) \right),
$$
and $g$ is the genus of $X$. Moreover, for each point $x$ of $D$, the valuation $v_x(f')$ of $f'$ at $x$ is prime to $p$ since $f'$ is tamely ramified at this point (cf. \ref{chap3specialcover}), and the restriction $u_!\F_{|X_{(x)}}$ is tamely ramified as well, since $\F$ is $f'$-special. Hence we can apply Proposition \ref{chap3computation}, which implies
$$
\varepsilon_{\overline{x}}(X_{(x)},u_! \G_1 \otimes \Lc_{\psi} \lbrace -f' \rbrace_{|X_{(x)}}, df' ) = \varepsilon_{\overline{x}}(X_{(x)}, \Lc_{\psi} \lbrace -f' \rbrace_{|X_{(x)}}, df' ).
$$
Thus we obtain the required formula with
\begin{align}\label{chap3formulalanglands}
\lambda_{f,\pi}(d\pi) = \chi_{\cyc}^{g-1} \prod_{x \in D} \delta_{x/s}^{1 - v_x(f')} \Ver_{x/s} \left( \varepsilon_{\overline{x}}(X_{(x)}, \Lc_{\psi} \lbrace -f' \rbrace_{|X_{(x)}}, df' ) \right),
\end{align}
which depends indeed only on $f$ and $\pi$. When $f = \mathrm{id}$, we have $X = \mathbb{P}^1_s$ and $D$ is the closed point $\infty$, so that Proposition \ref{chap3computation} applies with $n=1$ and yields
$$
\varepsilon_{\overline{s}}(\mathbb{P}^1_{s,(\infty)}, \Lc_{\psi} \lbrace -t \rbrace_{|\mathbb{P}^1_{s,(\infty)}}, dt ) = \chi_{\cyc}.
$$
By inserting this identity in (\ref{chap3formulalanglands}), in which $g=0$, we obtain $\lambda_{\mathrm{id},\pi}(d\pi) = 1$, which in turn implies $\lambda_{\mathrm{id},\pi}(\omega) = 1$ for any $\omega$; this concludes the proof of Proposition \ref{chap3sanity54}.

\begin{rema} The formula (\ref{chap3formulalanglands}), combined with Proposition \ref{chap3computation}, yields an expression of $\lambda_{f,\pi}(d\pi)$ as a product of certain quadratic Gauss sums and of a power of the cyclotomic character. More precisely, with notation as in (\ref{chap3formulalanglands}), let $D^+$ (resp. $D^-$) be the closed subset of $X$ consisting of the points of $D$ such that $v_x(f')$ is odd (resp. even), endowed with its reduced scheme structure. Then we have, for any totally ramified extension $f$,
$$
\lambda_{f,\pi}(d\pi) = \chi_{\cyc}^{g-1 + \deg(D^-) +\frac{1}{2}(\deg(f) + \deg(D^+))} \prod_{x \in D^-} \delta_{x/s} \Ver_{x/s} \left( \gamma_{\psi}(v_x(f') h_x) \right),
$$
with notation as in \ref{chap3computation}, where $h_x$ is an element of $k(x)^{\times}$ such that $h_x f'$ is a square in the field of fractions of $\Ow_{X_{(x)}}$. This formula simplifies greatly when $p=2$, since $D^-$ is then empty.
\end{rema}

\begin{cor}\label{chap3rank1comp} Let $\F$ be a $\mu$-twisted $\Lambda$-sheaf on $T$ such that $j^{-1} \F$ is of rank at most $1$. Then for any element $\omega$ of $\Omega^{1,\times}_{\eta}$, we have
$$
\varepsilon_{\pi,\overline{s}}(T,\F,\omega) = \varepsilon_{\overline{s}}(T,\F,\omega),
$$
where $\varepsilon_{\overline{s}}(T,\F,\omega)$ is the $\varepsilon$-factor defined in \ref{chap3localfact1}.
\end{cor}

Indeed, this results from Proposition \ref{chap3onapoint} if $\F$ is supported on $s$, or from Proposition \ref{chap3sanity54} with $\eta' = \eta$ and $s'=s$ if $\F$ is supported on $\eta$. 

\begin{cor}\label{chap3indepunif} Let $\pi'$ be an other uniformizer of $k(\eta)$. Then for any element $\omega$ of $\Omega^{1,\times}_{\eta}$ and any $\mu$-twisted $\Lambda$-sheaf $\F$ on $T$, we have
$$
\varepsilon_{\pi,\overline{s}}(T,\F,\omega) = \varepsilon_{\pi',\overline{s}}(T,\F,\omega).
$$
\end{cor}

Indeed, by multiplicativity of $\varepsilon_{\pi}$ and $\varepsilon_{\pi'}$ in short exact sequences (cf. \ref{chap3multiplicativity}), the maps $\F \mapsto \varepsilon_{\pi,\overline{s}}(T,\F,\omega)$ and $\F \mapsto \varepsilon_{\pi',\overline{s}}(T,\F,\omega)$ extend to homomorphisms $\varepsilon_{\pi,\overline{s}}(T,\omega)$ and $\varepsilon_{\pi',\overline{s}}(T,\omega)$ respectively from the abelian group
$$
\bigoplus_{\nu} K_0(T,\nu,\Lambda),
$$
where the sum runs over all unitary $\Lambda$-admissible multipliers on $G_k$ (cf. \ref{chap35.5}), to $\Lambda^{\times}$. It is therefore sufficient to prove that $\varepsilon_{\pi',k}(T,\omega)$ and $\varepsilon_{\pi',k}(T,\omega)$ agree on the three types of generators described in Proposition \ref{chap3generators}. For generators of type $(1)$ or $(2)$ this follows from Corollary \ref{chap3rank1comp}, while this follows from Proposition \ref{chap3sanity54} for generators of type $(3)$.

\begin{notat}\label{chap3localepsfact2} Let $\F$ be a $\mu$-twisted $\Lambda$-sheaf on $T$ and let $\omega$ be an element of $\Omega^{1,\times}_{\eta}$. We denote by $\varepsilon_{\overline{s}}(T,\F,\omega)$ the $\Lambda$-admissible map $\varepsilon_{\pi,\overline{s}}(T,\F,\omega)$ (cf. \ref{chap3localepsfact}), which does not depend on the uniformizer $\pi$ by Corollary \ref{chap3indepunif}.
\end{notat}

By Corollary \ref{chap3indepunif}, this definition does not depend on the choice of the uniformizer $\pi$, and by Corollary \ref{chap3rank1comp} this does not conflict with the definition \ref{chap3localfact2} when $j^{-1}\F$ is of rank at most $1$.

\begin{prop}\label{chap3inductionformula} Let $ \eta' \rightarrow \eta$ be a connected finite \'etale cover of $\eta$, with normalization $f : T' \rightarrow T$. Let $\overline{s}' : \Spec(\bk) \rightarrow T'$ be an $s$-morphism over the closed point $s'$ of $T'$ such that $f(\overline{s}') = \overline{s}$, and let $j' : \eta' \rightarrow T'$ be the natural inclusion. Then for any element $\omega$ of $\Omega^{1,\times}_{\eta}$, there exists a (unique) $\Lambda$-admissible homomorphism $\lambda_{f}(\omega)$ from $G_{\overline{s}}$ to $\Lambda^{\times}$, depending only on $f$ and $\omega$ with the following property: for any $\mu$-twisted $\Lambda$-sheaf $\F$ on $T'$, we have
$$
\varepsilon_{\overline{s}}(T,f_*\F,\omega) = \lambda_{f}(\omega)^{\rk(j'^{-1}\F)} \delta_{s'/s}^{a(T',\F,f^*\omega)} \Ver_{s'/s} \left( \varepsilon_{\overline{s}'}(T',\F,f^*\omega) \right),
$$
where $\delta_{s'/s}$ and $\Ver_{s'/s}$ are defined in \ref{chap31.4.6} and \ref{chap3verlagerung} respectively.
\end{prop}

Indeed, Corollary \ref{chap3indepunif} implies that the group homomorphism $\lambda_{f,\pi}(\omega)$ from Proposition \ref{chap3sanity54} does not depend on $\pi$. Let us denote it by $\lambda_{f}(\omega)$. The maps
\begin{align*}
A : \F &\mapsto \varepsilon_{\overline{s}}(T,f_*\F,\omega) \\
B : \F &\mapsto \lambda_{f}(\omega)^{\rk(j'^{-1}\F)} \delta_{s'/s}^{a(T',\F,f^*\omega)} \Ver_{s'/s} \left( \varepsilon_{\overline{s}'}(T',\F,f^*\omega) \right),
\end{align*}
both extend by multiplicativity (cf. \ref{chap3multiplicativity}) to homomorphisms from the abelian group
$$
\bigoplus_{\nu} K_0(T',\nu,\Lambda),
$$
where the sum runs over all unitary $\Lambda$-admissible multipliers on $G_s$ (cf. \ref{chap35.5}), to $\Lambda^{\times}$. Consequently, it is sufficient to check that the homomorphisms $A$ and $B$ coincide on the three type of generators described in Proposition \ref{chap3generators}. For the generator of type $(1)$, this follows from Proposition \ref{chap3sanity54} with $\nu_1 = \nu_2 = 1$ and $\G_1 = \G_2 = \Lambda$, while it holds as well for generators of type $(2)$ by Proposition \ref{chap3onapoint}. It remains to handle generators of the third type described in Proposition \ref{chap3generators}.

Let $ \eta'' \rightarrow \eta'$ be a connected finite \'etale cover of $\eta'$, with normalization $f' : T'' \rightarrow 'T$. Let $\overline{s}'' : \Spec(\bk) \rightarrow T''$ be an $s'$-morphism over the closed point $s''$ of $T''$ such that $f(\overline{s}'') = \overline{s}'$, and let $j'' : \eta' \rightarrow T'$ be the natural inclusion. Let $\mu_1$ and $\mu_2$ be unitary $\Lambda$-admissible multipliers on $G_{s''}$ such that $\mu_1 \mu_2 = \mu_{| G_{s''}}$, let $\F_1$ be a $\mu_1$-twisted $\Lambda$-sheaf of rank $1$ on $\eta''$, with finite geometric monodromy, and let $\F_2$ be an unramified $\mu_2$-twisted $\Lambda$-sheaf on $\eta''$. We must prove that the homomorphisms $A$ and $B$ associate the same map to the class
$$
[f'_* j''_! (\F_1 \otimes \F_2)] - \rk(\F_2) [f'_* j''_! \Lambda].
$$
By applying Proposition \ref{chap3sanity54} to the extension $\eta'' \rightarrow \eta'$, we obtain an equality
\begin{align*}
 A( f'_* j''_! (\F_1 \otimes \F_2) ) =\lambda_{ff'}(\omega)^{\rk(\F_2)} \delta_{s''/s}^{a} \Ver_{s''/s} \left( \Lc \right),
\end{align*}
with $a = a(T'',j''_! \F_1, (ff')^{*} \omega)\rk(\F_2)$ and
$$
\Lc = \det \left( \F_{2,\overline{s}''} \right)^{a(T'',j''_! \F_1, (ff')^{*} \omega)} \varepsilon_{\overline{s}''}(T'',j''_! \F_1, (ff')^{*} \omega)^{\rk(\F_2)}.
$$ 
A second application of Proposition \ref{chap3sanity54} yields, together with \ref{chap32.11998},
\begin{align*}
B( f'_* j''_! (\F_1 \otimes \F_2) ) = \lambda_{f}(\omega)^{\rk(\F_2) \mathrm{deg}(f')} \delta_{s'/s}^{[s'':s']a} \Ver_{s'/s} \left( \lambda_{f'}(f^* \omega)^{\rk(\F_2)}\delta_{s''/s'}^{a} \Ver_{s''/s'} \left( \Lc \right)\right),
\end{align*}
hence by Proposition \ref{chap3transferhomo} we have
$$
B( f'_* j''_! (\F_1 \otimes \F_2) ) =\left( \lambda_{f}(\omega) \Ver_{s'/s} \left( \lambda_{f'}(f^* \omega) \right) \right)^{\rk(\F_2)} \delta_{s'/s}^{[s'':s']a} \Ver_{s'/s} \left( \delta_{s''/s'} \right)^{a} \Ver_{s'/s} \Ver_{s''/s'}(\Lc),
$$
By Corollary \ref{chap3compositionverlag}, this yields
$$
(A / B)( f'_* j''_!(\F_1 \otimes \F_2)) = \left( \lambda_{ff'}(\omega) \lambda_{f}(\omega)^{-1}\Ver_{s'/s} \left( \lambda_{f'}(f^* \omega) \right)^{-1} \right)^{\rk(\F_2)}.
$$
By applying this formula to $\F_1 = \F_2 = \Lambda$, we obtain
$$
(A / B)( f'_* j''_!(\F_1 \otimes \F_2)) = (A / B)( f'_* j''_!\Lambda)^{\rk(\F_2)},
$$
hence the equality
$$
A( [f'_* j''_!(\F_1 \otimes \F_2)] - \rk(\F_2) [f'_* j''_!\Lambda] ) = B( [f'_* j''_!(\F_1 \otimes \F_2)] - \rk(\F_2) [f'_* j''_!\Lambda] ),
$$
which concludes our proof of Corollary \ref{chap3inductionformula}.



\begin{cor}\label{chap3transitivityformula} Let us consider a tower $\eta'' \rightarrow \eta' \rightarrow \eta$ of connected finite separable extensions, and let
$$
T'' \xrightarrow[]{f'} T' \xrightarrow[]{f} T,
$$
be the normalizations of $T$ in $\eta'$ and $\eta''$, with closed points $s''$ and $s'$. Then for any element $\omega$ of $\Omega^{1,\times}_{\eta}$, we have
$$
\lambda_{ff'}(\omega) = \lambda_{f}(\omega) \Ver_{s'/s} \left( \lambda_{f'}(f^* \omega) \right),
$$
with notation as in Proposition \ref{chap3inductionformula}.
\end{cor}

This is an immediate consequence of Propositions \ref{chap3inductionformula} and \ref{chap3compositionverlag}.

\subsection{\label{chap35.8}} We can now prove Theorem \ref{chap3teo2}. The rule $\varepsilon$ defined in \ref{chap3localepsfact2} and \ref{chap3localepsfact} clearly satisfies the properties $(i)$ and $(ii)$ from \ref{chap310.2}. It also satisfies the properties $(iii),(iv),(v),(vi),(vii)$ from \ref{chap310.2} by \ref{chap3multiplicativity}, \ref{chap3onapoint}, \ref{chap3inductionformula}, \ref{chap3rank1comp} and \ref{chap3unramtwist} respectively. This proves the existence statement in Theorem \ref{chap3teo2}, while the uniqueness is an immediate consequence of Proposition \ref{chap3generators}. The property $(viii)$ in Theorem \ref{chap3teo2} follows from the uniqueness, while the property $(ix)$ follows from Proposition \ref{chap3sanity5}.

\subsection{\label{chap35.9}} Let us now prove Theorem \ref{chap3teo4}. We assume that $k$ is a finite field. By Theorem \ref{chap3teo2} and by Proposition \ref{chap3compafini}, the rule which associates the quantity
$$
(-1)^{a(T,\F)} \varepsilon_{\overline{s}}(T,\F,\omega)(\Frob_s),
$$
to any quadruple $(T,\F,\omega, \overline{s})$, where $T$ is a henselian trait over $k$, with closed point $s$ finite over $k$, where $\overline{s} : \Spec(\bk) \rightarrow T$ is a $k$-morphism, where $\F$ is a $\Lambda$-sheaf on $T$ and where $\omega$ is a non zero meromorphic $1$-form on $T$, satisfies all the properties listed in (\cite{La87}, Th. 3.1.5.4), hence must coincide with the local $\varepsilon$-factor considered there.

\section{The product formula: preliminaries\label{chap3productsection2}}

 Let $\Lambda$ be an $\ell$-adic coefficient ring (cf. \ref{chap3conv}, \ref{chap30.0.0.1}) which is an algebraically closed field, and let $\psi : \mathbb{F}_{p} \rightarrow \Lambda^{\times}$ be a non trivial homomorphism. We fix a unitary $\Lambda$-admissible mutiplier $\mu$ on the topological group $G_k$ (cf. \ref{chap30.0}, \ref{chap3unitary}). Let $X$ be a connected smooth projective curve of genus $g$ over $k$ and let $\omega$ be a non zero global meromorphic differential $1$-form on $X$. We denote by $|X|$ the set of closed points of $X$. For any closed point $x$ of $X$, we denote by $X_{(x)}$ the henselization of $X$ at $x$.
 
 \begin{defi}\label{defiproductform} Let $\F$ be a $\mu$-twisted $\Lambda$-sheaf on $X$, or a class in $\bigoplus_{\nu} K_0(X,\nu,\Lambda)$, where the sum runs over all unitary $\Lambda$-admissible mutipliers on $G_k$. We denote by $\tilde{\varepsilon}_{\bk}(X, \F,\omega)$ the quantity
 $$ 
 \tilde{\varepsilon}_{\bk}(X, \F,\omega) = \chi_{\cyc}^{N(g-1) \rk(\F)} \prod_{x \in |X|} \delta_{x/k}^{a(X_{(x)},\F_{|X_{(x)}})} \Ver_{x/k} \left( \varepsilon_{\overline{x}}(X_{(x)},\F_{|X_{(x)}},\omega_{|X_{(x)}}) \right),
 $$
  where $N$ is the number of connected components of $X_{\bk}$, and where $\delta_{x/k}$ and $\Ver_{x/k}$ are defined in \ref{chap31.4.6}. We say that $(X,\F, \omega)$ \textit{satisfies the product formula over} $k$ if the $\Lambda$-admissible maps $\varepsilon_{\bk}(X, \F,\omega)$ and $ \tilde{\varepsilon}_{\bk}(X, \F)$ coincide.

 \end{defi}
 
In particular, Theorem \ref{chap3productform1} can be stated as establishing that for any $\mu$-twisted $\Lambda$-sheaf of generic rank at most $1$, the triple $(X,\F,\omega)$ satisfies the product formula over $k$. It is usually convenient to introduce the character
\begin{align}\label{chap3delta}
\delta_{X/k}(\omega)= \chi_{\cyc}^{N(g-1)} \prod_{x \in |X|} \delta_{x/k}^{v_x(\omega)},
\end{align} 
so that we have, by noting that $\delta_{x/k}^{v_x(\omega)} = \delta_{x/k}^{-v_x(\omega)}$, the following formula for $\tilde{\varepsilon}_{\bk}(X, \F,\omega)$:
\begin{align}\label{chap3productform5}
\tilde{\varepsilon}_{\bk}(X, \F,\omega)= \delta_{X/k}(\omega)^{\rk(\F)} \prod_{x \in |X|} \delta_{x/k}^{a(X_{(x)},\F_{|X_{(x)}},\omega_{|X_{(x)}})} \Ver_{x/k} \left( \varepsilon_{\overline{x}}(X_{(x)},\F_{|X_{(x)}},\omega_{|X_{(x)}}) \right).
\end{align}

\subsection{\label{36.0}} Let $s$ be the spectum of a finite extension of $k$, such that $X$ is a geometrically connected $s$-scheme, and let us fix a $k$-morphism $\overline{s} : \Spec(\bk) \rightarrow s$.

\begin{prop} If $(X,\F,\omega)$ satisfies the product formula over $s$, then it  satisfies the product formula over $k$.
\end{prop}

Indeed, we have
$$
\varepsilon_{\bk}(X, \F) = \delta_{s/k}^{-\chi(X_{\overline{s}},\F)} \Ver_{s/k}(\varepsilon_{\overline{s}}(X,\F)),
$$
by \ref{chap3transferhomo}, while Corollary \ref{chap3compositionverlag} yields
$$
\tilde{\varepsilon}_{\bk}(X, \F,\omega) = \delta_{s/k}^{\tau} \Ver_{s/k}(\tilde{\varepsilon}_{\overline{s}}(X,\F,\omega)),
$$
 where $\tau = \sum_{x \in |X|} [k(x):k(s)] a(X_{(x)},\F_{|X_{(x)}})$. The Grothendieck-Ogg-Shafarevich formula implies that $\tau$ has the same parity as $-\chi(X_{\overline{s}},\F)$, hence the result.
 
 \subsection{\label{36.1}} We henceforth assume that the smooth projective curve $X$ is geometrically connected over $k$. 
 
  \begin{prop}\label{unramtwistglobal} Let $\nu$ be a $\Lambda$-admissible mutiplier on $G_k$ and let $\G$ be a geometrically constant $\nu$-twisted $\Lambda$-sheaf on $X$. If $(X,\F,\omega)$ satisfies the product formula, then so does $(X,\F \otimes \G,\omega)$.
 \end{prop}

Indeed, if we denote by $V$ the $\Lambda$-admissible $G_k$-representation associated to $\G$, then we have
$$
\varepsilon_{\bk}(X, \F \otimes \G) = \varepsilon_{\bk}(X, \F)^{\rk(V)} \det(V)^{-\chi(X_{\bk},\F)}, 
$$
while Proposition \ref{chap3unramtwist} yields
$$
\tilde{\varepsilon}_{\bk}(X, \F \otimes \G,\omega) = \tilde{\varepsilon}_{\bk}(X, \F,\omega)^{\rk(V)} \det(V)^{\tau},  
$$
where $\tau = \sum_{x \in |X|} [k(x):k] a(X_{(x)},\F_{|X_{(x)}}, \omega_{|X_{(x)}})$. The result then follows from the Grothendieck-Ogg-Shafarevich formula, which states that $\tau$ coincides with $-\chi(X_{k},\F)$.

 \begin{prop}\label{inductionlemma} Let $f : X' \rightarrow X$ be a finite generically \'etale morphism of smooth connected projective curves over $k$, and let $\F$ be a class in $\bigoplus_{\nu} K_0(X,\nu,\Lambda)$ of generic rank $0$. Then $(X',\F,f^* \omega)$ satisfies the product formula over $k$, if and only if $(X,f_* \F, \omega)$ satisfies the product formula over $k$.
 \end{prop}
 
 Indeed, we have 
 $$
 \varepsilon_{\bk}(X, f_*\F) = \varepsilon_{\bk}(X', \F),
 $$
 while Corollary \ref{chap3compositionverlag} yields that $\tilde{\varepsilon}_{\bk}(X', \F,f^*\omega) $ can be written as
\begin{align*}
 & \prod_{x' \in |X'|} \delta_{x'/k}^{a(X'_{(x')},\F, f^* \omega)} \Ver_{x'/k} \left( \varepsilon_{\overline{x'}}(X'_{(x')},\F_{|X'_{(x')}},f^*\omega_{|X'_{(x')}}) \right),\\
=& \prod_{x \in |X|} \delta_{x/k}^{a(X_{(x)},f_* \F, \omega)} \Ver_{x/k} \left( \prod_{\substack{x' \in |X'| \\ f(x') = x}} \delta_{x'/x}^{a(X'_{(x')},\F, f^* \omega)} \Ver_{x'/x}\left( \varepsilon_{\overline{x'}}(X'_{(x')},\F_{|X'_{(x')}},f^*\omega_{|X'_{(x')}}) \right)\right)
\end{align*}
The induction formula \ref{chap3inductionformula} then yields
\begin{align*}
\tilde{\varepsilon}_{\bk}(X', \F,f^*\omega)  &= \prod_{x \in |X|} \delta_{x/k}^{a(X_{(x)},f_* \F, \omega)} \Ver_{x/k} \left( \varepsilon_{\overline{x}}(X_{(x)},f_*\F_{|X_{(x)}},\omega_{|X_{(x)}}) \right) \\
&= \tilde{\varepsilon}_{\bk}(X, f_*\F,\omega),
\end{align*}
hence the conclusion of Proposition \ref{inductionlemma}.
 
 \subsection{\label{chap36.3}} Let $\eta$ be the generic point of $X$ and let $\overline{\eta}$ be the spectrum of a separable closure of $k(\eta)$. For any non empty open subscheme $U$ in $X$, we have an exact sequence
 $$
 1 \rightarrow \pi_1(U_{\overline{k}}, \overline{\eta}) \rightarrow \pi_1(U, \overline{\eta}) \rightarrow G_k \rightarrow 1,
 $$
 of profinite groups.

 \begin{defi}\label{chap3finitemonodromy} A $\mu$-twisted $\Lambda$-sheaf $\F$ on $X$ has \textit{finite geometric monodromy} if there exists a non empty open subscheme $U$ of $X$ such that $\F_{|U}$ is a $\mu$-twisted $\Lambda$-local system such that $\pi_1(U_{\overline{k}}, \overline{\eta})$ acts through a finite quotient on the $\Lambda$-admissible representation $\F_{\overline{\eta}}$ of $( \pi_1(U, \overline{\eta}),\mu)$.
 \end{defi}

Let $K_0^{\fin}(X,\mu,\Lambda)$ be the Grothendieck group of the full subcategory of $\Sh(X,\mu,\Lambda)$ consisting of $\mu$-twisted $\Lambda$-sheaves on $X$ with finite geometric monodromy (cf. \ref{chap3finitemonodromy}). Thus any $\mu$-twisted $\Lambda$-sheaf $\F$ on $X$ with finite geometric monodromy has a well defined class $[\F]$ in $K_0^{\fin}(X,\mu,\Lambda)$, and the latter is generated by such classes with relations $[\F] = [\F'] + [\F'']$ for each short exact sequence
$$
0 \rightarrow \F' \rightarrow \F \rightarrow \F'' \rightarrow 0,
$$
of $\mu$-twisted $\Lambda$-sheaves on $X$, with finite geometric monodromy.

\begin{prop}\label{chap3generators2} Let $X$ be a connected smooth projective curve over $k$, let $\overline{s}$ and $K_0^{\fin}(X,\mu,\Lambda)$ be as in \ref{chap36.3}. The abelian group
$$
\bigoplus_{\nu} K_0^{\fin}(X,\nu,\Lambda),
$$
where the sum runs over all unitary $\Lambda$-admissible multipliers on $G_k$, is generated by its subset of elements of the following three types:
\begin{enumerate}
\item the class $[\Lambda]$ in $K_0^{\fin}(X,1,\Lambda)$,
\item for any unitary $\Lambda$-admissible multiplier $\mu$ on $G_k$, any closed point $i_x : x \rightarrow X$ of $X$ and any $\mu$-twisted $\Lambda$-sheaf $\G$ on $x$, the class $[i_{x*} \G]$ in $K_0^{\fin}(X,\mu,\Lambda)$,
\item for any finite extension $s' \rightarrow \Spec(k)$, together with a $k$-morphism $\overline{s}' : \Spec(\bk) \rightarrow s'$, for any unitary $\Lambda$-admissible multiplier $\mu$ on $G_k$, for any non empty open subscheme $j : U \rightarrow X$, any finite \'etale morphism $f : U' \rightarrow U$ such that $U'$ is a geometrically connected $s'$-scheme, any unitary $\Lambda$-admissible multipliers $\mu_1$ and $\mu_2$ on $G_{s'}$ such that $\mu_1 \mu_2 = \mu_{| G_{s'}}$, any $\mu_1$-twisted $\Lambda$-sheaf $\F_1$ of rank $1$ over $U'$ with finite geometric monodromy, and any geometrically constant $\mu_2$-twisted $\Lambda$-sheaf $\F_2$ over $U'$, the class $[j_! f_*(\F_1 \otimes \F_2)] - \rk(\F_2) [j_!f_*\Lambda] $ in the sum of $K_0^{\fin}(X,\mu,\Lambda)$ and $K_0^{\fin}(X,1,\Lambda)$.

\end{enumerate}
\end{prop}

Let $\F$ be a $\mu$-twisted $\Lambda$-sheaf on $X$ with finite geometric monodromy. Let $j : U \rightarrow X$ be a non empty open subscheme of $X$ such that $j^{-1} \F$ is a $\mu$-twisted $\Lambda$-local system. We have an exact sequence
$$
0 \rightarrow j_! j^{-1} \F \rightarrow \F \rightarrow \bigoplus_{x \in |X \setminus U|}i_{x*} i_x^{-1} \F \rightarrow 0,
$$
hence we can assume (and we do) that $\F = j_! j^{-1} \F$. For $\Lambda = C$, the result then follows from Proposition \ref{chap3generators5} by using the dictionary between $\mu$-twisted $C$-local systems on $U$ and $C$-admissible representations of $(\pi_1(U, \overline{\eta}),\mu)$, for a fixed geometric point $\overline{\eta}$ over the generic point of $X$, cf. \ref{chap31.4.4} and \ref{chap31.4.5}. The general case follows from the case $\Lambda = C$ by \ref{CGTRES2}.

\begin{cor}\label{chap3productform2} If $\F$ is a $\mu$-twisted $\Lambda$-sheaf on $X$ with finite geometric monodromy, cf. \ref{chap3finitemonodromy}, then $(X,\F,\omega)$ satisfies the product formula over $k$.
\end{cor}

This follows from Propositions \ref{chap3generators2}, \ref{inductionlemma}, \ref{unramtwistglobal} and from Theorem \ref{chap3productform1}.

\begin{cor}\label{chap3productform3} Let $f : X' \rightarrow X$ be a finite generically \'etale morphism of connected smooth projective curves over $k$, of respective genera $g'$ and $g$. Then for any non zero global meromorphic differential $1$-form $\omega$ on $X$, we have
$$
\prod_{x' \in |X'|} \Ver_{f(x')/k}\left( \lambda_{f_{x'}}(\omega_{|X_{(f(x'))}}) \right)= \delta_{X'/k}(f^*\omega) \delta_{X/k}(\omega)^{-1},
$$
where $f_{x'} : X'_{(x)} \rightarrow X_{(f(x'))}$ is for each closed point $x'$ of $X'$ the morphism induced by $f$ on the henselizations of $X'$ and $X$ at $x'$ and $f(x')$ respectively, where $\lambda_{f_{x'}}(\omega_{|X_{(f(x'))}})$ is the homomorphism defined in \ref{chap3inductionformula}, and where $\delta_{X'/k}(f^*\omega)$, $ \delta_{X/k}(\omega)$ are as in \ref{chap3delta}.
\end{cor}

Indeed, this follows from Proposition \ref{chap3inductionformula} and from Corollary \ref{chap3productform2} for the triples $(X',\Lambda,f^*\omega)$ and $(X,f_*\Lambda,\omega)$, since $\varepsilon(X, f_* \Lambda)$ is equal to $\varepsilon(X', \Lambda)$.

 \subsection{\label{chap36.9}} We conclude this section by giving a variant of Corollary \ref{chap3productform2}.
 
 \begin{prop}\label{chap3productform20} Let $\F$ is a $\mu$-twisted $\Lambda$-sheaf on $X$. Let us assume that there exists a non empty open subscheme $U$ of $X$ and a closed normal subgroup $I$ of $\pi_1(U_{\overline{k}}, \overline{\eta})$ with the following properties: the restriction $\F_{|U}$ is a $\mu$-twisted $\Lambda$-local system, the group $I$ acts through a finite quotient on the $\Lambda$-admissible representation $\F_{\overline{\eta}}$ of $( \pi_1(U, \overline{\eta}),\mu)$, and $\pi_1(U_{\overline{k}},\overline{\eta})/I$ is topologically of finite type with open centralizer in $\pi_1(U, \overline{\eta})/I$. Then $(X,\F,\omega)$ satisfies the product formula over $k$.
 \end{prop} 

We can assume (and we henceforth do) that $\Lambda =C$, since otherwise $F$ would have finite geometric monodromy and Corollary \ref{chap3productform2} would apply. By Proposition \ref{chap3dec}, we can assume (and we henceforth do) that there exists a smooth geometricall connected curve $X'$ over a finite extension $k'$ of $k$ contained in $\bk$, a finite generically \'etale $k$-morphism $f : X' \rightarrow X$, some $\Lambda$-admissible unitary $2$-cocycles $\mu_1,\mu_2$ on $G_{k'}$ such that $\mu_1 \mu_2 = \mu_{|G_{k'}}$, some geometrically constant $\mu_2$-twisted $\Lambda$-sheaf $\F_2$ and some $\mu_1$-twisted $\Lambda$-sheaf $\F_1$, such that 
$$
\F = f_*(\F_1 \otimes \F_2),
$$
such that $F_{1|f^{-1}(U)}$ is a $\Lambda$-local system, and such that the $\Lambda$-admissible representation $\F_{1,\overline{\eta}}$ of $( \pi_1(U, \overline{\eta}),\mu)$ has finite projective image. By Theorem \ref{chap3brauer2}, we can further assume (and we do) that $\F_1$ is generically of rank $1$. The product formula for $f_*(\F_1 \otimes \F_2)$ then follows from Propositions \ref{inductionlemma}, \ref{unramtwistglobal} and from Theorem \ref{chap3productform1}.

\section{The product formula: proof \label{chap3productsection3}}

 
 Let $\Lambda$ be an $\ell$-adic coefficient ring (cf. \ref{chap3conv}, \ref{chap30.0.0.1}) which is an algebraically closed field, and let $\psi : \mathbb{F}_{p} \rightarrow \Lambda^{\times}$ be a non trivial homomorphism. We fix a unitary $\Lambda$-admissible mutiplier $\mu$ on the topological group $G_k$ (cf. \ref{chap30.0}, \ref{chap3unitary}). We prove in this section the following product formula.

\begin{teo}\label{chap3xteo3} Let $X$ be a connected smooth projective curve over $k$, let $\omega$ be a non zero global meromorphic differential $1$-form on $X$ and let $\F$ be a $\mu$-twisted $\Lambda$-sheaf on $X$. Then $(X,\F,\omega)$ satisfies the product formula over $k$ (cf. \ref{defiproductform}).
\end{teo}

We postpone the proof of Theorem \ref{chap3xteo3} to \ref{chap315.8} below.
When $k$ is finite, this result is due to Laumon (\cite{La87}, 3.2.1.1). For the general case, we follow closely Laumon's proof, or rather its exposition by Katz in \cite{katz3}. The main ingredient we use is the $\ell$-adic stationary phase method, of which we only use the special case stated in Theorem \ref{chap3xteostat} below, and which was already established by Laumon (\cite{La87}, 2.3.3.1) in the case of an arbitrary perfect base field of positive characteristic. The only innovation in our proof lies in the treatment of Theorem \ref{chap3xteolaumon} below. Laumon's proof of the latter result when $k$ is finite (\cite{La87}, 3.5.1.1) starts with a reduction to the tamely ramified case (\cite{La87}, 3.5.3.1), and then resorts to a computation in the latter case (\cite{La87}, 2.5.3.1). We give instead a direct proof in the general case by using geometric local class field theory.

\subsection{\label{chap315.0}} Throughout this section, we consider two copies $\mathbb{A} = \Spec(k[t])$ and $\mathbb{A}' = \Spec(k[t'])$ of the affine line over $k$, with natural compactifications $\mathbb{P}$ and $\mathbb{P}'$ respectively, and we denote by $\pr$ and $\pr'$ the projections of $\mathbb{A} \times_k \mathbb{A}'$ onto its first and second factors respectively. For any $\mu$-twisted $\Lambda$-sheaf $\F$ on $\mathbb{A}$, we denote by $\mathrm{F}_{\psi}(\F)$ its Fourier transform defined as follows:
$$\mathrm{F}_{\psi}(\F) = R^1 \pr'_! \left( \pr^{-1} \F \otimes \Lc_{\psi} \lbrace t t' \rbrace \right),
$$
which is a $\mu$-twisted $\Lambda$-sheaf on $\mathbb{A}'$. This functor $\mathrm{F}_{\psi}$ would be denoted $\mathcal{H}^0(\mathcal{F}_{\psi})$ in Laumon's notation (\cite{La87}, 1.2.1.1), and is called the ``naive Fourier transform'' by Katz in (\cite{katz3}, p.112).

\subsection{\label{chap315.1}} For any closed points $s,s'$ of $\mathbb{P}$ and $\mathbb{P}'$ respectively, with respective separable closures $\overline{s}$ and $\overline{s}'$, we denote by $\eta_s$ (resp. $\eta_{s'}$) the generic point of the henselisation $\mathbb{P}_{(s)}$ (resp. $\mathbb{P}'_{(s')}$) of $\mathbb{P}$ at $s$ (resp. of $\mathbb{P}'$ at $s'$), and by $\overline{\eta}_{s}$ (resp. $\overline{\eta}_{s'}$) a separable closure thereof. For any $\mu$-twisted $\Lambda$-sheaf $\F$ on $\mathbb{A}_{(s)}$ with vanishing fiber at $s$, we denote by $\mathrm{F}_{\psi}^{(s,s')}(\F)$ its local Fourier transform defined as follows:
$$
\mathrm{F}_{\psi}^{(s,s')}(\F) = H^1 \left( \left( (\mathbb{P} \times_k \mathbb{P}')_{(\overline{s},\overline{s}')} \right)_{\overline{\eta}_{s'}},u_!( \pr^{-1} \F \otimes \Lc_{\psi} \lbrace t t' \rbrace ) \right),
$$
where $u$ is the natural open immersion of $\mathbb{A} \times_k \mathbb{A}'$ into $\mathbb{P} \times_k \mathbb{P}'$. Thus $\mathrm{F}_{\psi}^{(s,s')}(\F) $ is a $\Lambda$-admissible representation of $(G_{\eta_{s'}}, \mu)$, where $G_{\eta_{s'}} = \Gal(\overline{\eta}_{s'}/\eta_{s'})$.

%

\subsection{\label{chap315.3}} We now state a special case of Laumon's $\ell$-adic stationary phase method. 

\begin{teo}[$\ell$-adic stationary phase]\label{chap3xteostat} Let $\F$ be a $\mu$-twisted $\Lambda$-sheaf on $\mathbb{P}$, whose fibers at $0$ and $\infty$ vanish, whose restriction to $\mathbb{A}\setminus \{ 0 \}$ is a $\Lambda$-local system, and whose ramification at $\infty$ is bounded by $1$, i.e. the ramification slopes of $\F_{|\eta_{\infty}}$ are strictly less than $1$. Then the $\Lambda$-sheaf $\mathrm{F}_{\psi}(\F_{|\mathbb{A}})$ on $\mathbb{A}'$ (cf. \ref{chap315.0}) has the following properties:
\begin{itemize}
\item[$(i)$] the restriction of $\mathrm{F}_{\psi}(\F_{|\mathbb{A}})$ to $\mathbb{A}' \setminus \{ 0\}$ is a $\Lambda$-local system;
\item[$(ii)$] there is a functorial isomorphism 
$$
\mathrm{F}_{\psi}(\F_{|\mathbb{A}})_{| \overline{\eta}_{\infty'}} \cong \mathrm{F}_{\psi}^{(0,\infty')}(\F_{|\mathbb{A}_{(0)}}),
$$
of $\Lambda$-admissible representations of $(G_{\eta_{\infty'}}, \mu)$ (cf. \ref{chap315.1});
\item[$(iii)$] the restriction of $\mathrm{F}_{\psi}(\F)$ to $\eta_{0'}$ fits into a functorial exact sequence,
$$
0 \rightarrow H^{1}(\mathbb{P}_{\bk}, \F) \rightarrow \mathrm{F}_{\psi}(\F_{|\mathbb{A}})_{| \overline{\eta}_{0'}} \rightarrow \mathrm{F}_{\psi}^{(\infty,0')}(\F_{|\mathbb{P}_{(\infty)}}) \rightarrow H^{2}(\mathbb{P}_{\bk}, \F) \rightarrow 0,
$$
of $\Lambda$-admissible representations of $(G_{\eta_{0'}}, \mu)$ (cf. \ref{chap315.1}), where $H^{\nu}(\mathbb{P}_{\bk}, \F)$ is considered as an unramified representation of $(G_{\eta_{0'}}, \mu)$ by the natural homomorphism $G_{\eta_{0'}} \rightarrow G_k$.

\end{itemize}
\end{teo}

By functoriality, Theorem \ref{chap3xteostat} follows from its untwisted special case, i.e. when $\mu = 1$, which follows itself from (\cite{katz3}, Th. 3 and 10) or from (\cite{La87}, 2.3.3.1, 2.3.2), the latter being applied to the perverse sheaf $\F[1]$.

\begin{cor}\label{chap3detcor} Let $\F$ be a $\mu$-twisted $\Lambda$-sheaf on $\mathbb{P}$, whose fibers at $0$ and $\infty$ vanish, whose restriction to $\mathbb{A}\setminus \{ 0 \}$ is a $\Lambda$-local system, and whose ramification at $\infty$ is bounded by $1$. Then we have an equality
$$
\varepsilon_{\bk}(\mathbb{P}, \F) \langle \chi_{\det(\mathrm{F}_{\psi}^{(\infty,0')}(\F_{|\mathbb{P}_{(\infty)}}))} \rangle (t') = \langle \chi_{\det(\mathrm{F}_{\psi}^{(0,\infty')}(\F_{|\mathbb{A}_{(0)}}))} \rangle (t^{\prime-1}),
$$
of $\Lambda$-admissible maps on $G_k$, with notation as in \ref{chap32.30}.
\end{cor}

Indeed, Theorem \ref{chap3xteostat}(iii) yields
$$
\varepsilon_{\bk}(\mathbb{P}, \F) \langle \chi_{\det(\mathrm{F}_{\psi}^{(\infty,0')}(\F_{|\mathbb{P}_{(\infty)}}))} \rangle (t') = \langle \chi_{\det( \mathrm{F}_{\psi}(\F_{|\mathbb{A}})_{| \overline{\eta}_{0'}})} \rangle (t').
$$
Let $D$ be an effective Cartier divisor on $\mathbb{P}'$, supported on $0'$ and $\infty'$, such that the $\mu$-twisted $\Lambda$-local system $\det(\mathrm{F}_{\psi}(\F_{|\mathbb{A}}))_{|\mathbb{A}'\setminus \{ 0' \}}$ of rank $1$ (cf. \ref{chap3xteostat}(i)) has ramification bounded by $D$, and let $\chi_{\det(\mathrm{F}_{\psi}(\F_{|\mathbb{A}}))}$ be the $\mu$-twisted multiplicative $\Lambda$-local system on $\Pic_k(\mathbb{P}',D)_k$ associated to $\det(\mathrm{F}_{\psi}(\F_{|\mathbb{A}}))$ by geometric class field theory, cf. \ref{chap3ggcfttwisted}. If $x_0$ is the $k$-point of $\Pic_k(\mathbb{P}',D)_k$ corresponding to the line bundle $\Ow([0'])$ trivialized by $t^{\prime-1}$ at $0'$ and by $1$ at $\infty'$, then $\langle \chi_{\det( \mathrm{F}_{\psi}(\F_{|\mathbb{A}})_{| \overline{\eta}_{0'}})} \rangle (t')$ is the trace function of the stalk of $\chi_{\det( \mathrm{F}_{\psi}(\F_{|\mathbb{A}})_{| \overline{\eta}_{0'}})}$ at $t^{\prime -1}$ (cf. \ref{chap32.30}), or alternatively of the stalk of $\chi_{\det(\mathrm{F}_{\psi}(\F_{|\mathbb{A}}))}$ at $\overline{x}_0$, by local-global compatibility (cf. \ref{chap32.6}).

Likewise, Theorem \ref{chap3xteostat}(ii) yields
$$
\langle \chi_{\det(\mathrm{F}_{\psi}^{(0,\infty')}(\F_{|\mathbb{A}_{(0)}}))} \rangle (t^{\prime-1}) = \langle \chi_{\det(\mathrm{F}_{\psi}(\F_{|\mathbb{A}})_{| \overline{\eta}_{\infty'}} )} \rangle (t^{\prime -1}),
$$
and the latter coincides, by local-global compatibility (cf. \ref{chap32.6}), with the trace function of the stalk of $\chi_{\det(\mathrm{F}_{\psi}(\F_{|\mathbb{A}}))}$ at $\overline{x}_{\infty}$, where $x_{\infty}$ is the $k$-point of $\Pic_k(\mathbb{P}',D)_k$ corresponding to the line bundle $\Ow([\infty'])$ trivialized by $1$ at $0'$ and by $t'$ at $\infty'$. The conclusion of Proposition \ref{chap3detcor} then follows from the fact that $x_{0} = x_{\infty}$ in $\Pic_k(\mathbb{P}',D)_k$.

\subsection{\label{chap315.2}} Let $T$ be the spectrum of a $k$-algebra, which is a henselian discrete valuation ring $\Ow_T$ with residue field $k$, and let $i : s \rightarrow T$ be its closed point. Let $\pi$ be a uniformizer of $k(\eta)$, and let 
$$
\pi : T \rightarrow \mathbb{A}_{(0)},
$$
be the $k$-morphism sending $\pi$ to the $t$. Laumon's cohomological formula for local $\varepsilon$-factors (\cite{La87}, 3.5.1.1) admits the following extension to the case of an arbitrary perfect base field of positive characteristic $p$.

\begin{teo}\label{chap3xteolaumon} Let $\F$ be a $\mu$-twisted $\Lambda$-sheaf on $T$, with vanishing fiber at $s$. Then we have
$$
\varepsilon_{\overline{k}}(T, \F, d \pi) = \langle \chi_{\det(\mathrm{F}^{(0,\infty')}(\pi_* \F))} \rangle (t^{\prime-1}),
$$
with notation as in \ref{chap32.30}.
\end{teo}

Let us prove Theorem \ref{chap3xteolaumon}. We can assume (and we do) that $T$ is the henselization $\mathbb{A}_{(0)}$ of $\mathbb{A}$ at $0$, and that $\pi = t$ (cf. \ref{chap31.6.0}). Let $\F$ be a $\mu$-twisted $\Lambda$-sheaf on $\mathbb{A}_{(0)}$, with vanishing fiber at $0$, and let $t_{\diamondsuit} \F$ be its Gabber-Katz extension to $\mathbb{A}$ with respect to the uniformizer $t$ (cf. \ref{chap3GK3}). Thus $t_{\diamondsuit} \F$ is tamely ramified at $\infty$, and its fiber at $0$ vanishes. By Theorem \ref{chap3xteostat}(ii), we have an isomorphism
$$
 \mathrm{F}_{\psi}^{(0,\infty')}( \F) \cong \mathrm{F}_{\psi}(t_{\diamondsuit}\F)_{| \overline{\eta}_{\infty'}},
$$
of $\Lambda$-admissible representations of $(G_{\eta_{\infty'}}, \mu)$. Let $\eta_{\infty'}^{\mathrm{perf}}$ (resp. $\overline{\eta}^{\mathrm{perf}}_{\infty'}$) be the perfection of $\eta_{\infty'}$ (resp. $\overline{\eta}_{\infty'}$), so that $\overline{\eta}^{\mathrm{perf}}_{\infty'}$ is an algebraic closure of $\eta_{\infty'}^{\mathrm{perf}}$. By the proper base change theorem, we have
\begin{align*}
 \mathrm{F}_{\psi}(t_{\diamondsuit}\F)_{| \overline{\eta}_{\infty'}} &\cong H_c^1( \mathbb{A}_{\overline{\eta}_{\infty'}}, t_{\diamondsuit}\F \otimes \Lc_{\psi} \lbrace t t' \rbrace) \\
 &\cong H_c^1( \mathbb{A}_{\overline{\eta}^{\mathrm{perf}}_{\infty'}}, t_{\diamondsuit}\F \otimes \Lc_{\psi} \lbrace t t' \rbrace).
\end{align*}
Since the complex $R\Gamma_c( \mathbb{A}_{\overline{\eta}^{\mathrm{perf}}_{\infty'}}, t_{\diamondsuit}\F \otimes \Lc_{\psi} \lbrace t t' \rbrace)$ is concentrated in degree $1$, we obtain
\begin{align}\label{chap3eq897}
\det\left( \mathrm{F}_{\psi}^{(0,\infty')}( \F) \right) \cong \det\left(\mathrm{F}_{\psi}(t_{\diamondsuit}\F)_{| \overline{\eta}_{\infty'}} \right) \cong \det \left(  R\Gamma_c( \mathbb{A}_{\overline{\eta}^{\mathrm{perf}}_{\infty'}}, t_{\diamondsuit}\F \otimes \Lc_{\psi} \lbrace t t' \rbrace)\right)^{-1}.
\end{align}
Let us consider the uniformizer $\widetilde{t} = - t' t$ on $(\mathbb{A}_{\eta^{\mathrm{perf}}_{\infty'}})_{(0)}$, and the isomorphism $\theta$ from $\mathbb{A}_{\eta^{\mathrm{perf}}_{\infty'}}$ to itself which sends $t$ to $-t^{\prime -1} t$. The natural morphism $t : (\mathbb{A}_{\eta^{\mathrm{perf}}_{\infty'}})_{(0)} \rightarrow \mathbb{A}_{\eta^{\mathrm{perf}}_{\infty'}} $ factors as a composition
$$
(\mathbb{A}_{\eta^{\mathrm{perf}}_{\infty'}})_{(0)} \xrightarrow[]{\widetilde{t}} \mathbb{A}_{\eta^{\mathrm{perf}}_{\infty'}} \xrightarrow[]{\theta} \mathbb{A}_{\eta^{\mathrm{perf}}_{\infty'}},
$$
hence we have isomorphisms
\begin{align*}
\theta^{-1} t_{\diamondsuit}\F  &\cong \widetilde{t}_{\diamondsuit}\F\\
\theta^{-1} \Lc_{\psi} \lbrace t t' \rbrace  &\cong \Lc_{\psi}^{-1},
\end{align*}
where $\widetilde{t}_{\diamondsuit} \F$ is the Gabber-Katz extension of (the restriction to $(\mathbb{A}_{\eta^{\mathrm{perf}}_{\infty'}})_{(0)}$ of) $\F$ to $\mathbb{A}$ with respect to the uniformizer $\widetilde{t}$ (cf. \ref{chap3GK3}). Hence (\ref{chap3eq897}) yields
$$
\det\left( \mathrm{F}_{\psi}^{(0,\infty')}( \F) \right) \cong \det \left(  R\Gamma_c( \mathbb{A}_{\overline{\eta}^{\mathrm{perf}}_{\infty'}}, \widetilde{t}_{\diamondsuit}\F \otimes \Lc_{\psi}^{-1})\right)^{-1}.
$$
By Definition \ref{chap3localepsfact} with $\pi = \widetilde{t}$, we obtain that the composition of $\varepsilon_{\overline{k}}(\mathbb{A}_{(0)}, \F, d t)$  with the canonical surjective homomorphism $r$ from $ G_{\eta_{\infty'}^{\mathrm{perf}}} = G_{\eta_{\infty'}}$ to $G_k$ is given by
\begin{align}
\begin{split}\label{chap3eq852}
\varepsilon_{\overline{k}}(\mathbb{A}_{(0)}, \F, d t) \circ r &=   \varepsilon_{\overline{\eta}^{\mathrm{perf}}_{\infty'}}((\mathbb{A}_{\eta^{\mathrm{perf}}_{\infty'}})_{(0)}, \F, d t) \\
&= \langle \chi_{\det(j^{-1} \F)} \rangle \left( \frac{dt}{d\widetilde{t}} \right) \det\left( \mathrm{F}_{\psi}^{(0,\infty')}( \F) \right) \\
&= \langle \chi_{\det(j^{-1} \F)} \rangle \left( -t^{\prime -1} \right) \det\left( \mathrm{F}_{\psi}^{(0,\infty')}( \F) \right),
\end{split}
\end{align}
where $\chi_{\det(j^{-1} \F)}$ is the multiplicative $\Lambda$-local system on $\Pic(\mathbb{A}_{(0)}, \nu [0])_k$,  for some integer $\nu$ such that $\det(j^{-1} \F)$ has ramification bounded by $\nu$, associated to $\det(j^{-1} \F)$ by geometric local class field theory, cf. \ref{chap3lgcft2twisted}, and $ \langle \chi_{\det(j^{-1} \F)} \rangle \left(  -t^{\prime -1} \right)$ is the trace function of the stalk of $\chi_{\det(j^{-1} \F)}$ at the $\eta_{\infty'}$-point of $\Pic^0(\mathbb{A}_{(0)}, \nu [0])_k$ corresponding to the unit $-t'$ (cf. \ref{chap32.5}). Let us rewrite (\ref{chap3eq852}) as
\begin{align}\label{chap3eq5469}
 \det\left( \mathrm{F}_{\psi}^{(0,\infty')}( \F) \right) = \langle \chi_{\det(j^{-1} \F)} \rangle \left( -t'\right) \left( \varepsilon_{\overline{k}}(\mathbb{A}_{(0)}, \F, d t) \circ r \right).
\end{align}
We now use (\ref{chap3eq5469}) in order to compute the multiplicative $\Lambda$-local system associated to the determinant of $ \mathrm{F}_{\psi}^{(0,\infty')}( \F)$ by geometric local class field theory, cf. \ref{chap3lgcft2twisted}. Let us consider the local Abel-Jacobi morphism $\Phi_{\eta_{\infty'}, t^{\prime -1}}$ corresponding to the divisor $[\infty']$ on $\mathbb{P}'_{(\infty')}$ (cf. \ref{chap3lgcft2twisted}). The morphism $\Phi_{\eta_{\infty'}, t^{\prime -1}}$ factors as the composition
$$
\eta_{\infty'} \xrightarrow[]{t^{\prime -1}} \mathbb{G}_{m,k} \xrightarrow[]{t \mapsto 1 - t t'} \Pic^1( \mathbb{P}'_{(\infty')}, [\infty'] )_k.
$$
Let us consider the $k$-isomorphism $\tau$ from $\Pic^1( \mathbb{P}'_{(\infty')}, [\infty'] )_k$ to $\mathbb{G}_{m,k}$ which sends a section $u$ (cf. \ref{chap32.5}) to the section $t^{\prime -1} u$ of $\Pic^0( \mathbb{P}'_{(\infty')}, [\infty'] )_k = \mathbb{G}_{m,k}$. For any section $t$ of $\mathbb{G}_{m,k}$, the sections $1 - t t'$ and $-tt'$ of $\Pic^1( \mathbb{P}'_{(\infty')}, [\infty'] )_k$ coincides, hence the following commutative diagram.
\begin{center}
 \begin{tikzpicture}[scale=1]

\node (A) at (0,0) {$\eta_{\infty'}$};
\node (B) at (4,0) {$\Pic^1( \mathbb{P}'_{(\infty')}, [\infty'] )_k$};
\node (C) at (8,0) {$\mathbb{G}_{m,k}$};

\path[->,font=\scriptsize]
(A) edge node[above]{$\Phi_{\eta_{\infty'}, t^{\prime -1}}$} (B)
(B) edge node[above]{$\tau$} (C)
(A) edge[bend right] node[above]{$-t^{\prime -1}$} (C);
\end{tikzpicture} 
\end{center}
Let $\chi_0$ be the restriction of $\chi_{\det(j^{-1} \F)}$ to the subgroup $\mathbb{G}_{m,k}$ of $\Pic^0(\mathbb{A}_{(0)}, \nu [0])_k$, and let $\widetilde{\chi}$ be the pullback of $\chi_0 \otimes  \varepsilon_{\overline{k}}(\mathbb{A}_{(0)}, \F, d t) $ by $\tau$, where $ \varepsilon_{\overline{k}}(\mathbb{A}_{(0)}, \F, d t) $ is considered as a $\Lambda$-local system on $\Spec(k)$, pulled back to $\mathbb{G}_{m,k}$. Then the commutative diagram above, together with (\ref{chap3eq5469}), shows that $\Phi_{\eta_{\infty'}, t^{\prime -1}}^{-1} \widetilde{\chi}$ is isomorphic to $ \det\left( \mathrm{F}_{\psi}^{(0,\infty')}( \F) \right)$.

Moreover, $\widetilde{\chi}$ is the restriction to $\Pic^1( \mathbb{P}'_{(\infty')}, [\infty'] )_k$ of the unique (up to isomorphism) multiplicative $\Lambda$-local system on $\Pic( \mathbb{P}'_{(\infty')}, [\infty'] )_k$, still denoted by $\widetilde{\chi}$, whose restriction to $\Pic^0( \mathbb{P}'_{(\infty')}, [\infty'] )_k = \mathbb{G}_{m,k}$ is given by $\chi_0$, and whose stalk at the $k$-point $t'$ of $\Pic^1( \mathbb{P}'_{(\infty')}, [\infty'] )_k$ is given by $\varepsilon_{\overline{k}}(\mathbb{A}_{(0)}, \F, d t)$. Thus the multiplicative $\Lambda$-local system associated to $\det\left( \mathrm{F}_{\psi}^{(0,\infty')}( \F) \right)$ by geometric local class field theory (cf. \ref{chap3lgcft2twisted}) is $\widetilde{\chi}$ (up to isomorphism). We therefore obtain the equality
$$
\langle \chi_{\det\left( \mathrm{F}_{\psi}^{(0,\infty')}( \F) \right)} \rangle (t^{\prime -1}) = \langle \widetilde{\chi}\rangle(t^{\prime -1}) = \varepsilon_{\overline{k}}(\mathbb{A}_{(0)}, \F, d t),
$$
which concludes the proof of Theorem \ref{chap3xteolaumon}.

\begin{rema} This proof of Theorem \ref{chap3xteolaumon} incidentally shows that $\det\left( \mathrm{F}_{\psi}^{(0,\infty')}( \F) \right)$ is tamely ramified. The latter fact alternatively follows from the Hasse-Arf theorem and from (\cite{La87}, 2.4.3(i)(b)), which asserts that the ramification breaks of $ \mathrm{F}_{\psi}^{(0,\infty')}( \F)$ are strictly less than $1$.
\end{rema}

\begin{cor}\label{chap3infzero} Let $\F$ be a $\mu$-twisted $\Lambda$-sheaf on $\mathbb{P}_{(\infty)}$, with vanishing fiber at $\infty$ and with ramification bounded by $1$. Then we have
$$
\varepsilon_{\overline{k}}(\mathbb{P}_{(\infty)}, \F, d t) \langle \chi_{\det(\mathrm{F}_{\psi}^{(\infty,0')}(\F))} \rangle (t')= \chi_{\cyc}^{\rk(\F)}.
$$
\end{cor}

By Theorem \ref{chap3GK3} there exists a $\Lambda$-sheaf $\G$ on $\mathbb{A}$, with vanishing fiber at $0$, which is tamely ramified at $0$, and such that $\G_{| \mathbb{P}_{(\infty)}}$ is isomorphic to $\F$. Then Corollary \ref{chap3detcor} and Theorem \ref{chap3xteolaumon} yield
\begin{align}\label{chap3eq213}
\varepsilon_{\bk}(\mathbb{A}, \G) \langle \chi_{\det(\mathrm{F}_{\psi}^{(\infty,0')}(\F))} \rangle (t') = \varepsilon_{\overline{k}}(\mathbb{A}_{(0)}, \G, d t),
\end{align}
while Proposition \ref{chap3productform20}, whose hypotheses are satisfied by \ref{2.1.0} and \ref{localstructure}, yields
\begin{align}\label{chap3eq214}
\varepsilon_{\bk}(\mathbb{A}, \G) = \chi_{\cyc}^{-\rk(\F)} \varepsilon_{\overline{k}}(\mathbb{A}_{(0)}, \G, d t) \varepsilon_{\overline{k}}(\mathbb{P}_{(\infty)}, \F, d t).
\end{align}
The conclusion of Corollary \ref{chap3infzero} then follows by combining (\ref{chap3eq213}) with (\ref{chap3eq214}).


\subsection{\label{chap315.8}} We now prove Theorem \ref{chap3xteo3}. Its conclusion holds when $\F$ is generically of rank at most $1$ by \ref{chap3productform2}. In particular, it holds for the constant sheaf $\Lambda$. Thus Theorem \ref{chap3xteo3}
 is equivalent to the formula
$$
\frac{\varepsilon_{\bk}(X, \F_1)}{\varepsilon_{\bk}(X, \F_2)} = \frac{\tilde{\varepsilon}_{\bk}(X, \F_1, \omega)}{\tilde{\varepsilon}_{\bk}(X, \F_2, \omega)},
$$
cf. \ref{defiproductform}, for any $\F_1,\F_2$ satisfying the assumptions of Theorem \ref{chap3xteo3} (twisted by possibly different cocycles), with the same generic rank. If $f : X \rightarrow \mathbb{P}$ is a finite generically \'etale $k$-morphism, then by Proposition \ref{inductionlemma} the latter formula holds for $(X,\F_1,\F_2)$ if and only if it does for $(\mathbb{P}, f_* \F_1, f_* \F_2)$. Thus the conclusion of Theorem \ref{chap3xteo3} holds in general if and only if it holds for $X = \mathbb{P}$. If $X = \mathbb{P}$, then there exists a non empty open subscheme $U$ of $\mathbb{A}$ such that $\F_{1|U}$ and $\F_{2|U}$ are $\Lambda$-local systems, and we can find a polynomial $h$ in $k[t]$ whose vanishing locus in $\mathbb{A}$ is the complement of $U$ in $\mathbb{A}$. The $k$-morphism $\theta : \mathbb{P} \rightarrow \mathbb{P}$ which sends $t$ to $(t - h(t)^{-p})^{-1}$ is finite, and induces a finite \'etale morphism from $U$ onto $\mathbb{P} \setminus \{ 0 \}$. By replacing $(\mathbb{P},\F_1,\F_2)$ with $(\mathbb{P}, \theta_* \F_1, \theta_* \F_2)$, we can thus assume that the restrictions of $\F_1$ and $\F_2$ to $\mathbb{P} \setminus \{ 0 \}$ are $\Lambda$-local systems.

\begin{rema} The last reduction to $\Lambda$-sheaves on $\mathbb{P}$ with ramification concentrated on a single point is due to Katz, cf. (\cite{katz3}, Lemma 16).
\end{rema}

In order to prove Theorem \ref{chap3xteo3}, we can thus assume (and we do) that $X = \mathbb{P}$ and that $\F$ is a $\Lambda$-sheaf on $\mathbb{P}$, whose restriction to $\mathbb{A}\setminus \{ 0 \}$ is a $\Lambda$-local system, which is unramified at $\infty$. We can further assume (and we do) that the fibers of $\F$ at $0$ and $\infty$ vanish. The formula to be proved is then
$$
\varepsilon_{\bk}(\mathbb{P}, \F) = \chi_{\cyc}^{-\rk(\F)} \varepsilon_{\bk}(\mathbb{A}_{(0)}, \F, dt) \varepsilon_{\bk}(\mathbb{P}_{(\infty)}, \F, dt),
$$ 
By Corollary \ref{chap3detcor}, we have
$$
\varepsilon_{\bk}(\mathbb{P}, \F) = \langle \chi_{\det(\mathrm{F}_{\psi}^{(0,\infty')}(\F_{|\mathbb{A}_{(0)}}))} \rangle (t^{\prime-1}) \langle \chi_{\det(\mathrm{F}_{\psi}^{(\infty,0')}(\F_{|\mathbb{P}_{(\infty)}}))} \rangle (t')^{-1},
$$
and the conclusion then results from the formulas
\begin{align*}
\langle\chi_{\det(\mathrm{F}_{\psi}^{(0,\infty')}(\F_{|\mathbb{A}_{(0)}}))} \rangle (t^{\prime-1}) &= \varepsilon_{\bk}(\mathbb{A}_{(0)}, \F, dt), \\ 
\langle \chi_{\det(\mathrm{F}_{\psi}^{(\infty,0')}(\F_{|\mathbb{P}_{(\infty)}}))} \rangle (t')^{-1} &= \chi_{\cyc}^{-\rk(\F)} \varepsilon_{\bk}(\mathbb{P}_{(\infty)}, \F, dt)
,
\end{align*}
which follow respectively from Theorem \ref{chap3xteolaumon} and from Corollary \ref{chap3infzero}.

\bibliographystyle{amsalpha}

\end{document}